\documentclass[12pt]{article}
\usepackage{amsmath,amssymb,amsfonts}
\usepackage{graphicx}
\usepackage{epsfig}
\usepackage{amsfonts}
\usepackage{amssymb, latexsym, amsmath, pb-diagram}
\usepackage{pstricks-add}
\usepackage{multicol}
\usepackage{bm}
\usepackage{mathrsfs}
\usepackage{xy}
\usepackage{epic,eepic}
\DeclareMathAlphabet{\mathpzc}{OT1}{pzc}{m}{it}
\xyoption{all}
\begin{document}
\title{The Cohomology Algebra of Polyhedral Product Objects}
\author {
Qibing Zheng\\School of Mathematical Science and LPMC, Nankai University\\
Tianjin 300071, China\\
zhengqb@nankai.edu.cn\footnote{Project Supported by Natural Science
Foundation of China, grant No. 11071125 and No. 11671154\newline
\hspace*{5.5mm}Key words and phrases: polyhedral product, diagonal tensor product.\newline
\hspace*{5.5mm}Mathematics subject classification: 55N10}}\maketitle
\input amssym.def
%\newsymbol\leqslant 1336
%\newsymbol\geqslant 133E
\baselineskip=20pt
\def\w{\widetilde}
\def\makm{\mathbf k{\scriptstyle[m]}}
\def\makn{\mathbf k{\scriptstyle[n]}}
\def\mak{\mathbf k}
\def\AA{\mathscr A}
\def\BB{\mathscr B}
\def\CC{\mathscr C}
\def\DD{\mathscr D}
\def\RR{\mathscr R}
\def\SS{\mathscr S}
\def\XX{\mathscr X}
\def\TT{\mathscr T}
\def\NN{\mathscr N}
\def\LL{\mathscr L}
\def\ZZ{\mathscr Z}

\begin{abstract} In this paper, we compute the homology group and cohomology algebra
of various polyhedral product objects uniformly from the point of view of  diagonal tensor product.
As applications, we introduce the polyhedral product method into commutative algebra
and show that the homotopy types of polyhedral product spaces
depend on not only the homotopy type of each summand pair but also on
the character coproduct of the pair.
\end{abstract}

\hspace*{40mm}${\displaystyle{\bf Table\,\,of\,\,Contents}}$

Section 1\, Introduction

Section 2\, Diagonal Tensor Product of Indexed Groups and Complexes

Section 3\, (Co)homology Group of Polyhedral Product Chain Complexes

Section 4\, (Co)homology Group of Polyhedral Product Objects

Section 5\, Duality Isomorphisms

Section 6\, Diagonal Tensor Product of (Co)algebras

Section 7\, Local (Co)products of Total Objects of Simplicial Complexes

Section 8\, Cohomology Algebra of Polyhedral Product Chain Complexes

Section 9\, Cohomology Algebra of Polyhedral Product Objects
\vspace{3mm}

\section{Introduction}\vspace{3mm}

The polyhedral product functor is a relatively new construction with its origins
in toric topology and therefore closely related to toric objects coming from
algebraic and symplectic geometry. Since its formal appearance in work of
Buchstaber and Panov \cite{BP} and Grbi\'{c} and Theriault \cite{GR} in 2004,
the theory (especially the homotopical characteristics) of polyhedral products
has been developing rapidly. Due to its combinatorial nature coming from the underlying
simplicial complex and being a product space, the polyhedral product functor was quickly
recognized  as a complex construction but at the same time approachable.
As a result, polyhedral products are nowadays used not only in topology and geometry
but also group theory (abstract and geometric) as a collection of spaces on which to test
major conjectures and to form a rather delicate insight in building new theories.

In this paper, we compute the cohomology algebra of various polyhedral product objects
uniformly from the point of view of diagonal tensor product.

In section 2, we give the basic definitions of the paper, indexed groups and indexed complexes in Definition~2.1 and
the tensor product and diagonal tensor product of indexed objects in Definition~2.4 and Definition~2.5.
A $\Lambda$-indexed group is just a graded group simultaneously graded by the set $\Lambda$.
To distinguish the two different gradations, we use the word index.
The diagonal tensor product is always defined with respect to an index set but has nothing to do with the usual grade.

In section 3, we give the definition of polyhedral product chain complex
${\cal Z}(K;\underline{D_*},\underline{C_*})$ in Definition~3.1, which is the most general model for the computation
of the (co)homology of polyhedral product objects.
From the point of view of diagonal tensor product,
we compute the (co)homology group of split polyhedral product chain complexes in Theorem~3.12.
Precisely, the (co)homology group is a diagonal tensor product group
$$H_*({\cal Z}(K;\underline{D_*},\underline{C_*}))\cong H_*^{\XX_m}(K)\widehat\otimes H_*^{\XX_m}(\underline{D_*},\underline{C_*}),$$
$$H^*({\cal Z}(K;\underline{D_*},\underline{C_*}))\cong H^{\,*}_{\!\XX_m}(K)\widehat\otimes H^{\,*}_{\!\XX_m}(\underline{D_*},\underline{C_*}),$$
where $H_-^-(K)$ is the total (co)homology group of $K$ in Definition~3.7
and $H_-^-(\underline{D_*},\underline{C_*})$ is the indexed (co)homology group of $(\underline{D_*},\underline{C_*})$ in Definition~3.11.

In section 4, we compute the (co)homology group of various polyhedral product objects uniformly by Theorem~3.12,
precisely, the cohomology group of
polyhedral product space (simplicial complex) ${\cal Z}(K,\underline{X},\underline{A})$,
polyhedral smash product space ${\cal Z}^\wedge(K,\underline{X},\underline{A})$ and
polyhedral join space (simplicial complex) ${\cal Z}^*(K,\underline{X},\underline{A})$ in Theorem~4.6,
the total (co)homology group of polyhedral  join simplicial complex ${\cal Z}^*(K,\underline{X},\underline{A})$ in Theorem~4.13,
the total and right total (co)homology group of composition complex ${\cal Z}^*(K;L_1,{\cdots},L_m)$ in Example~4.15,
the Tor group of polyhedral tensor product module ${\cal Z}^\otimes(K,\underline{X},\underline{A})$ in Theorem~4.18,
the Tor group of composition ideals ${\cal Z}^\otimes(K;I_1,{\cdots},I_m)$ in Theorem~4.19.

In section 5, we prove duality isomorphisms between various polyhedral product objects.
In Theorem~5.2, we prove that for a polyhedral product space ${\cal Z}(K;\underline{X},\underline{A})$,
the complement space $(X_1{\times}{\cdots}{\times}X_m){\setminus}{\cal Z}(K;\underline{X},\underline{A})$ is the polyhedral product space
${\cal Z}(K^\circ;\underline{X},\underline{A}^c)$,
where $K^\circ$ is the Alexander dual of $K$ and $A_k^c=X_k{\setminus}A_k$.
In Theorem~5.4, we prove the duality isomorphisms
$$\gamma_{K,\sigma\!,\,\omega}\colon H_*^{\sigma\!,\,\omega}(K)\to
H^{|\omega|-*-1}_{\sigma'\!,\,\omega}(K^\circ),\quad
\gamma^\circ_{K,\sigma\!,\,\omega}\colon H^*_{\sigma\!,\,\omega}(K)
\to H_{|\omega|-*-1}^{\sigma'\!,\,\omega}(K^\circ),$$
that are very important for other duality isomorphisms.
In Definition~5.5, we define the duality isomorphisms
$$\gamma_K\colon H_*^{\XX;\XX_m}(\Delta\!^{[m]}{\!,}K)\to H^{\,*}_{\!\XX;\XX_m}(\Delta\!^{[m]}{\!,}K^\circ),$$
$$\gamma_K^{\,\circ}\colon H^{\,*}_{\!\XX;\XX_m}(\Delta\!^{[m]}{\!,}K)\to H_*^{\XX;\XX_m}(\Delta\!^{[m]}{\!,}K^\circ)\,$$
and their restriction isomorphisms
$$\overline\gamma_K\colon H_*^{\LL_m}(K)\to H^{\,*}_{\!\LL_m}(K^\circ),\quad
\overline\gamma_K^{\,\circ}\colon H^{\,*}_{\!\LL_m}(K)\to H_*^{\LL_m}(K^\circ).$$
In Theorem~5.6, we prove the duality equality for composition complexes
$${\cal Z}^*(K;L_1,{\cdots},L_m)^\circ={\cal Z}^*(K^\circ;L_1^\circ,{\cdots},L_m^\circ).$$
In Theorem~5.7, we prove
$$\gamma_{{\cal Z}^*(K;L_1,\cdots,L_m),\tilde\sigma,\tilde\omega}
=\gamma_{K,\sigma,\omega}{\otimes}\gamma_{L_1,\sigma_1,\omega_1}{\otimes}{\cdots}{\otimes}\gamma_{L_m,\sigma_m,\omega_m}$$
In Theorem~5.9, we generalize Hoschter Theorem from face ideal $I_K$ to monomial ideal $I_{(K;\underline{r})}$.
In Theorem~5.10, we prove the ideal equality
$${\cal Z}^\otimes(K;I_{(L_1;\underline{r_1})},{\cdots},I_{(L_m;\underline{r_m})})
=I_{({\cal Z}^*(K^\circ;L_1,{\cdots},L_m);(\underline{r_1},{\cdots},\underline{r_m}))}$$
and compute the group ${\rm Tor}^{\mak[n]}_*({\cal Z}^\otimes(K;I_{(L_1;\underline{r_1})},{\cdots},I_{(L_m;\underline{r_m})}),\mak)$.
In Definition~5.11, we define the complement duality isomorphisms
$$\gamma_{(\underline{X},\underline{A})}\colon H_*^{\XX_m}(\underline{X},\underline{A})\to H^{\,*}_{\!\XX_m}(\underline{X},\underline{A}^c),$$
$$\gamma_{(\underline{X},\underline{A})}^\circ\colon H^{\,*}_{\!\XX_m}(\underline{X},\underline{A})\to H_*^{\XX_m}(\underline{X},\underline{A}^c)\,$$
based on the Alexander duality isomorphism,
where each $A_k$ is a polyhedron subspace of the manifold $X_k$ and then in Theorem~5.14, we prove that
the Alexander duality isomorphisms
$$\alpha\colon H_*({\cal Z}(K;\underline{X},\underline{A}))\to H^*(X_1{\times}{\cdots}{\times}X_m,{\cal Z}(K;\underline{X},\underline{A})^c)$$
$$\alpha^\circ\colon H^*({\cal Z}(K;\underline{X},\underline{A}))\to H_*(X_1{\times}{\cdots}{\times}X_m,{\cal Z}(K;\underline{X},\underline{A})^c)$$
have decompositions $\alpha=\hat\alpha\oplus(\overline\gamma_K{\widehat\otimes}\overline\gamma_{(\underline{X},\underline{A})})$,
$\alpha^\circ=\hat\alpha^\circ\oplus(\overline\gamma_K^{\,\circ}{\widehat\otimes}\overline\gamma_{(\underline{X},\underline{A})}^{\,\circ})$.

In section 6, we generalize the results of section 2 from indexed groups and complexes to
indexed (co)algebras, chain coalgebras and cochain algebras.
We stress that in this paper, a (co)algebra may not be (co)associative or degree-preserving.
When $A_*^\DD$ and $B_*^\DD$ are indexed groups,
the diagonal tensor product group $A_*^\DD\widehat\otimes B_*^\DD$ is naturally a subgroup of $A_*^\DD{\otimes}B_*^\DD$.
When $A_*^\DD$ and $B_*^\DD$ are (co)algebras,
the (co)algebra $A_*^\DD\widehat\otimes B_*^\DD$ is in general not a sub(co)algebra or a quotient (co)algebra of  $A_*^\DD{\otimes}B_*^\DD$.
The elements of $A_*^\DD\widehat\otimes B_*^\DD$ have to be written in the form $a\widehat\otimes b$.

In section 7, we compute the local coproducts of the total homology coalgebra
and the local products of the total cohomology algebra of simplicial complexes in Theorem~7.7.
We correct a mistake in Definition~7.1 of \cite{Z}.

In  section 8, we prove that the cohomology group isomorphisms in Theorem~3.12
are algebra isomorphisms if each $(D_{k\,*},C_{k\,*})$ is a split coalgebra pair
in the polyhedral product chain  complex ${\cal Z}(K;\underline{D_*},\underline{C_*})$.

In section 9, we compute the cohomology algebra of polyhedral product spaces and simplicial complexes ${\cal Z}(K;\underline{X},\underline{A})$
in Theorem~9.3 and the total cohomology algebra of polyhedral join simplicial complex ${\cal Z}^*(K;\underline{X},\underline{A}))$ in Theorem~9.7.
As an application, we give two polyhedral product spaces
${\cal Z}(K;X_1,A_1)$ and ${\cal Z}(K;X_2,A_2)$ such that
$A_1\simeq A_2$, $X_1\simeq X_2$ and the two cohomology homomorphisms
$H^*(X_i)\to H^*(A_i)$ induced by inclusion are the same,
but the cohomology rings $H^*({\cal Z}(K;X_1,A_1))$ and $H^*({\cal Z}(K;X_2,A_2))$ are not isomorphic.
\vspace{3mm}

\section{Diagonal Tensor Product of Indexed Groups and Complexes}\vspace{3mm}

\hspace{5.5mm}{\bf Definition 2.1} An {\it indexed group},  or a {\it group indexed by} $\Lambda$, or a {\it $\Lambda$-group}
is a direct sum over $\Lambda$ of graded groups $A_*^\Lambda=\oplus_{\alpha\in\Lambda}\,A_*^{\alpha}$.
Each $A_*^\alpha$ is called a {\it local group} of $A_*^\Lambda$.

An {\it indexed chain complex}, or a {\it chain complex indexed by} $\Lambda$, or a {\it chain $\Lambda$-complex}
is a direct sum over $\Lambda$ of complexes $(C_*^\Lambda,d)=\oplus_{\alpha\in\Lambda}(C_*^{\alpha},d)$.
Each $(C_*^\alpha,d)$ is called a {\it local complex} of $(C_*^\Lambda,d)$.

An {\it indexed cochain complex}, or a {\it cochain complex indexed by} $\Lambda$, or a {\it cochain $\Lambda$-complex}
is a direct sum over $\Lambda$ of complexes $(C^*_\Lambda,\delta)=\oplus_{\alpha\in\Lambda}(C^*_{\alpha},\delta)$.
Each $(C^*_\alpha,\delta)$ is called a {\it local complex} of $(C^*_\Lambda,\delta)$.

The set $\Lambda$ is called the index set of the above indexed objects.
$\Lambda$ is called {\it trivial} if  it has only one element.
\vspace{3mm}

A $\Lambda$-group $A_*^\Lambda$ is just the graded group $A_*$ that is also graded by the set $\Lambda$.
To distinguish the two different gradations, we use the word index.
When considering (co)algebra structures, we have to use such definitions as
graded chain $\Lambda$-coalgebras. So the distinction of grade and index is inevitable.

All the basic properties of graded groups and (co)chain complexes can be generalized to indexed groups
and indexed (co)chain complexes.
\vspace{2mm}

{\bf Conventions.} For simplicity, we often have to abbreviate the index set of an indexed group.
When the index set is denoted by Greek letters such as $\Lambda,\Gamma,\cdots$,
the indexed groups $A_*^\Lambda,A_*^\Gamma,\cdots$ are often simply  denoted by $A_*$.
When the index set is denoted by $\DD,\SS,\TT,\XX,\RR,\cdots$,
the indexed groups $A_*^\DD,A_*^\SS,A_*^\TT,A_*^\XX,A_*^\RR,\cdots$ can not be simplified.
So the indexed group $A_*^{\DD{\times}\Lambda}$ is often simply denoted by $A_*^\DD$.

A group $A$ is always regarded as a graded group $A_*$ such that $A_0=A$ and $A_k=0$ if $k\neq 0$.
Specifically, the integer group $\Bbb Z$ is regarded as such a graded group.
A graded group $A_*$ is always regarded as a $\Lambda$-group $A_*$ such that $\Lambda$ is trivial.

Analogue conventions hold for (co)chain complexes.
\vspace{3mm}

{\bf Lemma 2.2} {\it Every subgroup or quotient group of a $\Lambda$-group
is also a $\Lambda$-group.

For a free $\Lambda$-group $A_*=\oplus_{\alpha\in\Lambda}\,A_*^\alpha$,
its dual group $A^*=\oplus_{\alpha\in\Lambda}\,A^*_\alpha$ with
$A^*_\alpha={\rm Hom}(A_*^\alpha,\Bbb Z)$ ($A^*_\alpha={\rm Hom}(A_*^\alpha,\mak)$
if $A_*$ is a vector space over $\mak$) is also a $\Lambda$-group.
For a free subgroup $B_*$ of $A_*$ such that $A_*/B_*$ is also free,
the dual group $B^*$ is the quotient group $A^*/I^*$, where $I^*=\{f\in A^*\,|\, f(a)=0\}$.

Analogue conclusions hold for indexed chain and cochain complexes.
\vspace{2mm}

Proof}\, Suppose $B_*$ is a subgroup of $A_*$.
Then $B_*=\oplus_{\alpha\in\Lambda}\,B_*{\cap}A_*^{\alpha}$ and
$A_*/B_*=\oplus_{\alpha\in\Lambda}\,A_*^{\alpha}/(B{\cap}A_*^{\alpha})$ are naturally $\Lambda$-groups.

If $B_*$ is a subgroup of the gree group $A_*$ such that $A_*/B_*$ is also free, then
we have a decomposition $A_*=B_*\oplus(A_*/B_*)$ and its dual decomposition $A^*=B^*\oplus(A_*/B_*)^*$.
For a $f\in  I^*$, there is a quotient $\overline f\colon A_*/B_*\to \Bbb Z$ given by $\overline f([a])=[f(a)]$ for $[a]\in A_*/B_*$.
Conversely, for a $\overline f\in (A_*/B_*)^*$, $\overline f$ can naturally be extended to $f\colon A_*\to\Bbb Z$
by defining $f(a)=0$ for all $a\in B_*$. So the correspondence $f\to\overline f$
is an isomorphism from $I^*$ to $(A_*/B_*)^*$.
\hfill$\Box$\vspace{3mm}

{\bf Definition 2.3} A {\it  $\Lambda$-group homomorphism} from $\Lambda$-group $A_*=\oplus_{\alpha\in\Lambda}A_*^\alpha$
to $\Lambda$-group $B_*=\oplus_{\alpha\in\Lambda}B_*^\alpha$ is the direct sum $f=\oplus_{\alpha\in\Lambda}f_\alpha\colon A_*\to B_*$
such that each $f_\alpha\colon A_*^{\alpha}\to B_*^{\alpha}$ is a graded group homomorphism
and is called a {\it local homomorphism} of $f$.
By Lemma~2.2, the kernel, cokernel, coimage and image of $f$ are all $\Lambda$-groups as follows.
$$\begin{array}{ll}
{\rm ker}\,f=\oplus_{\alpha\in\Lambda}{\rm ker}\,f_\alpha,&{\rm coker}\,f=\oplus_{\alpha\in\Lambda}{\rm coker}\,f_\alpha,\vspace{2mm}\\
{\rm coim}\,f\,=\oplus_{\alpha\in\Lambda}\,{\rm coim}\,f_\alpha,&{\rm im}\,f\,=\oplus_{\alpha\in\Lambda}\,{\rm im}\,f_\alpha.
\end{array}$$

We have analogue definitions for indexed (co)chain complexes.
\vspace{3mm}

{\bf Definition 2.4} Let $A_*$ be a $\Lambda$-group and $B_*$ be a $\Gamma$-group.
Their tensor product $A_*\otimes B_*=\oplus_{(\alpha,\beta)\in\Lambda\times\Gamma}\,A_*^{\alpha}\otimes B_*^{\beta}$
is naturally a $(\Lambda{\times}\Gamma)$-group.

For $\Lambda$-group homomorphism $f\colon A_*\to B_*$ and $\Gamma$-group homomorphism $g\colon C_*\to D_*$,
their tensor product $f{\otimes}g\colon A_*{\otimes}C_*\to B_*{\otimes}D_*$
is naturally a $(\Lambda{\times}\Gamma)$-group homomorphism with
$f{\otimes}g=\oplus_{(\alpha,\beta)\in\Lambda{\times}\Gamma}f_\alpha{\otimes}g_\beta$.

We have analogue definitions for indexed (co)chain complexes.
\vspace{3mm}

{\bf Definition~2.5} Let $A_*^\DD$ be a $(\DD{\times}\Lambda)$-group
and $B_*^\DD$ be a $(\DD{\times}\Gamma)$-group.
Precisely, $A_*^\DD=A_*^{\DD;\Lambda}=\oplus_{s\in\DD,\alpha\in\Lambda}\,A_*^{s;\alpha}$ and $B_*^\DD=B_*^{\DD;\Gamma}=\oplus_{s\in\DD,\beta\in\Gamma}\,B_*^{s;\beta}$.
Their {\it diagonal tensor product group} with respect to $\DD$ is the $(\DD{\times}\Lambda{\times}\Gamma)$-group
$$A_*^\DD\,\widehat\otimes\, B_*^\DD=\oplus_{s\in\DD,\,\alpha\in\Lambda,\,\beta\in\Gamma}\,A_*^{s;\alpha}{\otimes}B_*^{s;\beta}.$$

For $(\DD{\times}\Lambda)$-group homomorphism $f\colon A_*^\DD\to C_*^\DD$ and $(\DD{\times}\Gamma)$-group homomorphism $g\colon B_*^\DD\to D_*^\DD$,
their {\it diagonal tensor product homomorphism} with respect to $\DD$ is the $(\DD{\times}\Lambda{\times}\Gamma)$-group homomorphism
$$f\,\widehat\otimes\,g=\oplus_{s\in\DD,\,\alpha\in\Lambda,\,\beta\in\Gamma}\,f_{s;\alpha}{\otimes}g_{s;\beta}
\colon A_*^\DD\widehat\otimes B_*^\DD\to C_*^\DD\widehat\otimes D_*^\DD.$$

We have analogue definitions for indexed (co)chain complexes.
\vspace{3mm}

{\bf Theorem~2.6} {\it We have indexed group (complex) isomorphism and
indexed group (complex) homomorphism equality
$$(A_*^{\DD_1}\widehat\otimes B_*^{\DD_1}){\otimes}{\cdots}{\otimes}
(A_*^{\DD_m}\widehat\otimes B_*^{\DD_m})
\cong(A_*^{\DD_1}{\otimes}{\cdots}{\otimes}A_*^{\DD_m})\widehat\otimes
(B_*^{\DD_1}{\otimes}{\cdots}{\otimes}B_*^{\DD_m}),
$$
$$(f_1\widehat\otimes g_1){\otimes}{\cdots}{\otimes}
(f_m\widehat\otimes g_m)
=(f_1{\otimes}{\cdots}{\otimes}f_m)\widehat\otimes
(g_1{\otimes}{\cdots}{\otimes}g_m),
$$
where the right side diagonal tensor product is with respect to $\DD_1{\times}{\cdots}{\times}\DD_m$.
\vspace{2mm}

Proof}\, Let the factor-permuting isomorphism
$$\phi\colon A_*^{\DD_1}{\otimes}B_*^{\DD_1}{\otimes}{\cdots}{\otimes}
A_*^{\DD_m}{\otimes}B_*^{\DD_m}\stackrel{\cong}{\longrightarrow}
A_*^{\DD_1}{\otimes}{\cdots}{\otimes}A_*^{\DD_m}{\otimes}
B_*^{\DD_1}{\otimes}{\cdots}{\otimes}B_*^{\DD_m}\vspace{-2mm}$$
be defined as follows.
For $a_i\in A_*^{\DD_i}$, $b_i\in B_*^{\DD_i}$,
$$\phi((a_1{\otimes}b_1){\otimes}{\cdots}{\otimes}(a_m{\otimes}b_m))=
(-1)^s(a_1{\otimes}{\cdots}{\otimes}a_m){\otimes}(b_1{\otimes}{\cdots}{\otimes}b_m),\vspace{-2mm}$$
where $s=\Sigma_{i=2}^m(|b_1|{+}{\cdots}{+}|b_{i-1}|)|a_i|$. Then restriction of $\phi$
on the subgroup $(A_*^{\DD_1}\widehat\otimes B_*^{\DD_1})\otimes\cdots\otimes(A_*^{\DD_m}\widehat\otimes B_*^{\DD_m})$
is just the isomorphism of the theorem.
\hfill$\Box$\vspace{3mm}

{\bf Conventions} For a  $(\DD{\times}\Lambda)$-group $A_*^\DD=\oplus_{s\in\DD,\alpha\in\Lambda}A_*^{s;\alpha}$,
the local group $A_*^s$ satisfies certain property implies $A_*^{s,\alpha}$ satisfies the property for all $\alpha\in\Lambda$.
For example, $A_*^s=0$ implies $A_*^{s;\alpha}=0$ for all $\alpha\in\Lambda$.
$A_*^s\neq 0$ implies $A_*^{s;\alpha}\neq 0$ for some $\alpha\in\Lambda$.
\vspace{3mm}

{\bf Definition~2.7} For a $(\DD{\times}\Lambda)$-group $A_*^\DD=\oplus_{s\in\DD}\,A_*^s$ and a subset $\TT$ of $\DD$,
the {\it restriction group} of $A_*^\DD$ on $\TT$ is the $(\TT{\times}\Lambda)$-group $A_*^\TT=\oplus_{t\in\TT}\,A_*^t$.

The {\it support index set} of $A_*^\DD$ is the set $\SS=\{s\in\DD\,|\, A_*^s\neq 0\}$.
The {\it support group} $A_*^\SS$ of $A_*^\DD$ is its restriction group on the support index set $\SS$.

We have analogue definitions for indexed (co)chain complexes.
\vspace{3mm}

{\bf Theorem~2.8} {\it Suppose $A_{i\,*}^{\SS_i}$ is the support group of the $(\DD{\times}\Lambda_i)$-group $A_{i\,*}^{\DD}$ for $i=1,2$.
Then the support group  of $A_{1\,*}^\DD\widehat\otimes A_{2\,*}^\DD$
is the restriction group  $A_{1\,*}^\SS{\widehat\otimes} A_{2\,*}^\SS$ such that $\SS=\SS_1{\cap}\SS_2$.
For any set $\TT$ such that $\SS\subset\TT\subset\DD$, we have $(\Lambda_1{\times}\Lambda_2)$-group isomorphisms (index $\DD,\TT$ neglected)
$$A_{1\,*}^\DD{\widehat\otimes}A_{2\,*}^\DD=A_{1\,*}^\TT{\widehat\otimes}A_{2\,*}^\TT.$$

Analogue conclusions hold for indexed (co)chain complexes.
\vspace{2mm}

Proof}\, By definition.
\hfill$\Box$\vspace{3mm}

Note that a restriction group (complex) is always a subgroup (subcomplex), but the converse is not true.
The diagonal tensor product group (complex) $A_*^\DD\widehat\otimes B_*^\DD$ is always a restriction group
(complex) of $A_*^\DD{\otimes}B_*^\DD$.
\vspace{3mm}

\section{(Co)homology Group of Polyhedral Product Chain Complexes}\vspace{3mm}

\hspace*{5.5mm}{\bf Conventions and Notations.} All definitions and proofs in this paper have their dual analogues.
We omit the dual definitions and proofs when necessary.

$\Bbb Z(x_1,{\cdots},x_n)$ is the free group generated by $x_1,{\cdots},x_n$.

A simplicial complex $K$ on $[m]$ means
the vertex set of $K$ is a subset of $[m]=\{1,{\cdots},m\}$ and so ghost vertex $\{i\}\notin K$ is allowed.
The void complex $\{\,\}$ (not the empty complex $\{\emptyset\}$) is inevitable in this paper.
For $m<n$, a simplicial complex $K$ on $[m]$ is also a simplicial complex on $[n]$.
We always regard them as two different simplicial complexes.

For a finite set $S$,
$\Delta\!^S$ is the full simplicial complex of all subsets of $S$ (including $\emptyset$).
Define $\partial\Delta\!^S=\Delta\!^S{\setminus}\{S\}$. Specifically, for $S=[1]$, $\partial\Delta\!^{[1]}=\{\emptyset\}$.
For $S=\emptyset$, $\Delta\!^{\emptyset}=\{\emptyset\}$ and so $\partial\Delta\!^\emptyset=\{\,\}$.
We use $|S|$ to denote the cardinality of $S$, i.e., the number of elements of $S$.
\vspace{3mm}

{\bf Definition~3.1} Let $K$ be a simplicial complex on $[m]$
and $(\underline{D_*},\underline{C_*})=\{(D_{k\,*},C_{k\,*})\}_{k=1}^m$ be a sequence of indexed chain complex pairs,
i.e., each $(C_{k\,*},d)$ is a subcomplex of the $\Lambda_k$-complex $(D_{k\,*},d)$.
The {\it polyhedral product chain complex} ${\cal Z}(K;\underline{D_*},\underline{C_*})$
indexed by $\underline{\Lambda}=\Lambda_1{\times}{\cdots}{\times}\Lambda_m$ is the subcomplex of
$(D_{1\,*}{\otimes}{\cdots}{\otimes}D_{m\,*},d)$ defined as follows.
For a subset $\tau$ of $[m]$, define
$$(E_*(\tau),d)=(E_1{\otimes}\cdots{\otimes}E_m,d),\quad (E_k,d)=\left\{\begin{array}{cl}
(D_{k\,*},d)&{\rm if}\,\,k\in \tau,\vspace{1mm}\\
(C_{k\,*},d)&{\rm if}\,\,k\not\in \tau.
\end{array}
\right.$$
Then $({\cal Z}(K;\underline{D_*},\underline{C_*}),d)=(+_{\tau\in K}\,E_*(\tau),d)$.
Define $({\cal Z}(\{\,\};\underline{D_*},\underline{C_*}),d)=0$.
\vspace{3mm}

The (co)homology group of ${\cal Z}(K;\underline{D_*},\underline{C_*})$ can be computed when
each pair $(D_{k\,*},C_{k\,*})$ is split as in the following definition.
\vspace{3mm}

{\bf Definition~3.2} A $\Lambda$-group homomorphism $\theta\colon A_*\to B_*$ is called {\it split}
if ${\rm ker}\,\theta$, ${\rm coker}\,\theta$ and ${\rm coim}\,\theta$ are all free groups.
Specifically, if $A_*$ and $B_*$ are vector spaces over a field, then $\theta$ is always split .

A chain $\Lambda$-complex pair $(D_*,C_*)$ is called {\it homology split} if it satisfies the following conditions.

(1) $(C_*,d)$ is a chain subcomplex of the free complex $(D_*,d)$.

(2) The homology $\Lambda$-group homomorphism $\theta\colon H_*(C_*)\to H_*(D_*)$ induced by inclusion is split.

A polyhedral product chain complex ${\cal Z}(K;\underline{D_*},\underline{C_*})$ is called {\it homology split}
if each pair $(D_{k\,*},C_{k\,*})$ is homology split.
\vspace{3mm}

{\bf Definition~3.3} Let $(D_*,C_*)$ be a homology split pair with
$$\theta\colon H_*(C_*)\to H_*(D_*)\quad{\rm and}\quad \theta^\circ\colon H^*(D_*)\to H^*(C_*)$$
the homology and cohomology $\Lambda$-group homomorphism induced by inclusion.
The index set $\XX$ is $\{\,{\mathpzc n},{\mathpzc e},{\mathpzc i}\,\}$.

The {\it indexed homology group} $H_*^{\XX}\!(D_*,C_*)$ and
{\it the indexed cohomology group} $H^{\,*}_{\!\XX}(D_*,C_*)$ of $(D_*,C_*)$ are $\XX\!$-groups with local groups given by
$$H_*^{\XX}\!(D_*,C_*)=\oplus_{s\in\XX}\,H_*^{s}(D_*,C_*),\quad
H_*^s(D_*,C_*)=\left\{\begin{array}{ll}
{\rm coker}\,\theta&{\rm if}\,\,s={\mathpzc e},\\
{\rm ker}\,\theta&{\rm if}\,\,s={\mathpzc n},\\
{\rm coim}\,\theta&{\rm if}\,\,s={\mathpzc i},
\end{array}\right.$$
$$H^{\,*}_{\!\XX}(D_*,C_*)=\oplus_{s\in\XX}\,H^*_{s}(D_*,C_*),\quad
H^*_s(D_*,C_*)=\left\{\begin{array}{ll}
{\rm ker}\,\theta^\circ&{\rm if}\,\,s={\mathpzc e},\\
{\rm coker}\,\theta^\circ&{\rm if}\,\,s={\mathpzc n},\\
{\rm im}\,\theta^\circ&{\rm if}\,\,s={\mathpzc i}.
\end{array}\right.$$

The {\it character chain complex} $(C_*^{\XX}\!(D_*,C_*),d)$ of $(D_*,C_*)$ is a chain $\XX\!$-complex with local complexes given by
$$C_*^\XX(D_*,C_*)=\oplus_{s\in\XX}\,C_*^s(D_*,C_*),\quad
C_*^s(D_*,C_*)=\left\{\begin{array}{ll}
{\rm coker}\,\theta&{\rm if}\,\,s={\mathpzc e},\\
{\rm ker}\,\theta{\oplus}\Sigma{\rm ker}\,\theta&{\rm if}\,\,s={\mathpzc n},\\
{\rm coim}\,\theta&{\rm if}\,\,s={\mathpzc i},
\end{array}\right.\quad\quad$$
where $d$ is trivial on $C_*^{{\mathpzc i}}(D_*,C_*)$ and $C_*^{{\mathpzc e}}(D_*,C_*)$
and is the desuspension isomorphism on $C_*^{{\mathpzc n}}(D_*,C_*)$.
Precisely, for $x\in{\rm ker}\,\theta$, we always denote by $\overline x$ the unique element in $\Sigma{\rm ker}\,\theta$
such that $d\,\overline x=x$.

The {\it support index set} of $(D_*,C_*)$ is $\SS=\{{\mathpzc s}\in\XX\,|\,C_*^{\mathpzc s}(D_*,C_*)\neq 0\}$.
\vspace{3mm}

Note that when $\theta$ is an epimorphism, we have $H_*^\XX(D_*,C_*)=H_*(C_*)$ and $H^{\,*}_{\!\XX}(D_*,C_*)=H^*(C_*)$
by forgetting the $\XX$-group structure.
\vspace{3mm}

{\bf Remarks} We use the symbols ${\mathpzc n}$, ${\mathpzc e}$, ${\mathpzc i}$ to denote the elements of $\XX$
for the following reasons.
We have $H_*(C_*^{\mathpzc n}(D_*,C_*))=0$ and $H^*(C^{\,*}_{\!\mathpzc n}(D_*,C_*))=0$.
Since both the groups are null, we use the symbol ${\mathpzc n}$ to denote the index.
We have $H_*(C_*^{\mathpzc e}(D_*,C_*))={\rm coker}\,\theta$ and $H^*(C^{\,*}_{\!\mathpzc e}(D_*,C_*))={\rm ker}\,\theta^\circ$.
Since both ``coker" and ``ker" have the vowel letter ``e", we use the symbol ${\mathpzc e}$ to denote the index.
We have $H_*(C_*^{\mathpzc i}(D_*,C_*))={\rm coim}\,\theta$ and $H^*(C^{\,*}_{\!\mathpzc i}(D_*,C_*))={\rm im}\,\theta^\circ$.
Since both ``coim" and ``im" have the vowel letter ``i", we use the symbol ${\mathpzc i}$ to denote the index.
\vspace{3mm}

{\bf Definition~3.4} The {\it atom chain complex} $(T_*^\XX\!,d)$ is a chain $\XX\!$-complex with local complexes given by
$$(T_*^{{\mathpzc i}},d)=\Bbb Z({\mathpzc i}),\quad(T_*^{{\mathpzc n}},d)=(\Bbb Z({\mathpzc n},\overline{\mathpzc n}),d),\quad
(T_*^{{\mathpzc e}},d)=\Bbb Z({\mathpzc e}),$$
where $|{\mathpzc n}|=|{\mathpzc e}|=|{\mathpzc i}|=0$, $|\overline{\mathpzc n}|=1$, $d\,\overline{\mathpzc n}={\mathpzc n}$.

The {\it atom cochain complex} $(T^{\,*}_{\!\XX},\delta)$ is the dual of $(T_*^\XX\!,d)$ given by
$$(T^*_{{\mathpzc i}},\delta)=\Bbb Z({\mathpzc i}),\quad(T^*_{{\mathpzc n}},\delta)=(\Bbb Z({\mathpzc n},\overline{\mathpzc n}),\delta),\quad
(T^*_{{\mathpzc e}},\delta)=\Bbb Z({\mathpzc e}),\vspace{-2mm}$$
where we use the same symbol to denote a generator and its dual generator. So $\delta{\mathpzc n}=\overline{\mathpzc n}$.

The chain $\XX\!$-complex $(S_*^\XX,d)$ is the subcomplex of $(T_*^\XX\!,d)$ such that
$S_*^\XX=\Bbb Z({\mathpzc n},{\mathpzc i})$.
\vspace{3mm}

{\bf Theorem~3.5}\, {\it For a homology split pair $(D_*,C_*)$,
there is a quotient chain homotopy equivalence $q$
satisfying the following commutative diagram
$$\begin{array}{ccc}
(C_*,d)&\stackrel{q'}{\longrightarrow}&H_*(C_*)\,\vspace{1mm}\\
\cap&&\cap\\
(D_*,d)&\stackrel{q}{\longrightarrow}&(C_*^{\XX}\!(D_*,C_*),d),
\end{array}
$$
where the homotopy equivalence $q'$ is the restriction of $q$ and the index $\XX$ is neglected.
There are also chain $(\XX{\times}\Lambda)$-complex isomorphisms $\phi$ and its restriction $\phi'$ satisfying the following commutative diagram
$$\begin{array}{ccc}
H_*(C_*)&\stackrel{\phi'}{\cong}&S_*^\XX\,\widehat\otimes\, H_*^{\XX}\!(D_*,C_*)\vspace{1mm}\\
\cap&&\cap\\
(C_*^{\XX}\!(D_*,C_*),d)&\stackrel{\phi}{\cong}&(T_*^\XX\,\widehat\otimes\, H_*^{\XX}\!(D_*,C_*),d).
\end{array}$$

Let $\SS$ be the support index set of $(D_*,C_*)$.
Then for any index set $\TT$ such that $\SS\subset\TT\subset\XX$,
there are $\Lambda$-group isomorphisms ($\XX$ and $\TT$ neglected) satisfying the following commutative diagram
$$\begin{array}{ccc}
S_*^\XX\,\widehat\otimes\, H_*^{\XX}\!(D_*,C_*)&=&S_*^\TT\,\widehat\otimes\, H_*^{\TT}\!(D_*,C_*)\vspace{1mm}\\
\cap&&\cap\\
(T_*^\XX\,\widehat\otimes\, H_*^{\XX}\!(D_*,C_*),d)&=&(T_*^\TT\,\widehat\otimes\, H_*^{\TT}\!(D_*,C_*),d),
\end{array}$$
where $(-)^\TT$ means the restriction group (complex) of $(-)^\XX$ on $\TT$.
\vspace{2mm}

Proof}\, Take a representative $n_{i}$ in $C_*$ for every generator of ${\rm ker}\,\theta$
and let $\overline n_{i}\in D_*$ be any element such that $d\,\overline n_{i}= n_{i}$.
Take a representative $i_{j}$ in $C_*$ for every generator of ${\rm coim}\,\theta$.
Take a representative $e_k$ in $D_*$ for every generator of ${\rm coker}\,\theta$.
So we may regard $H_*(C_*)$ as the chain subcomplex of $C_*$ freely generated by all $n_{i}$'s and $i_{j}$'s
and regard $(C_*^{\XX\!}(D_*,C_*),d)$ as the chain subcomplex of $D_*$ freely generated by all
$n_{i}$'s, $\overline n_{i}$'s, $i_{j}$'s and $e_k$'s.
Then we have the following commutative diagram of short exact sequences of chain complexes
$$\begin{array}{ccccccc}
 0\to& H_*(C_*)&\stackrel{f'}{\longrightarrow}& C_*&\stackrel{g'}{\longrightarrow}&C_*/H_*(C_*)&\to0\vspace{1mm}\\
 &\cap&&\cap&&\downarrow&\\
0\to& C_*^{\XX\!}(D_*,C_*)&\stackrel{f}{\longrightarrow}& D_*&\stackrel{g}{\longrightarrow}&D_*/C_*^{\XX\!}(D_*,C_*)&\to0.
  \end{array}
$$
Since all the groups are free, the inclusions $f,f'$ have inverse group homomorphism.
From $H_*(C_*/H_*(C_*))=0$ and $H_*(D_*/C_*^{\XX\!}(D_*,C_*))=0$ we have that the inverse homomorphisms of $f,f'$ are chain homomorphisms.
So we may take $q,q'$ to be the inverse homomorphisms of $f,f'$.

$\phi$ is defined as shown in the following table.
\begin{center}
\begin{tabular}{|c|c|c|c|c|}
\hline
{\rule[-2mm]{0mm}{6mm}$\quad x\in\quad$}&${\rm coker}\,\theta$&$\Sigma\,{\rm ker}\,\theta$
&${\rm ker}\,\theta$&${\rm coim}\,\theta$\\
\hline
{\rule[-2mm]{0mm}{7mm}$\phi(x)=$}&${\mathpzc e}\otimes x$&$\overline{\mathpzc n}\otimes dx$
&${\mathpzc n}\otimes x$&${\mathpzc i}\otimes x$\\
\hline
\end{tabular}
\end{center}

The equality of the theorem is by Theorem~2.8.
\hfill$\Box$\vspace{3mm}

{\bf Definition~3.6} The index set $\XX_m$ is defined as follows.
$$\XX_m=\{(\sigma,\omega)\,|\,\sigma,\omega\subset[m],\,\sigma{\cap}\omega=\emptyset\}.$$
$\XX_m$ is always identified with $\XX\!{\times}{\cdots}{\times}\XX$ ($m$-fold) by the 1-1 correspondence
$$(\sigma\!,\,\omega)\in\XX_m\to(s_1,{\cdots},s_m)\in\XX\!{\times}{\cdots}{\times}\XX\!,$$
where $s_k={\mathpzc e}$ if $k\in\sigma$, $s_k={\mathpzc n}$ if $k\in\omega$ and $s_k={\mathpzc i}$ otherwise.

The index set $\RR$ is the subset $\{{\mathpzc i},{\mathpzc n}\}$ of $\XX$.
The index set $\RR_m$ is the subset $\{(\emptyset,\omega)\,|\,\omega\subset[m]\}$ of $\XX_m$.
$\RR_m$ is always identified with $\RR{\times}{\cdots}{\times}\RR$ ($m$-fold) by the 1-1 correspondence
$(\emptyset,\,\omega)\in\RR_m\to(s_1,{\cdots},s_m)\in\RR{\times}{\cdots}{\times}\RR\!$,
where $s_k={\mathpzc n}$ if $k\in\omega$ and $s_k={\mathpzc i}$ otherwise.
\vspace{3mm}

Let $(A_k)_*^\XX$ be a $(\XX{\times}\Lambda_k)$-group for $k=1,{\cdots},m$ and $\underline{\Lambda}=\Lambda_1{\times}{\cdots}{\times}\Lambda_m$.
Then the group $(A_1)_*^\XX{\otimes}{\cdots}{\otimes}(A_m)_*^\XX$
indexed by $(\XX\!{\times}\Lambda_1){\times}{\cdots}{\times}(\XX\!{\times}\Lambda_m)$ is also a group indexed by $\XX_m{\times}\underline{\Lambda}$.
So we may write $(\underline{A})_*^{\XX_m}=\oplus_{(\sigma\!,\,\omega)\in\XX_m}(\underline{A})_*^{\sigma\!,\,\omega}$
for $(A_1)_*^\XX{\otimes}{\cdots}{\otimes}(A_m)_*^\XX$ with local groups
$$(\underline{A})_*^{\sigma\!,\,\omega}=H_1{\otimes}{\cdots}{\otimes}H_m,\quad
H_k=\left\{\begin{array}{ll}
(A_k)_*^{{\mathpzc e}}&{\rm if}\,k\in\sigma,\\
(A_k)_*^{{\mathpzc n}}&{\rm if}\,k\in\omega,\\
(A_k)_*^{{\mathpzc i}}&{\rm othetwise}.\\
\end{array}\right.$$
Moreover, suppose $(A_k)_*^\XX=(A_k)_*^{\XX;\Lambda_k}=\oplus_{s\in\XX\!,\,\alpha\in\Lambda_k}(A_k)_*^{s,\alpha}$.
Then
$$(\underline{A})_*^{\XX_m}=(\underline{A})_*^{\XX_m;\underline{\Lambda}}=
\oplus_{(\sigma\!,\,\omega)\in\XX_m,\,\alpha_k\in\Lambda_k}\,(\underline{A})_*^{\sigma\!,\,\omega;\alpha_1,\cdots,\alpha_m}$$
with local groups
$$(\underline{A})_*^{\sigma\!,\,\omega;\alpha_1,{\cdots},\alpha_m}=H_1^{\alpha_1}{\otimes}{\cdots}{\otimes}H_m^{\alpha_m},\quad
H_k^{\alpha_k}=\left\{\begin{array}{ll}
(A_k)_*^{{\mathpzc e};\alpha_k}&{\rm if}\,k\in\sigma,\\
(A_k)_*^{{\mathpzc n};\alpha_k}&{\rm if}\,k\in\omega,\\
(A_k)_*^{{\mathpzc i};\,\alpha_k}&{\rm otherwise}.\\
\end{array}\right.$$

Analogue conventions hold for $\RR_m$ by replacing $\XX$ and $\sigma$ in the above formulas respectively by $\RR$ and $\emptyset$.
\vspace{3mm}

{\bf Definition~3.7} The {\it total chain  complex} $(T_*^{\XX_m},d)$ is the chain $\XX_m$-complex
$(T_*^\XX{\otimes}{\cdots}{\otimes}T_*^\XX\!,d)$ ($m$-fold).
Dually, the {\it total cochain  complex} $(T^{\,*}_{\!\XX_m},\delta)$ is the cochain $\XX_m$-complex
$(T^{\,*}_{\!\XX}{\otimes}{\cdots}{\otimes}T^{\,*}_{\!\XX},\delta)$

Let $K$ be a simplicial complex on $[m]$.

The {\it total chain complex} $(T_*^{\XX_m}(K),d)$ of $K$ is the chain subcomplex of
$(T_*^{\XX_m},d)$ defined as follows.
For a subset $\tau$ of $[m]$, define
$$(U_*(\tau),d)=(U_1{\otimes}\cdots{\otimes}U_m,d),\quad (U_k,d)=\left\{\begin{array}{cl}
(T_*^\XX\!,d)&{\rm if}\,\,k\in \tau,\vspace{1mm}\\
S_*^\XX&{\rm if}\,\,k\not\in \tau.
\end{array}
\right.$$
Then $(T_*^{\XX_m}(K),d)=(+_{\tau\in K}\,U_*(\tau),d)$.
Define $(T_*^{\XX_m}(\{\,\}),d)=0$.

The {\it total homology group} of $K$ is $H_*^{\XX_m}(K)=H_*(T_*^{\XX_m}(K))$.

The {\it total cochain complex} $(T^{\,*}_{\!\XX_m}(K),\delta)$ of $K$ is the dual of $(T_*^{\XX_m}(K),d)$.
Since $T_*^{\XX_m}/T_*^{\XX_m}(K)$ is free, by Lemma~2.2, $(T^{\,*}_{\!\XX_m}(K),\delta)$ is a quotient complex of $(T^{\,*}_{\!\XX_m},\delta)$.

The {\it total cohomology group} of $K$ is $H^{\,*}_{\!\XX_m}(K)=H^*(T^{\,*}_{\!\XX_m}(K))$.

For an index set $\DD\subset\XX_m$, the {\it total object} $(T_*^{\DD}(K),d)$, $H_*^{\DD}(K)$, $(T^{\,*}_{\!\DD}(K),d)$, $H^{\,*}_{\!\DD}(K)$
of $K$ on $\DD$ is the restriction group (complex) on $\DD$ of the corresponding total object.
Specifically, if $\DD=\RR_m$, we have the {\it right total object}
$(T_*^{\RR_m}(K),d)$, $H_*^{\RR_m}(K)$, $(T^{\,*}_{\!\RR_m}(K),d)$, $H^{\,*}_{\!\RR_m}(K)$ of $K$.
\vspace{3mm}

{\bf Theorem~3.8} {\it Denote the local groups (complexes) of the total objects of $K$ on $\DD$ as follows.
$$(T_*^{\DD}(K),d)=\oplus_{(\sigma,\omega)\in\DD}\,(T_*^{\sigma\!,\,\omega}(K),d),\quad
H_*^{\DD}(K)=\oplus_{(\sigma,\omega)\in\DD}\,H_*^{\sigma\!,\,\omega}(K),\vspace{-2mm}$$
$$(T^{\,*}_{\!\DD}(K),\delta)=\oplus_{(\sigma,\omega)\in\DD}\,(T^*_{\sigma\!,\,\omega}(K),\delta),\quad
H^{\,*}_{\!\DD}(K)=\oplus_{(\sigma,\omega)\in\DD}\,H^*_{\sigma\!,\,\omega}(K).$$
Then for $(\sigma\!,\,\omega)\in\DD$,\vspace{-2mm}
$$(T_*^{\sigma\!,\,\omega}(K),d)\cong(\Sigma\w C_*(K_{\sigma\!,\,\omega}),d),\quad
H_*^{\sigma\!,\,\omega}(K)\cong\w H_{*-1}(K_{\sigma\!,\,\omega}),\vspace{-2mm}$$
$$(T^*_{\sigma\!,\,\omega}(K),\delta)\cong(\Sigma\w C^*(K_{\sigma\!,\,\omega}),\delta),\quad
H^*_{\sigma\!,\,\omega}(K)\cong\w H^{*-1}(K_{\sigma\!,\,\omega}),$$
where $\Sigma\w C$ means the suspension augmented simplicial (co)chain complex,
$K_{\sigma\!,\,\omega}=({\rm link}_{_K}\sigma)|_\omega=\{\tau\,\,|\,\,\tau\subset\omega,\,\tau{\cup}\sigma\in K\}$
if $\sigma\in K$ and $K_{\sigma\!,\,\omega}=\{\,\}$ if $\sigma\notin K$.
\vspace{2mm}

Proof}\, Write $t_{E,\overline N,N,I}$ for the generator $t_1{\otimes}{\cdots}{\otimes}t_m$ of $T_*^{\XX_m}$
such that
$$E=\{k\,|\,t_k\!=\!{\mathpzc e}\},\,\,\overline N=\{k\,|\,t_k\!=\!\overline{\mathpzc n}\},\,\,N=\{k\,|\,t_k\!=\!{\mathpzc n}\},\,\,I=\{k\,|\,t_k\!=\!{\mathpzc i}\}$$
Then $T_*^{\XX_m}=\oplus_{(\sigma,\omega)\in\XX_m}\,T_*^{\sigma,\omega}$ with
$T_*^{\sigma,\omega}=\Bbb Z(\{t_{E,\overline N,N,I}\}_{E=\sigma,\,\overline N\cup N=\omega})$.
For the $U(\tau)$ in Definition~3.7,
we have $U(\tau)=\Bbb Z(\{t_{E,\overline N,N,I}\}_{E\cup\overline N\subset\,\tau})$. So
$$T_*^{\sigma\!,\,\omega}(K)=+_{\tau\in K}\,T_*^{\sigma\!,\,\omega}\cap U_*(\tau)
=\Bbb Z(\{t_{\sigma,\overline N,N,I}\}
_{\sigma\cup \overline N\,\in\, K,\,\,\overline N\cup N\,=\,\omega}).$$

For $\sigma\in K$, the 1-1 correspondence $\epsilon(\overline N)=t_{\sigma,\,\overline N,\,\omega{\setminus}\overline N,\,[m]{\setminus}(\sigma{\cup}\omega)}$
for all $\overline N\in K_{\sigma\!,\,\omega}$ induces an isomorphism
$\epsilon\colon(\Sigma\w C_*(K_{\sigma\!,\,\omega}),d)\to(T_*^{\sigma\!,\,\omega}(K),d)$.
For $\sigma\notin K$,  $T_*^{\sigma,\omega}(K)=\Sigma\w C_*(\{\,\})=0$.
\hfill$\Box$\vspace{3mm}

{\bf Definition~3.9} Let $K$ a simplicial complex on $[m]$, $(\sigma,\omega)\in\XX_m$ and everything else be as in Theorem~3.8.

The simplicial complex  $K_{\sigma,\omega}$ is called a {\it local complex} of $K$ with respect to $[m]$.
The (co)chain complexes $(T_*^{\sigma,\omega}(K),d),(T^*_{\sigma,\omega}(K),\delta)$ are called
a {\it local (co)chain complex} of $K$ with respect to $[m]$.
The (co)homology groups  $H_*^{\sigma,\omega}(K),H^*_{\sigma,\omega}(K)$ are called a {\it local (co)homology group} of $K$
with respect to $[m]$.
\vspace{3mm}

{\bf Example~3.10} We compute $H_*^{\XX_m}(K)$ for $K=\{\,\},\{\emptyset\},\Delta\!^{[m]},\Delta\!^S,\partial\Delta\!^S$
regarded as simplicial complexes on $[m]$, where $S\neq\emptyset$ and $S\subsetneq[m]$.

$H_*^{\sigma,\omega}(\{\,\})=0$ for all $(\sigma,\omega)\in\XX_m$.

$H_*^{\emptyset,\omega}(\{\emptyset\})=\Bbb Z$ for all $\omega\subset[m]$. $H_*^{\sigma\!,\,\omega}(\{\emptyset\})=0$ if $\sigma\neq\emptyset$.

$H_*^{\sigma,\emptyset}(\Delta\!^{[m]})=\Bbb Z$ for all $\sigma\subset[m]$. $H_*^{\sigma\!,\,\omega}(\{\emptyset\})=0$ if $\omega\neq\emptyset$.

$H_*^{\sigma\!,\,\omega}(\Delta\!^S)=\Bbb Z$ if $\sigma\subset S$, $\omega{\cap}S=\emptyset$.
$H_*^{\sigma\!,\,\omega}(\Delta\!^S)=0$ otherwise.

$H_*^{\sigma\!,\,\omega}(\partial\Delta\!^S)=\Bbb Z$ if $\sigma\subset S$, $\sigma\neq S$, $\omega{\cap}S=\emptyset$.
$H_*^{\sigma\!,\,\omega}(\partial\Delta\!^S)=\Sigma^{|S|-|\sigma|-1}\Bbb Z$ if $\sigma\subset S$, $\sigma\neq S$, $S{\setminus}\sigma\subset \omega$.
$H_*^{\sigma\!,\,\omega}(\partial\Delta\!^S)=0$ otherwise.
\vspace{3mm}

{\bf Definition~3.11} Let $(\underline{D_*},\underline{C_*})=\{(D_{k\,*},C_{k\,*})\}_{k=1}^m$ be
a sequence such that each $(D_{k\,*},C_{k\,*})$ is a homology split pair by Definition~3.2 with
$$\theta_k\colon H_*(C_{k\,*})\to H_*(D_{k\,*})\,\,{\rm and}\,\,\theta_k^\circ\colon H^*(D_{k\,*})\to H^*(C_{k\,*})$$
the (co)homology $\Lambda_k$-group homomorphism induced by inclusion.
Denote by $\underline{\Lambda}=\Lambda_1{\times}{\cdots}{\times}\Lambda_m$.

The {\it indexed homology group} $H_*^{\XX_m}(\underline{D_*},\underline{C_*})$ and
the {\it indexed cohomology group} $H^{\,*}_{\!\XX_m}(\underline{D_*},\underline{C_*})$ of $(\underline{D_*},\underline{C_*})$
are the $(\XX_m{\times}\underline{\Lambda})$-groups given by
$$H_*^{\XX_m}(\underline{D_*},\underline{C_*})=H_*^{\XX}\!(D_{1\,*},C_{1\,*}){\otimes}{\cdots}{\otimes}H_*^{\XX}\!(D_{m\,*},C_{m\,*}),$$
$$H^{\,*}_{\!\XX_m}(\underline{D_*},\underline{C_*})=H^{\,*}_{\!\XX}(D_{1\,*},C_{1\,*}){\otimes}{\cdots}{\otimes}H^{\,*}_{\!\XX}(D_{m\,*},C_{m\,*}).$$
For an index set $\DD\subset\XX_m$, the {\it indexed homology group} $H_*^\DD(\underline{D_*},\underline{C_*})$ and
{\it indexed cohomology group} $H^{\,*}_{\!\DD}(\underline{D_*},\underline{C_*})$ of $(\underline{D_*},\underline{C_*})$ on $\DD$
are respectively the restriction  group of $H_*^{\XX_m}(\underline{D_*},\underline{C_*})$ and $H^{\,*}_{\!\XX_m}(\underline{D_*},\underline{C_*})$
on $\DD{\times}\underline{\Lambda}$.

The {\it character chain complex} $(C_*^{\XX_m}(\underline{D_*},\underline{C_*}),d)$ of $(\underline{D_*},\underline{C_*})$
is the chain $(\XX_m{\times}\underline{\Lambda})$-complex given by
$$(C_*^{\XX_m}(\underline{D_*},\underline{C_*}),d)=(C_*^{\XX}\!(D_{1\,*},C_{1\,*}){\otimes}{\cdots}{\otimes}C_*^{\XX}\!(D_{m\,*},C_{m\,*}),d).$$

Let $K$ be a simplicial complex on $[m]$.

The {\it polyhedral product character chain complex} $({\cal Z}^{\XX_m}(K;\underline{D_*},\underline{C_*}),d)$ is the subcomplex of
$(C_*^{\XX_m}(\underline{D_*},\underline{C_*}),d)$ defined as follows.
For a subset $\tau$ of $[m]$, define
$$(P_*(\tau),d)=(P_1{\otimes}\cdots{\otimes}P_m,d),\quad (P_k,d)=\left\{\begin{array}{cl}
(C_*^{\XX}\!(D_{k\,*},C_{k\,*}),d)&{\rm if}\,\,k\in \tau,\vspace{1mm}\\
H_*(C_{k\,*})&{\rm if}\,\,k\not\in \tau.
\end{array}\right.$$
Then ${\cal Z}^{\XX_m}(K;\underline{D_*},\underline{C_*})=+_{\tau\in K}\,P_*(\tau)$.
Define ${\cal Z}^{\XX_m}(\{\,\};\underline{D_*},\underline{C_*})=0$.
\vspace{3mm}

{\bf Theorem~3.12}\, {\it For a homology split ${\cal Z}(K;\underline{D_*},\underline{C_*})$ by Definition~3.2,
there is a quotient chain homotopy equivalence (index $\XX_m$ neglected)
$$q_{(K;\underline{D},\underline{C})}\colon({\cal Z}(K;\underline{D_*},\underline{C_*}),d)
\stackrel{\simeq}{\longrightarrow}
({\cal Z}^{\XX_m}(K;\underline{D_*},\underline{C_*}),d)$$
and a chain $(\XX_m{\times}\underline{\Lambda})$-complex isomorphism
$$\phi_{(K;\underline{D},\underline{C})}\colon({\cal Z}^{\XX_m}(K;\underline{D_*},\underline{C_*}),d)
\stackrel{\cong}{\longrightarrow}
(T_*^{\XX_m}(K)\,\widehat\otimes\,H_*^{\XX_m}(\underline{D_*},\underline{C_*}),d).$$

Let $\underline{\SS}=\SS_1{\times}{\cdots}{\times}\SS_m$, where each $\SS_k$ is the support index set of $(D_{k\,*},C_{k\,*})$.
Then for any index set $\DD$ such that $\underline{\SS}\subset\DD\subset\XX_m$, we have $\underline{\Lambda}$-group isomorphisms
(index $\XX_m,\DD$ neglected)
$$H_*({\cal Z}(K;\underline{D_*},\underline{C_*}))\cong H_*^{\XX_m}(K)\,\widehat\otimes\, H_*^{\XX_m}(\underline{D_*},\underline{C_*})
=H_*^{\DD}(K)\,\widehat\otimes\, H_*^{\DD}(\underline{D_*},\underline{C_*}),$$
$$H^*({\cal Z}(K;\underline{D_*},\underline{C_*}))\cong H^{\,*}_{\XX_m}(K)\,\widehat\otimes\, H^{\,*}_{\XX_m}(\underline{D_*},\underline{C_*})
=H^{\,*}_{\DD}(K)\,\widehat\otimes\, H^{\,*}_{\!\DD}(\underline{D_*},\underline{C_*}).$$

Precisely, suppose
$C_{k\,*}=\oplus_{\alpha_k\in\Lambda_k}C_{k\,*}^{\alpha_k}$,
$D_{k\,*}=\oplus_{\alpha_k\in\Lambda_k}D_{k\,*}^{\alpha_k}$
and
$$\theta_k=\oplus_{\alpha_k\in\Lambda_k}\,\theta_{k,\alpha_k},\quad
\theta_{k,\alpha_k}\colon H_*(C_{k\,*}^{\alpha_k})\to H_*(D_{k\,*}^{\alpha_k}),$$
$$\theta_k^\circ=\oplus_{\alpha_k\in\Lambda_k}\,\theta_{k,\alpha_k}^\circ,\quad
\theta_{k,\alpha_k}^\circ\colon H^*(D_{k\,*}^{\alpha_k})\to H^*(C_{k\,*}^{\alpha_k}).$$
Then we have
$$H_*({\cal Z}(K;\underline{D_*},\underline{C_*}))
\cong\oplus_{(\sigma\!,\,\omega)\in\DD,\,\alpha_k\in\Lambda_k}H_{*}^{\sigma\!,\,\omega}(K)
\otimes H_1^{\alpha_1}{\otimes}{\cdots}{\otimes}H_m^{\alpha_m},
$$
$$H^*({\cal Z}(K;\underline{D_*},\underline{C_*}))
\cong\oplus_{(\sigma\!,\,\omega)\in\DD,\,\alpha_k\in\Lambda_k}H^{*}_{\sigma\!,\,\omega}(K)
\otimes H^1_{\alpha_1}{\otimes}{\cdots}{\otimes}H^m_{\alpha_m},
$$
where
$$H_k^{\alpha_k}=\left\{
\begin{array}{ll}
{\rm coker}\,\theta_{k,\alpha_k}& {\rm if}\,\, k\in\sigma,\vspace{1mm}\\
{\rm ker}\,\theta_{k,\alpha_k}& {\rm if}\,\, k\in\omega,\vspace{1mm}\\
{\rm coim}\,\theta_{k,\alpha_k}& {\rm otherwise},
\end{array}\right. \,\,
H^k_{\alpha_k}=\left\{
\begin{array}{ll}
{\rm ker}\,\theta_{k,\alpha_k}^\circ& {\rm if}\,\, k\in\sigma,\vspace{1mm}\\
{\rm coker}\,\theta_{k,\alpha_k}^\circ& {\rm if}\,\, k\in\omega,\vspace{1mm}\\
{\rm im}\,\theta_{k,\alpha_k}^\circ& {\rm otherwise}.
\end{array}\right.
$$

If each $\theta_k$ is an epimorphism, then we have
$$H_*({\cal Z}(K;\underline{D_*},\underline{C_*}))
\cong\oplus_{\omega\notin K,\,\alpha_k\in\Lambda_k}H_*^{\emptyset,\omega}(K)
\otimes H_1^{\alpha_1}{\otimes}{\cdots}{\otimes}H_m^{\alpha_m},$$
$$H^*({\cal Z}(K;\underline{D_*},\underline{C_*}))
\cong\oplus_{\omega\notin K,\,\alpha_k\in\Lambda_k}H^*_{\emptyset,\omega}(K)
\otimes H^1_{\alpha_1}{\otimes}{\cdots}{\otimes}H^m_{\alpha_m},
\vspace{2mm}$$
$$H_k^{\alpha_k}=\left\{\begin{array}{ll}
{\rm ker}\,\theta_{k,\alpha_k}&{\rm if}\,\,k\in\omega,\vspace{1mm}\\
{\rm coim}\,\theta_{k,\alpha_k}&{\rm otherwise},
\end{array}\right.\,\,
H^k_{\alpha_k}=\left\{\begin{array}{ll}
{\rm coker}\,\theta_{k,\alpha_k}^\circ&{\rm if}\,\,k\in\omega,\vspace{1mm}\\
{\rm im}\,\theta_{k,\alpha_k}^\circ&{\rm otherwise}.
\end{array}\right.$$

If each $H_*(D_{k\,*})=0$, then we have
$$H_*({\cal Z}(K;\underline{D_*},\underline{C_*}))
\cong\w H_{*-1}(K)\otimes H_*(C_{1\,*}){\otimes}{\cdots}{\otimes}H_*(C_{m\,*}),$$
$$H^*({\cal Z}(K;\underline{D_*},\underline{C_*}))
\cong\w H^{*-1}(K)\otimes H^*(C_{1\,*}){\otimes}{\cdots}{\otimes}H^*(C_{m\,*}).$$

Proof}\, Write $q_k,q'_k,\phi_k,\phi'_k$ for the $q,q',\phi,\phi'$ in Theorem~3.5 for $(D_*,C_*)=(D_{k\,*},C_{k\,*})$.

For $\tau\in K$, define $q_\tau=p_1{\otimes}{\cdots}{\otimes}p_m$, where $p_k=q_k$ if $k\in\tau$ and $p_k=q'_k$ if $k\notin K$.
Then each $q_\tau$ is a chain homotpy equivalence.
So $q_{(K;\underline{D_*},\underline{C_*})}=+_{\tau\in K}\,q_\tau$ is also a chain homotopy equivalence.

For $\tau\subset[m]$, define $\phi_\tau=\lambda_1{\otimes}{\cdots}{\otimes}\lambda_m$,
where $\lambda_k=\phi_k$ if $k\in\tau$ and $\lambda_k=\phi'_k$ if $k\notin\tau$.
Then
$$\phi_\tau\colon H_*(\tau)\to
\big(U_1\widehat\otimes H_*^\XX\!(D_{1\,*},C_{1\,*})\big){\otimes}{\cdots}{\otimes}
\big(U_m\widehat\otimes H_*^\XX\!(D_{m\,*},C_{m\,*})\big)$$
is an isomorphism, where $U(\tau)=U_1{\otimes}{\cdots}{\otimes}U_m$ is as in Definition~3.7.
By Theorem~2.6,
$$\big(U_1\widehat\otimes H_*^\XX\!(D_{1\,*},C_{1\,*})\big){\otimes}{\cdots}{\otimes}
\big(U_m\widehat\otimes H_*^\XX\!(D_{m\,*},C_{m\,*})\big)
\cong U_*(\tau)\widehat\otimes H_*^{\XX_m}(\underline{D_*},\underline{C_*}).$$
Identify the above two chain complexes. Then
$\phi_{(K;\underline{D_*},\underline{C_*})}=+_{\tau\in K}\phi_\tau$ is also an isomorphism
from ${\cal Z}_*^{\XX_m}(K;\underline{D_*},\underline{C_*})$ to $T_*^{\XX_m}(K)\widehat\otimes H_*^{\XX_m}(\underline{D_*},\underline{C_*})$.

Since $H_*^{\XX_m}(\underline{D_*},\underline{C_*})$ is free, we have by K\"{u}nneth Theorem
$$H_*\big(T_*^{\XX_m}(K)\,\widehat\otimes\, H_*^{\XX_m}(\underline{D_*},\underline{C_*})\big)
\cong H_*^{\XX_m}(K)\,\widehat\otimes\, H_*^{\XX_m}(\underline{D_*},\underline{C_*}).$$

The equality
$T_*^{\XX_m}(K)\,\widehat\otimes\, H_*^{\XX_m}(\underline{D_*},\underline{C_*})=T_*^{\DD}(K)\,\widehat\otimes\, H_*^{\DD}(\underline{D_*},\underline{C_*})$
holds by Theorem~2.8.

If each $\theta_k$ is an epimorphism, then the support index set $\SS_k$ of $\theta_k$ satisfies $\SS_k\subset\RR$.
So we may take $\DD$ to be $\RR_m$.
By definition, $H_*^{\DD}(K)=H_*^{\RR_m}(K)=\oplus_{\omega\subset[m]}\,H_{*}^{\emptyset,\omega}(K)
=\oplus_{\omega\notin K}\,H_{*}^{\emptyset,\omega}(K)$.

If each $H_*(D_{k\,*})=0$, then the support index set $\SS_k$ of $\theta_k$ is $\{{\mathpzc n}\}$.
So we may take $\DD$ to be the trivial index set $\{(\emptyset,[m])\}$.
By definition, $H_*^{\DD}(K)=H_*^{\emptyset,[m]}(K)=\w H_{*-1}(K)$.
\hfill$\Box$\vspace{3mm}

{\bf Example~3.13} We compute the homology group $H_*({\cal Z}(K;\underline{D_*},\underline{C_*}))$ in Theorem~3.12
for $K=\{\,\},\{\emptyset\},\Delta\!^{[m]},\Delta\!^S,\partial\Delta\!^S$ as in Example~3.10.
Denote the local group $H_*^{\sigma,\omega}(\underline{D_*},\underline{C_*})$ simply by $H_*^{\sigma,\omega}$.
Then $H_*^{\sigma,\omega}=H_1{\otimes}{\cdots}{\otimes}H_m$, where $H_k={\rm coker}\,\theta_k$ if $k\in\sigma$,
$H_k={\rm ker}\,\theta_k$ if $k\in\omega$ and $H_k={\rm coim}\,\theta_k$ otherwise.

By definition, $H_*(\{\,\};\underline{D_*},\underline{C_*}))=0$.

By K\"{u}nneth Theorem, $H_*(\{\emptyset\};\underline{D_*},\underline{C_*}))=H_*(C_{1\,*}){\otimes}{\cdots}{\otimes}H_*(C_{m\,*})$
and $H_*(\Delta\!^{[m]};\underline{D_*},\underline{C_*}))=H_*(D_{1\,*}){\otimes}{\cdots}{\otimes}H_*(D_{m\,*})$.
By Theorem~3.12,
$$\begin{array}{l}
\quad H_*({\cal Z}(\{\emptyset\};\underline{D_*},\underline{C_*}))\vspace{2mm}\\
\cong\oplus_{\omega\subset[m]}\,H_*^{\emptyset,\omega}(\{\emptyset\}){\otimes}H_*^{\emptyset,\omega}\vspace{2mm}\\
\cong\oplus_{\omega\subset[m]}\,H_*^{\emptyset,\omega}\vspace{2mm}\\
\cong H_*(C_{1\,*}){\otimes}{\cdots}{\otimes}H_*(C_{m\,*}),
\end{array}\quad
\begin{array}{l}
\quad H_*({\cal Z}(\Delta\!^{[m]};\underline{D_*},\underline{C_*}))\vspace{2mm}\\
\cong\oplus_{\sigma\subset[m]}\,H_*^{\sigma,\emptyset}(\Delta\!^{[m]}){\otimes}H_*^{\sigma,\emptyset}\vspace{2mm}\\
\cong\oplus_{\sigma\subset[m]}\,H_*^{\sigma,\emptyset}\vspace{2mm}\\
\cong H_*(D_{1\,*}){\otimes}{\cdots}{\otimes}H_*(D_{m\,*}).
\end{array}$$

By K\"{u}nneth Theorem,
$$H_*({\cal Z}(\Delta\!^S;\underline{D_*},\underline{C_*}))=U_1{\otimes}{\cdots}{\otimes}U_m,\quad U_k=\left\{\begin{array}{ll}
H_*(D_{k\,*})&{\rm if}\,k\in S,\\
H_*(C_{k\,*})&{\rm if}\,k\notin S.
\end{array}\right.$$
By Theorem~3.12 ,
$$\begin{array}{l}
\quad H_*({\cal Z}(\Delta\!^S;\underline{D_*},\underline{C_*}))\vspace{2mm}\\
\cong\oplus_{\sigma\subset S,\,\omega{\cap}S=
\emptyset}\,H_*^{\sigma,\omega}(\Delta\!^S){\otimes}H_*^{\sigma,\omega}\vspace{2mm}\\
\cong\oplus_{\sigma\subset S,\,\omega{\cap}S=\emptyset}\,H_*^{\sigma,\omega}\vspace{2mm}\\
\cong U_1{\otimes}{\cdots}{\otimes}U_m.
\end{array}$$

Let $E_{k\,*}=D_{k\,*}/C_{k\,*}$ if $k\in S$ and $E_{k\,*}=C_{k\,*}$ if $k\notin S$. Since we have a short exact sequence of chain complexes
$$0\to {\cal Z}(\partial\Delta\!^S;\underline{D_*},\underline{C_*})\stackrel{i}{\to} D_{1\,*}{\otimes}{\cdots}{\otimes}D_{m\,*}\stackrel{j}{\to}
E_{1\,*}{\otimes}{\cdots}{\otimes}E_{m\,*}\to 0,$$
we have $H_*({\cal Z}(\partial\Delta\!^S;\underline{D_*},\underline{C_*}))\cong\Sigma^{-1}{\rm coker}\,j_*\oplus{\rm ker}\,j_*$,
where $j_*$ is induced by $j$.
By Theorem~3.12,
$$\begin{array}{l}
\quad H_*({\cal Z}(\partial\Delta\!^S;\underline{D_*},\underline{C_*}))\vspace{2mm}\\
\cong\oplus_{\sigma'\subsetneq S,\,\omega'{\cap}S=\emptyset}\,H_*^{\sigma'\!,\omega'}(\partial\Delta\!^S){\otimes}H_*^{\sigma',\omega'}
\oplus_{\sigma''\subsetneq S,\,(S\setminus\sigma'')\subset\omega''}\,H_*^{\sigma''\!,\omega''}(\partial\Delta\!^S){\otimes}H_*^{\sigma'',\omega''}\vspace{2mm}\\
\cong\oplus_{\sigma'\subsetneq S,\,\omega'{\cap}S=\emptyset}\,H_*^{\sigma'\!,\omega'}
\oplus_{\sigma''\subsetneq S,\,(S\setminus\sigma'')\subset\omega''}\,\Sigma^{|S|-|\sigma|-1}H_*^{\sigma''\!,\omega''},
\end{array}$$
Since $H_*(D_{k\,*}/C_{k\,*})={\rm coker}\,\theta_k\oplus\Sigma{\rm ker}\,\theta_k$, we have
$\oplus_{\sigma'\subsetneq S,\,\omega'{\cap}S=\emptyset}\,H_*^{\sigma'\!,\omega'}\cong {\rm ker}\,j_*$,
$\oplus_{\sigma''\subsetneq S,\,(S\setminus\sigma'')\subset\omega''}\,\Sigma^{|S|-|\sigma|-1}H_*^{\sigma''\!,\omega''}
\cong\Sigma^{-1}{\rm coker}\,j_*$.
\vspace{3mm}

\section{(Co)homology Group of Polyhedral Product Objects}\vspace{3mm}
\vspace{3mm}

\hspace*{5.5mm} We need the following well-known definitions.

For topological spaces $X$ and $Y$, $X{\times}Y$ is the cartesian product space.
Empty space is allowed and define $\emptyset{\times}Z=\emptyset$ for all space $Z$.

For simplicial complexes $X$ and $Y$, $X{\times}Y$ is the simplicial complex defined as follows.
Suppose the vertex set of $X$ and $Y$ are respectively $S$ and $T$ and both sets are given a total order.
Then the vertex set of $X{\times}Y$ is $S{\times}T$.
$$X{\times}Y=\{\,\,\{(i_1,j_1),{\cdots},(i_n,j_n)\}\,|\,i_k\in S,\,j_k\in T,\,i_k\leqslant i_l\,{\rm implies}\,j_k\leqslant j_l\,\}.$$
Empty complex $\{\emptyset\}$ is allowed and define $\{\emptyset\}{\times}Z=\{\emptyset\}$ for all simplicial complex $Z$.
The void complex $\{\,\}$ is not allowed in a product of simplicial complexes.

For pointed topological spaces $X$ and $Y$ with base point denoted by $*$,
$X{\wedge}Y$ is the pointed quotient space $(X{\times}Y)/(*{\times}Y\cup X{\times}*)$ with base point the quotient class $*{\times}Y\cup X{\times}*$.
The one point space $*$ with base point $*$ satisfies $*{\wedge}Z=*$ for all pointed spaces $Z$.

For topological spaces $X$ and $Y$, the join $X*Y$ is the topological space defined as follows.
$X*Y=(I{\times}X{\times}Y)/\!\sim$ with $(0,x,y)\sim(0,x',y)$ and $(1,x,y)\sim(1,x,y')$ for all
$x,x'\in X$ and $y,y'\in Y$.
Empty space $\emptyset$ is allowed and define $\emptyset*Z=Z$ for all spaces $Z$.

For simplicial complexes $X$ and $Y$, $X*Y$ is the join complex
$$X*Y=\{\sigma\sqcup\tau\,|\,\sigma\in X,\,\tau\in Y\},$$
where $\sqcup$ means disjoint union. Void complex $\{\,\}$ is allowed and define $\{\,\}*Z=\{\,\}$ for all simplicial complex $Z$.

For a topological space $Z\neq\emptyset$ without base point, the suspension space $\Sigma Z$ of $Z$ is a pointed space defined as follows.
$\Sigma Z=([0,1]{\times}Z)/\sim$ with $(0,z)\sim (0,z')$ and $(1,z)\sim(1,z')$ for all $z,z'\in Z$.
The base point of $\Sigma Z$ is the equivalent class $\{0\}{\times}Z$.
Specifically, define $\Sigma\emptyset=S^0$, the two point space with one of them base point.
Then we have $\Sigma(X{*}Y)\cong (\Sigma X){\wedge}(\Sigma Y)$ for all spaces $X,Y$ without base point.
\vspace{2mm}

{\bf Definition~4.1} Let $K$ be a simplicial complex on $[m]$
and $(\underline{X},\underline{A})$ be a sequence of topological space or simplicial complex pairs $\{(X_k,A_k)\}_{k=1}^m$.

For topological pairs $(X_k,A_k)$, the {\it polyhedral product space} ${\cal Z}(K;\underline{X},\underline{A})$
is the subspace of $X_1{\times}{\cdots}{\times}X_m$ defined as follows.
For a subset $\tau$ of $[m]$, define
$$D(\tau)= Y_1{\times}\cdots{\times}Y_m,\quad Y_k=\left\{\begin{array}{cl}
X_k&{\rm if}\,\,k\in \tau, \\
A_k&{\rm if}\,\,k\not\in \tau.
\end{array}
\right.$$
Then ${\cal Z}(K;\underline{X},\underline{A})=\cup_{\tau\in K}\,D(\tau)$.
Specifically, define ${\cal Z}(\{\,\};\underline{X},\underline{A})=\emptyset$.

For simplicial pairs $(X_k,A_k)$, the {\it polyhedral product simplicial complex} ${\cal Z}(K;\underline{X},\underline{A})$
is similarly defined.
Specifically, define ${\cal Z}(\{\,\};\underline{X},\underline{A})=\{\emptyset\}$.

For pointed topological pairs $(X_k,A_k)$, the {\it polyhedral smash product space} ${\cal Z}^\wedge(K;\underline{X},\underline{A})$
is similarly defined by replacing $\times$ by $\wedge$.
Specifically, define ${\cal Z}^\wedge(\{\,\};\underline{X},\underline{A})=*$.

For topological pairs $(X_k,A_k)$, the {\it polyhedral join space} ${\cal Z}^*(K;\underline{X},\underline{A})$
is similarly defined by replacing $\times$ by $*$.
Specifically, define ${\cal Z}^*(\{\,\};\underline{X},\underline{A})=\emptyset$.

For \,simplicial\, pairs\, $(X_k,A_k)$, \,the\, {\it polyhedral\, join\, simplicial\, complex}\\ ${\cal Z}^*(K;\underline{X},\underline{A})$
is similarly defined by replacing $\times$ by $*$.
Specifically, define ${\cal Z}^*(\{\,\};\underline{X},\underline{A})=\{\,\}$.

If all $(X_k,A_k)$ are the same pair $(X,A)$,
then ${\cal Z}^-(K;\underline{X},\underline{A})$ is denoted by ${\cal Z}^-(K;X,A)$.
\vspace{3mm}

{\bf Example~4.2} For a simplicial complex $K$ on $[m]$ and a subset $S$ of $[m]$
($S\notin K$ is allowed), the simplicial complex ${\rm link}_K S$ is the usual link
${\rm link}_KS=\{\tau\subset [m]{\setminus}S\,|\,S{\cup}\tau\in K\}$ regarded as a simplicial complex on $[m]{\setminus}S$
when $S\in K$.
Define ${\rm link}_KS=\{\,\}$ if $S\notin K$.

For a polyhedral product space (simplicial complex) ${\cal Z}(K;\underline{X},\underline{A})$,
let $S=\{k\,|\,A_k=\emptyset\,\,\}$ ($=\{k\,|\,A_k=\{\emptyset\}\,\}$).
For a polyhedral smash product space ${\cal Z}^\wedge(K;\underline{X},\underline{A})$,
let $S=\{k\,|\,A_k=*\,\}\neq\emptyset$.
For a polyhedral join simplicial complex ${\cal Z}^*(K;\underline{X},\underline{A})$ (no topological analogue),
let $S=\{k\,|\,A_k=\{\,\}\,\}$.
Then for any product ${\scriptstyle\lozenge}$ in Definition~4.1, we have
$${\cal Z}^\lozenge(K;\underline{X},\underline{A})={\cal Z}^\lozenge({\rm link}_KS;\underline{X'},\underline{A'})
{\scriptstyle\lozenge}({\scriptstyle\lozenge}_{k\in S}X_k),$$
where $(\underline{X'},\underline{A'})=\{(X_k,A_k)\}_{k\notin S}$.

For topological space pair $(*,\emptyset)$, ${\cal Z}^*(K;*,\emptyset)=|K|$, where $|K|$ is the geometrical realization space of $K$.
For simplicial complex pair $([1],\{\emptyset\})$, ${\cal Z}^*(K;[1],\{\emptyset\})=K$.
\vspace{3mm}

{\bf Definition~4.3} A (pointed) topological pair $(X,A)$ is called {\it homology split}
if the (reduced) singular homology group homomorphism induced by the inclusion from $A$ to $X$ is a split homomorphism by Definition~3.2.

A simplicial pair $(X,A)$ is called {\it homology split}
if the (reduced) simplicial homology homomorphism induced by the inclusion from $A$ to $X$ is a split homomorphism.

A polyhedral product object ${\cal Z}^-(K;\underline{X},\underline{A})$ in Definition~4.1
is called {\it homology split} if each pair $(X_k,A_k)$ is homology split.
\vspace{3mm}

{\bf Definition~4.4} For a homology split topological pair $(X,A)$,
the {\it indexed homology group} $H_*^\XX(X,A)$ and {\it indexed cohomology group} $H^{\,*}_{\!\XX}(X,A)$ of $(X,A)$
are the $H_*^\XX(D_*,C_*)$ and $H^{\,*}_{\!\XX}(D_*,C_*)$ in Definition~3.3 by taking $(D_*,C_*)$
to be the pair $(S_*(X),S_*(A))$ (the index set $\Lambda$ is trivial),
where $S_*$ means the singular chain complex.
The {\it support index set} of $(X,A)$ is that of $(D_*,C_*)$.

For a homology split pointed topological pair $(X,A)$,
the {\it reduced indexed homology group} $\w H_*^\XX(X,A)$ and {\it reduced indexed cohomology group} $\w H^{\,*}_{\!\XX}(X,A)$ of $(X,A)$
are the $H_*^\XX(D_*,C_*)$ and $H^{\,*}_{\!\XX}(D_*,C_*)$ in Definition~3.3 by taking $(D_*,C_*)$
to be the pair $(\w S_*(X),\w S_*(A))$ (the index set $\Lambda$ is trivial),
where $\w S_*$ means the augmented singular chain complex.
The {\it reduced support index set} of $(X,A)$ is that of $(D_*,C_*)$.

For homology split topological pair $(X,A)$,
the {\it suspension reduced indexed homology group} $\w H_*^\XX(\Sigma X,\Sigma A)$ and
{\it suspension reduced indexed cohomology group} $\w H^{\,*}_{\!\XX}(\Sigma X,\Sigma A)$ of $(X,A)$
are the $H_*^\XX(D_*,C_*)$ and $H^{\,*}_{\!\XX}(D_*,C_*)$ in Definition~3.3 by taking $(D_*,C_*)$
to be the pair $(\w S_*(\Sigma X),\w S_*(\Sigma A))$ (the index set $\Lambda$ is trivial),
where $\w S_*$ means the augmented singular chain complex.
The {\it suspension reduced support index set} of $(X,A)$ is that of $(D_*,C_*)$.

We have all analogue definitions for simplicial complex pairs by replacing the (augmented) singular chain complex
by (augmented) simplicial chain complex and for a simplicial complex $K$, we define $\w C_*(\Sigma K)=\Sigma\w C_*(K)=\w C_{*-1}(K)$ and
$\w C^*(\Sigma K)=\Sigma\w C^*(K)=\w C^{*-1}(K)$.
\vspace{3mm}

{\bf Definition~4.5} Let $(\underline{X},\underline{A})=\{(X_k,A_k)\}_{k=1}^m$ be
a sequence of homology split (space, pointed space or simplicial complex) pairs.

The {\it (reduced, suspension reduced) indexed homology group} and {\it (reduced, suspension reduced) indexed cohomology group}
of $(\underline{X},\underline{A})$ are
$$\begin{array}{rcl}
H_*^{\XX_m}(\underline{X},\underline{A})\!\!&=&\!\!H_*^\XX(X_1,A_1){\otimes}{\cdots}{\otimes}H_*^\XX(X_m,A_m),\vspace{3mm}\\
H^{\,*}_{\!\XX_m}(\underline{X},\underline{A})\!\!&=&\!\!H^{\,*}_{\!\XX}(X_1,A_1){\otimes}{\cdots}{\otimes}H^{\,*}_{\!\XX}(X_m,A_m),\vspace{3mm}\\
\w H_*^{\XX_m}(\underline{X},\underline{A})\!\!&=&\!\!\w H_*^\XX(X_1,A_1){\otimes}{\cdots}{\otimes}\w H_*^\XX(X_m,A_m),\vspace{3mm}\\
\w H^{\,*}_{\!\XX_m}(\underline{X},\underline{A})\!\!&=&\!\!\w H^{\,*}_{\!\XX}(X_1,A_1){\otimes}{\cdots}{\otimes}\w H^{\,*}_{\!\XX}(X_m,A_m).\vspace{3mm}\\
\w H_*^{\XX_m}(\underline{\Sigma X},\underline{\Sigma A})
\!\!&=&\!\!\w H_*^\XX(\Sigma X_1,\Sigma A_1){\otimes}{\cdots}{\otimes}\w H_*^\XX(\Sigma X_m,\Sigma A_m),\vspace{3mm}\\
\w H^{\,*}_{\!\XX_m}(\underline{\Sigma X},\underline{\Sigma A})
\!\!&=&\!\!\w H^{\,*}_{\!\XX}(\Sigma X_1,\Sigma A_1){\otimes}{\cdots}{\otimes}\w H^{\,*}_{\!\XX}(\Sigma X_m,\Sigma A_m).
\end{array}$$

For an index set $\DD\subset\XX_m$, we have all analogue groups $(-)_*^\DD(-)$ and $(-)^{\,*}_{\!\DD}(-)$ of $(X,A)$ on $\DD$
defined to be respectively the restriction group of $(-)_*^{\XX_m}(-)$ and $(-)^{\,*}_{\!\XX_m}(-)$ on $\DD$.
\vspace{3mm}

{\bf Theorem~4.6} {\it Suppose all the following polyhedral product objects are homology split (no smash product for simplicial complex)
and everything be as in Definition~4.5.
Let $\underline{\SS}=\SS_1{\times}{\cdots}{\times}\SS_m$, where each $\SS_k$ is the (reduced, suspension reduced) support index set of $(X_k,A_k)$.
Denote by $\theta_k$ the $\theta$ in Definition~4.4 for $(X_k,A_k)$.
Then for any index set $\DD$ such that $\underline{\SS}\subset\DD\subset\XX_m$, we have
$$\begin{array}{l}
\quad H_*({\cal Z}(K;\underline{X},\underline{A}))\vspace{2mm}\\
\cong H_*^{\DD}(K)\,\widehat\otimes\,H_*^{\DD}(\underline{X},\underline{A})
=\oplus_{(\sigma,\omega)\in\DD}\,H_*^{\sigma,\omega}(K){\otimes}H_*^{\sigma,\omega}(\underline{X},\underline{A}),\hspace{11mm}
  \end{array}
$$
$$\begin{array}{l}
\quad H^*({\cal Z}(K;\underline{X},\underline{A}))\vspace{2mm}\\
\cong H^{\,*}_{\DD}(K)\,\widehat\otimes\,H^{\,*}_{\DD}(\underline{X},\underline{A})
=\oplus_{(\sigma,\omega)\in\DD}\,H^{\,*}_{\sigma,\omega}(K){\otimes}H^{\,*}_{\sigma,\omega}(\underline{X},\underline{A}),\hspace{11mm}
  \end{array}
$$
$$\begin{array}{l}
\quad H_*({\cal Z}^\wedge(K;\underline{X},\underline{A}))\vspace{2mm}\\
\cong H_*^{\DD}(K)\,\widehat\otimes\,\w H_*^{\DD}(\underline{X},\underline{A})
=\oplus_{(\sigma,\omega)\in\DD}\,H_*^{\sigma,\omega}(K){\otimes}\w H_*^{\sigma,\omega}(\underline{X},\underline{A}),\hspace{11mm}
  \end{array}
$$
$$\begin{array}{l}
\quad H^*({\cal Z}^\wedge(K;\underline{X},\underline{A}))\vspace{2mm}\\
\cong H^{\,*}_{\DD}(K)\,\widehat\otimes\,\w H^{\,*}_{\DD}(\underline{X},\underline{A})
=\oplus_{(\sigma,\omega)\in\DD}\,H^{\,*}_{\sigma,\omega}(K){\otimes}\w H^{\,*}_{\sigma,\omega}(\underline{X},\underline{A}),\hspace{11mm}
  \end{array}
$$
$$\begin{array}{l}
\quad H_{*-1}({\cal Z}^*(K;\underline{X},\underline{A}))\vspace{2mm}\\
\cong\w H_*^{\DD}(K)\,\widehat\otimes\,\w H_*^{\DD}(\underline{\Sigma X},\underline{\Sigma A})
=\oplus_{(\sigma,\omega)\in\DD}\,H_*^{\sigma,\omega}(K){\otimes}\w H_*^{\sigma,\omega}(\underline{\Sigma X},\underline{\Sigma A}),
  \end{array}
$$
$$\begin{array}{l}
\quad H^{*-1}({\cal Z}^*(K;\underline{X},\underline{A}))\vspace{2mm}\\
\cong\w H^{\,*}_{\DD}(K)\,\widehat\otimes\,\w H^{\,*}_{\DD}(\underline{\Sigma X},\underline{\Sigma A})
=\oplus_{(\sigma,\omega)\in\DD}\,H^{\,*}_{\sigma,\omega}(K){\otimes}\w H^{\,*}_{\sigma,\omega}(\underline{\Sigma X},\underline{\Sigma A}),
  \end{array}
$$
where the local groups of the indexed (co)homology groups satisfy
$$H_*^{\sigma,\omega}(-),\w H_*^{\sigma,\omega}(-)=H_1{\otimes}{\cdots}{\otimes}H_m,\quad
H_k=\left\{\begin{array}{ll}
{\rm coker}\,\theta_k&{\rm if}\,k\in\sigma,\vspace{1mm}\\
{\rm ker}\,\theta_k&{\rm if}\,k\in\omega,\vspace{1mm}\\
{\rm coim}\,\theta_k&{\rm otherwise},
\end{array}\right.$$
$$H^*_{\sigma,\omega}(-),\w H^*_{\sigma,\omega}(-)=H^1{\otimes}{\cdots}{\otimes}H^m,\quad
H^k=\left\{\begin{array}{ll}
{\rm ker}\,\theta_k^\circ&{\rm if}\,k\in\sigma,\vspace{1mm}\\
{\rm coker}\,\theta_k^\circ&{\rm if}\,k\in\omega,\vspace{1mm}\\
{\rm im}\,\theta_k^\circ&{\rm otherwise}.
\end{array}\right.$$

If each $\theta_k$ is an epimorphism, then we may take $\DD=\RR_m$
and have all analogue equalities by replacing $\DD$ and $\sigma$ respectively by $\RR_m$ and $\emptyset$.

If each $\w H_*(X_k)=0$, then we have
$$\w H_{*-1}({\cal Z}^*(K;\underline{X},\underline{A}))
\cong\w H_{*-1}(K)\otimes\w H_{*-1}(A_1){\otimes}{\cdots}{\otimes}\w H_{*-1}(A_m),$$
$$H^{*-1}({\cal Z}^*(K;\underline{X},\underline{A}))
\cong\w H^{*-1}(K)\otimes\w H^{*-1}(A_1){\otimes}{\cdots}{\otimes}\w H^{*-1}(A_m).$$

Proof} \, Take the $(D_{k\,*},C_{k\,*})$ in Theorem~3.12 to be $(C_*(X_k),C_*(A_k))$,
where $C_*$ means the singular (simplicial) chain complex.
Then we have the equalities for ${\cal Z}(K;\underline{X},\underline{A})$.

Take the $(D_{k\,*},C_{k\,*})$ in Theorem~3.12 to be $(\w C_*(X_k),\w C_*(A_k))$,
where $\w C_*$ means the augmented singular chain complex.
Then we have the equalities for ${\cal Z}^\wedge(K;\underline{X},\underline{A})$.

Since $\Sigma{\cal Z}^*(K;\underline{X},\underline{A})
={\cal Z}^\wedge(K;\underline{\Sigma X},\underline{\Sigma A})$,
the equalities for topological ${\cal Z}^*(K;\underline{X},\underline{A})$ is
by the equalities for topological ${\cal Z}^\wedge(K;\underline{\Sigma X},\underline{\Sigma A})$.

Take the $(D_{k\,*},C_{k\,*})$ in Theorem~3.12 to be $(\Sigma\w C_*(X_k),\Sigma\w C_*(A_k))$,
where $\Sigma\w C_*$ means the suspension augmented simplicial chain complex.
Note that $\Sigma\w  C_*(X*Y)\cong\Sigma\w  C_*(X){\otimes}\Sigma\w  C_*(Y)$.
So the polyhedral product chain complex ${\cal Z}(K;\underline{D_*},\underline{C_*})$
is isomorphic to $\Sigma\w C_*({\cal Z}^*(K;\underline{X},\underline{A}))$.
\hfill$\Box$\vspace{3mm}

{\bf Example~4.7} We compute the (co)homology groups in Theorem~4.6 for the topological pair $(S^r\!,S^p)$,
where the sphere $S^p$ is a subspace of the sphere $S^r$ ($p<r$) by any tame embedding.
Since all the $H_k$ in Theorem~4.6 satisfy that either $H_k=0$ or $H_k=\Sigma^{t_k}\Bbb Z$,
we may identify $H_*^{\sigma\!,\,\omega}(K){\otimes}H_1{\otimes}{\cdots}{\otimes}H_m$
with $\Sigma^{t_1+\cdots+t_m}H_*^{\sigma\!,\,\omega}(K)$ if all $H_k\neq 0$. Then
$$H_*({\cal Z}(K;S^r,S^p))=\oplus_{(\sigma\!,\,\omega)\in\XX_m}\,\w H_{*-r|\sigma|-p|\omega|-1}(K_{\sigma\!,\,\omega}),$$
$$H^*({\cal Z}(K;S^r,S^p))=\oplus_{(\sigma\!,\,\omega)\in\XX_m}\,\w H^{*-r|\sigma|-p|\omega|-1}(K_{\sigma\!,\,\omega}),$$
$$\w H_{*}({\cal Z}^\wedge(K;S^r,S^p))=\oplus_{\sigma\in K}\,\w H_{*-r|\sigma|-p(m{-}|\sigma|)-1}({\rm link}_K\sigma),$$
$$\w H^{*}({\cal Z}^\wedge(K;S^r,S^p))=\oplus_{\sigma\in K}\,\w H^{*-r|\sigma|-p(m{-}|\sigma|)-1}({\rm link}_K\sigma),$$
$$\w H_{*}({\cal Z}^*(K;S^r,S^p))=\oplus_{\sigma\in K}\,\w H_{*-r|\sigma|-p(m-|\sigma|)-m}({\rm link}_K\sigma),$$
$$\w H^{*}({\cal Z}^*(K;S^r,S^p))=\oplus_{\sigma\in K}\,\w H^{*-r|\sigma|-p(m-|\sigma|)-m}({\rm link}_K\sigma).$$

Analogue conclusions (no polyhedral smash product) hold for simplicial sphere pair $(S^r,S^p)$.
\vspace{3mm}

{\bf Example 4.8} We compute the (co)homology groups in Theorem~4.6 for the topological pair $(D^n\!,S^{n-1})$,
where $D^n$ is a ball with boundary $S^{n-1}$.
With similar identification as in Example~4.7, we have
$$H_*({\cal Z}(K;D^n,S^{n-1}))\cong\oplus_{\omega\notin K}\,\w H_{*-(n-1)|\omega|-1}(K|_{\omega}),$$
$$H^*({\cal Z}(K;D^n,S^{n-1}))\cong\oplus_{\omega\notin K}\,\w H^{*-(n-1)|\omega|-1}(K|_{\omega}),$$
$$\w H_{*}({\cal Z}^\wedge(K;D^n,S^{n-1}))\cong\w H_{*-(n-1)m-1}(K),$$
$$\w H^{*}({\cal Z}^\wedge(K;D^n,S^{n-1}))\cong\w H^{*-(n-1)m-1}(K),$$
$$\w H_{*}({\cal Z}^*(K;D^n,S^{n-1}))\cong\w H_{*-nm}(K),$$
$$\w H^{*}({\cal Z}^*(K;D^n,S^{n-1}))\cong\w H^{*-nm}(K).$$

Analogue conclusion (no polyhedra smash product) holds for simplicial pair $(D^n\!,S^{n-1})$,
where $D^n$ is a simplicial ball with boundary $S^{n-1}$.
\vspace{3mm}

The polyhedral join simplicial complexes play an important role in polyhedral product theory as in the following theorem.
\vspace{2mm}

{\bf Theorem~4.9} {\it Let ${\cal Z}^-(K;-,-)$ be any of the polyhedral product object defined before,
for example, polyhedral product (character) chain complex or those in Definition~4.1.
Suppose ${\cal Z}^-(K;\underline{Y},\underline{B})$,
$(\underline{Y},\underline{B})=\{(Y_k,B_k)\}_{k=1}^m$ is a polyhedral product object
defined as follows. For each $k$, $(Y_k,B_k)$ is a pair of polyhedral product objects given by
$$(Y_k,B_k)=\big({\cal Z}^-(X_k;\underline{U_k},\underline{C_k}),\,{\cal Z}^-(A_k;\underline{U_k},\underline{C_k})\big),\,\,
(\underline{U_k},\underline{C_k})=\{(U_i,C_i)\}_{i\in[n_k]},$$
where $(X_k,A_k)$ is a simplicial complex pair on $[n_k]$.
Then
$${\cal Z}^-(K;\underline{Y},\underline{B})={\cal Z}^-({\cal Z}^*(K;\underline{X},\underline{A});\underline{U},\underline{C}),$$
where $(\underline{U}{,}\underline{C})=\{(U_k{,}C_k)\}_{k=1}^n$, $[n]=[n_1]\sqcup{\cdots}\sqcup[n_m]$
(\,$\sqcup$ means disjoint union and so $n=n_1{+}{\cdots}{+}n_m$\,).
\vspace{2mm}

Proof}\, For $\sigma\subset[m]$, let $P_k=X_k$ if $k\in\sigma$ and $P_k=A_k$ if $k\notin\sigma$. Then
$$\begin{array}{l}
\quad {\cal Z}^-(\Delta\!^\sigma;\underline{Y},\underline{B})\vspace{2mm}\\
={\cal Z}^-(P_1;\underline{U_1},\underline{C_1}){\times}{\cdots}{\times}{\cal Z}^-(P_m;\underline{U_m},\underline{C_m})\vspace{2mm}\\
={\cal Z}^-(P_1{*}{\cdots}{*}P_m;\underline{U},\underline{C})\vspace{2mm}\\
={\cal Z}^-({\cal Z}^*(\Delta\!^\sigma;\underline{X},\underline{A});\underline{U},\underline{C}).
\end{array}$$
So for $K\neq\{\,\}$,
$$\begin{array}{l}
\quad {\cal Z}^-(K;\underline{Y},\underline{B})\vspace{2mm}\\
=\cup_{\sigma\in K}{\cal Z}^-(\Delta\!^\sigma;\underline{Y},\underline{B})\vspace{2mm}\\
=\cup_{\sigma\in K}{\cal Z}^-({\cal Z}^*(\Delta\!^\sigma;\underline{X},\underline{A});\underline{U},\underline{C})\vspace{2mm}\\
={\cal Z}^-({\cal Z}^*(K;\underline{X},\underline{A});\underline{U},\underline{C}).
\end{array}$$

For $K=\{\,\}$, ${\cal Z}^-(K;\underline{Y},\underline{B})={\cal Z}^-({\cal Z}^*(K;\underline{X},\underline{A});\underline{U},\underline{C})=\emptyset$.
\hfill$\Box$\vspace{3mm}

The above theorem is a trivial extension of Proposition 5.1 of \cite{AY}.
\vspace{3mm}

{\bf Theorem~4.10} {\it \, Let ${\cal Z}^*(K;\underline{X},\underline{A})$ be such that
each $(X_k,A_k)$ is a simplicial complex pair on $[n_k]$ and $[n]=[n_1]\sqcup{\cdots}\sqcup[n_m]$.
For an index $(\w\sigma\!,\,\w\omega)\in\XX_n$, denote by $\sigma_k\!=\!\w\sigma{\cap}[n_k]$, $\omega_k\!=\!\w\omega{\cap}[n_k]$.
Then we have local complex with respect to $[n]$ (as in Definition~3.9) equalities
$${\cal Z}^*(K;\underline{X},\underline{A})_{\tilde\sigma,\,\tilde\omega}
={\cal Z}^*(K;\underline{X}\,_{\tilde\sigma,\,\tilde\omega},\underline{A}\,_{\tilde\sigma,\,\tilde\omega}),$$
where $(\underline{X}\,_{\tilde\sigma,\,\tilde\omega},\underline{A}\,_{\tilde\sigma,\,\tilde\omega})\!=\!
\{(\,(X_k)_{\sigma_k,\,\omega_k},(A_k)_{\sigma_k,\,\omega_k})\}_{k=1}^m$.
We also have
$${\rm link}_{{\cal Z}^*(K;\underline{X},\underline{A})}\w\sigma=
{\cal Z}^*(K;{\rm link}_{(\underline{X},\underline{A})}\w\sigma),$$
where ${\rm link}_{(\underline{X},\underline{A})}\w\sigma{=}\{({\rm link}_{X_k}\sigma_k,{\rm link}_{A_k}\sigma_k)\}_{k=1}^m$
and link is as in Example~4.2.
\vspace{2mm}

Proof}\, Let $Y^\tau_k=X_k$ if $k\in\tau$ and $Y^\tau_k=A_k$ if $k\notin\tau$.
Then
\begin{eqnarray*}&&{\cal Z}^*(K;\underline{X},\underline{A})_{\tilde\sigma\!,\,\tilde\omega}\vspace{1mm}\\
&=&\cup_{\tau\in K}(Y^\tau_1*{\cdots}*Y^\tau_m)_{\tilde\sigma\!,\,\tilde\omega}\vspace{1mm}\\
&=&\cup_{\tau\in K}(Y^\tau_1)_{\sigma_1,\,\omega_1}*\cdots*(Y^\tau_m)_{\sigma_m,\,\omega_m}\quad
((Y^\tau_k)_{\sigma_k,\omega_k}=\{\,\}\,{\rm allowed})\vspace{1mm}\\
&=&{\cal Z}^*(K;\underline{X}\,_{\tilde\sigma,\,\tilde\omega},\underline{A}\,_{\tilde\sigma,\,\tilde\omega}).
\end{eqnarray*}

For a simplicial  complex $L$ on $[l]$ and $\sigma\in L$, ${\rm link}_L\sigma=L_{\sigma,[l]\setminus\sigma}$.
So the second equality of the theorem is the special case of the first equality for $\w\omega=[n]{\setminus}\w\sigma$.
\vspace{3mm}\hfill$\Box$

{\bf Definition~4.11} Let $(X,A)$ be a simplicial complex pair on $[m]$ and $\TT$ be a subset of $\XX_m$.
The pair $(X,A)$ is called {\it densely split} on $\TT$ if the total (co)homology $\TT$-group (by Definition~3.7) homomorphisms
$$\theta\colon H_*^{\TT}(A)\to H_*^{\TT}(X),\,\,\theta^\circ\colon H^{\,*}_{\!\TT}(X)\to H^{\,*}_{\!\TT}(A)$$
induced by inclusion are split.
Equivalently, for each $(\sigma,\omega)\in\TT$, the local homology homomorphism
$\theta_{\sigma,\omega}\colon H_*^{\sigma,\omega}(A)\to H_*^{\sigma,\omega}(X)$ induced by inclusion
satisfies ${\rm coker}\,\theta_{\sigma,\omega}$, ${\rm ker}\,\theta_{\sigma,\omega}$, ${\rm coim}\,\theta_{\sigma,\omega}$
are all free groups.

For a pair $(X,A)$ densely split on $\TT$,
the {\it densely indexed homology group} $H_*^{\XX;\TT}(X,A)$ and {\it densely indexed cohomology group} $H^{\,*}_{\!\XX;\TT}(X,A)$
of $(X,A)$ on $\TT$ are the $H_*^\XX(D_*,C_*)$ and $H^{\,*}_{\!\XX}(D_*,C_*)$ in Definition~3.3
by taking $(D_*,C_*)$ to be the pair $(T_*^\TT(X),T_*^\TT(A))$ (the index set $\Lambda=\TT$),
where $T_*^\TT(-)$ means the total chain complex on $\TT$ in Definition~3.7.
$\TT$ is called the {\it reference index set} of $(X,A)$.
The {\it dense support index set} of $(X,A)$ on $\TT$ is the support of $(D_*,C_*)$ in Definition~3.3.
\vspace{3mm}

{\bf Definition~4.12} Let $(\underline{X},\underline{A})=\{(X_k,A_k)\}_{k=1}^m$ be a sequence
such that each $(X_k,A_k)$ is a simplicial complex pair on $[n_k]$ and is densely split on $\TT_k$ with
$$\theta_{\sigma_k,\omega_k}\colon H_*^{\sigma_k,\omega_k}(A_k)\to H_*^{\sigma_k,\omega_k}(X_k),\,\,
\theta^\circ_{\sigma_k,\omega_k}\colon H^*_{\sigma_k,\omega_k}(X_k)\to H^*_{\sigma_k,\omega_k}(A_k)$$
the local (co)homology homomorphism induced by inclusion.
$[n]$ is the disjoint union $[n_1]\sqcup{\cdots}\sqcup[n_m]$.

The {\it reference index set} $\underline{\TT}$ of $(\underline{X},\underline{A})$ is the subset of $\XX_n$ given by
$$\underline{\TT}=\{(\w\sigma,\w\omega)\in\XX_n\,|\,(\w\sigma{\cap}[n_k],\w\omega{\cap}[n_k])\in\TT_k\,\,{\rm for}\,\,k=1,{\cdots},m\,\}.$$
It is obvious that $\underline{\TT}=\TT_1{\times}{\cdots}{\times}\TT_m$ by identifying the subset $(\w\sigma,\w\omega)\in\underline{\TT}$
with the subset
$(\w\sigma{\cap}[n_1],\w\omega{\cap}[n_1],{\cdots},\w\sigma{\cap}[n_m],\w\omega{\cap}[n_m])\in\TT_1{\times}{\cdots}{\times}\TT_m$.

The {\it densely indexed homology group} and {\it densely indexed cohomology group} of $(\underline{X},\underline{A})$ on $\underline{\TT}$ are
$$H_*^{\XX_m;\underline{\TT}}(\underline{X},\underline{A})=H_*^{\XX;\TT_1}(X_1,A_1){\otimes}{\cdots}{\otimes}H_*^{\XX;\TT_m}(X_m,A_m),$$
$$H^{\,*}_{\!\XX_m;\underline{\TT}}(\underline{X},\underline{A})=H^{\,*}_{\!\XX;\TT_1}(X_1,A_1){\otimes}{\cdots}{\otimes}H^{\,*}_{\!\XX;\TT_m}(X_m,A_m).$$
Precisely, for $(\sigma,\omega)\in\XX_m$ and $(\w\sigma,\w\omega)\in\underline{\TT}$ with $\sigma_k=\w\sigma{\cap}[n_k]$, $\omega_k=\w\omega{\cap}[n_k]$,
$$H_*^{\sigma\!,\,\omega;\tilde\sigma\!,\,\tilde\omega}(\underline{X},\underline{A})=
H_*^{\sigma\!,\,\omega;\sigma_1,\omega_1,{\cdots},\sigma_m,\omega_m}(\underline{X},\underline{A})
=H_1^{\sigma_1,\omega_1}{\otimes}{\cdots}{\otimes}H_m^{\sigma_m,\omega_m},
$$
$$H^*_{\sigma\!,\,\omega;\tilde\sigma\!,\,\tilde\omega}(\underline{X},\underline{A})=
H^*_{\sigma\!,\,\omega;\sigma_1,\omega_1,{\cdots},\sigma_m,\omega_m}(\underline{X},\underline{A})
=H^1_{\sigma_1,\omega_1}{\otimes}{\cdots}{\otimes}H^m_{\sigma_m,\omega_m},
$$
where
$$H_k^{\sigma_k,\omega_k}=\left\{
\begin{array}{ll}
{\rm coker}\,\theta_{\sigma_k,\omega_k}& {\rm if}\,\, k\in\sigma,\vspace{1mm}\\
{\rm ker}\,\theta_{\sigma_k,\omega_k}& {\rm if}\,\, k\in\omega,\vspace{1mm}\\
{\rm coim}\,\theta_{\sigma_k,\omega_k}& {\rm otherwise},
\end{array}\right.\,\,
H^k_{\sigma_k,\omega_k}=\left\{
\begin{array}{ll}
{\rm ker}\,\theta_{\sigma_k,\omega_k}^\circ& {\rm if}\,\, k\in\sigma,\vspace{1mm}\\
{\rm coker}\,\theta_{\sigma_k,\omega_k}^\circ& {\rm if}\,\, k\in\omega,\vspace{1mm}\\
{\rm im}\,\theta_{\sigma_k,\omega_k}^\circ& {\rm otherwise}.
\end{array}\right.$$
\vspace{3mm}

{\bf Theorem~4.13} {\it Let ${\cal Z}^*(K;\underline{X},\underline{A})$ be a polyhedral join simplicial complex
such that each $(X_k,A_k)$ is densely split on $\TT_k$.
Let $\underline{\SS}=\SS_1{\times}{\cdots}{\times}\SS_m$, where each $\SS_k$ is the dense support index set of $(X_k,A_k)$ on $\TT_k$.
Then for any index set $\DD$ such that $\underline{\SS}\subset\DD\subset\XX_m$
and any index set $\TT\subset\underline{\TT}=\TT_1{\times}{\cdots}{\times}\TT_m$,
the total (co)homology group of ${\cal Z}^*(K;\underline{X},\underline{A})$ on $\TT$ satisfies
$$H_*^{\TT}({\cal Z}^*(K;\underline{X},\underline{A}))\cong
H_*^{\DD}(K)\,\widehat\otimes\,H_*^{\DD;\TT}(\underline{X},\underline{A}),$$
$$H^*_{\TT}({\cal Z}^*(K;\underline{X},\underline{A}))\cong
H^*_{\DD}(K)\,\widehat\otimes\,H^*_{\DD;\TT}(\underline{X},\underline{A}).$$
Precisely, for $(\w\sigma,\w\omega)\in\TT$ with $\sigma_k=\w\sigma{\cap}[n_k]$, $\omega_k=\w\omega{\cap}[n_k]$,
$$H_{*}^{\tilde\sigma,\tilde\omega}({\cal Z}^*(K;\underline{X},\underline{A}))
=\oplus_{(\sigma,\omega)\in\DD}\,H_*^{\sigma,\omega}(K){\otimes}
H_1^{\sigma_1,\omega_1}{\otimes}{\cdots}{\otimes}H_m^{\sigma_m,\omega_m},
$$
$$H^*_{\tilde\sigma,\tilde\omega}({\cal Z}^*(K;\underline{X},\underline{A}))
=\oplus_{(\sigma,\omega)\in\DD}\,H^*_{\sigma,\omega}(K){\otimes}
H^1_{\sigma_1,\omega_1}{\otimes}{\cdots}{\otimes}H^m_{\sigma_m,\omega_m},
$$
where
$$H_k^{\sigma_k,\omega_k}=\left\{
\begin{array}{ll}
{\rm coker}\,\theta_{\sigma_k,\omega_k}& {\rm if}\,\, k\in\sigma,\vspace{1mm}\\
{\rm ker}\,\theta_{\sigma_k,\omega_k}& {\rm if}\,\, k\in\omega,\vspace{1mm}\\
{\rm coim}\,\theta_{\sigma_k,\omega_k}& {\rm otherwise},
\end{array}\right.\,\,
H^k_{\sigma_k,\omega_k}=\left\{
\begin{array}{ll}
{\rm ker}\,\theta_{\sigma_k,\omega_k}^\circ& {\rm if}\,\, k\in\sigma,\vspace{1mm}\\
{\rm coker}\,\theta_{\sigma_k,\omega_k}^\circ& {\rm if}\,\, k\in\omega,\vspace{1mm}\\
{\rm im}\,\theta_{\sigma_k,\omega_k}^\circ& {\rm otherwise}.
\end{array}\right.$$

If all $\theta_{\sigma_k,\omega_k}$ are epimorphisms, then we may take $\DD=\RR_m$.
\vspace{2mm}

Proof}\, Take the $(D_{k\,*},C_{k\,*})$ in Theorem~3.12 to be the total chain complex pair $(T_*^{\TT_k}(X_k),T_*^{\TT_k}(A_k))$.
Then the polyhedral product chain complex ${\cal Z}^*(K;\underline{D_*},\underline{C_*})$
is just $T_*^{\underline{\TT}}({\cal Z}^*(K;\underline{X},\underline{A}))$. So we have
$$H_*^{\underline{\TT}}({\cal Z}^*(K;\underline{X},\underline{A}))\cong
H_*^{\DD}(K)\,\widehat\otimes\,H_*^{\DD;\underline{\TT}}(\underline{X},\underline{A}).$$
So the equality also holds for the restriction group on $\TT$.
\hfill$\Box$\vspace{3mm}

The (co)homology groups of the above theorem can be computed in another way. By Theorem~4.10,
$$H_*^{\tilde\sigma\!,\,\tilde\omega}({\cal Z}^*(K;\underline{X},\underline{A}))
=\w H_{*-1}({\cal Z}^*(K;\underline{X},\underline{A})_{\,\tilde\sigma\!,\,\tilde\omega}))
=\w H_{*-1}({\cal Z}^*(K;\underline{X}_{\,\tilde\sigma\!,\,\tilde\omega},\underline{A}_{\,\tilde\sigma\!,\,\tilde\omega})).
$$
Apply Theorem~4.6 by taking $(\underline{X},\underline{A})=\{(X_{\sigma_k,\omega_k},A_{\sigma_k,\omega_k})\}_{k=1}^m$, we have
$$\w H_{*-1}({\cal Z}^*(K;\underline{X}_{\,\tilde\sigma\!,\,\tilde\omega},\underline{A}_{\,\tilde\sigma\!,\,\tilde\omega}))
=\oplus_{(\sigma\!,\,\omega)\in\DD}\,H_*^{\sigma\!,\,\omega}(K){\otimes}H_1^{\sigma_1,\omega_1}{\otimes}{\cdots}{\otimes}H_m^{\sigma_m,\omega_m}.
\vspace{3mm}$$

{\bf Definition~4.14} The {\it composition complex} ${\cal Z}^*(K;L_1,{\cdots},L_m)$ is
the polyhedral join simplicial complex ${\cal Z}^*(K;\underline{X},\underline{A})$
such that each pair $(X_k,A_k)$ is $(\Delta\!^{[n_k]},L_k)$ on $[n_k]$, where $L_{k}=\{\,\}$ or $\Delta\!^{[n_k]}$ is not allowed.
\vspace{3mm}

The name composition complex comes from Definition 4.5 of \cite{AY}.
\vspace{3mm}

{\bf Example~4.15} We compute the total (co)homology group of the composition complex ${\cal Z}^*(K;L_1,{\cdots},L_k)$.
Suppose that either all $H_*^{\XX_{n_k}}(L_k)$ are free groups, or the (co)homology is taken over a field.
So each pair $(\Delta\!^{[n_k]},L_k)$ is densely split on $\XX_{n_k}$.

Denote by $H_*^{\XX\!;\XX_{n_k}}(\Delta\!^{[n_k]},L_k)=\oplus_{s\in\XX\!,(\sigma_k,\omega_k)\in\XX_{n_k}}\,H_*^{s;\sigma_k,\omega_k}$
the densely indexed homology group of $(\Delta\!^{[n_k]},L_k)$. We have the following by definition.

(1) ${\rm ker}\,\theta_{\sigma_k,\omega_k}= H_*^{\sigma_k,\omega_k}(L_k)$ for $\omega_k\neq\emptyset$
and ${\rm ker}\,\theta_{\sigma_k,\omega_k}= 0$ if $\omega_k=\emptyset$.
This implies $H_*^{{\mathpzc n};\sigma_k,\omega_k}= H_*^{\sigma_k,\omega_k}(L_k)$ for $\omega_k\neq\emptyset$
and $H_*^{{\mathpzc n};\sigma_k,\omega_k}= 0$ if $\omega_k=\emptyset$.

(2) ${\rm coim}\,\theta_{\sigma_k,\emptyset}= H_0^{\sigma_k,\emptyset}(L_k)=\Bbb Z$
for $\sigma_k\in L_k$ with generator denoted by ${\mathpzc i}_{\,\sigma_k}$ and ${\rm coim}\,\theta_{\sigma_k,\omega_k}=0$ if
$\sigma_k\notin K$ or $\omega_k\neq\emptyset$.
This implies $H_*^{{\mathpzc i};\sigma_k,\emptyset}=\Bbb Z({\mathpzc i}_{\,\sigma_k})$ for $\sigma_k\in L_k$
and $H_*^{{\mathpzc i};\sigma_k,\omega_k}=0$ if $\sigma_k\notin L_k$ or $\omega_k\neq\emptyset$.

(3) ${\rm coker}\,\theta_{\sigma_k,\emptyset}=H_0^{\sigma_k,\emptyset}(\Delta\!^{[n_k]})=\Bbb Z$ for $\sigma_k\notin L_k$
with generator denoted by ${\mathpzc e}_{\,\sigma_k}$ and ${\rm coker}\,\theta_{\sigma_k,\omega_k}=0$ if
$\sigma_k\in K$ or $\omega_k\neq\emptyset$.
This implies $H_*^{{\mathpzc e};\sigma_k,\emptyset}=\Bbb Z({\mathpzc e}_{\,\sigma_k})$ for $\sigma_k\notin L_k$
and $H_*^{{\mathpzc e};\sigma_k,\omega_k}=0$ if $\sigma_k\in K$ or $\omega_k\neq\emptyset$.

Apply Theorem~4.13 by taking $\DD=\DD_1{\times}{\cdots}{\times}\DD_m$ with $\DD_k=\XX_{n_k}$. Then we have
$$\begin{array}{l}
\quad H_*^{\XX_n}\big({\cal Z}^*(K;L_1,{\cdots},L_m)\big)\vspace{2mm}\\
\cong\oplus_{(\sigma\!,\,\omega)\in\XX_m}\,H_*^{\sigma\!,\,\omega}(K){\otimes}
\big(\otimes_{i\notin\sigma\cup\omega}\Bbb Z({\mathpzc i}_{\,\sigma_i})
\otimes_{j\in\sigma}\Bbb Z({\mathpzc e}_{\,\sigma_j})\otimes_{k\in\omega} H_*^{\LL_{n_k}}(L_k)\big)\vspace{2mm}\\
\cong\oplus_{(\sigma\!,\,\omega)\in\XX_m}\,H_*^{\sigma\!,\,\omega}(K){\otimes}
\big(\otimes_{k\in\omega} H_*^{\LL_{n_k}}(L_k)\big),
\end{array}$$
$$\begin{array}{l}
\quad H^{\,*}_{\!\XX_n}\big({\cal Z}^*(K;L_1,{\cdots},L_m)\big)\vspace{2mm}\\
\cong\oplus_{(\sigma\!,\,\omega)\in\XX_m}\,H^*_{\sigma\!,\,\omega}(K){\otimes}
\big(\otimes_{i\notin\sigma\cup\omega}\Bbb Z({\mathpzc i}_{\,\sigma_i})
\otimes_{j\in\sigma}\Bbb Z({\mathpzc e}_{\,\sigma_j})\otimes_{k\in\omega} H^{\,*}_{\LL_{n_k}}(L_k)\big)\vspace{2mm}\\
\cong\oplus_{(\sigma\!,\,\omega)\in\XX_m}\,H^*_{\sigma\!,\,\omega}(K){\otimes}
\big(\otimes_{k\in\omega} H^{\,*}_{\LL_{n_k}}(L_k)\big),
\end{array}$$
where the index set $\LL_t=\{(\sigma,\omega)\in\XX_t\,|\,\omega\neq\emptyset\}$.

So the total (co)homology group of ${\cal Z}^*(K;L_1,{\cdots},L_m)$ on any index set $\TT\subset\XX_n$ can be computed.
Specifically, take $\TT=\{(\emptyset,[n])\}$, then we have
$$\begin{array}{l}
\quad\w H_{*-1}\big({\cal Z}^*(K;L_1,{\cdots},L_m)\big)\vspace{2mm}\\
\cong H_*^{\emptyset,[n]}\big({\cal Z}^*(K;L_1,{\cdots},L_m)\big)\vspace{2mm}\\
\cong \w H_{*-1}(K){\otimes}\w H_{*-1}(L_1){\otimes}{\cdots}{\otimes}\w H_{*-1}(L_m),
\end{array}$$
This result coincides with Theorem 4.6 for the case when each $H_*(X_k)=0$.

Similarly, take $\TT=\RR_n$, then we have
$$\begin{array}{l}
\quad H_*^{\RR_n}\big({\cal Z}^*(K;L_1,{\cdots},L_m)\big)\vspace{2mm}\\
\cong\oplus_{(\emptyset,\omega)\in\RR_m}\, H_*^{\emptyset,\omega}(K){\otimes}
\big(\otimes_{i\notin\omega}\Bbb Z({\mathpzc i}_{\,\emptyset})
\otimes_{k\in\omega} H_*^{\overline\RR_{n_k}}(L_k)\big),
\end{array}$$
$$\begin{array}{l}
\quad H^{\,*}_{\!\RR_n}\big({\cal Z}^*(K;L_1,{\cdots},L_m)\big)\vspace{2mm}\\
\cong\oplus_{(\emptyset,\omega)\in\RR_m}\, H^*_{\emptyset,\omega}(K){\otimes}
\big(\otimes_{i\notin\omega}\Bbb Z({\mathpzc i}_{\,\emptyset})
\otimes_{k\in\omega} H^{\,*}_{\overline\RR_{n_k}}(L_k)\big),
\end{array}$$
where $\overline\RR_{n_k}=\{(\emptyset,\omega)\in\RR_{n_k}\,|\,\omega\neq\emptyset\}$.

We compute the above right total homology group directly by Theorem~3.12 by taking $(D_{k\,*},C_{k\,*})$
to be $(T_*^{\RR_{n_k}}(\Delta\!^{[n_k]}),T_*^{\RR_{n_k}}(L_k))$, where $T_*^{\RR_{n_k}}$ means the right total chain complex in Definition~3.7.
In this case, each $\theta_k\colon H_*^{\RR_{n_k}}(L_k)\to H_*^{\RR_{n_k}}(\Delta\!^{[n_k]})$
is an epimorphism. So by Theorem~3.12,
$$H_*^{\RR_n}\big({\cal Z}^*(K;L_1,{\cdots},L_m)\big)\cong
H_*^{\RR_m}(K)\,\widehat\otimes\,\big(H_*^{\RR_{n_1}}(L_1){\otimes}{\cdots}{\otimes}H_*^{\RR_{n_m}}(L_m)\big),$$
where the $(\RR_m{\times}\RR_n)$-group structure (denoted by $H_*^{\RR_m;\RR_n}$) of the right side tensor product group is defined as follows.
For generators $a_k\in H_*^{\emptyset,\omega_k}(L_k)$, $a_1{\otimes}{\cdots}{\otimes}a_m\in H_*^{\emptyset,\omega;\emptyset,\tilde\omega}$,
where $\w\omega{\cap}[n_k]=\omega_k$, $\omega=\{k\,|\,\omega_k\neq\emptyset\}$,
$a_k\in H_*^{\overline{\RR}_{n_k}}(L_k)$ if $\omega_k\neq\emptyset$ and $a_k={\mathpzc i}_{\,\emptyset}$ otherwise.
So
$$\begin{array}{l}
\quad H_*^{\RR_m}(K)\,\widehat\otimes\,\big(H_*^{\RR_{n_1}}(L_1){\otimes}{\cdots}{\otimes}H_*^{\RR_{n_m}}(L_m)\big)\vspace{2mm}\\
\cong\oplus_{(\emptyset,\omega)\in\RR_m}\, H_*^{\emptyset,\omega}(K){\otimes}
\big(\otimes_{i\notin\omega}\Bbb Z({\mathpzc i}_{\,\emptyset})
\otimes_{k\in\omega} H_*^{\overline\RR_{n_k}}(L_k)\big).
\end{array}$$
This result coincides with the previous computation.
\vspace{3mm}

{\bf Conventions} All rings (always commutative) and modules in this paper are vector spaces over the field $\mak$.
$\otimes$ means $\otimes_\mak$.
\vspace{3mm}

{\bf Definition~4.16} For a simplicial complex $K$ on $[m]$ and a sequence of module pairs $(\underline{X},\underline{A})=\{(X_k,A_k)\}_{k=1}^m$,
where each $A_k$ is a submodule of the module $X_k$ over the ring $R_k$,
the {\it polyhedral tensor product module} ${\cal Z}^\otimes(K;\underline{X},\underline{A})$
is a module over $R=R_1{\otimes}{\cdots}{\otimes}R_m$ defined as follows.
For a subset $\tau$ of $[m]$, define module over $R$
$$M(\tau)= Y_1{\otimes}{\cdots}{\otimes}Y_m,\quad Y_k=\left\{\begin{array}{cl}
X_k&{\rm if}\,\,k\in \tau, \\
A_k&{\rm if}\,\,k\not\in \tau.
\end{array}
\right.$$
Then ${\cal Z}^\otimes(K;\underline{X},\underline{A})=+_{\tau\in K}\,M(\tau)$.
Specifically, define ${\cal Z}^\otimes(\{\,\};\underline{X},\underline{A})=0$.

If each $(X_k,A_k)$ is a pair of ideals of $R_k$, the polyhedral tensor product module
is called a {\it polyhedral tensor product ideal} of $R$.
If each $(X_k,A_k)=(R_k,I_k)$ with $I_k$ a proper ideal $\neq 0$,
then ${\cal Z}^\otimes(K;\underline{X},\underline{A})$ is called a
{\it composition ideal} and is denoted by ${\cal Z}^\otimes(K;I_1,{\cdots},I_m)$.
\vspace{3mm}

{\bf Example~4.17} For the polyhedral tensor product module ${\cal Z}^\otimes(K;\underline{X},\underline{A})$,
let $S=\{k\,|\,A_k=0\,\,\}$. Then
$${\cal Z}^\otimes(K;\underline{X},\underline{A})={\cal Z}^\otimes({\rm link}_KS;\underline{X'},\underline{A'}){\otimes}(\otimes_{k\in S}X_k),$$
where $(\underline{X'},\underline{A'})=\{(X_k,A_k)\}_{k\notin S}$.
\vspace{3mm}

{\bf Theorem~4.18} {\it Suppose $N_k$ is a module over $R_k$  for $k=1,{\cdots},m$ and $N=N_1{\otimes}{\cdots}{\otimes}N_m$.
Let $\theta_k\colon{\rm Tor}_*^{R_k}(A_k,N_k)\to{\rm Tor}_*^{R_k}(X_k,N_k)$ be the {\rm Tor} group homomorphism induced by inclusion.
Then
$${\rm Tor}_*^R({\cal Z}^\otimes(K;\underline{X},\underline{A}),N)
\cong\oplus_{(\sigma\!,\,\omega)\in\XX_m} H_{*}^{\sigma\!,\,\omega}(K)\otimes
(H_1{\otimes}{\cdots}{\otimes}H_m),
$$
where $H_k\!=\!{\rm coker}\,\theta_k$ if $k\in\sigma$;
$H_k\!=\!{\rm ker}\,\theta_k$ if $k\in\omega$;
$H_k\!=\!{\rm coim}\,\theta_k$ otherwise.
\vspace{2mm}

Proof}\, For a module $M$, denote by $(F_*(M),d)$ the free resolution of $M$ as follows.
$$\cdots\stackrel{d}{\to} F_2(M)\stackrel{d}{\to} F_1(M)\stackrel{d}{\to} F_0(M)\stackrel{\varepsilon}{\to}M\to 0.$$
Since $A_k$ is a submodule of $X_k$, we may suppose $F_*(A_k)$ is a chain subcomplex of $F_*(X_k)$
(for example, we may take $F_*(X_k)=F_*(A_k)\oplus F_*(X_k/A_k)$).

Let $M(\tau)=Y_1{\otimes}{\cdots}{\otimes}Y_m$ be as in Definition~4.16.
Then $(F_*(M(\tau)),d)=(F_*(Y_1){\otimes}{\cdots}{\otimes}F_*(Y_m),d)$.
So the free resolution $(F_*({\cal Z}^\otimes(K;\underline{X},\underline{A})),d)$ can be taken to be
the polyhedral product chain complex $({\cal Z}(K;\underline{F_*(X)},\underline{F_*(A)}),d)$
with $(\underline{F_*(X)},\underline{F_*(A)})=\{(F_*(X_k),F_*(A_k))\}_{k=1}^m$.
Take the $(D_{k\,*},C_{k\,*})$ in Theorem~3.12 to be
$(F(X_k){\otimes}_{R_k}N_k,F(A_k){\otimes}_{R_k}N_k)$. Then
$$({\cal Z}(K;\underline{D_*},\underline{C_*}),d)=
({\cal Z}(K;\underline{F_*(X)},\underline{F_*(A)}){\otimes}_RN,d)$$
and $H_*({\cal Z}(K;\underline{D_*},\underline{C_*}))$ is just the Tor group of the theorem.
\hfill$\Box$\vspace{3mm}

{\bf Theorem~4.19} {\it Suppose $\mak$ is a trivial module over $R_k$, i.e., there is an algebra homomorphism
$\varepsilon_k\colon R_k\to\mak$.
Then
$${\rm Tor}^{R}_*({\cal Z}^\otimes(K;I_1,{\cdots},I_m),\mak)
\cong\oplus_{\sigma\in K}\, H_{*}^{\sigma,[m]\setminus\sigma}(K)\otimes
\big(\otimes_{k\notin\sigma}{\rm  Tor}^{R_k}_*(I_k,\mak)\big).
\vspace{2mm}$$

Proof}\, The long exact sequence
$${\scriptstyle\cdots\,\to\,{\rm Tor}_1^{R_k}(R_k/I_k,\mak)\,\to\,
{\rm Tor}_0^{R_k}(I_k,\mak)\,\to\,{\rm Tor}_0^{R_k}(R_k,\mak)\,\to\,{\rm Tor}_0^{R_k}(R_k/I_k,\mak)\,\to\, 0}$$
satisfies ${\rm Tor}_0^{R_k}(I_k,\mak)=0$, ${\rm Tor}_0^{R_k}(R_k,\mak)={\rm Tor}_0^{R_k}(R_k/I_k,\mak)=\mak$,
${\rm Tor}_s^{R_k}(R_k,\mak)=0$ if $s>0$. So the homomorphism
$$\theta_k\colon {\rm Tor}_*^{R_k}(I_k,\mak)\to{\rm Tor}_*^{R_k}(R_k,\mak)$$
induced by inclusion satisfies ${\rm coker}\,\theta_k={\rm Tor}_0^{R_k}(R_k,\mak)\cong\mak$,
${\rm coim}\,\theta_k=0$ and
${\rm ker}\,\theta_k={\rm Tor}_*^{R_k}(I_k,\mak)$.
So by identifying ${\rm coker}\,\theta_k{\otimes}X$ with $X$, we have the equality of the theorem by Theorem~4.18.
\hfill$\Box$\vspace{3mm}

\section{Duality Isomorphisms}\vspace{3mm}

\hspace*{5.5mm} {\bf Definition~5.1} Let $K$ be a simplicial complex with vertex set a subset of $S\neq\emptyset$.
The Alexander {\it dual of $K$ relative to} $S$ is the simplicial
complex
$$K^\circ=\{\,S{\setminus}\sigma\,\,|\,\,\sigma\subset S,\,\sigma\!\notin\! K\,\}.$$

It is obvious that
$(K^\circ)^\circ=K$, $(K_1{\cup}K_2)^\circ=(K_1)^\circ{\cap}(K_2)^\circ$ and
$(K_1{\cap}K_2)^\circ=(K_1)^\circ{\cup}(K_2)^\circ$. Specifically, $(\Delta\!^S)^\circ=\{\,\}$
and $(\partial\Delta\!^S)^\circ=\{\emptyset\}$.
\vspace{3mm}

{\bf Theorem~5.2} {\it For a polyhedral product space ${\cal Z}(K;\underline{X},\underline{A})$, we have the
complement space duality isomorphism
$$(X_1{\times}{\cdots}{\times}X_m){\setminus}{\cal Z}(K;\underline{X},\underline{A})
={\cal Z}(K;\underline{X},\underline{A})^c,$$
where $(\underline{X},\underline{A}^c)=\{(X_k,A_k^c)\}_{k=1}^m$ with $A_k^c=X_k{\setminus}A_k$
and $K^\circ$ is the dual of $K$ relative to $[m]$.
}\vspace{2mm}

{\it Proof}\, For $\sigma\subset[m]$ but $\sigma\neq[m]$ ($\sigma=\emptyset$ is allowed),
$$\begin{array}{l}
\quad(X_1{\times}{\cdots}{\times}X_m)\setminus{\cal Z}(\Delta\!^\sigma;\underline{X},\underline{A})\vspace{2mm}\\
=\cup_{j\notin\sigma}\,X_1{\times}{\cdots}{\times}(X_j{\setminus}A_j){\times}{\cdots}{\times}X_m\vspace{2mm}\\
=\cup_{j\in[m]\setminus\sigma}\,{\cal Z}(\Delta\!^{[m]\setminus\{j\}};\underline{X},\underline{A}^c)\vspace{2mm}\\
={\cal Z}((\Delta\!^\sigma)^\circ;\underline{X},\underline{A}^c).
\end{array}$$
So for $K\neq\Delta\!^{[m]}$ or $\{\,\}$,
$$\begin{array}{l}
\quad{\cal Z}(\Delta\!^{[m]};\underline{X},\underline{A})\setminus{\cal Z}(K;\underline{X},\underline{A})\vspace{2mm}\\
={\cal Z}(\Delta\!^{[m]};\underline{X},\underline{A})\setminus\big(\cup_{\sigma\in K}{\cal Z}(\Delta\!^\sigma;\underline{X},\underline{A})\big)\vspace{2mm}\\
=\cap_{\sigma\in K}\big({\cal Z}(\Delta\!^{[m]};\underline{X},\underline{A})\setminus{\cal Z}(\Delta\!^\sigma;\underline{X},\underline{A})\big)\vspace{2mm}\\
=\cap_{\sigma\in K}{\cal Z}((\Delta\!^\sigma)^\circ;\underline{X},\underline{A}^c)\vspace{2mm}\\
={\cal Z}(K^\circ;\underline{X},\underline{A}^c)\end{array}$$

For $K=\Delta\!^{[m]}$ or $\{\,\}$, the above equality holds naturally. \hfill $\Box$
\vspace{3mm}

{\bf Example~5.3} Let $\Bbb F$ be a field and $V$ be a linear space over $\Bbb F$ with base $e_1,{\cdots},e_m$.
For a subset $\sigma=\{i_1,{\cdots},i_s\}\subset[m]$, denote by $\Bbb F(\sigma)$ the subspace of $V$ with base $e_{i_1},{\cdots},e_{i_s}$.
Then for $\Bbb F=\Bbb R$ or $\Bbb C$ and a simplicial complex $K$ on $[m]$, we have
$$V\setminus(\cup_{\sigma\in K}\Bbb R(\sigma))=\Bbb R^m\setminus{\cal Z}(K;\Bbb R,\{0\})={\cal Z}(K^\circ;\Bbb R,\Bbb R{\setminus}\{0\})
\simeq {\cal Z}(K^\circ;D^1,S^0),$$
$$V\setminus(\cup_{\sigma\in K}\Bbb C(\sigma))=\Bbb C^m\setminus{\cal Z}(K;\Bbb C,\{0\})={\cal Z}(K^\circ;\Bbb C,\Bbb C{\setminus}\{0\})
\simeq {\cal Z}(K^\circ;D^2,S^1).$$
This example is applied in Lemma 2.4 of \cite{GW}.
\vspace{3mm}

{\bf Theorem~5.4} {\it Let $K$ and $K^\circ$ be the dual of each other relative to $[m]$.
Then for any $(\sigma\!,\,\omega)\in\XX_m$ such that $\omega\neq\emptyset$, we have local complex with respect to $\omega$ (Definition~3.9) duality isomorphism
$$(K_{\sigma\!,\,\omega})^\circ=(K^\circ)_{\sigma'\!,\,\omega},\quad\sigma'=[m]{\setminus}(\sigma{\cup}\omega),$$
where $(K_{\sigma\!,\,\omega})^\circ$ is the dual of $K_{\sigma\!,\,\omega}$ relative to $\omega$.
We also have local (co)homology group duality isomorphisms
$$\gamma_{K,\sigma\!,\,\omega}\colon H_*^{\sigma\!,\,\omega}(K)\to
H^{|\omega|-*-1}_{\sigma'\!,\,\omega}(K^\circ),\quad
\gamma^\circ_{K,\sigma\!,\,\omega}\colon H^*_{\sigma\!,\,\omega}(K)
\to H_{|\omega|-*-1}^{\sigma'\!,\,\omega}(K^\circ)$$
such that $\gamma_{K,\sigma,\omega}^\circ=\gamma_{K^\circ\!,\sigma',\omega}^{-1}$.
\vspace{2mm}

Proof}\, Suppose $\sigma\in K$. Then
$$\begin{array}{l}
\quad (K^\circ)_{\sigma'\!,\,\omega}\vspace{1mm}\\
=\{\eta\,\,|\,\,\eta\subset\omega,\,\,[m]{\setminus}(\sigma'{\cup}\eta)=\sigma\!\cup\!(\omega\!\setminus\!\eta)\notin K\,\}\vspace{1mm}\\
=\{\omega\!\setminus\!\tau\,\,|\,\,\tau\subset\omega,\,\,\sigma\!\cup\!\tau\notin K\,\}
\quad(\tau=\omega\!\setminus\!\eta)\vspace{1mm}\\
=(K_{\sigma\!,\,\omega})^\circ.
\end{array}$$

If $\sigma\notin K$, then $(K_{\sigma\!,\omega})^\circ=\Delta\!^\omega=(K^\circ)_{\sigma'\!,\omega}$.

Let $(C_*(\Delta\!^{\omega},K_{\sigma\!,\,\omega}),d)$ be the relative simplicial chain complex.
Since $\w H_*(\Delta\!^{\omega})=0$, we have a boundary isomorphism\vspace{1mm}\\
\hspace*{30mm}$\partial\colon H_*(\Delta\!^\omega,K_{\sigma\!,\,\omega})\stackrel{\cong}{\longrightarrow}
\w H_{*-1}(K_{\sigma\!,\,\omega})=H_*^{\sigma\!,\,\omega}(K)$.\vspace{1mm}\\
$C_*(\Delta\!^{\omega},K_{\sigma\!,\,\omega})$ has a set of generators consisting of all non-simplices of $K_{\sigma\!,\,\omega}$, i.e., $K^c_{\sigma\!,\,\omega}=\{\eta\subset\omega\,|\,\eta\not\in K_{\sigma\!,\,\omega}\}$ is a set of generators of $C_*(\Delta\!^{\omega},K_{\sigma\!,\,\omega})$.
So we may denote $(C_*(\Delta\!^\omega,K_{\sigma\!,\,\omega}),d)$ by $(C_*(K^c_{\sigma\!,\,\omega}),d)$, where
$\eta\in K^c_{\sigma\!,\,\omega}$ has degree $|\eta|{-}1$.
The correspondence $\eta\to \omega{\setminus}\eta$ for all $\eta\in K^c_{\sigma\!,\,\omega}$ induces a dual complex isomorphism\vspace{1mm}\\
\hspace*{32mm}$\psi\colon(C_*(K^c_{\sigma\!,\,\omega}),d)\to
(\w C^{*}((K_{\sigma\!,\,\omega})^\circ),\delta)$.\vspace{1mm}\\
Since $(K_{\sigma\!,\,\omega})^\circ=(K^\circ)_{\sigma'\!,\,\omega}$, we have induced homology group isomorphism
$\bar\psi\colon H_*(\Delta\!^\omega,K_{\sigma\!,\,\omega})\to H^{|\omega|-*-1}_{\sigma'\!,\,\omega}(K^\circ)$.
Define $\gamma_{K,\sigma\!,\,\omega}=\bar\psi\partial^{-1}$.
\hfill $\Box$\vspace{3mm}

{\bf Definition~5.5} For a simplicial complex $K$ on $[m]$,
the densely indexed (co)homology group $H_*^{\XX;\XX_m}(\Delta\!^{[m]},K)$ and $H^{\,*}_{\!\XX;\XX_m}(\Delta\!^{[m]},K)$
are defined as in Example~4.15 even if $H_*^{\XX_m}(K)$ is not free.
Precisely, the local groups of the densely indexed homology (dual case is similar) are defined as follows.
\begin{eqnarray*}H_*^{{\mathpzc n};\sigma,\omega}(\Delta\!^{[m]},K)&=&\left\{\begin{array}{ll}
 H_*^{\sigma,\omega}(K)&{\rm if}\,\omega\neq\emptyset,\vspace{2mm}\\
0&{\rm  otherwise},
\end{array}\right.\\
H_*^{{\mathpzc i};\sigma,\omega}(\Delta\!^{[m]},K)&=&\left\{\begin{array}{ll}
 H_*^{\sigma,\omega}(K)\cong\Bbb Z&{\rm if}\,\sigma\in K,\,\omega=\emptyset,\vspace{2mm}\\
0&{\rm  otherwise},
\end{array}\right.\\
H_*^{{\mathpzc n};\sigma,\omega}(\Delta\!^{[m]},K)&=&\left\{\begin{array}{ll}
 H_*^{\sigma,\omega}(\Delta\!^{[m]})\cong\Bbb Z&{\rm if}\,\sigma\notin K,\,\omega=\emptyset,\vspace{2mm}\\
0&{\rm  otherwise}.
\end{array}\right.\end{eqnarray*}
The generators of $H_0^{\sigma,\emptyset}(\Delta\!^{[m]})$ and $H_0^{\sigma,\emptyset}(K)$ are
respectively denoted by ${\mathpzc e}_{\,\sigma}$ and ${\mathpzc i}_{\,\sigma}$.
Then the densely indexed (co)homology group duality isomorphisms
$$\gamma_K\colon H_*^{\XX;\XX_m}(\Delta\!^{[m]}{\!,}K)\to H^{\,*}_{\!\XX;\XX_m}(\Delta\!^{[m]}{\!,}K^\circ),$$
$$\gamma_K^{\,\circ}\colon H^{\,*}_{\!\XX;\XX_m}(\Delta\!^{[m]}{\!,}K)\to H_*^{\XX;\XX_m}(\Delta\!^{[m]}{\!,}K^\circ)\,$$
are defined as follows.
For $x\in H_*^{{\mathpzc n};\sigma,\omega}(\Delta\!^{[m]}{,}K)=H_*^{\sigma,\omega}(K)$ with $\omega\neq\emptyset$,
$\gamma_K(x)=\gamma_{K,\sigma,\omega}(x)$ with $\gamma_{K,\sigma,\omega}$ as in Theorem~5.4.
Define $\gamma_K({\mathpzc e}_{\,\sigma})={\mathpzc i}_{\,\sigma'}$ and $\gamma_K({\mathpzc i}_{\,\sigma})={\mathpzc e}_{\,\sigma'}$,
where $\sigma'=[m]{\setminus}(\sigma{\cup}\omega)$.
We also denote the local homomorphisms of $\gamma_K,\gamma_K^\circ$ by $\gamma_{K,\sigma,\omega},\gamma_{K,\sigma,\omega}^\circ$
and have the equality  $\gamma_{K,\sigma,\omega}^\circ=\gamma_{K^\circ\!,\sigma',\omega}^{-1}$.

The duality isomorphisms
$$\overline\gamma_K\colon H_*^{\LL_m}(K)\!=\!H_*^{{\mathpzc n};\XX_m}(\Delta\!^{[m]},K)\to
H^{\,*}_{\!\LL_m}(K^\circ)\!=\!H^{\,*}_{\!{\mathpzc n};\XX_m}(\Delta\!^{[m]},K^\circ),$$
$$\overline\gamma_K^{\,\circ}\colon H^{\,*}_{\!\LL_m}(K)\!=\!H^{\,*}_{\!{\mathpzc n};\XX_m}(\Delta\!^{[m]},K)\to
H_*^{\LL_m}(K^\circ)\!=\!H_*^{{\mathpzc n};\XX_m}(\Delta\!^{[m]},K^\circ)\,$$
are the restriction of $\gamma_K$, $\gamma_K^\circ$, where $\LL_m=\{(\sigma,\omega)\in\XX_m\,|\,\omega\neq\emptyset\}$.
\vspace{3mm}

The dual of ${\cal Z}^*(K;\underline{X},\underline{A})$ relative to $[n]$ is in general not a polyhedral join complex.
But the composition complex in Definition~4.14 is.
\vspace{2mm}

{\bf Theorem~5.6} {\it Let ${\cal Z}^*(K;L_1,{\cdots},L_m)^\circ$ be the dual of ${\cal Z}^*(K;L_1,{\cdots},L_m)$
relative to $[n]=[n_1]\sqcup{\cdots}\sqcup[n_m]$ ($\sqcup$ disjoint union). Then
$${\cal Z}^*(K;L_1,{\cdots},L_m)^\circ={\cal Z}^*(K^\circ;L_1^\circ,{\cdots},L_m^\circ),$$
where $K^\circ$ is the dual of $K$ relative to $[m]$ and $L_k^\circ$ is the dual  of $L_k$ relative to $[n_k]$.
So if $K$ and all $L_k$ are self dual ($X\cong X^\circ$ relative to its non-empty vertex set),
then ${\cal Z}^*(K;L_1,{\cdots},L_m)$ is self dual.
\vspace{2mm}

Proof}\, For $\sigma\subset[m]$ but $\sigma\neq[m]$ ($\sigma=\emptyset$ is allowed),
$$\begin{array}{l}
\quad {\cal Z}^*(\Delta\!^\sigma;L_1,{\cdots},L_m)^\circ\vspace{2mm}\\
=\{[n]\!\setminus\!\tau\,|\,\tau\in\cup_{j\notin\sigma}\Delta\!^{n_1}*{\cdots}*(\Delta\!^{n_j}{\setminus}L_j)*{\cdots}*\Delta\!^{n_m}\}\vspace{2mm}\\
=\cup_{j\notin\sigma}\,\Delta\!^{n_1}*{\cdots}*L_j^\circ*{\cdots}*\Delta\!^{n_m}\vspace{2mm}\\
={\cal Z}^*((\Delta\!^\sigma)^\circ;L_1^\circ,{\cdots},L_m^\circ),
\end{array}$$
So for $K\neq[m]$ or $\{\,\}$,
$$\begin{array}{l}
\quad {\cal Z}^*(K;L_1,{\cdots},L_m)^\circ\vspace{2mm}\\
=(\cup_{\sigma\in K}{\cal Z}^*((\Delta\!^\sigma);L_1,{\cdots},L_m))^\circ\vspace{2mm}\\
={\cal Z}^*(\cap_{\sigma\in K}(\Delta\!^\sigma)^\circ;L_1^\circ,{\cdots},L_m^\circ)\vspace{2mm}\\
={\cal Z}^*(K^\circ;L_1^\circ,{\cdots},L_m^\circ).
\end{array}$$

For $K=\Delta\!^{[m]}$ or $\{\,\}$, the equality holds naturally.
\hfill$\Box$\vspace{3mm}

{\bf Theorem~5.7} {\it Let $L_k$ be a simplicial complex on $[n_k]$ such that the pair $(\Delta\!^{[n_k]},L_k)$ is densely split
on $\XX_{n_k}$ by Definition~4.11 for $k=1,{\cdots},m$.
Then the $\gamma_-$ and $\overline\gamma_-$ in Definition~5.5 satisfy the following commutative diagram
$$\begin{array}{ccc}
{\scriptstyle H_*^{\LL_n}({\cal Z}^*(K;L_1,\cdots,L_m))}&\cong&
{\scriptstyle  H_*^{\LL_m}(K){\otimes}H_*^{\LL_m;\LL_n}(\underline{\Delta},\underline{L})}\,\\
\downarrow\,\,{\scriptstyle\overline\gamma_{{\cal Z}^*(K;L_1,\cdots,L_m)}}&&
\downarrow\,\,{\scriptstyle\overline\gamma_K\otimes\overline\gamma_{(\underline{\Delta},\underline{L})}}\vspace{1mm}\\
{\scriptstyle H^{\,*}_{\!\LL_n}({\cal Z}^*(K^\circ;L_1^\circ,\cdots,L_m^\circ))}&\cong&
{\scriptstyle H^*_{\LL_m}(K^\circ){\otimes}H^*_{\LL_m;\LL_n}(\underline{\Delta},\underline{L})},
\end{array}$$
where $\overline\gamma_{(\underline{\Delta},\underline{L})}$ is the restriction of
$\gamma_{(\underline{\Delta},\underline{L})}=\gamma_{L_1}{\otimes}{\cdots}{\otimes}\gamma_{L_m}$.
Equivalently, for each $(\w\sigma,\w\omega)\in\LL_n$, we have
$$\gamma_{{\cal Z}^*(K;L_1,\cdots,L_m),\tilde\sigma,\tilde\omega}
=\gamma_{K,\sigma,\omega}{\otimes}\gamma_{L_1,\sigma_1,\omega_1}{\otimes}{\cdots}{\otimes}\gamma_{L_m,\sigma_m,\omega_m}$$
where  $\sigma_k=\w\sigma{\cap}[n_k]$, $\omega_k=\w\omega{\cap}[n_k]$, $\sigma=\{k\,|\,\sigma_k\notin L_k\}$, $\omega=\{k\,|\,\omega_k\neq\emptyset\}$,
$\gamma_{{\cal Z}^*(K;L_1,\cdots,L_m),\tilde\sigma,\tilde\omega}$ and $\gamma_{K,\sigma,\omega}$ are as defined in Theorem~5.4 and
$\gamma_{L_k,\sigma_k,\omega_k}$ is as defined in Definition~5.5.
\vspace{2mm}

Proof}\, Denote by $\overline{X}={\cal Z}^*(X;L_1,{\cdots},L_m)$.
Then for simplicial complexes $K_1,K_2$ on $[m]$, we have
the following commutative diagram of long exact sequences
\[\begin{array}{ccccccc}
{\scriptstyle\cdots\,\,\longrightarrow}\!\!&\!\!{\scriptstyle H_k^{\tilde\sigma\!,\,\tilde\omega}(\overline K_1{\cap}\overline K_2)}
\!\!&\!\!{\scriptstyle\longrightarrow}
\!\!&\!\!{\scriptstyle H_k^{\tilde\sigma\!,\,\tilde\omega}(\overline K_1){\oplus}H_k^{\tilde\sigma\!,\,\tilde\omega}(\overline K_2)}
\!\!&\!\!{\scriptstyle\longrightarrow}\!\!&\!\!
{\scriptstyle H_k^{\tilde\sigma\!,\,\tilde\omega}(\overline{K_1{\cup}K_2})}\!\!&\!\!{\scriptstyle\longrightarrow\,\,\cdots}\,\vspace{1mm}\\
\!\!&\!\!^{\gamma_{\,\overline{K_1\cap K_2}}}\downarrow\quad\quad\!\!&\!\!\!\!&\!\!
^{\gamma_{\,\overline K_1}\oplus\gamma_{\,\overline K_2}}\downarrow\quad\quad\quad\!\!&\!\!\!\!&\!\!
^{\gamma_{\,\overline{K_1\cup K_2}}}\downarrow\quad\quad\!\!&\!\!\,\\
{\scriptstyle\cdots\,\,\longrightarrow}\!\!&\!\!
{\scriptstyle H^{k'}_{\tilde\sigma'\!,\,\tilde\omega}(\overline K_1^\circ{\cup}\overline K_2^\circ)}
\!\!&\!\!{\scriptstyle\longrightarrow}\!\!&\!\!
{\scriptstyle H^{k'}_{\tilde\sigma'\!,\,\tilde\omega}(\overline K_1^\circ){\oplus}
H^{k'}_{\tilde\sigma'\!,\,\tilde\omega}(\overline K_2^\circ)}
\!\!&\!\!{\scriptstyle\longrightarrow}\!\!&\!\!{\scriptstyle H^{k'}_{\tilde\sigma'\!,\,\tilde\omega}
(\overline K_1^\circ\cap\overline K_2^\circ)}\!\!&\!\!{\scriptstyle\longrightarrow\,\,\cdots},
  \end{array}\vspace{1mm}\]
where $\gamma_{-}=\gamma_{-,\tilde\sigma,\tilde\omega}$, $k'=|\w\omega|{-}k{-}1$, $\w\sigma'=[n]{\setminus}(\w\sigma{\cup}\w\omega)$.
By Theorem~5.6 and the computation of Example~4.15,
$$H_*^{\tilde\sigma,\tilde\omega}(\overline K)
\cong\w H_{*-1}^{\sigma\!,\,\omega}(K){\otimes}
\big(\otimes_{i\notin\sigma\cup\omega}\Bbb Z({\mathpzc i}_{\,\sigma_i})
\otimes_{j\in\sigma}\Bbb Z({\mathpzc e}_{\,\sigma_j})\otimes_{k\in\omega} H_*^{\sigma_k,\omega_k}(L_k)\big),
$$
$$H^{|\tilde\omega|-*-1}_{\tilde\sigma'\!,\,\tilde\omega}(\overline K^\circ)\cong\,
H^{|\omega|-*-1}_{\sigma'\!,\,\omega}(K^\circ){\otimes}
\big(\otimes_{i\notin\sigma'\cup\omega}\Bbb Z({\mathpzc i}_{\,\sigma'_i})
\otimes_{j\in\sigma'}\Bbb Z({\mathpzc e}_{\,\sigma'_j})\otimes_{k\in\omega} H^{|\omega_k|-*-1}_{\sigma'_k,\omega_k}(L_k^\circ)\big),
$$
where $\sigma'_k=[n_k]{\setminus}(\sigma_k{\cup}\omega_k)$.
So we have
the following commutative diagram of long exact sequences
\[\begin{array}{ccccccc}
{\scriptstyle\cdots\,\,\longrightarrow}\!\!&\!\!{\scriptstyle H_l^{\sigma\!,\,\omega}( K_1{\cap} K_2)\otimes A}
\!\!&\!\!{\scriptstyle\longrightarrow}
\!\!&\!\!{\scriptstyle (H_l^{\sigma\!,\,\omega}( K_1){\oplus}H_l^{\sigma\!,\,\omega}( K_2))\otimes A}
\!\!&\!\!{\scriptstyle\longrightarrow}\!\!&\!\!
{\scriptstyle H_l^{\sigma\!,\,\omega}({K_1{\cup}K_2})\otimes A}\!\!&\!\!{\scriptstyle\longrightarrow\,\,\cdots}\,\vspace{1mm}\\
\!\!&\!\!^{\gamma_{{K_1\cap K_2}}\otimes\gamma}\downarrow\quad\quad\!\!&\!\!\!\!&\!\!
^{(\gamma_{K_1}\oplus\gamma_{K_2})\otimes\gamma}\downarrow\quad\quad\quad\!\!&\!\!\!\!&\!\!
^{\gamma_{{K_1\cup K_2}}\otimes\gamma}\downarrow\quad\quad\!\!&\!\!\,\\
{\scriptstyle\cdots\,\,\longrightarrow}\!\!&\!\!
{\scriptstyle H^{l'}_{\sigma'\!,\,\omega}(K_1^\circ{\cup} K_2^\circ)\otimes B}
\!\!&\!\!{\scriptstyle\longrightarrow}\!\!&\!\!
{\scriptstyle (H^{l'}_{\sigma'\!,\,\omega}(K_1^\circ){\oplus}
H^{l'}_{\sigma'\!,\,\omega}(K_2^\circ))\otimes B}
\!\!&\!\!{\scriptstyle\longrightarrow}\!\!&\!\!{\scriptstyle H^{l'}_{\sigma'\!,\,\omega}
( K_1^\circ\cap K_2^\circ)\otimes B}\!\!&\!\!{\scriptstyle\longrightarrow\,\,\cdots},
  \end{array}\vspace{1mm}\]
where $l'=|\omega|{-}l{-}1$, $\gamma_{-}=\gamma_{-,\sigma,\omega}$,
$\gamma=\gamma_{L_1,\sigma_1,\omega_1}{\otimes}{\cdots}{\otimes}\gamma_{L_m,\sigma_m,\omega_m}$,
$$A=\otimes_{i\notin\sigma\cup\omega}\Bbb Z({\mathpzc i}_{\,\sigma_i})
\otimes_{j\in\sigma}\Bbb Z({\mathpzc e}_{\,\sigma_j})\otimes_{k\in\omega} H_*^{\sigma_k,\omega_k}(L_k),
$$
$$B=\otimes_{i\notin\sigma'\cup\omega}\Bbb Z({\mathpzc i}_{\,\sigma'_i})
\otimes_{j\in\sigma'}\Bbb Z({\mathpzc e}_{\,\sigma'_j})\otimes_{k\in\omega} H^{|\omega_k|-*-1}_{\sigma'_k,\omega_k}(L_k^\circ).
$$
The above two commutative diagrams imply that if the theorem holds for $K_1{\cap}K_2$,
$K_1$ and $K_2$, then the theorem holds for $K_1{\cup}K_2$.
So by induction on the number of maximal simplexes, we only need prove the theorem for the case
that $K$ has only one maximal simplex.

Now suppose $K=\Delta\!^S$ with $S\neq[m]$ (the case $K=\Delta\!^{[m]}$ is trivial)
and $H_*^{\tilde\sigma,\tilde\omega}(\overline K)\neq 0$.
Then $H_*^{\sigma,\omega}(K)\neq 0$. This implies $\sigma\subset S$, $\omega\neq\emptyset$ and $\omega{\cap}S=\emptyset$.
So $H_*^{\tilde\sigma,\tilde\omega}(\overline K)$ is generated by homology classes $[x]$ such that
$x=(x_1,{\cdots},x_m)\in\Sigma\w C_*(\overline K)$ satisfying the following conditions.

(1) $x_k=\emptyset\in T_0^{\sigma_k,\emptyset}(\Delta\!^{[n_k]})$ for $k\in\sigma$.

(2) $[x_k]\in H_*^{\sigma_k,\omega_k}(L_k)$ for $k\in\omega$.

(3) $x_k=\emptyset\in T_0^{\sigma_k,\emptyset}(L_k)$ for $k\in[m]{\setminus}(\sigma{\cup}\omega)$.

From  the proof of Theorem~5.4 we know that $\gamma_{L_k,\sigma_k,\omega_k}([x_k])$ is defined as follows.
For a chain $\mu=\Sigma_i\,k_i\nu_i\in\Sigma\w C_*(\Delta\!^{\omega_k})$ with $k_i\in\Bbb Z$ and $\nu_i\subset\omega_k$,
denote by $\mu'$ the chain $\Sigma_i\,k_i\nu'_i\in\Sigma\w C^*(\Delta\!^{\omega_k})$, where $\nu'_i=\omega_k{\setminus}\nu_i$.
Since $\w H_*(\Delta\!^{\omega_k})=0$, we may take a $y_k\in\Sigma\w C_*(\Delta\!^{\omega_k})$ such that $dy_k=x_k$
and so $\delta y'_k=x'_k$ in $\Sigma\w C^*(\Delta\!^{\omega_k})$.
Then $\gamma_{L_k,\sigma_k,\omega_k}([x_k])=[y'_k]$.
Similarly, replace $L_k,\sigma_k,\omega_k$ in the above proof by $\overline K,\w\sigma\!,\,\w\omega$ and we have
$$\gamma_{\overline K,\tilde\sigma,\tilde\omega}([x])=[y]=(-1)^{s}[(x'_1,{\cdots},x'_{i-1},y'_i,x'_{i+1},{\cdots},x'_m)]
\in H^{n-*-1}_{\tilde\sigma'\!,\tilde\omega}(\overline K^\circ),$$
where $x'_k=\emptyset$ for $k\notin\omega$, $i$ is any given number in $\omega$ and $s=|x'_1|+\cdots+|x'_{i-1}|$.
Note that $K_{\sigma,\omega}=\{\emptyset\}$ and $(K^\circ)_{\sigma'\!,\omega}=\partial\Delta\!^\omega$.
Let $\phi=\phi_{(-;\underline{D},\underline{C})}\,q_{(-;\underline{D},\underline{C})}$ from
$T_*^{\XX_n}({\cal Z}(-;\underline{D_*},\underline{C_*}))$ to $T_*^{\XX_m}(-)\widehat\otimes(\otimes_k H_*^{\XX_m;\XX_{n_k}}(\Delta\!^{[n_k]}\!,L_k))$
be as defined in Theorem~3.12 with dual $\phi^\circ$. Then ($A,B$ as previously defined)
$$\phi(x)=(\emptyset){\otimes}([x_1]{\otimes}{\cdots}{\otimes}[x_m])\in
T_0^{\sigma,\omega}(K){\otimes} A,$$
$$(\phi^\circ)^{-1}(y)=(\omega{\setminus}\{i\}){\otimes}([y'_1]{\otimes}{\cdots}{\otimes}[y'_m])
\in T^{|\omega|-1}_{\sigma'\!,\omega}(K^\circ){\otimes}B,$$
where $y'_k=\emptyset$ for $k\notin\omega$.
So
$\gamma_{\overline K,\tilde\sigma,\tilde\omega}=\gamma_{K,\sigma,\omega}{\otimes}\gamma_{L_1,\sigma_1,\omega_1}{\otimes}{\cdots}{\otimes}
\gamma_{L_m,\sigma_m,\omega_m}$ for $K=\Delta\!^S$.
\hfill$\Box$\vspace{3mm}

{\bf Definition~5.8} For a simplicial complex $K$ on $[m]$ such that $K\neq\{\,\}$ or $\Delta\!^{[m]}$
and a sequence $\underline{r}=(r_1,{\cdots},r_m)$ of positive integers,
$I_{(K;\underline{r})}$ is the ideal of the polynomial ring $\mak[m]=\mak[x_1,{\cdots},x_m]$
generated by all monomials $x_{i_1}^{r_{i_1}}{\cdots}x_{i_s}^{r_{i_s}}$ such that $\{i_1,{\cdots},i_s\}\notin K$.
Specifically, for $\underline{1}=(1,{\cdots},1)$, $I_{(K;\underline{1})}$
is just the {\it Stanley-Reisner face ideal} $I_K$ of $K$ on $[m]$.
\vspace{2mm}

{\bf Theorem~5.9} (Hochster Theorem) {\it For a simplicial complex $K$ on $[m]$ such that $K\neq\{\,\}$ or $\Delta\!^{[m]}$
and $\underline{r}=(r_1,{\cdots},r_m)$, we have
$$I_{(K;\underline{r})}={\cal Z}^\otimes(K^\circ;(x^{r_1}),{\cdots},(x^{r_m})),$$
$${\rm Tor}^{\mak[m]}_*(I_{(K;\underline{r})},\mak)\cong\oplus_{\sigma\in K^\circ}\, H_{*}^{\sigma,[m]\setminus\sigma}(K^\circ)
\cong\oplus_{\omega\notin K}\, H^{|\omega|-*-1}_{\emptyset,\omega}(K),
$$
where $K^\circ$ is the dual of $K$ relative to $[m]$, all the $R_k$ in the composition module is the same polynomial $\mak[x]$.
\vspace{2mm}

Proof}\, $I_{(K;\underline{r})}={\cal Z}^\otimes(K^\circ;(x^{r_1}),{\cdots},(x^{r_m}))$ is by definition.
The minimal resolution of $(x^{r_k})$ is $0\to\mak[x]\stackrel{{\cdot} x^{r_k}}{\longrightarrow}(x^{r_k})\to 0$.
So ${\rm Tor}_0^{\mak[x]}((x^{r_k}),\mak)=\mak$ and ${\rm Tor}_k^{\mak[x]}((x^{r_k}),\mak)=0$ if $k\neq 0$.
Then we have by Theorem~4.19
$$\begin{array}{l}
\quad{\rm Tor}^{\mak[m]}_*(I_{(K;\underline{r})},\mak)\vspace{2mm}\\
\cong\oplus_{\sigma\in K^\circ}\, H_{*}^{\sigma,[m]\setminus\sigma}(K^\circ)\otimes
\big(\otimes_{k\notin\sigma}{\rm  Tor}^{\mak[x]}_*((x^{r_k}),\mak)\big)\vspace{2mm}\\
\cong\oplus_{\sigma\in K^\circ}\, H_{*}^{\sigma,[m]\setminus\sigma}(K^\circ)\vspace{2mm}\\
\cong\oplus_{\omega\notin K}\, H^{|\omega|-*-1}_{\emptyset,\omega}(K),\quad({\rm by\,\,Theorem~5.4})
\end{array}$$
where $\omega=[m]{\setminus}\sigma$.
\hfill$\Box$\vspace{3mm}

{\bf Theorem~5.10} {\it Let $K$ be a simplicial complex on $[m]$, $L_k$ be a simplicial complex on $[n_k]$
and $\underline{r_k}=(r_{k,1},{\cdots},r_{k,n_k})$ for $k=1,{\cdots},m$.
Then we have
$$\begin{array}{l}
\quad{\cal Z}^\otimes(K;I_{(L_1;\underline{r_1})},{\cdots},I_{(L_m;\underline{r_m})})
=I_{({\cal Z}^*(K^\circ;L_1,{\cdots},L_m);(\underline{r_1},{\cdots},\underline{r_m}))},\vspace{4mm}\\
\quad{\rm Tor}^{\mak[n]}_*({\cal Z}^\otimes(K;I_{(L_1;\underline{r_1})},{\cdots},I_{(L_m;\underline{r_m})}),\mak)\vspace{2mm}\\
\cong\oplus_{\sigma\in K}\, H_{*}^{\sigma,[m]\setminus\sigma}(K)\otimes
\big(\otimes_{k\notin\sigma}(\oplus_{\sigma_k\in L_k^\circ} H_{*}^{\sigma_k,[n_k]\setminus\sigma_k}(L_k^\circ))\big)\vspace{2mm}\\
\cong\oplus_{\omega\notin K^\circ}\, H^{|\omega|-*-1}_{\emptyset,\omega}(K^\circ)\otimes
\big(\otimes_{k\in\omega}(\oplus_{\omega\notin L_k} H^{|\omega_k|-*-1}_{\emptyset,\omega_k}(L_k))\big),
\end{array}$$
where $[n]=[n_1]\sqcup{\cdots}\sqcup[n_m]$, $\omega=[m]{\setminus}\sigma$ and $\omega_k=[n_k]{\setminus}\sigma_k$.
\vspace{2mm}

Proof}\, Theorem~4.9 also holds for module pairs $(U_k,C_k)$. So we have
$$\begin{array}{l}
\quad{\cal Z}^\otimes(K;I_{(L_1;\underline{r_1})},{\cdots},I_{(L_m;\underline{r_m})})\vspace{2mm}\\
={\cal Z}^\otimes(K;{\cal Z}^\otimes(L_1^\circ;(x^{r_{1,1}}),{\cdots},(x^{r_{1,n_1}})),{\cdots},{\cal Z}^\otimes(L_m^\circ;(x^{r_{m,1}}),{\cdots},(x^{r_{m,n_m}})))\vspace{2mm}\\
={\cal Z}^\otimes({\cal Z}^*(K;L_1^\circ{\cdots},L_m^\circ);(x^{r_{1,1}}),{\cdots},(x^{r_{1,n_1}}),
{\cdots},(x^{r_{m,1}}),{\cdots},(x^{r_{m,n_m}}))\vspace{2mm}\\
=I_{({\cal Z}^*(K^\circ;L_1{\cdots},L_m);(\underline{r_1},{\cdots},\underline{r_m}))}.
\end{array}$$
So by Theorem~4.19.
$$\begin{array}{l}
\quad{\rm Tor}^{\mak[n]}_*({\cal Z}^\otimes(K;I_{(L_1;\underline{r_1})},{\cdots},I_{(L_m;\underline{r_m})}),\mak)\vspace{2mm}\\
\cong\oplus_{\sigma\in K}\, H_{*}^{\sigma,[m]\setminus\sigma}(K)\otimes
\big(\otimes_{k\notin\sigma} {\rm Tor}_{*}^{\mak[n_k]}(I_{(L_k;\underline{r_k})},\mak)\big)\vspace{2mm}\\
\cong\oplus_{\sigma\in K}\, H_{*}^{\sigma,[m]\setminus\sigma}(K)\otimes
\big(\otimes_{k\notin\sigma}(\oplus_{\sigma_k\in L_k^\circ} H_{*}^{\sigma_k,[n_k]\setminus\sigma_k}(L_k^\circ))\big)\vspace{2mm}\\
\cong\oplus_{\omega\notin K^\circ}\, H^{|\omega|-*-1}_{\emptyset,\omega}(K^\circ)\otimes
\big(\otimes_{k\in\omega}(\oplus_{\omega\notin L_k} H^{|\omega_k|-*-1}_{\emptyset,\omega_k}(L_k))\big).
\end{array}$$
From another point of view,
$$\begin{array}{l}
\quad{\rm Tor}^{\mak[n]}_*(I_{({\cal Z}^*(K^\circ;L_1,{\cdots},L_m);(\underline{r_1},{\cdots},\underline{r_m}))},\mak)\vspace{2mm}\\
\cong\oplus_{\tilde\sigma\in {\cal Z}^*(K;L_1^\circ,{\cdots},L_m^\circ)}\,
H_{*}^{\tilde\sigma,[n]\setminus\tilde\sigma}({\cal Z}^*(K;L_1^\circ,{\cdots},L_m^\circ))\vspace{2mm}\\
\cong\oplus_{\sigma\in K}\, H_{*}^{\sigma,[m]\setminus\sigma}(K)\otimes
\big(\otimes_{k\notin\sigma}(\oplus_{\sigma_k\in L_k^\circ} H_{*}^{\sigma_k,[n_k]\setminus\sigma_k}(L_k^\circ))\big).
\end{array}$$

The two results coincide.
\hfill$\Box$\vspace{3mm}

{\bf Definition~5.11} Let $(\underline{X},\underline{A})=\{(X_k,A_k)\}_{k=1}^m$ be a sequence of topological
pairs satisfying the following conditions.

1) Each homology homomorphism $i_k\colon H_*(A_k)\to H_*(X_k)$ induced by inclusion is split.
When the homology is taken over a field, this condition can be canceled.

2) Each $X_k$ is a closed orientable manifold of dimension $r_k$.

3) Each $A_k$ is a proper compact polyhedron subspace of $X_k$.

Denote by $(\underline{X},\underline{A}^c)=\{(X_k,A_k^c)\}_{k=1}^m$ with $A_k^c=X_k{\setminus}A_k$.
Then we have the {\it complement duality isomorphism}
$$\gamma_{(\underline{X},\underline{A})}\colon H_*^{\XX_m}(\underline{X},\underline{A})\to H^{\,*}_{\!\XX_m}(\underline{X},\underline{A}^c),$$
$$\gamma_{(\underline{X},\underline{A})}^\circ\colon H^{\,*}_{\!\XX_m}(\underline{X},\underline{A})\to H_*^{\XX_m}(\underline{X},\underline{A}^c),$$
defined as follows.
We have the following commutative diagram
\[\begin{array}{cccccccccc}
\cdots\longrightarrow \hspace{-1.5mm}&\hspace{-1.5mm} {\scriptstyle H_n(A_k)} \hspace{-1.5mm}&\hspace{-1.5mm} \stackrel{i_k}{\longrightarrow} \hspace{-1.5mm}&\hspace{-1.5mm} {\scriptstyle H_n(X_k)}
&\hspace{-1.5mm} \stackrel{j_k}{\longrightarrow} \hspace{-1.5mm}&\hspace{-1.5mm} {\scriptstyle H_n(X_k,A_k)} \hspace{-1.5mm}&\hspace{-1.5mm}\stackrel{\partial_k}{\longrightarrow} \hspace{-1.5mm}&\hspace{-1.5mm}
{\scriptstyle H_{n-1}(A_k)}\hspace{-1.5mm}&\hspace{-1.5mm} \longrightarrow \cdots\,\vspace{1mm}\\
&^{\alpha_k}\downarrow\quad&&^{\gamma_k}\downarrow\quad&&^{\beta_k}\downarrow\quad&&^{\alpha_k}\downarrow\quad\quad&&\\
\cdots\longrightarrow \hspace{-1.5mm}&\hspace{-1.5mm} {\scriptstyle H^{r_k-n}(X_k,A_k^c)} \hspace{-1.5mm}&\hspace{-1.5mm} \stackrel{q_k^\circ}{\longrightarrow} \hspace{-1.5mm}&\hspace{-1.5mm}{\scriptstyle H^{r_k-n}(X_k)}
\hspace{-1.5mm}&\hspace{-1.5mm} \stackrel{p_k^\circ}{\longrightarrow} \hspace{-1.5mm}&\hspace{-1.5mm}
{\scriptstyle H^{r_k-n}(A_k^c)} \hspace{-1.5mm}&\hspace{-1.5mm}\stackrel{\partial_k^\circ}{\longrightarrow} \hspace{-1.5mm}&\hspace{-1.5mm}
{\scriptstyle H^{r_k-n+1}(X_k,A^c_k)} \hspace{-1.5mm}&\hspace{-1.5mm}\longrightarrow\cdots,
  \end{array}\]
where $\alpha_k,\beta_k$ are the Alexander duality isomorphisms
and $\gamma_k$ is the Poncar\'{e} duality isomorphism.
So we have the following group isomorphisms
$$\begin{array}{rl}
(\partial_k^\circ)^{-1}\alpha_k\colon&\,\,\,{\rm ker}\,i_k\,\,\,\stackrel{\cong}
{\longrightarrow}\,\,{\rm coker}\,p_k^\circ,\vspace{1mm}\\
\gamma_ki_k\colon&\,\,\,{\rm coim}\,i_k\,\,\,\,\stackrel{\cong}{\longrightarrow}
\,\,\,{\rm ker}\,p_k^\circ,\vspace{1mm}\\
p_k^\circ\gamma_k\colon&{\rm coker}\,i_k\stackrel{\cong}{\longrightarrow}
\,\,\,\,{\rm im}\,p_k^\circ.
\end{array}$$
Define $\phi_k\colon H_*^{\XX}\!(X_k,A_k)\to H^{\,*}_{\!\XX}(X_k,A_k^c)$ (not degree keeping!)
to be the direct sum of the above three isomorphisms
and $\gamma_{(\underline{X},\underline{A})}$ to be $\phi_1{\otimes}{\cdots}{\otimes}\phi_m$.
The dual case is similar.

It is obvious that $\gamma_{(\underline{X},\underline{A})}=\oplus_{(\sigma\!,\,\omega)\in\XX_m}\,\gamma_{\sigma\!,\,\omega}$,
$\gamma_{(\underline{X},\underline{A})}^\circ=\oplus_{(\sigma\!,\,\omega)\in\XX_m}\,\gamma_{\sigma\!,\,\omega}^\circ$ with
$$\gamma_{\sigma\!,\,\omega}\colon H_*^{\sigma\!,\,\omega}(\underline{X},\underline{A})\to H^{r-|\omega|-*}_{\sigma'\!,\,\omega}(\underline{X},\underline{A}^c),$$
$$\gamma^\circ_{\sigma\!,\,\omega}\colon H^*_{\sigma\!,\,\omega}(\underline{X},\underline{A})\to H_{r-|\omega|-*}^{\sigma'\!,\,\omega}(\underline{X},\underline{A}^c),$$
where $\sigma'=[m]{\setminus}(\sigma{\cup}\omega)$, $r=r_1{+}{\cdots}{+}r_m$.

Define $\overline\gamma_{(\underline{X},\underline{A})}=\oplus_{(\sigma\!,\,\omega)\in\LL_m}\,\gamma_{\sigma\!,\,\omega}$,
$\overline\gamma_{(\underline{X},\underline{A})}^{\,\circ}=\oplus_{(\sigma\!,\,\omega)\in\LL_m}\,\gamma_{\sigma\!,\,\omega}^\circ$ with
$\LL_m$ as in Definition~5.5.
\vspace{3mm}

{\bf Definition~5.12} For a polyhedral product space $M={\cal Z}(K;\underline{X},\underline{A})$,
let $i\colon H_*(M)\to H_*(\w X)$ and $i^\circ\colon H^*(\w X)\to H^*(M)$
be the singular (co)homology homomorphism induced by the inclusion map from $M$ to
$\w X=X_1{\times}{\cdots}{\times}X_m$.
From the long exact exact sequences
$$\begin{array}{l}
\scriptstyle{\cdots\,\,\,\,\longrightarrow\,\,\,\, H_n(M)\,\,\,\,\stackrel{i}{\longrightarrow}\,\,\,\,H_n(\w X)
\,\,\,\,\stackrel{j}{\longrightarrow}\,\,\,\, H_n(\w X,M)\,\,\,\,\stackrel{\partial}{\longrightarrow}\,\,\,\, H_{n-1}(M)
\,\,\,\,\longrightarrow\,\,\,\,\cdots}\\\\
\scriptstyle{\cdots\,\,\,\,\longrightarrow\,\,\,\, H^{n-1}(M)\,\,\,\,\stackrel{\partial^\circ}{\longrightarrow}\,\,\,\, H^n(\w X,M)
\,\,\,\,\stackrel{j^\circ}{\longrightarrow}\,\,\,\, H^n(\w X)\,\,\,\,\stackrel{i^\circ}{\longrightarrow}\,\,\,\, H^{n}(M)
\,\,\,\, \longrightarrow\,\,\,\,\cdots}
  \end{array}$$
we define
$$\hat H_*(M)\!=\!{\rm coim}\,i,\,\overline H_*(M)\!=\!{\rm ker}\,i,\,
\hat H_*(\w X{,}M)\!=\!{\rm im}\,j,\,\overline H_*(\w X{,}M)\!=\!{\rm coker}\,j,$$
$$\hat H^*(M)\!=\!{\rm im}\,i^\circ,\,\overline H^*(M)\!=\!{\rm coker}\,i^\circ,\,
\hat H^*(\w X{,}M)\!=\!{\rm coim}\,j^\circ,\,H^*(\w X{,}M)\!=\!{\rm ker}\,j^\circ.$$

{\bf Theorem~5.13} {\it Suppose the $M$ in Definition~5.12 is homology split by Definition~4.3.
Then we have the following group decompositions
$$\begin{array}{cc}
H_*(M)=\hat H_*(M){\oplus}\overline H_*(M),&
H_*(\w X,M)=\hat H_*(\w X,M){\oplus}\overline H_*(\w X,M),\vspace{2mm}\\
H^*(M)=\hat H^*(M){\oplus}\overline H^*(M),&
H^*(\w X,M)=\hat H^*(\w X,M){\oplus}\overline H^*(\w X,M)
\end{array}$$
with
$$\begin{array}{c}
\overline H_{*+1}(\w X,M)\cong \overline H_*(M)\cong H_*^{\LL_m}(K)\widehat\otimes H_*^{\LL_m}(\underline{X},\underline{A}),\vspace{2mm}\\
\overline H^{*+1}(\w X,M)\cong \overline H^*(M)\cong H^{\,*}_{\!\LL_m}(K)\otimes H^{\,*}_{\!\LL_m}(\underline{X},\underline{A}),\vspace{2mm}\\
\hat H_*(M)\cong H_*^{K}(\underline{X},\underline{A}),\,\,\hat H_*(\w X,M)\cong H_*^{K^c}(\underline{X},\underline{A}),\vspace{2mm}\\
\hat H^*(M)\cong H^{\,*}_{\!K}(\underline{X},\underline{A}),\,\,\hat H^*(\w X,M)\cong H^{\,*}_{\!K^c}(\underline{X},\underline{A}),
\end{array}$$
where $\LL_m$ is as in Definition~5.5, $K$ and its complement $K^c$ are regarded as subsets of $\XX_m$ defined by
$K=\{(\sigma,\emptyset)\in\XX_m\,|\,\sigma\in K\}$, $K^c=\{(\sigma,\emptyset)\in\XX_m\,|\,\sigma\notin K\}$.
\vspace{2mm}

\it Proof}\, By definition, $i=\oplus_{(\sigma\!,\,\omega)\in\XX_m}\,i_{\sigma\!,\,\omega}$ with\vspace{2mm}\\
\hspace*{15mm}$i_{\sigma\!,\,\omega}\colon H_*^{\sigma\!,\,\omega}(K){\otimes}H_*^{\sigma\!,\,\omega}(\underline{X},\underline{A})
\stackrel{i{\otimes}1}{-\!\!\!\longrightarrow}
H_*^{\sigma\!,\,\omega}(\Delta\!^{[m]}){\otimes}H_*^{\sigma\!,\,\omega}(\underline{X},\underline{A})$,\vspace{1mm}\\
where $i$ is induced by inclusion and $1$ is the identity.
So
$$\hat H_*(M)=\oplus_{\sigma\in K}
H_*^{\sigma,\emptyset}(K){\otimes}H_*^{\sigma,\emptyset}(\underline{X},\underline{A})
\cong\oplus_{\sigma\in K} H_*^{\sigma,\emptyset}(\underline{X},\underline{A}),$$
$$\overline H_*(M)
=\oplus_{(\sigma\!,\,\omega)\in\LL_m}\,
H_*^{\sigma\!,\,\omega}(K)\otimes H_*^{\sigma\!,\,\omega}(\underline{X},\underline{A}).$$

The relative group case is similar.
\hfill $\Box$\vspace{3mm}

{\bf Theorem~5.14} {\it For $M={\cal Z}(K;\underline{X},\underline{A})$ such that $(\underline{X},\underline{A})$
satisfies the condition of Definition~5.11,
the Alexander duality isomorphisms
$$\alpha\colon H_*(M)\to H^{r-*}(\w X,M^c),\,\quad
\alpha^\circ\colon H^*(M)\to H_{r-*}(\w X,M^c)$$
have direct sum decomposition $\,\alpha=\hat\alpha\oplus\overline\alpha$,\,\,$\alpha^\circ=\hat\alpha^\circ\oplus\overline\alpha^\circ$,
where
$$\hat\alpha\colon \hat H_*(M)\to \hat H^{r-*}(\w X,M^c),\,\quad
\overline\alpha\colon \overline H_*(M)\to \overline H^{r-*}(\w X,M^c)\cong\overline H^{r-*-1}(M^c),$$
$$\hat\alpha^\circ\colon \hat H^*(M)\to \hat H_{r-*}(\w X,M^c),\,\quad
\overline\alpha^\circ\colon \overline H^*(M)\to \overline H_{r-*}(\w X,M^c)\cong\overline H_{r-*-1}(M^c)\,$$
are as follows. Identify all the above groups with the corresponding indexed groups in Theorem~5.13.
Then
$$\hat\alpha\,\,=\oplus_{\sigma\in K}\,\gamma_{\sigma,\emptyset},\,\quad
\overline\alpha\,\,=\overline\gamma_K\,\widehat\otimes\,\overline\gamma_{(\underline{X},\underline{A})}
=\oplus_{(\sigma\!,\omega)\in\LL_m}\,\,\gamma_{K,\sigma\!,\,\omega}{\otimes}\gamma_{\sigma\!,\,\omega},$$
$$\hat\alpha^\circ=\oplus_{\sigma\in K}\,\gamma_{\sigma,\emptyset}^\circ,\,\quad
\overline\alpha^\circ\,\,=\overline\gamma_K^{\,\circ}\,\widehat\otimes\,\overline\gamma_{(\underline{X},\underline{A})}^{\,\circ}
=\oplus_{(\sigma\!,\omega)\in\LL_m}\,\,\gamma_{K,\sigma\!,\,\omega}^\circ{\otimes}\gamma_{\sigma\!,\,\omega}^\circ,$$
where $\gamma_{-}$, $\gamma_{-}^\circ$ are as in Theorem~5.4 and Definition~5.11.
}\vspace{2mm}

{\it Proof}\, Denote by $\alpha=\alpha_M$, $\hat\alpha=\hat\alpha_M$,
$\overline\alpha=\overline\alpha_M$.
Then for $M={\cal Z}(K;\underline{X},\underline{A})$ and $N={\cal Z}(L;\underline{X},\underline{A})$,
we have the following commutative diagrams of exact sequences
\[\,\,\,\,\begin{array}{ccccccc}
{\scriptstyle\cdots\,\,\longrightarrow}\!\!&\!\!{\scriptstyle H_k(M{\cap}N)}\!\!&\!\!{\scriptstyle\longrightarrow}
\!\!&\!\!{\scriptstyle H_k(M){\oplus}H_k(N)}
\!\!&\!\!{\scriptstyle\longrightarrow}\!\!&\!\!{\scriptstyle H_k(M{\cup}N)}\!\!&\!\!{\scriptstyle\longrightarrow\,\,\cdots}\vspace{1mm}\\
\!\!&\!\!^{\alpha_{M\cap N}}\downarrow\quad\quad\!\!&\!\!\!\!&\!\!^{\alpha_M\oplus\alpha_N}\downarrow\quad\quad\quad\!\!&\!\!\!\!&\!\!^{\alpha_{M\cup N}}\downarrow\quad\quad\!\!&\!\!\\
{\scriptstyle\cdots\,\,\longrightarrow}\!\!&\!\!{\scriptstyle H^{r-k}(\w X,(M{\cap}N)^c)}\!\!&\!\!{\scriptstyle\longrightarrow}\!\!&\!\!
{\scriptstyle H^{r-k}(\w X,M^c){\oplus}H^{r-k}(\w X,N^c)}
\!\!&\!\!{\scriptstyle\longrightarrow}\!\!&\!\!{\scriptstyle H^{r-k}(\w X,(M{\cup}N)^c)}\!\!&\!\!{\scriptstyle\longrightarrow\,\,\cdots}
  \end{array}\quad\,\,\,\,(1)\vspace{1mm}\]
\[\begin{array}{ccccccc}
{\scriptstyle0\quad\longrightarrow}\!\!&\!\!{\scriptstyle \hat H_k(M{\cap}N)}\!\!&\!\!{\scriptstyle\longrightarrow}\!\!&\!\!
{\scriptstyle\hat H_k(M){\oplus}\hat H_k(N)}
\!\!&\!\!{\scriptstyle\longrightarrow}\!\!&\!\!
{\scriptstyle\hat H_k(M{\cup}N;\underline{X},\underline{A})}\!\!&\!\!{\scriptstyle\longrightarrow\quad 0}
\vspace{1mm}\\
\!\!&\!\!^{\hat\alpha_{M\cap N}}\downarrow\quad\quad\!\!&\!\!\!\!&\!\!^{\hat\alpha_M\oplus\hat\alpha_N}\downarrow\quad\quad\quad\!\!&\!\!\!\!&\!\!^{\hat\alpha_{M\cup N}}\downarrow\quad\quad\!\!&\!\!\\
{\scriptstyle0\quad\longrightarrow}\!\!&\!\!{\scriptstyle\hat H^{r-k}(\w X,(M{\cap}N)^c)}\!\!&\!\!{\scriptstyle\longrightarrow}\!\!&\!\!
{\scriptstyle\hat H^{r-k}(\w X,M^c){\oplus}\hat H^{r-k}(\w X,N^c)}
\!\!&\!\!{\scriptstyle\longrightarrow}\!\!&\!\!{\scriptstyle\hat H^{r-k}(\w X,(M{\cup}N)^c)}\!\!&\!\!{\scriptstyle\longrightarrow\quad 0}
  \end{array}\quad(2)\vspace{1mm}\]
For $(\sigma\!,\,\omega)\in\LL_m$,
$A=H_l^{\sigma\!,\,\omega}(\underline{X},\underline{A})$,
$B=H^{r-|\omega|-l}_{\sigma\!,\,\omega}(\underline{X},\underline{A}^c)$,
$\gamma_1=\gamma_{K{\cap}L,\sigma\!,\omega}$, $\gamma_2=\gamma_{K,\sigma\!,\omega}{\oplus}\gamma_{L,\sigma\!,\omega}$,
$\gamma_3=\gamma_{K\cup L,\sigma\!,\omega}$, we have the commutative diagram
\[\begin{array}{ccccccc}
{\scriptstyle\cdots\,\,\longrightarrow}\!\!&\!\!{\scriptstyle H_k^{\sigma\!,\omega}(K{\cap}L){\otimes}A}\!\!&\!\!{\scriptstyle\longrightarrow}\!\!&\!\!
{\scriptstyle (H_k^{\sigma\!,\omega}(K){\oplus} H_k^{\sigma\!,\omega}(L)){\otimes}A}\!\!&\!\!
{\scriptstyle\longrightarrow}\!\!&\!\!{\scriptstyle H_k^{\sigma\!,\omega}(K{\cup}L){\otimes}A}\!\!&\!\!{\scriptstyle\longrightarrow\,\,\cdots}\vspace{1mm}\\
\!\!&\!\!^{\gamma_1\otimes\gamma_{\sigma\!,\,\omega}}
\downarrow\quad\quad\!\!&\!\!\!\!&\!\!^{\gamma_2\otimes\gamma_{\sigma\!,\,\omega}}
\downarrow\quad\quad\quad\!\!&\!\!\!\!&\!\!
^{\gamma_{3}{\otimes}\gamma_{\sigma\!,\,\omega}}\downarrow\quad\quad\!\!&\!\!\\
{\scriptstyle\cdots\,\,\longrightarrow}\!\!&\!\!
{\scriptstyle H^{k'}_{\tilde\sigma\!,\omega}((K\cap L)^\circ){\otimes}B}\!\!&\!\!
{\scriptstyle\longrightarrow}\!\!&\!\!{\scriptstyle (H^{k'}_{\tilde\sigma\!,\omega}(K^\circ)
{\oplus}H^{k'}_{\tilde\sigma\!,\omega}(L^\circ)){\otimes}B}\!\!&\!\!{\scriptstyle\longrightarrow}\!\!&\!\!
{\scriptstyle H^{k'}_{\tilde\sigma\!,\omega}((K\cup L)^\circ){\otimes}B}\!\!&\!\!
{\scriptstyle\longrightarrow\,\,\cdots}
  \end{array}\]
where $k'=|\omega|{-}k{-}1$. The direct sum of all the above diagrams is the following diagram.
\[\begin{array}{ccccccc}
{\scriptstyle\cdots\,\,\longrightarrow}\!\!&\!\!{\scriptstyle\overline H_k(M{\cap}N)}\!\!&\!\!{\scriptstyle\longrightarrow}\!\!&\!\!
{\scriptstyle\overline H_k(M){\oplus}\overline H_k(N)}\!\!&\!\!{\scriptstyle\longrightarrow}\!\!&\!\!
{\scriptstyle\overline H_k(M{\cup}N)}\!\!&\!\!{\scriptstyle\longrightarrow\,\,\cdots}\vspace{1mm}\\
\!\!&\!\!^{\overline\alpha_{M\cap N}}\downarrow\quad\quad\!\!&\!\!\!\!&\!\!^{\overline\alpha_M\oplus\overline\alpha_N}\downarrow\quad\quad\quad\!\!&\!\!\!\!&\!\!^{\overline\alpha_{M\cup N}}\downarrow\quad\quad\!\!&\!\!\\
{\scriptstyle\cdots\,\,\longrightarrow}\!\!&\!\!{\scriptstyle\overline H^{r-k}(\w X,(M{\cap}N)^c)}\!\!&\!\!
{\scriptstyle\longrightarrow}\!\!&\!\!{\scriptstyle\overline H^{r-k}(\w X,M^c){\oplus}\overline H^{r-k}(\w X,N^c)}
\!\!&\!\!{\scriptstyle\longrightarrow}\!\!&\!\!{\scriptstyle\overline H^{r-k}(\w X,(M{\cup}N)^c)}\!\!&\!\!
{\scriptstyle\longrightarrow\,\,\cdots}
  \end{array}\quad(3)\]

(1), (2) and (3) imply that
if the theorem holds for $M$ and $N$ and $M{\cap}N$, then it holds for $M{\cup}N$.
So by induction on the number of maximal simplices of $K$,
we only need prove the theorem for the special case that $K$ has only one maximal simplex.

Now we prove the theorem for $M={\cal Z}(\Delta\!^S;\underline{X},\underline{A})$ with $S\subset[m]$.
Then
$$M=Y_1{\times}{\cdots}{\times}Y_m,\quad Y_k=\left\{\begin{array}{cc}
X_k & {\rm if}\,\,k\in S, \\
A_k & {\rm if}\,\,k\notin S.
\end{array}\right.$$
So $(\w X,M^c)=(X_1,Y_1^c){\times}{\cdots}{\times}(X_m,Y_m^c)$.

Construct a degree-keeping homomorphism $\w\phi_k$ similar to the $\phi_k$ in Definition~5.11.
From the commutative diagram
\[\begin{array}{cccccccccc}
\cdots\longrightarrow \hspace{-1.5mm}&\hspace{-1.5mm} {\scriptstyle H_n(A_k)} \hspace{-1.5mm}&\hspace{-1.5mm} \stackrel{i_k}{\longrightarrow} \hspace{-1.5mm}&\hspace{-1.5mm} {\scriptstyle H_n(X_k)}
&\hspace{-1.5mm} \stackrel{j_k}{\longrightarrow} \hspace{-1.5mm}&\hspace{-1.5mm} {\scriptstyle H_n(X_k,A_k)} \hspace{-1.5mm}&\hspace{-1.5mm}\stackrel{\partial_k}{\longrightarrow} \hspace{-1.5mm}&\hspace{-1.5mm}
{\scriptstyle H_{n-1}(A_k)}\hspace{-1.5mm}&\hspace{-1.5mm} \longrightarrow \cdots\,\vspace{1mm}\\
&^{\alpha_k}\downarrow\quad&&^{\gamma_k}\downarrow\quad&&^{\beta_k}\downarrow\quad&&^{\alpha_k}\downarrow\quad\quad&&\\
\cdots\longrightarrow \hspace{-1.5mm}&\hspace{-1.5mm} {\scriptstyle H^{r_k-n}(X_k,A_k^c)} \hspace{-1.5mm}&\hspace{-1.5mm} \stackrel{q_k^\circ}{\longrightarrow} \hspace{-1.5mm}&\hspace{-1.5mm}{\scriptstyle H^{r_k-n}(X_k)}
\hspace{-1.5mm}&\hspace{-1.5mm} \stackrel{p_k^\circ}{\longrightarrow} \hspace{-1.5mm}&\hspace{-1.5mm}
{\scriptstyle H^{r_k-n}(A_k^c)} \hspace{-1.5mm}&\hspace{-1.5mm}\stackrel{\partial_k^\circ}{\longrightarrow} \hspace{-1.5mm}&\hspace{-1.5mm}
{\scriptstyle H^{r_k-n+1}(X_k,A^c_k)} \hspace{-1.5mm}&\hspace{-1.5mm}\longrightarrow\cdots,
  \end{array}\]
we have group isomorphisms
$$\begin{array}{rl}
\alpha_k\colon&\,\,\,{\rm ker}\,i_k\,\,\,\stackrel{\cong}{\longrightarrow}\,\,{\rm im}\,\partial_k^\circ\subset H^*(X_k,A_k^c),\vspace{1mm}\\
\gamma_k\colon&\,{\rm coim}\,i_k\,\,\stackrel{\cong}{\longrightarrow}\,\,\,{\rm coim}\,q_k^\circ\subset H^*(X_k,A_k^c),\vspace{1mm}\\
p_k^\circ\gamma_k\colon&{\rm coker}\,i_k\stackrel{\cong}{\longrightarrow}\,\,\,\,{\rm coim}\,p_k^\circ\subset H^*(X_k).
\end{array}$$
Define $\w H^{\,*}_{\!\XX}(X_k,A_k^c)$ to be the direct sum of the above groups on the right side given by
$\w H^*_{\mathpzc i}(X_k,A_k^c)\cong H^{*}_{\mathpzc i}(X_k,A_k^c)$,
$\w H^*_{\mathpzc e}(X_k,A_k^c)\cong H^*_{\mathpzc e}(X_k,A_k^c)$ and
$\w H^{*-1}_{\mathpzc n}(X_k,A_k^c)\cong H^*_{\mathpzc n}(X_k,A_k^c)$.
Let $\w\phi_k\colon H_*^\XX\!(X_k,A_k)\to\w H^{\,*}_{\!\XX}(X_k,A_k^c)$ be the direct sum of the above three isomorphisms.
Then $\w\phi_1{\otimes}{\cdots}{\otimes}\w\phi_m=\oplus_{(\sigma\!,\omega)\in\XX_m}\,\w\gamma_{\sigma\!,\omega}$
and $\tilde\gamma_{\sigma\!,\,\emptyset}=\gamma_{\sigma\!,\,\emptyset}$.

We have the following commutative diagram
\[\begin{array}{ccc}
{\scriptstyle H_*(M)}
&\stackrel{\alpha_M}{-\!\!\!-\!\!\!-\!\!\!-\!\!\!-\!\!\!\longrightarrow}&
{\scriptstyle H^{r-*}(\w X,M^c)}\,\vspace{1mm}\\
\|&&\|\\
{\scriptstyle H_*(Y_1){\otimes}{\cdots}{\otimes}H_*(Y_m)}
&\stackrel{\alpha_M}{-\!\!\!-\!\!\!-\!\!\!-\!\!\!-\!\!\!\longrightarrow}&
{\scriptstyle H^{r_1-*}(X_1,Y_1^c){\otimes}{\cdots}{\otimes}H^{r_m-*}(X_m,Y_m^c)}\,\vspace{1mm}\\
\|\wr&&\|\wr\\
\oplus_{\sigma\subset S,\,\omega{\cap}S=\emptyset}\,{\scriptstyle H_*^{\sigma\!,\,\omega}(\underline{X},\underline{A})}
&\stackrel{\oplus\,\tilde\gamma_{\sigma\!,\omega}}{-\!\!\!-\!\!\!-\!\!\!-\!\!\!-\!\!\!\longrightarrow}&
\oplus_{\sigma'\subset S,\,\omega{\cap}S=\emptyset}\,{\scriptstyle
\w H^{r-*}_{\sigma'\!,\,\omega}(\underline{X},\underline{A}^c)}\vspace{1mm}\\
\cap&&\cap\\
{\scriptstyle H_*^{\XX}\!(X_1,A_1){\otimes}{\cdots}{\otimes}H_*^{\XX}\!(X_m,A_m)}
&\stackrel{\tilde\phi_1\otimes\cdots\otimes\tilde\phi_m}
{-\!\!\!-\!\!\!-\!\!\!-\!\!\!-\!\!\!\longrightarrow}&
{\scriptstyle \w H^{\,*}_{\!\XX}(X_1,A_1^c){\otimes}{\cdots}{\otimes}\w H^{\,*}_{\!\XX}(X_m,A_m^c)},
\end{array}\quad(4)\]
where ${\scriptstyle\w H^{\,*}_{\!\XX}(X_1,A_1^c){\otimes}{\cdots}{\otimes}\w H^{\,*}_{\!\XX}(X_m,A_m^c)=
\oplus_{(\sigma\!,\,\omega)\in\XX_m}\w H^{*}_{\sigma\!,\,\omega}(\underline{X},\underline{A}^c)}$ and $\sigma'=[m]{\setminus}(\sigma{\cup}\omega)$.

For $\sigma\subset S$, $\omega{\cap}S=\emptyset$, $\omega\neq\emptyset$,
we have
\[\begin{array}{ccc}
{\scriptstyle H_0^{\sigma\!,\,\omega}(\Delta\!^S){\otimes}H_*^{\sigma\!,\,\omega}(\underline{X},\underline{A})}
&\stackrel{\gamma_{\Delta\!^S\!,\sigma\!,\,\omega}\otimes\gamma_{\sigma\!,\,\omega}}
{-\!\!\!-\!\!\!-\!\!\!-\!\!\!-\!\!\!\longrightarrow}&
{\scriptstyle H^{|\omega|-1}_{\sigma'\!,\,\omega}((\Delta\!^S)^\circ){\otimes}
H^{r-|\omega|-*}_{\sigma'\!,\,\omega}(\underline{X},\underline{A}^c)\quad(\subset\overline H^{r-*-1}(M^c))}
\,\vspace{1mm}\\
\|\wr&&\hspace*{-20mm}\|\wr\\
{\scriptstyle H_*^{\sigma\!,\,\omega}(\underline{X},\underline{A})}
&\stackrel{\w\gamma_{\sigma\!,\,\omega}}{-\!\!\!-\!\!\!-\!\!\!-\!\!\!-\!\!\!\longrightarrow}&
{\scriptstyle
\w H^{r-*}_{\sigma'\!,\,\omega}(\underline{X},\underline{A}^c)\quad(\subset\overline H^{r-*}(\w X,M^c))}.
\end{array}\]
So with the identification of the theorem, the third row of (4) is the direct sum
of $\overline\alpha_M=\oplus_{\sigma\subset S,\,\omega{\cap}S=\emptyset,\,\omega\neq\emptyset}\,
\gamma_{\Delta\!^S\!,\sigma\!,\,\omega}{\otimes}\gamma_{\sigma\!,\,\omega}$ and
$\hat\alpha_M=\oplus_{\sigma\subset S}\,\gamma_{\sigma,\emptyset}$ ($\w\gamma_{\sigma,\emptyset}=\gamma_{\sigma,\emptyset}$).
So $\alpha_M=\hat\alpha_M{\oplus}\overline\alpha_M$.
\hfill $\Box$\vspace{3mm}

{\bf Example~5.15}\, Regard $S^{r+1}$ as one-point compactification of $\Bbb R^{r+1}$.
Then for $q\leqslant r$, the standard  space pair $(S^{r+1},S^q)$ is given by\vspace{1mm}\\
\hspace*{2mm}$S^q=\{(x_1,{\cdots},x_{r+1})\in\Bbb R^{r+1}\subset S^{r+1}
\,|\,x^2_1{+}{\cdots}{+}x_{q+1}^2=1,\,x_i=0,\,\,{\rm if}\,\,i>q{+}1\}.$

Let $M={\cal Z}_K\Big(\!
\begin{array}{ccc}
\scriptstyle{r_1{+}1} \!&\!\scriptstyle{\cdots}\!&\!\scriptstyle{r_m{+}1}\\
\scriptstyle{ q_1}   \!&\!\scriptstyle{\cdots}\!&\! \scriptstyle{q_m}\end{array}\!\Big)
={\cal Z}(K;\underline{X},\underline{A})$ be the polyhedral product space
such that $(X_k,A_k)=(S^{r_k+1},S^{q_k})$.
Since  $S^{r-q}$ is a deformation retract of $S^{r+1}{\setminus}S^q$,
the complement space $M^c={\cal Z}(K^\circ;\underline{X},\underline{A}^c)$
is homotopy equivalent to ${\cal Z}_{K^\circ}\Big(
\begin{array}{ccc}
\scriptstyle{r_1+\,1\,} &\scriptstyle{\cdots}&\scriptstyle{r_m+\,1\,}\\
\scriptstyle{r_1{-}q_1}   &\scriptstyle{\cdots}& \scriptstyle{r_m{-}q_m}\end{array}\Big)$.\vspace{1mm}

By the computation of Example~4.7, regardless of degree,
$\overline H_*(M)\cong H_*^{\LL_m}(K)
=\oplus_{(\sigma\!,\,\omega)\in\LL_m}H_*^{\sigma\!,\,\omega}(K)$ and
the direct sum $\oplus_{(\sigma\!,\,\omega)\in\LL_m}\gamma_{K,\sigma\!,\omega}$
is just the isomorphism
$\overline H_*(M)\cong\overline H\,^{r-*-1}(M^c)$.

Specifically, ${\cal Z}(K;S^{2n+1},S^n)={\cal Z}_K\Big(\!
\begin{array}{ccc}
\scriptstyle{2n{+}1} \hspace{-1mm}&\hspace{-1mm}\scriptstyle{\cdots}\hspace{-1mm}&\hspace{-1mm}\scriptstyle{2n{+}1}\\
\scriptstyle{n}   \hspace{-1mm}&\hspace{-1mm}\scriptstyle{\cdots}\hspace{-1mm}&\hspace{-1mm} \scriptstyle{n}\end{array}\!\Big)$.
Then we have $$\overline H_*({\cal Z}(K;S^{2n+1},S^n))\cong
\overline H\,^{(2n+1)m-*-1}({\cal Z}(K^\circ;S^{2n+1},S^n)).$$

\section{Diagonal Tensor Product of (Co)algebras}\vspace{3mm}

\hspace*{5.5mm}{\bf Definition~6.1} An {\it indexed coalgebra}, or a {\it coalgebra indexed by $\Lambda$},
or a {\it $\Lambda$-coalgebra} $(A_*,{\vartriangle}_A)$ is a pair satisfying the following conditions.

(1) $A_*$ is a $\Lambda$-group $A_*\!=\!\oplus_{\alpha\in\Lambda}\,A_*^\alpha$.

(2) The coproduct ${\vartriangle}_A\colon A_*\to A_*{\otimes}A_*$ is a group homomorphism
that may not be coassociative or keep degree,
where we neglect the indexed group structure of $A_*$ and $A_*{\otimes}A_*$.

The local coproducts of $\vartriangle_A$ with respect to $\Lambda$ are defined as follows.
For $a\in A_*^\alpha$ and each $\alpha'\!,\alpha''\in\Lambda$, there is a unique
$b_{\alpha'\!,\alpha''}\in A_*^{\alpha'}{\otimes}A_*^{\alpha''}$ such that ${\vartriangle}_A(a)=\Sigma_{\alpha',\alpha''\in\Lambda}b_{\alpha'\!,\alpha''}$.
The correspondence $a\to b_{\alpha'\!,\alpha''}$ is the group homomorphism
$$(\vartriangle_A)^\alpha_{\alpha',\alpha''}\colon A^\alpha_*\stackrel{i}{\to}A_*
\stackrel{\vartriangle_A}{\longrightarrow}A_*{\otimes}A_*\stackrel{p}{\to}
A^{\alpha'}_*{\otimes}A^{\alpha''}_*,$$
where $i$ is the inclusion and $p$ is the projection.
Each $(\vartriangle_A)^\alpha_{\alpha',\alpha''}$
is called a {\it local coproduct} of ${\vartriangle_A}$.
${\vartriangle_A}$ is defined if and only if all its local coproducts are defined.
$(A_*,\vartriangle_A)$ is a graded $\Lambda$-coalgebra, i.e., ${\vartriangle_A}$ keeps degree,
if and only if all its local coproducts keep degree.

A {\it subcoalgebra} $B_*$ of the $\Lambda$-coalgebra $(A_*,{\vartriangle}_A)$
is a subgroup of $A_*$ such that ${\vartriangle}_A(B_*)\subset B_*{\otimes}B_*$.
Denote by ${\vartriangle}_B$ the restriction of ${\vartriangle}_A$ on $B_*$.
By Lemma~2.2, $(B_*,{\vartriangle}_{B})$ is naturally a $\Lambda$-coalgebra
such that each local coproduct $(\vartriangle_{B})^\alpha_{\alpha',\alpha''}$ is the restriction of $(\vartriangle_A)^\alpha_{\alpha',\alpha''}$.
$(B_*,{\vartriangle}_{B})$ is also called the {\it subcoalgebra} of $(A_*,{\vartriangle}_A)$.

An {\it indexed coalgebra homomorphism}, or a {\it $\Lambda$-coalgebra homomorphism} $\theta\colon(A_*,{\vartriangle}_A)\to(B_*,{\vartriangle}_B)$
is a $\Lambda$-group homomorphism such that $(\theta{\otimes}\theta){\vartriangle}_A={\vartriangle}_B\theta$.
Precisely, if $\theta=\oplus_{\alpha\in\Lambda}\,\theta_\alpha$,
then $(\theta_{\alpha'}{\otimes}\theta_{\alpha''})({\vartriangle}_A)^\alpha_{\alpha'\!,\alpha''}
=({\vartriangle}_B)^\alpha_{\alpha'\!,\alpha''}\theta_\alpha$ for all local coproducts.
If $\theta$ is a $\Lambda$-group isomorphism,
then the two $\Lambda$-coalgebras are isomorphic and denoted by $(A_*,{\vartriangle}_A)\cong_\Lambda(B_*,{\vartriangle}_B)$,
or simply $A_*\cong_\Lambda B_*$.

Dually, an {\it indexed algebra}, or an {\it algebra indexed by $\Lambda$},
or a {\it $\Lambda$-algebra} $(A^*,{\triangledown}_{\!A})$ is a pair satisfying the following conditions.

(1) $A^*$ is a $\Lambda$-group $A^*\!=\!\oplus_{\alpha\in\Lambda}\,A^*_\alpha$.

(2) The product ${\triangledown}_{\!A}\colon A^*{\otimes}A^*\to A^*$ is a group homomorphism that may not be associative or keep degree,
where we neglect the indexed group structure of $A^*$ and $A^*{\otimes}A^*$.

The local products of ${\triangledown}_{\!A}$ with respect to $\Lambda$ are defined as follows.
For $b\in A^*_{\alpha'}$, $c\in A^*_{\alpha''}$ and each $\alpha\in\Lambda$, there is a unique
$a_\alpha\in A^*_\alpha$ such that
${\triangledown}_{\!A}(b{\otimes}c)=\Sigma_{\alpha\in\Lambda}\,a_\alpha$.
The correspondence $b{\otimes}c\to a_\alpha$ is the group homomorphism
$$(\triangledown_{\!A})_\alpha^{\alpha',\alpha''}\colon A_{\alpha'}^*{\otimes}A_{\alpha''}^*\stackrel{i}{\to}
A^*{\otimes}A^*\stackrel{\triangledown_{\!A}}{\longrightarrow}A^*\stackrel{p}{\to} A_\alpha^*,$$
where $i$ is the inclusion and $p$ is the projection.
Each \,$(\triangledown_{\!A})_\alpha^{\alpha',\alpha''}$\,
is called a {\it local product} of ${\triangledown}_{\!A}$.
${\triangledown}_{\!A}$ is defined if and only if all its local products are defined.
$(A^*,\triangledown_{\!A})$ is a graded $\Lambda$-algebra, i.e., ${\triangledown}_{\!A}$ keeps degree,
if and only if all its local products keep degree.

An {\it ideal} $I^*$ of the $\Lambda$-algebra $(A^*,{\triangledown}_{\!A})$
is a subgroup of $A^*$ such that ${\triangledown}_{\!A}(I^*{\otimes}A^*+A^*{\otimes}I^*)\subset I^*$.
The quotient group $A^*/I^*$ has an induced product ${\triangledown}_{\!A/I}\colon A^*/I^*{\otimes}A^*/I^*\to A^*/I^*$ defined by
${\triangledown}_{\!A/I}([x]{\otimes}[y])=[{\triangledown}_{\!A}(x{\otimes}y)]$.
By Lemma~2.2, $(A^*/I^*,{\triangledown}_{\!A/I})$ is naturally a $\Lambda$-algebra
such that each local product $(\triangledown_{\!A/I})^{\alpha',\alpha''}_\alpha$ is the quotient of $(\triangledown_{\!A})^{\alpha',\alpha''}_\alpha$.
$(A^*/I^*,{\triangledown}_{\!A/I})$ is called the {\it quotient algebra} of $(A^*,{\triangledown}_{\!A})$ over $I^*$.

An indexed algebra\, homomorphism, or a {\it $\Lambda$-algebra \,homomorphism} \,$\theta\colon$ \,$(A^*,\triangledown_{\!A})\to(B^*,\triangledown_{\!B})$
is a $\Lambda$-group homomorphism such that $\theta{\triangledown}_{\!A}={\triangledown}_{\!B}(\theta{\otimes}\theta)$.
Precisely, if $\theta=\oplus_{\alpha\in\Lambda}\,\theta_\alpha$,
then $\theta_\alpha({\triangledown}_{\!A})^{\alpha',\alpha''}_\alpha
=({\triangledown}_{\!B})^{\alpha',\alpha''}_\alpha(\theta_{\alpha'}{\otimes}\theta_{\alpha''})$ for all local products.
If $\theta$ is a $\Lambda$-group isomorphism,
then the two $\Lambda$-algebras are isomorphic and denoted by $(A^*,{\triangledown}_{\!A})\cong_\Lambda(B^*,{\triangledown}_{\!B})$,
or simply $A^*\cong_\Lambda B^*$.
\vspace{3mm}

In this paper, we often have to regard a $\Lambda$-coalgebra as a coalgebra without indexed group structure.
So for $\Lambda$-coalgebras $(A_*,{\vartriangle}_A)$ and $(B_*,{\vartriangle}_B)$ such that $A_*\not\cong_\Lambda B_*$,
there might be a coalgebra isomorphism $A_*\cong B_*$ when we neglect the indexed group structure.
For example, let $A_*^\XX$ be defined by $A_*^{\mathpzc i}=\Bbb Z(1)$, $A_*^{\mathpzc n}=\Bbb Z(n)$, $A_*^{\mathpzc c}=\Bbb Z(c)$.
Let ${\vartriangle}_1,{\vartriangle}_2,{\vartriangle}_3$ be coproducts on $A_*^\XX$ defined as follows.
$${\vartriangle}_1(1)={\vartriangle}_2(1)={\vartriangle}_3(1)=1{\otimes}1,
\quad {\vartriangle}_1(a)={\vartriangle}_2(a)={\vartriangle}_3(a)=a{\otimes}1{+}1{\otimes}a,$$
$${\vartriangle}_1(c)=c{\otimes}1{+}1{\otimes}c,\quad{\vartriangle}_2(c)={\vartriangle}_1(c){+}a{\otimes}1{+}1{\otimes}a
,\quad{\vartriangle}_3(c)={\vartriangle}_1(c){+}a{\otimes}a.$$
It is obvious that the three coalgebras $(A_*^\XX,\Delta_i)$ are not isomorphic $\XX$-coalgebras.
Define the group isomorphism $f\colon A_*^\XX\to A_*^\XX$ by $f(1)=1$, $f(a)=a$, $f(c)=c{-}a$.
Then $f$ is a coalgebra isomorphism from $(A_*^\XX,{\vartriangle}_1)$ to $(A_*^\XX,{\vartriangle}_2)$.
So we have $(A_*^\XX,{\vartriangle}_1)\cong(A_*^\XX,{\vartriangle}_2)$ but $(A_*^\XX,{\vartriangle}_1)\not\cong_\XX(A_*^\XX,{\vartriangle}_2)$.
It is obvious that $(A_*^\XX,{\vartriangle}_1)\not\cong(A_*^\XX,{\vartriangle}_3)$.
\vspace{3mm}

{\bf Lemma~6.2} {\it Let $(A_*,{\vartriangle}_A)$ be a free $\Lambda$-coalgebra with $A_*=\oplus_{\alpha\in\Lambda}A_*^\alpha$.
Then the dual group $A^*=\oplus_{\alpha\in\Lambda}A^*_\alpha$ with $A^*_\alpha={\rm Hom}(A_*^\alpha,\Bbb Z)$
is a $\Lambda$-algebra with product ${\triangledown}_{\!A}$ the dual of ${\vartriangle}_A$
defined by ${\triangledown}_{\!A}(f{\otimes}g)(a)=(f{\otimes}g)({\vartriangle}_A(a))$ for all $a\in A_*$ and $f,g\in A^*$.
Moreover, each local product $(\triangledown_{\!A})_{\alpha}^{\alpha'\!,\alpha''}$
is the dual of $(\vartriangle_A)^{\alpha}_{\alpha'\!,\alpha''}$.
$(A^*,{\triangledown}_{\!A})$ is called the dual algebra of $(A_*,{\vartriangle}_A)$.

For a free subcoalgebra $(B_*,{\vartriangle}_B)$ of the $\Lambda$-coalgebra $(A_*,{\vartriangle}_A)$ such that $A_*/B_*$ is also free,
the dual algebra $(B^*,{\triangledown}_{\!B})$ of $(B_*,{\vartriangle}_B)$
is a quotient algebra of $(A^*,{\triangledown}_{\!A})$.
\vspace{2mm}

Proof}\,
For a free subcoalgebra $B_*$ of $A_*$ such that $A_*/B_*$ is free,
the group $I^*=\{f\in A^*\,|\,f(B_*)=0\}$ is an ideal of $A^*$.
So we have algebra isomorphism $B^*\cong A^*/I^*$.
\hfill$\Box$\vspace{3mm}

The above definitions and lemma can be generalized to the (co)chain complex case.
To distinguish the two cases, we always use $\psi$ and $\pi$ to denote respectively
the coproduct of chain coalgebras and the product of cochain algebras.
\vspace{3mm}

{\bf Definition~6.3} An {\it indexed chain coalgebra}, or a {\it chain coalgebra indexed by $\Lambda$},
or a {\it chain $\Lambda$-coalgebra} is a triple $(C_*,\psi_C,d)$ satisfying the following conditions.

(1) $(C_*,d)$ is a chain $\Lambda$-complex.

(2) The coproduct $\psi_C\colon(C_*,d)\to(C_*{\otimes}C_*,d)$ is a chain homomorphism
if we forget the indexed group structure of $C_*$ and $C_*{\otimes}C_*$.
So a local coproduct $(\psi_C)^\alpha_{\alpha',\alpha''}\colon C_*^\alpha\to C_*^{\alpha'}{\otimes}C_*^{\alpha''}$
of $\psi_C$ may not be a chain homomorphism.

A {\it strong chain $\Lambda$-coalgebra} $(C_*,\psi_C,d)$ is a chain $\Lambda$-coalgebra such that
each local coproduct $(\psi_C)^\alpha_{\alpha',\alpha''}\colon C_*^\alpha\to
C_*^{\alpha'}{\otimes}C_*^{\alpha''}$ is a chain homomorphism.
If there is no confusion, a (strong) chain $\Lambda$-coalgebra $(C_*,\psi_C,d)$ is simply denoted by $(C_*,\psi_C)$ or $C_*$.

A {\it chain subcoalgebra} $(D_*,\psi_D,d)$ of the chain $\Lambda$-coalgebra $(C_*,\psi_C,d)$
is a triple such that $(D_*,\psi_D)$ is a subcoalgebra of $(C_*,\psi_C)$ and $(D_*,d)$
is a subcomplex of $(C_*,d)$.
$(D_*,\psi_D,d)$ is naturally a chain $\Lambda$-coalgebra.

A {\it chain $\Lambda$-coalgebra homomorphism} $\vartheta\colon(C_*,\psi_C,d)\to(D_*,\psi_D,d)$
is a $\Lambda$-group homomorphism such that $\vartheta d=d\,\vartheta$ and $(\vartheta{\otimes}\vartheta)\psi_C=\psi_D \vartheta$.
If $\vartheta$ is a group isomorphism,
then the two chain $\Lambda$-coalgebras are isomorphic and denoted by $(C_*,\psi_C,d)\cong_\Lambda(D_*,\psi_D,d)$,
or simply $C_*\cong_\Lambda D_*$.

A {\it chain $\Lambda$-coalgebra weak homomorphism} $\vartheta\colon(C_*,\psi_C,d)\to(D_*,\psi_D,d)$
is a $\Lambda$-group homomorphism such that $\vartheta d=d\,\vartheta$ and $(\vartheta{\otimes}\vartheta)\psi_C\simeq\psi_D \vartheta$.

Dually, an {\it indexed cochain algebra}, or a {\it cochain algebra indexed by $\Lambda$},
or a {\it cochain $\Lambda$-algebra} is a triple $(C^*,\pi_C,\delta)$ satisfying the following conditions.

(1) $(C^*,\delta)$ is a cochain $\Lambda$-complex.

(2) The product $\pi_C\colon(C^*{\otimes}C^*,\delta)\to(C^*,\delta)$ is a cochain homomorphism if we
forget the indexed group structure of $C^*$ and $C^*{\otimes}C^*$.
So a local product $(\pi_C)_\alpha^{\alpha',\alpha''}\colon C^*_{\alpha'}{\otimes}C^*_{\alpha''}\to C^*_\alpha$ of $\pi_C$
may not be a cochain homomorphism.

A {\it strong cochain $\Lambda$-algebra} $(C^*,\pi_C,\delta)$ is a cochain $\Lambda$-algebra such that
each local coproduct $(\pi_C)_\alpha^{\alpha',\alpha''}\colon C^*_{\alpha'}{\otimes}C^*_{\alpha''}\to C^*_\alpha$ is a cochain homomorphism.
If there is no confusion, a (strong) cochain $\Lambda$-algebra $(C^*,\pi_C,\delta)$ is simply denoted by $(C^*,\pi_C)$ or $C^*$.

A {\it cochain ideal} $(I^*,\delta)$ of the cochain $\Lambda$-algebra $(C^*,\pi_C,\delta)$
is a pair such that $I^*$ is an ideal of $C^*$ and $(I_*,\delta)$ is a subcomplex of $(C^*,\delta)$.
We have a quotient algebra $(C^*/I^*,\pi_{C/I})$ and a quotient complex $(C^*/I^*,\delta)$.
Then, $(C^*/I^*,\pi_{C/I},\delta)$ is naturally a cochain $\Lambda$-algebra and is called the
{\it quotient cochain algebra} of $(C^*,\pi_C,\delta)$ over the cochain ideal $(I^*,\delta)$.

A {\it cochain $\Lambda$-algebra homomorphism} $\vartheta\colon(C^*,\pi_C,\delta)\to(D^*,\pi_D,\delta)$
is a $\Lambda$-group homomorphism such that $\vartheta\delta=\delta\vartheta$ and $\vartheta\pi_C=\pi_D(\vartheta{\otimes}\vartheta)$.
If $\vartheta$ is a group isomorphism,
then the two (strong) cochain $\Lambda$-algebras are isomorphic and denoted by $(C^*,\pi_C,\delta)\cong_\Lambda(D^*,\pi_D,\delta)$,
or simply $C^*\cong_\Lambda D^*$.

A {\it cochain $\Lambda$-algebra weak homomorphism} $\vartheta\colon(C^*,\pi_C,\delta)\to(D^*,\pi_D,\delta)$
is a $\Lambda$-group homomorphism such that $\vartheta\delta=\delta\vartheta$ and $\vartheta\pi_C\simeq\pi_D(\vartheta{\otimes}\vartheta)$.
\vspace{3mm}

{\bf Lemma~6.4} {\it Let $(C_*,\psi_C,d)$ be a free chain $\Lambda$-coalgebra with $C_*=\oplus_{\alpha\in\Lambda}C_*^\alpha$.
Then the dual cochain complex $(C^*,\delta)=\oplus_{\alpha\in\Lambda}(C^*_\alpha,\delta)$ with $C^*_\alpha={\rm Hom}(C_*^\alpha,\Bbb Z)$
is a  cochain $\Lambda$-algebra with product $\pi_C$ the dual of $\psi_C$
defined by $\pi_C(f{\otimes}g)(a)=(f{\otimes}g)(\psi_C(a))$ for all $a\in C_*$ and $f,g\in C^*$.
Moreover, each local product $(\pi_C)_{\alpha}^{\alpha'\!,\alpha''}$
is the dual of the local coproduct $(\psi_C)^{\alpha}_{\alpha'\!,\alpha''}$.
$(C^*,\pi_C,\delta)$ is called the dual cochain algebra of $(C_*,\psi_C,d)$.

For a free chain subcoalgebra $(D_*,\psi_D,d)$ of the chain $\Lambda$-coalgebra $(C_*,\psi_C,d)$ such that $C_*/D_*$ is also free,
the dual cochain algebra $(D^*,\pi_D,\delta)$ of $(D_*,\psi_D,d)$ is a quotient cochain algebra of $(C^*,\pi_C,\delta)$.
\vspace{2mm}

Proof}\, The same as Lemma~6.2.
\hfill$\Box$\vspace{3mm}

We have four types of indexed coalgebras as follows.
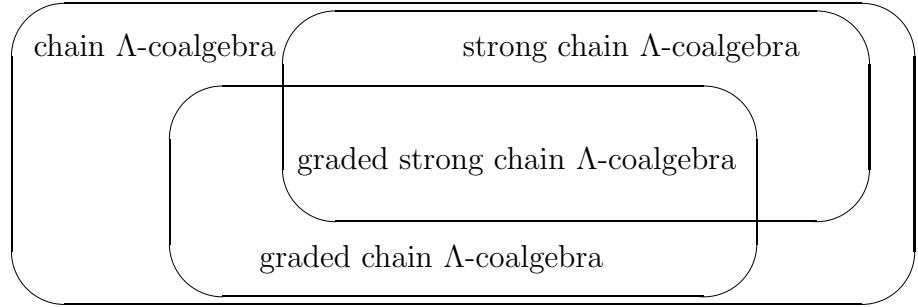
\begin{center}
\setlength{\unitlength}{1mm}
\begin{picture}(120,40)
\put(60,20){\oval(120,40)}\put(60,15){\oval(78,28)}\put(75,25){\oval(78,28)}
\put(3,33){chain $\Lambda$-coalgebra}\put(60,33){strong chain $\Lambda$-coalgebra}
\put(33,5){graded chain $\Lambda$-coalgebra}\put(38,18){graded strong chain $\Lambda$-coalgebra}
\end{picture}
\end{center}
The following two theorems and those in the remaining sections show that in the homology split condition (the only case concerned in this paper),
we may neglect the existence of graded, strong and graded strong chain $\Lambda$-coalgebras.
Moreover, we have to study chain $\Lambda$-coalgebra weak homomorphisms instead of homomorphisms.
This is because all the theorems for chain $\Lambda$-coalgebras with weak homomorphisms in this paper will have no better result
if the chain $\Lambda$-coalgebras or weak homomorphisms are replaced by
graded (strong, graded strong) chain $\Lambda$-coalgebras with homomorphisms.
\vspace{3mm}

{\bf Theorem~6.5} {\it Let $(C_*,\psi_C,d)$ be a free chain $\Lambda$-coalgebra
with dual cochain $\Lambda$-algebra $(C^*,\pi_C,\delta)$.

Suppose $H_*(C_*)$ is a free group.
By K\"{u}nneth Theorem, $\psi_C$ induces a coproduct ${\vartriangle}_C$ on $H_*(C_*)$
and so $(H_*(C_*),{\vartriangle}_C)$ is a $\Lambda$-coalgebra. The local coproduct
$$({\vartriangle}_C)^{\alpha}_{\alpha',\alpha''}\colon H_*(C_*^\alpha)\to H_*(C_*^{\alpha'}){\otimes}H_*(C_*^{\alpha''})$$
of ${\vartriangle}_C$ is defined as follows. For a homology class $[x]\in H_*(C_*^\alpha)$, $({\vartriangle}_C)_{\alpha}^{\alpha',\alpha''}([x])=[(\psi_C)_{\alpha}^{\alpha',\alpha''}(x)]$.
Dually, $\pi_C$ induces a product ${\triangledown}_{\!A}$ on $H^*(C^*)$
and so $(H^*(C^*),{\triangledown}\!_C)$ is a $\Lambda$-algebra. The local product
$$({\triangledown}\!_C)_{\alpha}^{\alpha',\alpha''}\colon H^{\,*}(C^*_{\alpha'}){\otimes}H^*(C^{\,*}_{\alpha''})\to H^*(C^{\,*}_\alpha)$$
of ${\triangledown}\!_C$ is defined as follows.
For cohomology classes $[x]\in H^*(C^{\,*}_{\alpha'})$ and $[y]\in H^*(C^{\,*}_{\alpha''})$,
$({\triangledown}\!_C)_{\alpha}^{\alpha',\alpha''}([x]{\otimes}[y])=[(\pi_C)_{\alpha}^{\alpha',\alpha''}(x{\otimes}y)]$.

If $H_*(C_*)$ is not a free group, then $\pi_C$ also induces a cup product
$$\cup_C\colon H^*(C^*){\otimes}H^*(C^*)\to H^*(C^*)$$ defined as follows.
For $[x]\in H^*(C^*)$ and $[y]\in H^*(C^*)$, $[x]\cup_C[y]=[\pi_C(x{\otimes}y)]$.
Then $(H^*(C^*),\cup_C)$ is a $\Lambda$-algebra with the local product
$$(\cup_C)_\alpha^{\alpha'\!,\alpha''}\colon H^{\,*}(C^*_{\alpha'}){\otimes}H^*(C^{\,*}_{\alpha''})\to H^*(C^{\,*}_\alpha)$$
of $\cup_C$ defined as follows. For cohomology classes $[x]\in H^*(C^{\,*}_{\alpha'})$ and $[y]\in H^*(C^{\,*}_{\alpha''})$,
$[x]{(\cup_C)_\alpha^{\alpha'\!,\alpha''}}[y]=[(\pi_C)_\alpha^{\alpha'\!,\alpha''}(x{\otimes}y)]$.
\vspace{2mm}

Proof}\, For $[x]\in H_*(C_*^\alpha)$, $[\psi_C(x)]=[\Sigma_{\alpha'\!,\alpha''}\,(\psi_C)^\alpha_{\alpha'\!,\alpha''}(x)]$
is a homology class in $H_*(C_*{\otimes}C_*)$.
By K\"{u}nneth Theorem, $H_*(C_*{\otimes}C_*)$ is the direct sum group
$\oplus_{\alpha'\!,\alpha''} H_*(C_*^{\alpha'}){\otimes}H_*(C_*^{\alpha''})$.
So $[(\psi_C)^\alpha_{\alpha'\!,\alpha''}(x)]$ is the summand group of ${\vartriangle}_C([x])$ in $H_*(C_*^{\alpha'}){\otimes}H_*(C_*^{\alpha''})$.
This implies that the equality $({\vartriangle}_C)^\alpha_{\alpha'\!,\alpha''}([x])=[(\psi_C)^\alpha_{\alpha'\!,\alpha''}(x)]$ holds even if $(\psi_C)^\alpha_{\alpha'\!,\alpha''}$ is not a chain homomorphism.
 \hfill$\Box$\vspace{3mm}

{\bf Conventions} As in the above definition, the notation $(A_*,{\vartriangle}_A)$ implies that $A_*$ is a free coalgebra.
Dually, the notation $(A^*,{\triangledown}_{\!A})$ implies that $A^*$ is a free algebra.
So we always denote the cup product $\cup_A$ by $a\cup_A b$ and denote the product ${\triangledown}_{\!A}$ by ${\triangledown}_{\!A}(a{\otimes}b)$.
\vspace{3mm}

{\bf Theorem~6.6} {\it Let $\vartheta\colon(C_*,\psi_C,d)\to(D_*,\psi_D,d)$ be a chain $\Lambda$-coalgebra
weak homomorphism with dual cochain $\Lambda$-algebra weak homomorphism $\vartheta^\circ$.
Then we have induced (co)homology $\Lambda$-(co)algebra homomorphisms
$$\theta\colon(H_*(C_*),{\vartriangle}_C)\to(H_*(D_*),{\vartriangle}_D),\quad
\theta^\circ\colon(H^*(D^*),{\triangledown}\!_D)\to(H^*(C^*),{\triangledown}\!_C),$$
$$\theta^\circ\colon(H^*(D^*),\cup_D)\to(H^*(C^*),\cup_C),$$
where ${\vartriangle}_-,{\triangledown}_-,\cup_-$ are as defined in Theorem~6.5.
\vspace{2mm}

Proof}\, For $[x]\in H_*(C_*)$, $\theta([x])=[\vartheta(x)]$ is a $\Lambda$-group homomorphism from $H_*(C_*)$ to $H_*(D_*)$.
The weak homomorphism implies
$${\vartriangle}_B(\theta([x]))=[\psi_B(\vartheta(x))]=[(\vartheta{\otimes}\vartheta)\psi_A(x)]=(\theta{\otimes}\theta)({\vartriangle}_A([x])).$$
So $\theta$ is a $\Lambda$-coalgebra homomorphism.
\hfill$\Box$\vspace{3mm}

{\bf Definition~6.7} Let $(A_*,{\vartriangle}_A)$ and $(B_*,{\vartriangle}_B)$ be coalgebras respectively indexed by $\Lambda$ and $\Gamma$.
The {\it tensor product coalgebra} $(A_*{\otimes} B_*,{\vartriangle}_A{\otimes}{\vartriangle}_B)$
is a $(\Lambda{\times}\Gamma)$-coalgebra defined as follows.
Suppose for $a\in A_*^\alpha$ and $b\in B_*^{\alpha'}$,
we have ${\vartriangle}_A(a)=\Sigma_i\,a'_i{\otimes}a''_i$
with $a'_i{\otimes}a''_i\in A_*^{\alpha'}{\otimes}A_*^{\alpha''}$
and ${\vartriangle}_B(b)=\Sigma_j\,b'_j{\otimes}b''_j$
with $b'_j{\otimes}b''_j\in B_*^{\alpha''}{\otimes}B_*^{\alpha'''}$.
Then
$$({\vartriangle}_A{\otimes}{\vartriangle}_B)(a{\otimes}b)=\Sigma_{i,j}\,
(-1)^{|a''_i||b'_j|}(a'_i{\otimes}b'_j){\otimes}(a''_i{\otimes}b''_j).$$
Equivalently, each local coproduct of ${\vartriangle}_A{\otimes}{\vartriangle}_B$ satisfies
$$({\vartriangle}_A{\otimes}{\vartriangle}_B)^{(\alpha,\alpha')}_{(\alpha'\!,\alpha''),(\alpha''\!,\alpha''')}=
({\vartriangle}_A)^{\alpha}_{\alpha'\!,\alpha''}\otimes({\vartriangle}_B)^{\alpha'}_{\alpha''\!,\alpha'''}.$$

Dually, let $(A^*,{\triangledown}\!_A)$ and $(B^*,{\triangledown}\!_B)$ be algebras respectively indexed by $\Lambda$ and $\Gamma$.
The {\it tensor product algebra} $(A^*{\otimes} B^*,{\triangledown}\!_A{\otimes}{\triangledown}\!_B)$ is
a $(\Lambda{\times}\Gamma)$-algebra defined as follows.
Suppose for $\alpha'\!\in A^*_{\alpha'}$,  $a''\in A^*_{\alpha''}$
and $b'\in B^*_{\beta''}$, $b''\in B^*_{\beta''}$, we have
${\triangledown}\!_A(a'{\otimes}a'')=\Sigma_i\,a_i$  with $a_i\in A^*_\alpha$ and
${\triangledown}\!_B(b'{\otimes}b'')=\Sigma_j\,b_j$ with $b_j\in B^*_{\beta}$. Then
$$({\triangledown}\!_A{\otimes}{\triangledown}\!_B)((a'{\otimes}b'){\otimes}(a''{\otimes}b''))
=(-1)^{|a''||b'|}(\Sigma_{i,j}\,a_i{\otimes}b_j).$$
Equivalently, each local product of ${\triangledown}\!_A{\otimes}{\triangledown}\!_B$ satisfies
$$({\triangledown}\!_A{\otimes}{\triangledown}\!_B)_{(\alpha,\beta)}^{(\alpha'\!,\beta'),(\alpha''\!,\beta'')}=
({\triangledown}\!_A)_{\alpha}^{\alpha'\!,\alpha''}\otimes
({\triangledown}\!_B)_{\beta}^{\beta'\!,\beta''}.$$

We have analogue definitions for indexed cochain algebras and indexed chain coalgebras
by neglecting the (co)chain complex structure.
\vspace{3mm}

When $A_*^\DD$ and $B_*^\DD$ are indexed groups, the diagonal tensor product group $A_*^\DD\widehat\otimes B_*^\DD$ is always regarded
as a subgroup of $A_*^\DD{\otimes}B_*^\DD$.
When $A_*^\DD$ and $B_*^\DD$ are (co)algebras,
the (co)algebra $A_*^\DD\widehat\otimes B_*^\DD$ is in general not a sub(co)algebra or a quotient (co)algebra of  $A_*^\DD{\otimes}B_*^\DD$.
So we have the following conventions.
\vspace{2mm}

{\bf Convention} For (co)algebras $A_*^\DD$ and $B_*^\DD$,
we use $a\widehat\otimes b$ to denote the element of $A_*^\DD\widehat\otimes B_*^\DD$
and $a{\otimes}b$ to denote the element of $A_*^\DD{\otimes}B_*^\DD$.
Precisely, for $a\in A_*^{s;\alpha}$ and $b\in B_*^{t,\beta}$,
define $a{\widehat\otimes}b= a{\otimes}b\in A_*^{s;\alpha}{\otimes}B_*^{s;\beta}\subset A_*^\DD\widehat\otimes B_*^\DD$
if $s=t$ and $a{\widehat\otimes}b=0$ if $s\neq t$.
\vspace{3mm}

{\bf Definition~6.8} Let $(A_*^\DD,{\vartriangle}_A)$ and $(B_*^\DD,{\vartriangle}_B)$
be  coalgebras respectively indexed by $\DD{\times}\Lambda$ and $\DD{\times}\Gamma$.
The {\it diagonal tensor product coalgebra} (with respect to $\DD$)
$(A_*^\DD\widehat\otimes B_*^\DD,{\vartriangle}_A{\widehat\otimes}{\vartriangle}_B)$
is a $(\DD{\times}\Lambda{\times}\Gamma)$-coalgebra defined as follows.
Suppose for $a\in A_*^{s;\alpha}$ and $b\in B_*^{t;\beta}$,
we have ${\vartriangle}_A(a)=\Sigma_i\,a'_i{\otimes}a''_i$ with $a'_i{\otimes}a''_i\in A_*^{s';\alpha'}{\otimes}A_*^{s'';\alpha''}$
and ${\vartriangle}_B(b)=\Sigma_j\,b'_j{\otimes}b''_j$ with $b'_j{\otimes}b''_j\in B_*^{t';\beta'}{\otimes}B_*^{t'';\beta''}$.
Define
$$({\vartriangle}_A{\widehat\otimes}{\vartriangle}_B)(a{\widehat\otimes}b)=\Sigma_{i,j}\,
(-1)^{|a''_i||b'_j|}(a'_i{\widehat\otimes}b'_j){\otimes}(a''_i{\widehat\otimes}b''_j).$$
Equivalently, each local coproduct of ${\vartriangle}_A\widehat\otimes{\vartriangle}_B$ satisfies
$$({\vartriangle}_A{\widehat\otimes}{\vartriangle}_B)^{(s;\alpha,\beta)}_{(s';\alpha'\!,\beta'),(s'';\alpha''\!,\beta'')}=
({\vartriangle}_A)^{(s;\alpha)}_{(s';\alpha'),(s'';\alpha'')}\otimes({\vartriangle}_B)^{(s;\beta)}_{(s';\beta'),(s'';\beta'')}.$$

Dually, let $(A^{\,*}_\DD,{\triangledown}\!_A)$ and $(B^{\,*}_\DD,{\triangledown}\!_B)$
be  algebras respectively indexed by $\DD{\times}\Lambda$ and $\DD{\times}\Gamma$.
The {\it diagonal tensor product algebra} (with respect of $\DD$) $(A^{\,*}_\DD{\widehat\otimes} B^{\,*}_\DD,{\triangledown}\!_A{\widehat\otimes}{\triangledown}\!_B)$
is a $(\DD{\times}\Lambda{\times}\Gamma)$-algebra defined as follows.
Suppose for $a'{\otimes}a''\in A^*_{s';\alpha'}{\otimes}A^*_{s'';\alpha''}$
and  $b'{\otimes}b''\in B^*_{t';\beta'}{\otimes}B^*_{t'';\beta''}$, we have
${\triangledown}\!_A(a'{\otimes}a'')=\Sigma_i\,a_i$  with $a_i\in A^*_{s;\alpha}$
and ${\triangledown}\!_B(b'{\otimes}b'')=\Sigma_j\,b_j$ with $b_j\in B^*_{t;\beta}$. Define
$$({\triangledown}\!_A{\widehat\otimes}{\triangledown}\!_B)((a'{\widehat\otimes}b'){\otimes}(a''{\widehat\otimes}b''))
=(-1)^{|a''||b'|}(\Sigma_{i,j}\,a_i{\widehat\otimes}b_j).$$
Equivalently, each local product of ${\triangledown}\!_A{\widehat\otimes}{\triangledown}\!_B$ satisfies
$$({\triangledown}\!_A{\widehat\otimes}{\triangledown}\!_B)_{(s;\alpha,\beta)}^{(s';\alpha'\!,\beta'),(s'';\alpha''\!,\beta'')}=
({\triangledown}\!_A)_{(s;\alpha)}^{(s';\alpha'),(s'';\alpha'')}\otimes
({\triangledown}\!_B)_{(s;\beta)}^{(s';\beta'),(s'';\beta'')}.$$

We have analogue definitions for indexed cochain algebras and indexed chain coalgebras by neglecting the (co)chain complex structure.
\vspace{3mm}

{\bf Theorem~6.9} {\it We have  indexed (chain) coalgebra,  indexed (cochain) algebra isomorphism and homomorphism equalities
$$(A_*^{\DD_1}\widehat\otimes B_*^{\DD_1}){\otimes}{\cdots}{\otimes}
(A_*^{\DD_m}\widehat\otimes B_*^{\DD_m})
\cong(A_*^{\DD_1}{\otimes}{\cdots}{\otimes}A_*^{\DD_m})\widehat\otimes
(B_*^{\DD_1}{\otimes}{\cdots}{\otimes}B_*^{\DD_m}),
$$
$$(A^{\,*}_{\DD_1}\widehat\otimes B^{\,*}_{\!\DD_1}){\otimes}{\cdots}{\otimes}
(A^{\,*}_{\DD_m}\widehat\otimes B^{\,*}_{\!\DD_m})
\cong(A^{\,*}_{\DD_1}{\otimes}{\cdots}{\otimes}A^{\,*}_{\DD_m})\widehat\otimes
(B^{\,*}_{\!\DD_1}{\otimes}{\cdots}{\otimes}B^{\,*}_{\!\DD_m}),
$$
$$(f_1\widehat\otimes g_1){\otimes}{\cdots}{\otimes}
(f_m\widehat\otimes g_m)
=(f_1{\otimes}{\cdots}{\otimes}f_m)\widehat\otimes
(g_1{\otimes}{\cdots}{\otimes}g_m),
$$
where the right side diagonal tensor product is with respect to $\DD_1{\times}{\cdots}{\times}\DD_m$.
\vspace{2mm}

Proof}\, The group isomorphism $\phi$ in Theorem~2.6 is naturally a (co)algebra isomorphism.
\hfill$\Box$\vspace{3mm}

The properties of diagonal tensor product differ greatly from that of tensor product.
For example, for algebras $A^*_\DD$ and $B^*_\DD$ that are not associative (degree-keeping),
the diagonal tensor product algebra $A^*_\DD{\widehat\otimes} B^*_\DD$ might be associative (degree-keeping).
For $A^*_\DD\not\cong_\DD A'^*_\DD$ and  $B^*_\DD\not\cong_\DD B'^*_\DD$,
there might be an isomorphism $A^*_\DD{\widehat\otimes} B^*_\DD\cong_{\DD} A'^*_\DD{\widehat\otimes} B'^*_\DD$.
This is because for $a{\otimes}b\neq 0$ in a tensor product group,
$a{\widehat\otimes}b$ may be $0$ in the diagonal tensor product group.
In this paper, we want to know the algebra isomorphism class of a diagonal tensor product algebra
$A^*_\DD{\widehat\otimes} B^*_\DD$. So even if $A^*_\DD{\widehat\otimes} B^*_\DD\not\cong_{\DD} A'^*_\DD{\widehat\otimes} B'^*_\DD$,
we may have $A^*_\DD{\widehat\otimes} B^*_\DD\cong A'^*_\DD{\widehat\otimes} B'^*_\DD$.
\vspace{3mm}

In the following definition, for a  $(\DD{\times}\Lambda)$-coalgebra $(A_*^\DD,{\vartriangle}_A)$,
the group summand $A_*^s$ satisfies certain property implies $A_*^{s,\alpha}$ satisfies the property for all $\alpha\in\Lambda$
as in the convention before Definition~2.7.
So the local coproduct $({\vartriangle}_A)^s_{s'\!,s''}$ of ${\vartriangle}_A$ satisfies certain property implies
that $({\vartriangle}_A)^{(s;\alpha)}_{(s';\alpha'),(s'';\alpha'')}$ satisfies the property for all $\alpha,\alpha',\alpha''\in\Lambda$.
\vspace{2mm}

{\bf Definition~6.10} For  $(\DD{\times}\Lambda)$-coalgebras $(A_*^\DD,{\vartriangle}_A)$ and $(A_*^\DD,{\vartriangle}'_A)$ on the same group $A_*^\DD$,
$(A_*^\DD,{\vartriangle}_A)$ is called a {\it partial coalgebra} of $(A_*^\DD,{\vartriangle}'_A)$,
equivalently, ${\vartriangle}_A$ is called a {\it partial coproduct} of ${\vartriangle}'_A$ and is denoted by ${\vartriangle}_A\prec{\vartriangle}'_A$,
if all local coproducts of ${\vartriangle}_A$ satisfy that either $({\vartriangle}_A)^s_{s',s''}=({\vartriangle}'_A)^s_{s',s''}$, or $({\vartriangle}_A)^s_{s',s''}=0$.

For a subset $\TT$ of $\DD$, the {\it restriction coalgebra} of $(A_*^\DD,{\vartriangle}_A)$ on $\TT$ is the
 coalgebra $(B_*^\TT,{\vartriangle}_B)$ such that $B_*^\TT$ is the restriction group of $A_*^\DD$ on $\TT$
and the local coproducts of ${\vartriangle}_B$ satisfy $({\vartriangle}_B)^s_{s',s''}=({\vartriangle}_A)^s_{s',s''}$ for all $s,s',s''\in\TT$.
We always denote the restriction coalgebra $(B_*^\TT,{\vartriangle}_B)$ by $(A_*^\TT,{\vartriangle}_A)$.

The {\it support coalgebra} of $(A_*^\DD,{\vartriangle}_A)$ is the restriction coalgebra $(A_*^\SS,{\vartriangle}_A)$
on  the support index set $\SS=\{s\in\DD\,|\,A_*^s\neq 0\}$ of $A_*^\DD$.

We have all dual analogues for indexed  algebras.
We also have analogue definitions for indexed cochain algebras and indexed chain coalgebras
by neglecting the (co)chain complex structure.
\vspace{3mm}

{\bf Theorem~6.11} {\it Suppose $(A_{i\,*}^{\SS_i},{\vartriangle}_i)$ is the support coalgebra of
the $(\DD{\times}\Lambda_i)$-coalgebra $(A_{i\,*}^{\DD},{\vartriangle}_i)$ for $i=1,2$.
Then the support coalgebra of the diagonal tensor product coalgebra $(A_{1\,*}^\DD\widehat\otimes A_{2\,*}^\DD,{\vartriangle}_1\widehat\otimes{\vartriangle}_2)$
is the restriction coalgebra $(A_{1\,*}^\SS{\widehat\otimes}A_{2\,*}^\SS,{\vartriangle}_1\widehat\otimes{\vartriangle}_2)$ such that $\SS=\SS_1{\cap}\SS_2$.
For any index set $\TT$ such that $\SS\subset\TT\subset\DD$ and any coproduct ${\vartriangle}'_i$ such that ${\vartriangle}_i\prec{\vartriangle}'_i$,
we have the same $(\Lambda_1{\times}\Lambda_2)$-coalgebra
$$(A_{1\,*}^\DD{\widehat\otimes}A_{2\,*}^\DD,{\vartriangle}_1\widehat\otimes{\vartriangle}_2)
=(A_{1\,*}^\TT{\widehat\otimes}A_{2\,*}^\TT,{\vartriangle}'_1\widehat\otimes{\vartriangle}'_2),
$$
where the index $\TT,\DD$ are neglected.

We have analogue conclusion for indexed algebras.
We also have analogue conclusions for indexed cochain algebras and indexed chain coalgebras
by neglecting the (co)chain complex structure.
\vspace{2mm}

Proof}\, For $a\widehat\otimes b\in A_*^\DD\widehat\otimes B_*^\DD$,
either $a\widehat\otimes b=0$, or it is the same element $a\widehat\otimes b\in A_*^\TT\widehat\otimes B_*^\TT$.
\hfill$\Box$\vspace{3mm}

Note that although a restriction group is always a subgroup, a restriction coalgebra is in general not a subcoalgebra.
Specifically, the diagonal tensor product coalgebra $A_*^\DD\widehat\otimes B_*^\DD$ is in general not a subcoalgebra of $A_*^\DD{\otimes}B_*^\DD$.
This is one of the reasons why we denote the element of a diagonal tensor product (co)algebra by $a{\widehat\otimes}b$.
\vspace{3mm}

\section{Local (Co)products of Total Objects of Simplicial Complexes}\vspace{3mm}

\hspace*{5.5mm} {\bf Definition~7.1} We have the following coproducts of chain $\XX$-coalgebras on the same
atom chain complex $(T_*^\XX\!,d)$ in Definition~3.4.

The {\it universal coproduct} $\widehat\psi\colon T_*^\XX\!\to T_*^\XX\!{\otimes}T_*^\XX\!$ is defined as follows.

$\widehat\psi({\mathpzc i})={\mathpzc i}{\scriptstyle\otimes}{\mathpzc i}+{\mathpzc i}{\scriptstyle\otimes}{\mathpzc n}
+{\mathpzc n}{\scriptstyle\otimes}{\mathpzc i}+{\mathpzc n}{\scriptstyle\otimes}{\mathpzc n}$.

$\widehat\psi({\mathpzc n})={\mathpzc n}{\scriptstyle\otimes}{\mathpzc n}+{\mathpzc n}{\scriptstyle\otimes}{\mathpzc i}
+{\mathpzc i}{\scriptstyle\otimes}{\mathpzc n}$.

$\widehat\psi({\overline{\mathpzc n}})={\overline{\mathpzc n}}{\scriptstyle\otimes}{\mathpzc n}+{\overline{\mathpzc n}}{\scriptstyle\otimes}{\mathpzc i}+{\mathpzc i}{\scriptstyle\otimes}{\overline{\mathpzc n}}+
{\mathpzc e}{\scriptstyle\otimes}{\mathpzc e}+{\mathpzc e}{\scriptstyle\otimes}{\mathpzc i}+{\mathpzc i}{\scriptstyle\otimes}{\mathpzc e}+{\mathpzc i}{\scriptstyle\otimes}{\mathpzc i}$.

$\widehat\psi({\mathpzc e})={\mathpzc e}{\scriptstyle\otimes}{\mathpzc e}+{\mathpzc e}{\scriptstyle\otimes}{\mathpzc i}+{\mathpzc i}{\scriptstyle\otimes}{\mathpzc e}
+{\mathpzc i}{\scriptstyle\otimes}{\mathpzc i}$.

The {\it normal coproduct} $\w\psi\colon T_*^\XX\!\to T_*^\XX\!{\scriptstyle\otimes}T_*^\XX\!$ is defined as follows.

$\w\psi({\mathpzc i})={\mathpzc i}{\scriptstyle\otimes}{\mathpzc i}+{\mathpzc i}{\scriptstyle\otimes}{\mathpzc n}+
{\mathpzc n}{\scriptstyle\otimes}{\mathpzc i}+{\mathpzc n}{\scriptstyle\otimes}{\mathpzc n}$.

$\w\psi({\mathpzc n})={\mathpzc n}{\scriptstyle\otimes}{\mathpzc n}+{\mathpzc n}{\scriptstyle\otimes}{\mathpzc i}+
{\mathpzc i}{\scriptstyle\otimes}{\mathpzc n}$.

$\w\psi({\overline{\mathpzc n}})={\overline{\mathpzc n}}{\scriptstyle\otimes}{\mathpzc n}+{\overline{\mathpzc n}}{\scriptstyle\otimes}{\mathpzc i}
+{\mathpzc i}{\scriptstyle\otimes}{\overline{\mathpzc n}}$.

$\w\psi({\mathpzc e})={\mathpzc e}{\scriptstyle\otimes}{\mathpzc e}+{\mathpzc e}{\scriptstyle\otimes}{\mathpzc i}
+{\mathpzc i}{\scriptstyle\otimes}{\mathpzc e}+{\mathpzc i}{\scriptstyle\otimes}{\mathpzc i}$.

The {\it strictly normal coproduct} $\tilde\psi\colon T_*^\XX\!\to T_*^\XX\!{\scriptstyle\otimes}T_*^\XX\!$ is defined as follows.

$\w\psi({\mathpzc i})={\mathpzc i}{\scriptstyle\otimes}{\mathpzc i}$.

$\w\psi({\mathpzc n})={\mathpzc n}{\scriptstyle\otimes}{\mathpzc n}+{\mathpzc n}{\scriptstyle\otimes}{\mathpzc i}+
{\mathpzc i}{\scriptstyle\otimes}{\mathpzc n}$.

$\w\psi({\overline{\mathpzc n}})={\overline{\mathpzc n}}{\scriptstyle\otimes}{\mathpzc n}+{\overline{\mathpzc n}}{\scriptstyle\otimes}{\mathpzc i}
+{\mathpzc i}{\scriptstyle\otimes}{\overline{\mathpzc n}}$.

$\w\psi({\mathpzc e})={\mathpzc e}{\scriptstyle\otimes}{\mathpzc e}+{\mathpzc e}{\scriptstyle\otimes}{\mathpzc i}
+{\mathpzc i}{\scriptstyle\otimes}{\mathpzc e}$.

The {\it special coproduct} $\overline\psi\colon T_*^\XX\!\to T_*^\XX\!{\scriptstyle\otimes}T_*^\XX\!$ is defined as follows.

$\overline\psi({\mathpzc i})={\mathpzc i}{\scriptstyle\otimes}{\mathpzc i}$.

$\overline\psi({\mathpzc n})={\mathpzc n}{\scriptstyle\otimes}{\mathpzc i}+{\mathpzc i}{\scriptstyle\otimes}{\mathpzc n}$.

$\overline\psi({\overline{\mathpzc n}})={\overline{\mathpzc n}}{\scriptstyle\otimes}{\mathpzc i}+{\mathpzc i}{\scriptstyle\otimes}{\overline{\mathpzc n}}$.

$\overline\psi({\mathpzc e})={\mathpzc e}{\scriptstyle\otimes}{\mathpzc i}+{\mathpzc i}{\scriptstyle\otimes}{\mathpzc e}$.

The {\it weakly special coproduct} $\bar\psi\colon T_*^\XX\!\to T_*^\XX\!{\scriptstyle\otimes}T_*^\XX\!$ is defined as follows.

$\bar\psi({\mathpzc i})={\mathpzc i}{\scriptstyle\otimes}{\mathpzc i}$.

$\bar\psi({\mathpzc n})={\mathpzc n}{\scriptstyle\otimes}{\mathpzc i}+{\mathpzc i}{\scriptstyle\otimes}{\mathpzc n}$.

$\bar\psi({\overline{\mathpzc n}})={\overline{\mathpzc n}}{\scriptstyle\otimes}{\mathpzc i}+{\mathpzc i}{\scriptstyle\otimes}{\overline{\mathpzc n}}
+{\mathpzc i}{\scriptstyle\otimes}{\mathpzc i}$.

$\bar\psi({\mathpzc e})={\mathpzc e}{\scriptstyle\otimes}{\mathpzc i}+{\mathpzc i}{\scriptstyle\otimes}{\mathpzc e}$.

The {\it right universal, right normal, right strictly normal, right special, right weakly special} coproduct
$\widehat\psi'$, $\w\psi'$, $\tilde\psi'$, $\overline{\psi'}$, $\bar\psi'$ are defined as follows.
For $\psi$ one of $\widehat\psi,\w\psi,\tilde\psi,\overline{\psi},\bar\psi$, we always have $\psi'({\mathpzc e})=0$.
For ${\mathpzc s}\in\{{\mathpzc i},{\mathpzc n},\overline{\mathpzc n}\}$,
$\psi'({\mathpzc s})$ is obtained from $\psi({\mathpzc s})$ by canceling all terms with factor ${\mathpzc e}$.

The dual product $T^{\,*}_{\!\XX}\otimes T^{\,*}_{\!\XX}\to T^{\,*}_{\!\XX}$ of the above coproduct
is called the {\it universal, right universal (normal, right normal $\cdots$)} product
and is denoted by $\widehat{\pi},\widehat\pi'$ ($\w\pi,\w\pi'\cdots$).
\vspace{3mm}

{\bf Definition~7.2} An {\it atom chain coalgebra} $(T_*^\XX,\psi)$ is a partial coalgebra (Definition~6.10)
of the universal chain coalgebra $(T_*^\XX,\widehat\psi)$.
Its dual algebra $(T^{\,*}_{\!\XX},\pi)$ is called an {\it atom cochain algebra} and is also a partial algebra of $(T^{\,*}_{\!\XX},\widehat\pi)$.

A {\it total chain coalgebra} is $(T_*^{\XX_m},\underline\psi)=(T_*^{\XX}{\otimes}{\cdots}{\otimes}T_*^{\XX},\psi_1{\otimes}{\cdots}{\otimes}\psi_m)$,
where each $(T_*^{\XX},\psi_k)$ is an atom chain coalgebra.
Its dual algebra
$(T^{\,*}_{\!\XX_m},\underline\pi)=(T^{\,*}_{\!\XX}{\otimes}{\cdots}{\otimes}T^{\,*}_{\!\XX},\pi_1{\otimes}{\cdots}{\otimes}\pi_m)$
is called a {\it total cochain algebra}.

The total chain coalgebra $(T_*^{\XX_m},\psi^{(m)})=(T_*^{\XX}{\otimes}{\cdots}{\otimes}T_*^{\XX},\psi{\otimes}{\cdots}{\otimes}\psi)$
is called the {\it universal (normal, $\cdots$) chain coalgebra}
if $\psi$ is the universal (normal, $\cdots$) coproduct $\widehat\psi$ ($\w\psi,\cdots$).
Its dual total cochain algebra $(T^{\,*}_{\!\XX_m},\pi^{(m)})$
is called the {\it universal (normal, $\cdots$) cochain algebra}
if $\pi$ is the universal (normal, $\cdots$) product $\widehat\pi$ ($\w\pi,\cdots$).
\vspace{3mm}

{\bf Lemma~7.3}{\it\, Let $K$ be a simplicial complex on $[m]$ and $(T_*^{\XX_m},\underline{\psi})$ be a total chain coalgebra
with dual total cochain algebra $(T^{\,*}_{\!\XX_m},\underline{\pi})$.

The total chain complex $(T_*^{\XX_m}(K),d)$ in Definition~3.7 is a chain subcoalgebra of $(T_*^{\XX_m},\underline{\psi})$
such that $T_*^{\XX_m}/T_*^{\XX_m}(K)$ is also free.
Denote this chain $\XX_m$-coalgebra by $(T_*^{\XX_m}(K),\psi_K)$.
Then by Lemma~6.4, its dual algebra $(T^{\,*}_{\!\XX_m}(K),\pi_K)$ is a quotient cochain algebra of $(T^{\,*}_{\!\XX_m},\underline{\pi})$.
\vspace{2mm}

\it Proof}\, The $(U_*(\tau),d)$ in Definition~3.7 is a chain subcoalgebra of $(T_*^{\XX_m},\underline{\psi})$.
Then $(T_*^{\XX_m}(K),d)=(+_{\tau\in  K}U(\tau),d)$ is also a chain subcoalgebra of $(T_*^{\XX_m},\underline{\psi})$.
By Lemma~6.4, $(T^{\,*}_{\!\XX_m}(K),\pi)$ is a quotient algebra of $(T^{\,*}_{\!\XX_m},\underline{\pi})$.
\hfill$\Box$\vspace{3mm}

{\bf Definition~7.4} Let $K$ be a simplicial complex on $[m]$ and $(T_*^{\XX_m},\underline{\psi})$ be
a total chain coalgebra.

The {\it total chain coalgebra} $(T_*^\XX(K),\psi_K)$ of $K$ with respect to $\underline{\psi}$ is the $(T_*^{\XX_m}(K),\psi_K)$ in Lemma~7.3.
The {\it total cochain algebra} $(T^{\,*}_{\!\XX}(K),\pi_K)$ of $K$ with respect to $\underline{\psi}$ is the $(T^{\,*}_{\!\XX_m}(K),\pi_K)$ in Lemma~7.3.

By Theorem~6.5, if $H_*^{\XX_m}(K)$ is a free group, then we have a $\XX_m$-coalgebra $(H_*^{\XX_m}(K),{\vartriangle}_K)$
called {\it the total homology coalgebra} of $K$ with respect to $\underline{\psi}$.
The dual $\XX_m$-algebra $(H^{\,*}_{\!\XX_m}(K),{\triangledown}_{\!K})$
is called the {\it total cohomology algebra} of $K$ with respect to $\underline{\psi}$.

We have a product ${\cup}_K$ induced by $\underline{\psi}$ even if $H_*^{\XX_m}(K)$ is not free.
The $\XX_m$-algebra $(H^{\,*}_{\!\XX_m}(K),{\cup}_K)$ is called
the {\it total cohomology algebra} of $K$ with respect to $\underline{\psi}$.
When $H_*^{\XX_m}(K)$ is free, ${\cup}_K={\triangledown}_{\!K}$, i.e., $a\cup_K b={\triangledown}_{\!K}(a{\otimes}b)$.

For an index set $\DD\subset\XX_m$, the {\it total objects} $(T_*^\DD(K),\psi_K)$, $(H_*^\DD(K),{\vartriangle}_K)$,
$(T^{\,*}_{\!\DD}(K),\pi_K)$, $(H^{\,*}_{\!\DD}(K),{\triangledown}_K)$, $(H^{\,*}_{\!\DD}(K),\cup_K)$ of $K$ on $\DD$ with respect to $\underline{\psi}$
are the restriction (co)algebra of the corresponding total object of $K$.
When $\DD=\RR_m$, we have the {\it right total objects} of $K$.
Denote them as in the following table.
\begin{center}
{\rm \begin{tabular}{|c|c|c|}
\hline
{\rule[-3mm]{0mm}{8mm}}
total object\,\,&\,\,right total object\,\,&\,\, total object on $\DD$\\
\hline
{\rule[-3mm]{0mm}{8mm}}
\,\,$(T_*^{\XX_m}(K),\psi_K)$\,\,&\,\,$(T_*^{\RR_m}(K),\psi_K)$\,\,&\,\,$(T_*^\DD(K),\psi_K)$\,\,\\
\hline
{\rule[-3mm]{0mm}{8mm}}
\,\,$(T^{\,*}_{\!\XX_m}(K),\pi_K)$\,\,&\,\,$(T^{\,*}_{\!\RR_m}(K),\pi_K)$
\,\,&\,\,$(T^{\,*}_{\!\DD}(K),\pi_K)$\,\,\\
\hline
{\rule[-3mm]{0mm}{8mm}}
\,\,$(H_*^{\XX_m}(K),\vartriangle_K)$\,\,&\,\,$(H_*^{\RR_m}(K),\vartriangle_K)$
\,\,&\,\,$(H_*^{\DD}(K),\vartriangle_K)$\,\,\\
\hline
{\rule[-3mm]{0mm}{8mm}}
\,\,$(H^{\,*}_{\!\XX_m}(K),\triangledown_K)$\,\,&\,\,$(H^{\,*}_{\!\RR_m}(K),\triangledown_K)$
\,\,&\,\,$(H^{\,*}_{\!\DD}(K),\triangledown_K)$\,\,\\
\hline
{\rule[-3mm]{0mm}{8mm}}
\,\,$(H^{\,*}_{\!\XX_m}(K),\cup_K)$\,\,&\,\,$(H^{\,*}_{\!\RR_m}(K),\cup_K)$\,\,&\,\,$(H^{\,*}_{\!\DD}(K),\cup_K)$\,\,\\
\hline
\end{tabular}}
\end{center}
\vspace{5mm}

The total object of $K$ with respect to $\underline{\psi}$ is called
the {\it universal, (normal, $\cdots$) object} of $K$ if $\underline{\psi}=\psi^{(m)}$
with $\psi$ the universal (normal, $\cdots$) coproduct in Definition~7.1.
The {\it right total object} of $K$ is the restriction (co)algebra of the total object of $K$ on $\RR_m$.
They are denoted as in the following table.\vspace{2mm}\\
{\rm \begin{tabular}{|c|c|}
\hline
{\rule[-3mm]{0mm}{8mm}}
universal,\hspace{10mm} normal, $\quad\cdots$&right universal, right normal, $\cdots$ \\
\hline
{\rule[-3mm]{0mm}{8mm}}
$(T_*^{\XX_m}(K),\widehat\psi_K)$, $(T_*^{\XX_m}(K),\w\psi_K)$, $\cdots$&
$(T_*^{\RR_m}(K),\widehat\psi_K)$, $(T_*^{\RR_m}(K),\w\psi_K)$, $\cdots$\\
\hline
{\rule[-3mm]{0mm}{8mm}}
$(T^{\,*}_{\!\XX_m}(K),\widehat\pi_{\!K})$, $(T^{\,*}_{\!\XX_m}(K),\w\pi_{\!K})$, $\cdots$&
$(T^{\,*}_{\!\RR_m}(K),\widehat\pi_{\!K})$, $(T^{\,*}_{\!\RR_m}(K),\w\pi_{\!K})$, $\cdots$\\
\hline
{\rule[-3mm]{0mm}{8mm}}
$(H_*^{\XX_m}(K),\widehat\vartriangle_K)$, $(H_*^{\XX_m}(K),\w\vartriangle_K)$, $\cdots$&
$(H_*^{\RR_m}(K),\widehat\vartriangle_K)$, $(H_*^{\RR_m}(K),\w\vartriangle_K)$, $\cdots$\\
\hline
{\rule[-3mm]{0mm}{8mm}}
$(H^{\,*}_{\!\XX_m}(K),\widehat\triangledown_{\!K})$, $(H^{\,*}_{\!\XX_m}(K),\w\triangledown_{\!K})$, $\cdots$&
$(H^{\,*}_{\!\RR_m}(K),\widehat\triangledown_{\!K})$, $(H^{\,*}_{\!\RR_m}(K),\w\triangledown_{\!K})$, $\cdots$\\
\hline
{\rule[-3mm]{0mm}{8mm}}
$(H^{\,*}_{\!\XX_m}(K),\widehat\cup_K)$, $(H^{\,*}_{\!\XX_m}(K),\w\cup_K)$, $\cdots$&
$(H^{\,*}_{\!\RR_m}(K),\widehat\cup_K)$, $(H^{\,*}_{\!\RR_m}(K),\w\cup_K)$, $\cdots$\\
\hline
\end{tabular}}
\vspace{3mm}

For $\psi'=\widehat\psi',\w\psi',\cdots$ in Definition~7.1 and $\DD\subset\XX_m$, denote by $(T_*^\DD(K),\psi'_K)$
the total chain coalgebra of $K$ on $\DD$ with respect to $\underline{\psi}=\psi'^{(m)}$.
Then we have $(T_*^{\RR_m}(K),\psi_K)=(T_*^{\RR_m}(K),\psi'_K)$, although $(T_*^{\XX_m},\psi^{(m)})\neq(T_*^{\XX_m},\psi'^{(m)})$.
So all the (co)product $-_K$ of the right objects in the above table can be replaced by $-'_K$.

The universal chain coalgebra $(T_*^{\XX_m}(K),\underline{\psi})$ is neither a graded chain $\XX_m$-coalgebra nor a strong chain $\XX_m$-coalgebra.
The normal chain coalgebra $(T_*^{\XX_m}(K),\w\psi^{(m)})$ is a graded strong chain $\XX_m$-coalgebra.
The special chain coalgebra $(T_*^{\XX_m}(K),\overline\psi^{(m)})$ is a coassociative, cocommutative, graded strong chain $\XX_m$-coalgebra.
\vspace{3mm}

{\bf Theorem~7.5} {\it For any total chain coalgebra $(T_*^{\XX_m},\underline{\psi})$,
$(T_*^{\XX_m}(K),\psi_K)$ is a partial coalgebra of $(T_*^{\XX_m}(K),\widehat\psi_K)$,
i.e., each local coproduct $(\psi_K)^s_{s',s''}$ of $\psi_K$ satisfies
that either $(\psi_K)^s_{s',s''}=(\widehat\psi_K)^s_{s',s''}$, or $(\psi_K)^s_{s',s''}=0$.
Dually, $(T^{\,*}_{\!\XX_m}(K),\pi_K)$ is a partial algebra of $(T^{\,*}_{\!\XX_m}(K),\widehat\pi_K)$,
i.e., each local product $(\pi_K)_s^{s',s''}$ of $\pi_K$ satisfies
that either $(\pi_K)_s^{s',s''}=(\widehat\pi_K)_s^{s',s''}$, or $(\pi_K)_s^{s',s''}=0$.
\vspace{2mm}

Proof} Since $(T_*^{\XX_m},\underline{\psi})$ is a partial coalgebra of $(T_*^{\XX_m},\widehat\psi^{(m)})$,
their local coproducts satisfy that either $(\underline{\psi})^s_{s',s''}=(\widehat\psi^{(m)})^s_{s',s''}$, or $(\underline{\psi})^s_{s',s''}=0$.
Since $\widehat\psi_K=\widehat\psi^{(m)}|_{T_*^{\XX_m}(K)}$ and $\psi_K=\underline{\psi}|_{T_*^{\XX_m}(K)}$,
their local coproducts also satisfy that either $(\psi_K)^s_{s',s''}=(\widehat\psi_K)^s_{s',s''}$, or $(\psi_K)^s_{s',s''}=0$.
\hfill$\Box$\vspace{3mm}

{\bf Theorem~7.6} {\it Let $(T_*^{\XX_m}(K),\widehat\psi_K)$ be the universal chain coalgebra of $K$ in Definition~7.4.
Identify $T_*^{\sigma,\omega}(K)$, $T^*_{\sigma,\omega}(K)$ with the suspension augmented simplicial (co)chain complex
$\Sigma\w C_*(K_{\sigma,\omega})$, $\Sigma\w C^*(K_{\sigma,\omega})$ as in Theorem~3.8.\

The local coproducts
$$(\widehat\psi_K)_l=(\widehat\psi_K)_{\sigma'\!,\,\omega';\sigma''\!,\,\omega''}^{\sigma,\omega}\colon
T_*^{\sigma,\,\omega}(K)\to T_*^{\sigma'\!,\,\omega'}(K)\otimes T_*^{\sigma''\!,\,\omega''}(K)$$
of $\widehat\psi_K$ is defined as follows.
If $(\sigma'{\cup}\sigma''){\setminus}\sigma\subset\omega{\setminus}(\omega'{\cup}\omega'')\in K$, then for
$\tau\in T_*^{\sigma,\omega}(K)$, we have
$$(\widehat\psi_K)_l(\tau)=\langle\tau',\tau''\rangle\,\tau'{\otimes}\tau'',$$
where $\tau'=\tau{\cap}(\omega'{\setminus}(\sigma'{\cup}\sigma''))$,
$\tau''=\tau{\cap}((\omega''{\setminus}\omega'){\setminus}(\sigma'{\cup}\sigma''))$, $\langle\{i_1,{\cdots}i_s\},\{j_1,{\cdots}j_t\}\rangle$
is the sign of the permutation
$\Big(
\begin{array}{cccccc}
\scriptstyle{i_1} \!&\!\scriptstyle\!{\cdots}\!\!&\!\scriptstyle{i_s}\!&\!\scriptstyle{j_{1}} \!&\!\scriptstyle\!{\cdots}\!\!&\!\scriptstyle{j_t}\\
\scriptstyle{k_1} \!&\!\scriptstyle\!{\cdots}\!\!&\!\scriptstyle{k_s}\!&\!\scriptstyle{k_{s+1}} \!&\!\scriptstyle\!{\cdots}\!\!&\!\scriptstyle{k_{s+t}}\end{array}\Big)$
with all three sets ordered.
$(\widehat\psi_K)_l=0$ otherwise.

Dually, the local products
$$(\widehat\pi_K)_l=(\widehat\pi_K)^{\sigma'\!,\,\omega';\sigma''\!,\,\omega''}_{\sigma,\,\omega}\colon
T^{\,*}_{\sigma'\!,\,\omega'}(K)\otimes T^{\,*}_{\sigma''\!,\,\omega''}(K)\to T^{\,*}_{\sigma,\omega}(K)$$
of $\widehat\pi_K$ is defined as follows.
If $(\sigma'{\cup}\sigma''){\setminus}\sigma\subset\omega{\setminus}(\omega'{\cup}\omega'')\in K$, then for
$\tau'\in T^{\,*}_{\sigma'\!,\,\omega'}(K)$ and $\tau''\in T^{\,*}_{\sigma''\!,\,\omega''}(K)$, we have
$$(\widehat\pi_K)_l(\tau'{\otimes}\tau'')=\langle\tau',\tau''\rangle\,\Sigma\,\tau.$$
where the sum is taken over all $\tau\in T^{\,*}_{\sigma\!,\,\omega}(K)$
such that $\tau'=\tau{\cap}(\omega'{\setminus}(\sigma'{\cup}\sigma''))$,
$\tau''=\tau{\cap}((\omega''{\setminus}\omega'){\setminus}(\sigma'{\cup}\sigma''))$
and the sum is $0$ if there is no such $\tau$.
$(\widehat\pi_K)_l=0$ otherwise.
\vspace{2mm}

Proof}\, Suppose $(\widehat\psi_K)_l\neq 0$.
Let $t_{E,\overline N,N,I}$ be as defined in the proof of Theorem~3.8.
Suppose for $t=t_{E,\overline N,N,I}\in T_*^{\sigma,\omega}(K)$, we have
$$\widehat\psi_K(t)=\Sigma\pm t_{E'\!,\overline N'\!,N'\!,I'}{\otimes}t_{E''\!,\overline N''\!,N''\!,I''}
=\Sigma\pm(t'_1{\otimes}{\cdots}{\otimes}t'_m){\otimes}(t''_1{\otimes}{\cdots}{\otimes}t''_m),$$
where $t_{E'\!,\overline N'\!,N'\!,I'}\in T_*^{\sigma'\!,\,\omega'}(K)$, $t_{E''\!,\overline N''\!,N''\!,I''}\in T_*^{\sigma''\!,\,\omega''}(K)$
and $\pm$ is the sign $\langle \overline N'\!,\overline N''\rangle$.
Then we have

(1) $N\subset N'{\cup}N''$, for if $t_k={\mathpzc n}$, then at least one of $t'_k$ and $t''_k$ is ${\mathpzc n}$.

(2) $\overline N{\setminus}(\overline N'{\cup}\overline N'')\subset(E'{\cup}I'){\cap}(E''{\cup}I'')$,
for if $t_k={\overline{\mathpzc n}}$ and $t'_k,t''_k\neq{\overline{\mathpzc n}}$, then $t'_k,t''_k={\mathpzc e}$ or ${\mathpzc i}$.

(3) $(E'{\cup}E''){\setminus}E\subset \overline N{\setminus}(\overline N'{\cup}\overline N'')$, for if $t_k\neq{\mathpzc e}$ and one of $t'_k$ and $t''_k$ is ${\mathpzc e}$,
then $t_k={\overline{\mathpzc n}}$ and neither of $t'_k$ and $t''_k$ is ${\overline{\mathpzc n}}$.

(2) implies $(\overline N{\setminus}(\overline N'{\cup}\overline N'')){\cap}(N'{\cup}N'')=\emptyset$.
So (1) and (2) imply\\
\hspace*{17mm}$\overline N{\setminus}(\overline N'{\cup}\overline N'')=
(\overline N{\cup}N){\setminus}(\overline N'{\cup}\overline N''{\cup}N'{\cup}N'')
=\omega{\setminus}(\omega'{\cup}\omega'')$.\\
So by (3) we have $(\sigma'{\cup}\sigma''){\setminus}\sigma\subset\omega{\setminus}(\omega'{\cup}\omega'')$.
$\overline N\in K$ implies $\omega{\setminus}(\omega'{\cup}\omega'')\in K$.

Now suppose $(\sigma'{\cup}\sigma''){\setminus}\sigma\subset\omega{\setminus}(\omega'{\cup}\omega'')\in K$
and we prove $(\widehat\psi_K)_l$ satisfies the formula of the theorem.

For free groups $G,G',G''$, a coproduct $\phi\colon G\to G'{\otimes}G''$ is called a base inclusion
if for every generator $g\in G$, there are unique generators $g'\in G'$ and
$g''\in G''$ such that $\phi(g)=\pm\, g'{\otimes}g''$.
It is easy to check that each local coproduct of the universal coproduct $\widehat\psi$
is either a base inclusion or $0$.
So as a restriction of $\widehat\psi^{(m)}$, each local coproduct of $\widehat\psi_K$
is either a base inclusion or $0$.

For the generator
$t=t_{\sigma\!,\,\overline N,\,\omega{\setminus}\overline N,\,[m]{\setminus}(\sigma{\cup}\omega)}
=t_1{\otimes}{\cdots}{\otimes}t_m\in T_*^{\sigma\!,\,\omega}(K)$,
$\widehat\psi_K(t)$ has a summand
$$\begin{array}{l}
\quad\pm\,t_{\sigma'\!,\,\overline N'\!\!,\,\omega'{\setminus}\overline N'\!\!,\,[m]{\setminus}(\sigma'{\cup}\omega')}
{\otimes}t_{\sigma''\!,\,\overline N''\!\!,\,\omega''{\setminus}\overline N''\!\!,\,[m]{\setminus}(\sigma''{\cup}\omega'')}\vspace{2mm}\\
=\pm(t'_1{\otimes}{\cdots}{\otimes}t'_m){\otimes}(t''_1{\otimes}{\cdots}{\otimes}t''_m)\vspace{2mm}\\
\in T_*^{\sigma'\!,\,\omega'}(K){\otimes}T_*^{\sigma''\!,\,\omega''}(K)
\end{array}$$
defined as follows.

(1) For $k\in\sigma$, $t_k={\mathpzc e}$. Take
$t'_k={\mathpzc e}$, $t''_k={\mathpzc i}$ if $k\in\sigma'{\setminus}\sigma''$;
$t'_k={\mathpzc i}$, $t''_k={\mathpzc e}$ if $k\in\sigma''{\setminus}\sigma'$;
$t'_k=t''_k={\mathpzc e}$ if $k\in\sigma'{\cap}\sigma''$; $t'_k=t''_k={\mathpzc i}$, otherwise.

(2) For $k\in(\sigma'{\cup}\sigma''){\setminus}\sigma$, $t_k={\overline{\mathpzc n}}$. Take
$t'_k={\mathpzc e}$, $t''_k={\mathpzc i}$ if $k\in\sigma'{\setminus}\sigma''$;
$t'_k={\mathpzc i}$, $t''_k={\mathpzc e}$ if $k\in\sigma''{\setminus}\sigma'$;
$t'_k=t''_k={\mathpzc e}$ if $k\in\sigma'{\cap}\sigma''$.

(3) For $k\in \overline N{\setminus}(\sigma'{\cup}\sigma'')$, $t_k={\overline{\mathpzc n}}$.
Take $t'_k={\overline{\mathpzc n}}$, $t''_k={\mathpzc n}$ if $k\in\omega'{\cap}\omega''$; $t'_k={\overline{\mathpzc n}}$, $t''_k={\mathpzc i}$ if $k\in\omega'{\setminus}\omega''$;
$t'_k={\mathpzc i}$, $t''_k={\overline{\mathpzc n}}$ if $k\in\omega''{\setminus}\omega'$; $t'_k=t''_k={\mathpzc i}$, otherwise.

(4) For $k\in(\omega{\setminus}\overline N){\setminus}(\sigma'{\cup}\sigma'')$, $t_k={\mathpzc n}$. Take
$t'_k={\mathpzc n}$, $t''_k={\mathpzc i}$ if $k\in\omega'{\setminus}\omega''$;
$t'_k={\mathpzc i}$, $t''_k={\mathpzc n}$ if $k\in\omega''{\setminus}\omega'$;
$t'_k=t''_k={\mathpzc n}$ if $k\in\omega'{\cap}\omega''$.

(5) For $k\in[m]{\setminus}(\sigma{\cup}\omega)$, $t_k={\mathpzc i}$. Take
$t'_k={\mathpzc n}$, $t''_k={\mathpzc i}$ if $k\in\omega'{\setminus}\omega''$;
$t'_k={\mathpzc i}$, $t''_k={\mathpzc n}$ if $k\in\omega''{\setminus}\omega'$;
$t'_k=t''_k={\mathpzc n}$ if $k\in\omega'{\cap}\omega''$; $t'_k=t''_k={\mathpzc i}$, otherwise.

It is obvious that $\overline N'= \overline N{\cap}(\omega'{\setminus}(\sigma'{\cup}\sigma''))$,
$\overline N''= \overline N{\cap}((\omega''{\setminus}\omega'){\setminus}(\sigma'{\cup}\sigma''))$.
Since $(\widehat\psi_K)_l$ is a base inclusion, we have
$$\begin{array}{l}
\quad(\widehat\psi_K)_l(t_{\sigma\!,\,\overline N,\,\omega{\setminus}\overline N,\,[m]{\setminus}(\sigma{\cap}\omega)})\vspace{2mm}\\
={\scriptstyle\langle \overline N'\!,\,\overline N''\rangle}\,
t_{\sigma'\!,\,\overline N'\!,\,[m]\setminus\overline N'\!,\,[m]\setminus(\sigma'\cup\omega')}{\otimes}
t_{\sigma''\!,\,\overline N''\!,\,[m]\setminus\,N''\!,\,[m]\setminus(\sigma''\cup\omega'')}
\end{array}$$
is just the formula of the theorem by the isomorphism $\epsilon$ in
the proof of Theorem~3.8.
\hfill$\Box$\vspace{3mm}

{\bf Remark} Note that the above local coproduct $(\widehat\psi_K)_l$ is in general not a chain homomorphism and so
we do not have an induced homology homomorphism from $H_*^{\sigma\!,\,\omega}(K)$
to $H_*^{\sigma'\!,\,\omega'}(K){\otimes}H_*^{\sigma''\!,\,\omega''}(K)$.
So the coproduct $\psi_\vartriangle$ in Definition~7.1 of \cite{Z} is in general not a chain homomorphism.
But we still have the formula $({\vartriangle}_K)_l([x])=[(\widehat\psi_K)_l(x)]$ by Theorem~6.5.
The computations in that paper neglect the existence of simplicial complex $K$ and should be as in the following theorem.\vspace{3mm}

{\bf Theorem~7.7} {\it For a simplicial complex $K$ on $[m]$, let $\lozenge_K$ be any of the
$\psi_K$, $\pi_K$, $\vartriangle_K$, $\triangledown_K$, $\cup_K$ in Definition~7.4.
Then for the universal (normal, $\cdots$) objects, the local (co)products satisfy the following table.\vspace{3mm}\\
\hspace*{3mm}{\rm \begin{tabular}{|c|c|c|}
\hline
{\rule[-3mm]{0mm}{8mm}}
&$=(\widehat\lozenge_K)_l$&$=0$\\
\hline
{\rule[-3mm]{0mm}{8mm}$(\widehat\lozenge_K)_l$}
&$(\sigma'{\cup}\sigma''){\setminus}\sigma\subset\omega{\setminus}(\omega'{\cup}\omega'')\in K$&otherwise\\
\hline
{\rule[-3mm]{0mm}{8mm}$(\w\lozenge_K)_l$}
&$\sigma'{\cup}\sigma''\subset\sigma$,\, $\omega\subset\omega'{\cup}\omega''$&otherwise\\
\hline
{\rule[-3mm]{0mm}{8mm}$(\tilde\lozenge_K)_l$}
&$\sigma'{\cup}\sigma''=\sigma$,\, $\omega=\omega'{\cup}\omega''$&otherwise\\
\hline
{\rule[-3mm]{0mm}{8mm}$(\overline\lozenge_K)_l$}
&$\sigma'{\cup}\sigma''=\sigma,\,\sigma'{\cap}\sigma''=\emptyset$,\, $\omega=\omega'{\cup}\omega'',\,\omega'{\cap}\omega''=\emptyset$&otherwise\\
\hline
{\rule[-3mm]{0mm}{8mm}$(\bar\lozenge_K)_l$}
&$\sigma'{\cup}\sigma''=\sigma,\,\sigma'{\cap}\sigma''=\emptyset$,\,
$\omega{\setminus}(\omega'{\cup}\omega'')\in K,\,\omega'{\cap}\omega''=\emptyset$&otherwise\\
\hline
\end{tabular}}\vspace{3mm}

For the right universal (right normal, $\cdots$) objects, the local (co)products satisfy the following table with
$\sigma=\sigma'=\sigma''=\emptyset$.\vspace{3mm}\\
\hspace*{28mm}{\rm \begin{tabular}{|c|c|c|}
\hline
{\rule[-3mm]{0mm}{8mm}}
&$=(\widehat\lozenge_K)_l$&$=0$\\
\hline
{\rule[-3mm]{0mm}{8mm}$(\widehat\lozenge_K)_l$}
&$\omega{\setminus}(\omega'{\cup}\omega'')\in K$&otherwise\\
\hline
{\rule[-3mm]{0mm}{8mm}$(\w\lozenge_K)_l$}
&$\omega\subset\omega'{\cup}\omega''$&otherwise\\
\hline
{\rule[-3mm]{0mm}{8mm}$(\tilde\lozenge_K)_l$}
&$\omega=\omega'{\cup}\omega''$&otherwise\\
\hline
{\rule[-3mm]{0mm}{8mm}$(\overline\lozenge_K)_l$}
&$\omega=\omega'{\cup}\omega'',\,\omega'{\cap}\omega''=\emptyset$&otherwise\\
\hline
{\rule[-3mm]{0mm}{8mm}$(\bar\lozenge_K)_l$}
&$\omega{\setminus}(\omega'{\cup}\omega'')\in K,\,\omega'{\cap}\omega''=\emptyset$&otherwise\\
\hline
\end{tabular}}\vspace{3mm}

Proof}\, Similar to the proof of Theorem~7.6.
\hfill$\Box$\vspace{3mm}

In general, it is very complicated to determine when the local coproduct $(\underline{\psi})_l$ of a total chain coalgebra is non-zero.
For a $\DD$-coalgebra $(A_*^{\DD},\psi)$, denote by $L(A_*^{\DD},\psi)$ the subset of $\DD_3=\DD{\times}\DD{\times}\DD$
such that the local coproduct $\psi^s_{s',s''}\colon A_*^s\to A_*^{s'}{\otimes}A_*^{s''}$ is non-zero
if and only if $(s,s',s'')\in L(A_*^{\DD},\psi)$.
Then we have
$$L(T_*^{\XX_m},\underline{\psi})=L(T_*^{\XX}\!,\psi_1){\times}{\cdots}{\times}L(T_*^{\XX}\!,\psi_m)\in
(\XX_3){\times}{\cdots}{\times}(\XX_3)\,\,(m\,{\rm fold}).$$
But $(\XX_3){\times}{\cdots}{\times}(\XX_3)=\XX_m{\times}\XX_m{\times}\XX_m$.
So we have another form $L'(T_*^{\XX_m},\underline{\psi})$ of  $L(T_*^{\XX_m},\underline{\psi})$ in $\XX_m{\times}\XX_m{\times}\XX_m$.
This new form is hard to get. For example, take all $\psi_k$ to be the same universal coproduct.
Then by taking $K=\Delta\!^{[m]}$ in Theorem~7.6, we have that $L'(T_*^{\XX_m},\widehat\psi^{(m)})$ is the set
$$\{(\sigma,\omega,\sigma',\omega',\sigma'',\omega'')\in\XX_m{\times}\XX_m{\times}\XX_m\,\,|\,\,
(\sigma'{\cup}\sigma''){\setminus}\sigma\subset\omega{\setminus}(\omega'{\cup}\omega'')\}.$$
So we only need the total chain coalgebras (cochain algebras) in the above theorem in actual computations.
\vspace{3mm}

{\bf Example~7.8} Let $K=\Delta\!^S$ with $S\subset[m]$ regarded as a simplicial complex on $[m]$.
Denote by $1_{\sigma,\omega}$ the generator of $H_0^{\sigma,\omega}(K)$.
Then we have group isomorphism
$$H^{\,*}_{\!\XX_m}(\Delta\!^S)=\oplus_{\sigma\subset S,\,S{\cap}\omega=\emptyset}\,\Bbb Z(1_{\sigma,\omega}).$$
By definition, the products $\widehat\cup_{K},\w\cup_{K},\bar\cup_{K},\overline\cup_{K},\bar\cup_{K}$ on $H^{\,*}_{\!\XX_m}(K)$ are given by
$$1_{\sigma'\!,\,\omega'}\,\widehat\cup_{\Delta\!^S}\,1_{\sigma''\!,\,\omega''}
=1_{\sigma'\!,\,\omega'}\,\w\cup_{\Delta\!^S}\,1_{\sigma''\!,\,\omega''}
=\Sigma_{\sigma'\cup\sigma''\subset\sigma\subset S,\,\omega\subset\omega'\cup\omega''}1_{\sigma\!,\,\omega},$$
$$1_{\sigma'\!,\,\omega'}\,\tilde\cup_{\Delta\!^S}\,1_{\sigma''\!,\,\omega''}=1_{\sigma'\cup\sigma'',\,\omega'\cup\omega''},$$
$$1_{\sigma'\!,\,\omega'}\,\overline\cup_{\Delta\!^S}\,1_{\sigma''\!,\,\omega''}=
\left\{\begin{array}{cl}
1_{\sigma'\cup\sigma'',\,\omega'\cup\omega''}&{\rm if\,\,}\sigma'{\cap}\sigma''=\emptyset,\,\omega'{\cap}\omega''=\emptyset,\\
0&{\rm otherwise},
\end{array}\right.$$
$$1_{\sigma'\!,\,\omega'}\,\bar\cup_{\Delta\!^S}\,1_{\sigma''\!,\,\omega''}=
\left\{\begin{array}{cl}
\Sigma_{\omega\subset\omega'\cup\omega''}\,
1_{\sigma'\cup\sigma'',\,\omega}&{\rm if\,\,}\sigma'{\cap}\sigma''=\emptyset,\,\omega'{\cap}\omega''=\emptyset,\\
0&{\rm otherwise}.
\end{array}\right.\vspace{3mm}$$

Take $S=\emptyset$, i.e., $K=\{\emptyset\}$ is regarded as a simplicial complex on $[m]$.
Then $H^{\,*}_{\!\XX_m}(\{\emptyset\})=\oplus_{\omega\subset[m]}\Bbb Z(1_{\emptyset,\omega})$ and we have
$$1_{\emptyset,\omega'}\,\widehat\cup_{\{\emptyset\}}\,1_{\emptyset,\omega''}
=1_{\emptyset,\omega'}\,\w\cup_{\{\emptyset\}}\,1_{\emptyset,\omega''}
=\Sigma_{\omega\subset\omega'\cup\omega''}1_{\emptyset,\omega},$$
$$1_{\emptyset,\omega'}\,\tilde\cup_{\{\emptyset\}}\,1_{\emptyset,\omega''}=1_{\emptyset,\,\omega'\cup\omega''},$$
$$1_{\emptyset,\omega'}\,\overline\cup_{\{\emptyset\}}\,1_{\emptyset,\omega''}=
\left\{\begin{array}{cl}
1_{\emptyset,\,\omega'\cup\omega''}&{\rm if\,\,}\omega'{\cap}\omega''=\emptyset,\\
0&{\rm otherwise},
\end{array}\right.$$
$$1_{\emptyset,\omega'}\,\bar\cup_{\{\emptyset\}}\,1_{\emptyset,\omega''}=
\left\{\begin{array}{cl}
\Sigma_{\omega\subset\omega'\cup\omega''}\,
1_{\emptyset,\,\omega}&{\rm if\,\,}\omega'{\cap}\omega''=\emptyset,\\
0&{\rm otherwise}.
\end{array}\right.\vspace{3mm}$$

Take $S=[m]$, i.e., $K=\Delta\!^{[m]}$ is regarded as a simplicial complex on $[m]$.
Then $H^{\,*}_{\!\XX_m}(\Delta\!^{[m]})=\oplus_{\sigma\subset[m]}\Bbb Z(1_{\sigma,\emptyset})$ and we have
$$1_{\sigma',\emptyset}\,\widehat\cup_{\Delta\!^{[m]}}\,1_{\sigma'',\emptyset}
=1_{\sigma',\emptyset}\,\w\cup_{\Delta\!^{[m]}}\,1_{\sigma'',\emptyset}
=\Sigma_{\sigma\subset\sigma'\cup\sigma''}1_{\sigma,\emptyset},$$
$$1_{\sigma',\emptyset}\,\tilde\cup_{\Delta\!^{[m]}}\,1_{\sigma'',\emptyset}=1_{\sigma'\cup\sigma'',\emptyset},$$
$$1_{\sigma',\emptyset}\,\overline\cup_{\Delta\!^{[m]}}\,1_{\sigma'',\emptyset}=
1_{\sigma',\emptyset}\,\bar\cup_{\Delta\!^{[m]}}\,1_{\sigma'',\emptyset}=
\left\{\begin{array}{cl}
1_{\sigma'\cup\sigma'',\emptyset}&{\rm if\,\,}\sigma'{\cap}\sigma''=\emptyset,\\
0&{\rm otherwise}.
\end{array}\right.\vspace{3mm}$$

{\bf Example~7.9} Let $K$ be the $m$-gon, $m>3$.
We compute the right universal cohomology algebra $(H^{\,*}_{\!\RR_m}(K),\widehat\cup_K)$.
The vertex set of $K$ is $[m]$ with edges $\{i,i{+}1\}$ for $i\in\Bbb Z_{m}$,
where $\Bbb Z_m$ is the group of integers modular $m$ regarded only as a set.

Denote by $H_*^{\omega}=H_*^{\emptyset,\omega}(K)$, $H^*_{\omega}=H^{*}_{\emptyset,\omega}(K)$.
We use the same symbol to denote the base element and its dual element.
Then all the non-zero $H_k^\omega$ and $H^k_\omega$ are the following.

(1) $H_0^{\emptyset}=\Bbb Z(1)$ with $1=[\emptyset]$. Dually,
$H^0_{\emptyset}=\Bbb Z(1)$ with $1=[\emptyset]$.

(2) Let $\omega$ be a subset of $[m]$ with connected component $\omega_1,\cdots,\omega_k$ ($k>1$).
For any $a\in\omega_i$ and $b\in\omega_j$, $[\{a\}{-}\{b\}]$ is a homology class in $H_1^\omega$ independent of the choice of $a$ and $b$.
The homology group $H_1^\omega$ is the group generated by all $\omega_i{-}\omega_j$ modulo the subgroup generated by
all $[\{a\}{-}\{b\}]+[\{b\}{-}\{c\}]-[\{a\}{-}\{c\}]$ for $a,b,c$ in different connected components.
Take the base of $H_1^\omega$ to be $h_{\omega,1},{\cdots},h_{\omega,{k-1}}$, where $h_{\omega,s}=[\{a\}{-}\{b\}]$ with $a\in\omega_s$ and $b\in\omega_k$.

Dually, $[\Sigma_{u\in\omega_i}\{u\}]$ represents a cohomology class in $H^1_{\omega}$.
$H^1_{\omega}$ is the group generated by $[\Sigma_{u\in\omega_1}\{u\}],\cdots,[\Sigma_{u\in\omega_k}\{u\}]$ modulo the
subgroup generated by $[\Sigma_{u\in\omega_1}\{u\}]{+}\cdots{+}[\Sigma_{u\in\omega_k}\{u\}]$.
Take the base of $H^1_\omega$ to be $h_{\omega,1},{\cdots},h_{\omega,k-1}$,
where $h_{\omega,s}=[\Sigma_{u\in\omega_s}\{u\}]$.

(3) $H_2^{[m]}=\Bbb Z(\kappa)$ with $\kappa=[\Sigma_{i\in\Bbb Z_m}\{i,i{+}1\}]$.
Dually, $H^2_{[m]}=\Bbb Z(\kappa)$ with $\kappa=[\{i,i{+}1\}]$ for any $i\in\Bbb Z_m$.

By definition, the local product
$(\widehat\pi_K)_{\,\omega}^{\,\omega'\!,\,\omega''}\colon T^1_{\emptyset,\omega'}(K){\otimes}T^1_{\emptyset,\omega''}(K)\to T^k_{\emptyset,\omega}(K)$ is $0$ if $k<2$.
Since $H^2_{\omega}=0$ except $\omega=[m]$,
we have that the local product $(\widehat\cup_K)_{\,\omega}^{\,\omega'\!,\,\omega''}$ of $\widehat\cup_K$ is $0$ if $\omega\neq[m]$.
The local product $(\widehat\pi_K)_{\,[m]}^{\,\omega'\!,\,\omega''}$ is as follows.
$(\widehat\pi_K)_{\,[m]}^{\,\omega'\!,\,\omega''}(\{i\}{\otimes}\{i{+}1\})=\{i,i{+}1\}$,
$(\widehat\pi_K)_{\,[m]}^{\,\omega'\!,\,\omega''}(\{i\}{\otimes}\{i{-}1\})=-\{i{-}1,i\}$,
$(\widehat\pi_K)_{\,[m]}^{\,\omega'\!,\,\omega''}(\{i\}{\otimes}\{j\})=0$ if $i{-}j\neq\pm 1$ mod$m$.
Define sign product $*\colon\Bbb Z_m\times\Bbb Z_m\to\Bbb Z$ by
$$i*j=\left\{\begin{array}{rl}
1&{\rm if}\,\, j\equiv i{+}1\,\,{\rm mod}\,m,\\
-1&{\rm if}\,\, j\equiv i{-}1\,\,{\rm mod}\,m,\\
0&{\rm otherwise}.
\end{array}\right.$$
For subsets $A,B$ of $[m]$, define $A*B=\Sigma_{i\in A,\,j\in B}\,i*j$.
Then for $\omega'\subset[m]$ with connected components $\omega'_1,{\cdots},\omega'_s$
and $\omega''\subset[m]$ with connected components $\omega''_1,{\cdots},\omega''_t$,
$$h_{\omega'\!,i}\,{\widehat\cup}_K\,h_{\omega''\!,j}=
h_{\omega'\!,i}\,(\widehat\cup_K)_{\,[m]}^{\,\omega'\!,\,\omega''}h_{\omega''\!,j}=
(\omega'_i*\omega''_j)\kappa.$$
It is easy to check that all other products are $1\,{\widehat\cup}_K\,x=x\,{\widehat\cup}_K\,1=x$,
$h_{\omega,i}\,{\widehat\cup}_K\,\kappa=\kappa\,{\widehat\cup}_K\,h_{\omega,i}=0$
and $\kappa\,{\widehat\cup}_K\,\kappa=0$.

Dually, we have
$$\widehat\vartriangle_K(1)=1{\otimes}1,\,\,\,\widehat\vartriangle_K(h_{\omega,i})=1{\otimes}h_{\omega,i}+h_{\omega,i}{\otimes}1,$$
$$\widehat\vartriangle_K(\kappa)=1{\otimes}\kappa+\kappa{\otimes}1+\Sigma(\omega'_i*\omega''_j)h_{\omega'\!,_i}{\otimes}h_{\omega''\!,j},$$
where the sum is taken over all $\omega',\omega''\subset[m]$ with more than one connected components.
It is complicated to compute the above equality by the formula
$(\widehat\vartriangle_K)_l(\kappa)=[(\widehat\psi_K)_l(\Sigma_{i=1}^{m}\{i,i{+}1\})]$.

Specifically, the right strictly normal cohomology algebra $(H^*_{\RR_m}(K),\tilde\cup_K)$
($\cong H^*({\cal Z}(K;D^1\!,S^0))$) satisfies
$h_{\omega'\!,i}\,{\tilde\cup}_K\,h_{\omega''\!,j}=(\omega'_i*\omega''_j)\kappa$ if $\omega'{\cup}\omega''=[m]$ and
$h_{\omega'\!,i}\,{\tilde\cup}_K\,h_{\omega''\!,j}=0$ otherwise.
\vspace{3mm}

{\bf Example~7.10} Let $A$ be a simplicial complex on $[s]$.
Then for any $t>0$, it is also a simplicial complex on $[s{+}t]$.
To distinguish the two cases, we denote by $A_s$ the simplicial complex $A$ on $[s]$
and denote by $A_{s+t}$ the one on $[s{+}t]$.
Then $A_{s+t}=A_s*\{\emptyset\}$, where $\{\emptyset\}$ is regarded as a simplicial complex on $[t]$ as in Example~7.8.
For total chain coalgebras $(T_*^{\XX_s},\underline{\psi'})$ and $(T_*^{\XX_t},\underline{\psi''})$, we have
a total chain coalgebra $(T_*^{\XX_{s+t}},\underline{\psi})=(T_*^{\XX_s}{\otimes}T_*^{\XX_t},\underline{\psi'}{\otimes}\underline{\psi''})$. Then
$$(T_*^{\XX_{s+t}}(A_{s+t}),\psi_{A_{s+t}})\cong
(T_*^{\XX_s}(A_s){\otimes}T_*^{\XX_t}(\{\emptyset\}),\psi'_{A_s}{\otimes}\psi''_{\{\emptyset\}}).$$
Since $H_*^{\XX_t}(\{\emptyset\})$ is a free group, we have by K\"{u}nneth Theorem,
$$(H^{\,*}_{\!\XX_{s+t}}(A_{s+t}),\cup_{A_{s+t}})\cong
(H^{\,*}_{\!\XX_s}(A_s){\otimes}H^{\,*}_{\!\XX_t}(\{\emptyset\}),\cup'_{A_s}{\otimes}\cup''_{\{\emptyset\}}).
\vspace{3mm}$$

\section{Cohomology Algebra of Polyhedral Product Chain Complexes}\vspace{3mm}

\hspace*{5.5mm}{\bf Definition~8.1} A {\it split coalgebra pair} $(D_*,C_*)=\big((D_*,\psi_C,d),(C_*,\psi_D,d)\big)$ is
a pair of chain $\Lambda$-coalgebras satisfying the following conditions.

(1) The pair $(D_*,C_*)$ is homology split by Definition~3.2.

(2) The inclusion $\vartheta\colon (C_*,\psi_C,d)\to (D_*,\psi_D,d)$ is a chain $\Lambda$-coalgebra weak homomorphism by Definition~6.3.
\vspace{3mm}

{\bf Definition~8.2} Suppose for the split coalgebra pair $(D_*,C_*)$,
a given quotient homotopy equivalence $q$ (and so its restriction $q'$) as constructed in Theorem~3.5 is chosen.
Then all the groups and complexes in Definition~3.3 have (co)algebra structures defined as follows.

Let $\theta\colon(H_*(C_*),{\vartriangle}_C)\to(H_*(D_*),{\vartriangle}_D)$ be the homology coalgebra homomorphism induced by inclusion
and denote by
$${\mathpzc i}_{\,*}={\rm coim}\,\theta,\quad{\mathpzc n}_{\,*}={\rm ker}\,\theta,\quad
\overline{\mathpzc n}_{\,*}=\Sigma{\rm ker}\,\theta,\quad{\mathpzc e}_{\,*}={\rm coker}\,\theta.$$
Then we have group equalities $H_*(C_*)={\mathpzc i}_{\,*}\oplus{\mathpzc n}_{\,*}$,
$H_*^\XX(D_*,C_*)={\mathpzc i}_{\,*}\oplus{\mathpzc n}_{\,*}\oplus{\mathpzc e}_{\,*}$ and
$C_*^\XX(D_*,C_*)={\mathpzc i}_{\,*}\oplus{\mathpzc n}_{\,*}\oplus\overline{\mathpzc n}_{\,*}\oplus{\mathpzc e}_{\,*}$.

The {\it normal homology coalgebra} $(H_*^\XX(D_*,C_*),\vartriangle_{(D,C)})$ of $(D_*,C_*)$ (irrelevant to $q$) is given by
$${\vartriangle}_{(D,C)}(x)=\left\{\begin{array}{cl}
{\vartriangle}_C(x)&{\rm if}\,\,x\in {\mathpzc i}_{\,*},\vspace{1mm}\\
{\vartriangle}_C(x)&{\rm if}\,\,x\in {\mathpzc n}_{\,*},\vspace{1mm}\\
{\vartriangle}_D(x)&{\rm if}\,\,x\in {\mathpzc e}_{\,*}.
\end{array}\right.$$

The {\it normal cohomology algebra} $(H^{\,*}_{\!\XX}(D_*,C_*),{\triangledown}_{(D,C)})$ of $(D_*,C_*)$ is
the dual algebra of $(H_*^\XX(D_*,C_*),\vartriangle_{(D,C)})$.

The {\it character chain coalgebra} $(C_*^\XX(D_*,C_*),{\vartriangle}_q^{\!\XX},d)$ of $(D_*,C_*)$ with respect to $q$
is the chain $(\XX\!{\times}\Lambda)$-coalgebra defined as follows.
The restriction of ${\vartriangle}_q^{\!\XX}$ on $H_*^\XX(D_*,C_*)$ is ${\vartriangle}_{(D,C)}$.
For $x\in{\mathpzc n}_{\,*}$ with ${\vartriangle}_C(x)=\Sigma x'_i{\otimes}x''_i+\Sigma y'_j{\otimes}y''_j$,
where each $x'_i\in{\mathpzc n}_{\,*}$ and each $y'_j\notin{\mathpzc n}_{\,*}$ but $y''_j\in{\mathpzc n}_{\,*}$, define
$${\vartriangle}_q^{\!\XX}(\overline x)=\Sigma\overline x'_i{\otimes}x''_i+\Sigma(-1)^{|y'_j|}y'_j{\otimes}\overline y''_j
+\xi(q{\otimes}q)\psi_D(\overline x),$$
where $\xi$ is defined as follows. Let $C_*^{\XX}=C_*^\XX(D_*,C_*)$ and $N={\mathpzc n}_{\,*}{\otimes}C_*^\XX+C_*^\XX\!{\otimes}{\mathpzc n}_{\,*}$
and $N\oplus H=C_*^\XX\!{\otimes}C_*^\XX$. So $H_*(N)=0$ and $H\cong H_*(C_*^\XX\!{\otimes}C_*^\XX)$.
Then $\xi$ is the projection from $C_*^\XX\!{\otimes}C_*^\XX$ to $H$.

The {\it indexed homology coalgebra} $(H_*^\XX(D_*,C_*),{\vartriangle}_q)$ of $(D_*,C_*)$ with respect to $q$ is
the $(\XX\!{\times}\Lambda)$-coalgebra given by
$${\vartriangle}_q(x)=\left\{\begin{array}{ll}
{\vartriangle}_C(x)&{\rm if}\,\,x\in{\mathpzc i}_{\,*},\vspace{1mm}\\
{\vartriangle}_C(x){+}\xi(q{\otimes}q)\psi_D(\overline x)&{\rm if}\,\,x\in {\mathpzc n}_{\,*},\vspace{1mm}\\
{\vartriangle}_D(x)&{\rm if}\,\,x\in {\mathpzc e}_{\,*}.
\end{array}\right.$$

The {\it indexed cohomology algebra} $(H^{\,*}_{\!\XX}(D_*,C_*),{\triangledown}_{\!q})$ of $(D_*,C_*)$ with respect to $q$ is
the dual algebra of $(H_*^\XX(D_*,C_*),{\vartriangle}_q)$.

$(D_*,C_*)$ is called {\it normal} if there is a homotopy equivalence $q$ such that ${\vartriangle}_q={\vartriangle}_{(D,C)}$
and ${\triangledown}_{\!q}={\triangledown}_{\!(D,C)}$.

The {\it atom coproduct} $\psi_q$ and the {\it atom chain coalgebra $(T_*^\XX\!,\psi_q)$} with respect to $q$ are defined as follows.
The generators ${\mathpzc e},{\bar{\mathpzc n}},{\mathpzc n},{\mathpzc i}$ of $T_*^\XX$
and the group summands ${\mathpzc e}_{\,*},\bar{\mathpzc n}_{\,*},{\mathpzc n}_{\,*},{\mathpzc i}_{\,*}$ of $C_*^{\XX}\!(D_*,C_*)$
have a 1-1 correspondence ${\mathpzc s}\to{\mathpzc s}_{\,*}$.
For ${\mathpzc s}\in T_*^\XX$, suppose ${\vartriangle}_q^{\!\XX}({\mathpzc s}_*)\subset\oplus_i\,({\mathpzc s}'_i)_*{\otimes}({\mathpzc s}''_i)_*$
but ${\vartriangle}_q^{\!\XX}({\mathpzc s}_*){\cap}(({\mathpzc s}'_i)_*{\otimes}({\mathpzc s}''_i)_*)\neq 0$ for each $i$,
then define $\psi_q({\mathpzc s})=\Sigma_i\,{\mathpzc s}'_i{\otimes}{\mathpzc s}''_i$.

The index set of $(D_*,C_*)$ is $\SS=\{{\mathpzc s}\in T_*^\XX\,|\,{\mathpzc s}_{\,*}\neq 0\}$.
By definition, $\psi_q({\mathpzc s})=0$ if ${\mathpzc s}\notin\SS$.
\vspace{3mm}

Note that even if both $H_*(C_*)$ and $H_*(D_*)$ are coassociative coalgebras,
the normal homology coalgebra $(H_*^{\XX}(D_*,C_*),\vartriangle_{(D,C)})$ may not be coassociative.
If $\theta\colon H_*(C_*)\to H_*(D_*)$ is an epimorphism and $(D_*,C_*)$ is normal,
then the normal homology coalgebra $(H_*^\XX(D_*,C_*),{\vartriangle}_{(D,C)})$ is just $(H_*(C_*),{\vartriangle}_C)$ when the index $\XX$ is neglected.
\vspace{3mm}

{\bf Theorem~8.3} {\it For a split coalgebra pair $(D_*,C_*)$,
all the chain complex homomorphisms in Theorem~3.5 are chain coalgebra homomorphisms when the quotient homotopy equivalence $q$ is chosen.
Precisely, we have the following commutative diagram
$$\begin{array}{ccc}
(C_*,\psi_C)&\stackrel{q'}{\longrightarrow}&(H_*(C_*),{\vartriangle}_C)\vspace{1mm}\\
\cap&&\cap\\
(D_*,\psi_D)&\stackrel{q}{\longrightarrow}&(C_*^{\XX}(D_*,C_*),{\vartriangle}_q^{\!\XX})
\end{array}
$$
of chain $\Lambda$-coalgebra weak homomorphisms (index $\XX$ neglected) such that
$q$ and $q'$ induce homology coalgebra isomorphisms.
We have the commutative diagram
$$\begin{array}{ccc}
(H_*(C_*),{\vartriangle}_C)&\stackrel{\phi'}{\cong}&(S_*^\XX\,\widehat\otimes\,
H_*^{\XX}\!(D_*,C_*),\,\widehat\psi\,\widehat\otimes\,{\vartriangle}_q)\vspace{1mm}\\
\cap&&\cap\\
(C_*^{\XX}(D_*,C_*),{\vartriangle}_q^{\!\XX})&\stackrel{\phi}{\cong}&
(T_*^\XX\,\widehat\otimes\, H_*^{\XX}\!(D_*,C_*),\,\widehat\psi\,\widehat\otimes\,{\vartriangle}_q)
\end{array}
$$
of chain $(\XX{\times}\Lambda)$-coalgebra isomorphisms, where $\widehat\psi$ is the universal coproduct in Definition~7.1.
For any coproduct $\psi$ such that $\psi_q\prec\psi\prec\widehat\psi$ (by Definition~6.10),
we have the commutative diagram
$$\begin{array}{ccc}
(S_*^\XX\,\widehat\otimes\,H_*^{\XX}\!(D_*,C_*),\,\widehat\psi\,\widehat\otimes\,{\vartriangle}_q)&=&
(S_*^\TT\,\widehat\otimes\,H_*^{\TT}\!(D_*,C_*),\,\psi\,\widehat\otimes\,{\vartriangle}_q)\vspace{1mm}\\
\cap&&\cap\\
(T_*^\XX\,\widehat\otimes\, H_*^{\XX}\!(D_*,C_*),\,\widehat\psi\,\widehat\otimes\,{\vartriangle}_q)&=&
(T_*^\TT\,\widehat\otimes\, H_*^{\TT}\!(D_*,C_*),\,\psi\,\widehat\otimes\,{\vartriangle}_q)
\end{array}$$
of chain $\Lambda$-coalgebra isomorphisms (index $\XX,\TT$ neglected),
where $(-)^\TT$ means the restriction coalgebra of $(-)^\XX$.
\vspace{2mm}

Proof}\, Consider the diagram
$$\begin{array}{ccc}
(C_*,\psi_C)&\stackrel{q'}{\longrightarrow}&(H_*(C_*),{\vartriangle}_C)\,\vspace{1mm}\\
^{\vartheta}\downarrow\quad&&^{\vartheta'}\downarrow\quad\\
(D_*,\psi_D)&\stackrel{q}{\longrightarrow}&(C_*^{\XX}(D_*,C_*),{\vartriangle}_q^{\!\XX}),
\end{array}$$
where $\vartheta'$ is the restriction of $q\vartheta$ on $H_*(C_*)$.
The inclusion $\vartheta$ is a chain $\Lambda$-coalgebra weak homomorphism by (2) of Definition~8.1.
$q'$ is naturally a chain $\Lambda$-coalgebra homomorphism that induces homology coalgebra isomorphism.
If $q$ is a chain $\Lambda$-coalgebra weak homomorphism, then as a chain homotopy equivalence,
it induces a homology coalgebra isomorphism and as a restriction of $q\vartheta$, $\vartheta'$ is also a chain $\Lambda$-coalgebra weak homomorphism.
So we only need prove that $q$ is a chain $\Lambda$-coalgebra weak homomorphism.

Let everything be as in Definition~8.2. Construct homotopy $s\colon D_*\to\Sigma(C_*^\XX{\otimes}C_*^\XX)$ such that $ds+sd={\vartriangle}_q^{\!\XX}q-(q{\otimes}q)\psi_D$ as follows.

For $x\in {\mathpzc i}_{\,*}$, we have ${\vartriangle}_q^{\!\XX}(q(x)){-}(q{\otimes}q)(\psi_D(x))={\vartriangle}_C(x){-}{\vartriangle}_D(x)\in N$.
So there is $y\in N$ such that $dy={\vartriangle}_C(x){-}{\vartriangle}_D(x)$. Define $s(x)=y$.
Then $(ds{+}sd)(x)=({\vartriangle}_q^{\!\XX}q{-}(q{\otimes}q)\psi_D)(x)$.

For $x\in {\mathpzc e}_{\,*}$, we have ${\vartriangle}_q^{\!\XX}(q(x)){-}(q{\otimes}q)(\psi_D(x))={\vartriangle}_D(x){-}{\vartriangle}_D(x)=0$.
Define $s(x)=0$. Then $(ds+sd)(x)=({\vartriangle}_q^{\!\XX}q-(q{\otimes}q)\psi_D)(x)$.

For $x\in{\mathpzc n}_{\,*}$ and $\overline x\in\overline{\mathpzc n}_{\,*}$,
we have ${\vartriangle}_q^{\!\XX}(q(x))=(q{\otimes}q)(\psi_D(x))={\vartriangle}_C(x)$ and
$w={\vartriangle}_q^{\!\XX}(q(\overline x)){-}\xi((q{\otimes}q)(\psi_D(\overline x)))\in N$ with $dw={\vartriangle}_C(x)$.
We also have  $v=(q{\otimes}q)(\psi_D(\overline x)){-}\xi((q{\otimes}q)(\psi_D(\overline x)))\in N$ with $dv={\vartriangle}_C(x)$.
So there is $y\in N$ such that $dy=w{-}v$. Define $s(x)=0$ and $s(\overline x)=y$.
Then $(ds+sd)(x)=({\vartriangle}_q^{\!\XX}q-(q{\otimes}q)\psi_D)(x)$ and
$(ds+sd)(\overline x)=({\vartriangle}_q^{\!\XX}q-(q{\otimes}q)\psi_D)(\overline x)$.

For $x\in{\rm ker}\,q$ such that $dx=0$, denote by $\overline x$ the unique element in ${\rm ker}\,q$ such that $d\overline x=x$.
Define $s(x)={-}(q{\otimes}q)(\psi_D(\overline x))$ and $s(\overline x)=0$.
Then $(ds+sd)(x)=({\vartriangle}_q^{\!\XX}q-(q{\otimes}q)\psi_D)(x)$ and
$(ds+sd)(\overline x)=({\vartriangle}_q^{\!\XX}q-(q{\otimes}q)\psi_D)(\overline x)$.

So $q$ is a chain $\Lambda$-coalgebra weak homomorphism.

Now we prove  $\phi\colon(C_*^{\XX},{\vartriangle}_q^{\!\XX})\to
(T_*^\XX\,\widehat\otimes\, H_*^{\XX},\,\widehat\psi\,\widehat\otimes\,{\vartriangle}_q)$
is a $(\XX\!{\times}\Lambda)$-coalgebra isomorphism.
We use the specific symbols to denote elements of the corresponding groups as in the following table.
\vspace{2mm}\\
\hspace*{5mm}
\begin{tabular}{|c|c|c|c|c|}
\hline
{\rule[-2mm]{0mm}{6mm}elements\,\,of}&${\mathpzc e}_{\,*}{=}{\rm coker}\theta$&$\bar{\mathpzc n}_{\,*}{=}\Sigma\,{\rm ker}\theta$
&${\mathpzc n}_{\,*}{=}{\rm ker}\theta$&${\mathpzc i}_{\,*}{=}{\rm coim}\theta$\\
\hline
{\rule[-2mm]{0mm}{7mm}symbols}&$e,e'_1,e''_1,\cdots$&${\overline n},{\overline n}'_1,{\overline n}''_1,\cdots$
&$n,n'_1,n''_1,\cdots$&$i,i'_1,i''_1,\cdots$\\
\hline
\end{tabular}
\vspace{3mm}\\
Then we have

${\vartriangle}_q^{\!\XX}({\mathpzc i}_{\,*})\subset{\mathpzc i}_{\,*}{\scriptstyle\otimes}{\mathpzc i}_{\,*}\oplus{\mathpzc i}_{\,*}{\scriptstyle\otimes}{\mathpzc n}_{\,*}
\oplus{\mathpzc n}_{\,*}{\scriptstyle\otimes}{\mathpzc i}_{\,*}\oplus{\mathpzc n}_{\,*}{\scriptstyle\otimes}{\mathpzc n}_{\,*}$.

${\vartriangle}_q^{\!\XX}({\mathpzc n}_{\,*})\subset{\mathpzc n}_{\,*}{\scriptstyle\otimes}{\mathpzc n}_{\,*}\oplus{\mathpzc n}_{\,*}{\scriptstyle\otimes}{\mathpzc i}_{\,*}
\oplus{\mathpzc i}_{\,*}{\scriptstyle\otimes}{\mathpzc n}_{\,*}$.

${\vartriangle}_q^{\!\XX}({\overline{\mathpzc n}_{\,*}})\subset{\overline{\mathpzc n}_{\,*}}{\scriptstyle\otimes}{\mathpzc n}_{\,*}\oplus{\overline{\mathpzc n}_{\,*}}{\scriptstyle\otimes}{\mathpzc i}_{\,*}\oplus{\mathpzc i}_{\,*}{\scriptstyle\otimes}{\overline{\mathpzc n}_{\,*}}\oplus
{\mathpzc e}_{\,*}{\scriptstyle\otimes}{\mathpzc e}_{\,*}\oplus{\mathpzc e}_{\,*}{\scriptstyle\otimes}{\mathpzc i}_{\,*}\oplus{\mathpzc i}_{\,*}{\scriptstyle\otimes}{\mathpzc e}_{\,*}\oplus{\mathpzc i}_{\,*}{\scriptstyle\otimes}{\mathpzc i}_{\,*}$.

${\vartriangle}_q^{\!\XX}({\mathpzc e}_{\,*})\subset{\mathpzc e}_{\,*}{\scriptstyle\otimes}{\mathpzc e}_{\,*}\oplus{\mathpzc e}_{\,*}{\scriptstyle\otimes}{\mathpzc i}_{\,*}\oplus{\mathpzc i}_{\,*}{\scriptstyle\otimes}{\mathpzc e}_{\,*}\oplus{\mathpzc i}_{\,*}{\scriptstyle\otimes}{\mathpzc i}_{\,*}$.

In the following computations, all $\Sigma$ are omitted, i.e., $x{\otimes}y$ means $\Sigma\, x{\otimes}y$.
Suppose for $e,{\overline n},n,i\in C_*^\XX$,

${\vartriangle}_q^{\!\XX}(i)= i'_1{\otimes}i''_1+i'_2{\otimes}n''_2+n'_3{\otimes}i''_3+n'_4{\otimes}n''_4$,

${\vartriangle}_q^{\!\XX}(n)= n'_1{\otimes}n''_1+n'_2{\otimes}i''_2+i'_3{\otimes}n''_3$,

${\vartriangle}_q^{\!\XX}({\overline n})= {\overline n}'_1{\otimes}n''_1+{\overline n}'_2{\otimes}i''_2+
{\scriptstyle(-1)}^{|i'_3|}i'_3{\otimes}{\overline n}''_3\,+\,
e'_4{\otimes}e''_4\,+\,e'_5{\otimes}i''_5\,+\,i'_{6}{\otimes}e''_{6}\,+\,i'_{7}{\otimes}i''_{7}$,

${\vartriangle}_q^{\!\XX}(e)= e'_{1}{\otimes}e''_{1}+e'_{2}{\otimes}i''_{2}+i'_{3}{\otimes}e''_{3}+i'_{4}{\otimes}i''_{4}$,\\
where $d{\overline x}=x$.
Then

${\vartriangle}_q(i)= i'_1{\otimes}i''_1+i'_2{\otimes}n''_2+n'_3{\otimes}i''_3+n'_4{\otimes}n''_4$,

${\vartriangle}_q(n)= n'_1{\otimes}n''_1+n'_2{\otimes}i''_2+i'_3{\otimes}n''_3+\,
e'_4{\otimes}e''_4\,+\,e'_5{\otimes}i''_5\,+\,i'_{6}{\otimes}e''_{6}\,+\,i'_{7}{\otimes}i''_{7}$,

${\vartriangle}_q(e)= e'_{1}{\otimes}e''_{1}+e'_{2}{\otimes}i''_{2}+i'_{3}{\otimes}e''_{3}+i'_{4}{\otimes}i''_{4}$.

For simplicity, $x{\otimes}y$ is abbreviated to $xy$ and $x\widehat\otimes y$ is abbreviated to
$x{\scriptstyle\wedge}y$ in the following computation.

\hspace*{20mm}$(\widehat\psi\,\widehat\otimes\,{\vartriangle}_q)(\phi(i))=(\widehat\psi\,\widehat\otimes\,{\vartriangle}_q)({\mathpzc i}{\scriptstyle\wedge} i)\\
\hspace*{20mm}=({\mathpzc i}\hspace{0.1mm}{\mathpzc i}{+}{\mathpzc i}\hspace{0.1mm}{\mathpzc n}
{+}{\mathpzc n}\hspace{0.1mm}{\mathpzc i}{+}{\mathpzc n}\hspace{0.1mm}{\mathpzc n}){\scriptstyle\wedge}
(i'_1i''_1{+}i'_2n''_2{+}n'_3i''_3{+}n'_4n''_4)\\
\hspace*{20mm}=({\mathpzc i}{\scriptstyle\wedge} i'_1)({\mathpzc i}{\scriptstyle\wedge} i''_1)
{+}({\mathpzc i}{\scriptstyle\wedge} i'_2)({\mathpzc n}{\scriptstyle\wedge} n''_2)
{+}({\mathpzc n}{\scriptstyle\wedge} n'_3)({\mathpzc i}{\scriptstyle\wedge} i''_3)
{+}({\mathpzc n}{\scriptstyle\wedge} n'_4)({\mathpzc n}{\scriptstyle\wedge} n''_4)\\
\hspace*{20mm}=\phi(i'_1)\phi(i''_1)
{+}\phi(i'_2)\phi(n''_2){+}\phi(n'_3)\phi(i''_3){+}\phi(n'_4)\phi(n''_4)\\
\hspace*{20mm}=(\phi{\otimes}\phi)({\vartriangle}_q^{\!\XX}(i))$,

\hspace*{20mm}$(\widehat\psi\,\widehat\otimes\,{\vartriangle}_q)(\phi(n))=(\widehat\psi\,\widehat\otimes\,{\vartriangle}_q)({\mathpzc n}{\scriptstyle\wedge} n)\\
\hspace*{20mm}=({\mathpzc n}\hspace{0.1mm}{\mathpzc n}{+}{\mathpzc n}\hspace{0.1mm}{\mathpzc i}{+}{\mathpzc i}\hspace{0.1mm}{\mathpzc n}){\scriptstyle\wedge}
(n'_1n''_1{+}n'_2i''_2{+}i'_3n''_3{+}e'_4e''_4{+}e'_5i''_5{+}i'_{6}e''_{6}{+}i'_{7}i''_{7})\\
\hspace*{20mm}=({\mathpzc n}{\scriptstyle\wedge} n'_1)({\mathpzc n}{\scriptstyle\wedge} n''_1){+}
({\mathpzc n}{\scriptstyle\wedge} n'_2)({\mathpzc i}{\scriptstyle\wedge} i''_2){+}
({\mathpzc i}{\scriptstyle\wedge} i'_3)({\mathpzc n}{\scriptstyle\wedge} n''_3)\\
\hspace*{20mm}=\phi(n'_1)\phi(n''_1){+}\phi(n'_2)\phi(i''_2){+}\phi(i'_3)\phi(n''_3)\\
\hspace*{20mm}=(\phi{\otimes}\phi)({\vartriangle}_q^{\!\XX}(n))$,

\hspace*{20mm}$(\widehat\psi\hspace{0.1mm}\widehat\otimes\hspace{0.1mm}{\vartriangle}_q)(\phi({\overline n}))=
(\widehat\psi\hspace{0.1mm}\widehat\otimes\hspace{0.1mm}{\vartriangle}_q)({\bar{\mathpzc n}}{\scriptstyle\wedge} n)\\
\hspace*{20mm}=({\bar{\mathpzc n}}\hspace{0.1mm}{\mathpzc n}{+}{\bar{\mathpzc n}}\hspace{0.1mm}{\mathpzc i}{+}{\mathpzc i}\hspace{0.1mm}{\bar{\mathpzc n}}
{+}{\mathpzc e}\hspace{0.1mm}{\mathpzc e}{+}{\mathpzc e}\hspace{0.1mm}{\mathpzc i}{+}{\mathpzc i}\hspace{0.1mm}{\mathpzc e}{+}{\mathpzc i}\hspace{0.1mm}{\mathpzc i})$$\\
\hspace*{30mm}{\scriptstyle\wedge}
(n'_1n''_1{+}n'_2i''_2{+}i'_3n''_3{+}e'_4e''_4{+}e'_5i''_5{+}i'_{6}e''_{6}{+}i'_{7}i''_{7})\\
\hspace*{20mm}=({\bar{\mathpzc n}}{\scriptstyle\wedge} n'_1)({\mathpzc n}{\scriptstyle\wedge} n''_1){+}
({\bar{\mathpzc n}}{\scriptstyle\wedge} n'_2)({\mathpzc i}{\scriptstyle\wedge} i''_2){+}
(-1)^{|i'_3||{\bar{\mathpzc n}}|}({\mathpzc i}{\scriptstyle\wedge} i'_3)({\bar{\mathpzc n}}{\scriptstyle\wedge} n''_3)\\
\hspace*{23mm}{+}({\mathpzc e}{\scriptstyle\wedge} e'_4)({\mathpzc e}{\scriptstyle\wedge} e''_4){+}
({\mathpzc e}{\scriptstyle\wedge} e'_5)({\mathpzc i}{\scriptstyle\wedge} i''_5){+}
({\mathpzc i}{\scriptstyle\wedge} i'_{6})({\mathpzc e}{\scriptstyle\wedge} e''_{6}){+}
({\mathpzc i}{\scriptstyle\wedge} i'_{7})({\mathpzc i}{\scriptstyle\wedge} i''_{7})\\
\hspace*{20mm}=\phi({\overline n}'_1)\phi(n''_1){+}\phi({\overline n}'_2)\phi(i''_2){+}
(-1)^{|i'_3|}\phi(i'_3)\phi({\overline n}''_3)\\
\hspace*{23mm}{+}\phi(e'_4)\phi(e''_4){+}
\phi(e'_5)\phi(i''_5){+}\phi(i'_{6})\phi(e''_{6}){+}\phi(i'_{7})\phi(i''_{7})\\
\hspace*{20mm}=(\phi{\otimes}\phi)({\vartriangle}_q^{\!\XX}({\overline n}))$,

$\hspace*{20mm}(\widehat\psi\,\widehat\otimes\,{\vartriangle}_q)(\phi(e))=(\widehat\psi\,\widehat\otimes\,{\vartriangle}_q)({\mathpzc e}{\scriptstyle\wedge} e)\\
\hspace*{20mm}=({\mathpzc e}\hspace{0.1mm}{\mathpzc e}{+}{\mathpzc e}\hspace{0.1mm}{\mathpzc i}{+}
{\mathpzc i}\hspace{0.1mm}{\mathpzc e}{+}{\mathpzc i}\hspace{0.1mm}{\mathpzc i}){\scriptstyle\wedge}
(e'_{1}e''_{1}{+}e'_{2}i''_{2}{+}i'_{3}e''_{3}{+}i'_{4}i''_{4})\\
\hspace*{20mm}=({\mathpzc e}{\scriptstyle\wedge} e'_{1})({\mathpzc e}{\scriptstyle\wedge} e''_{1})
{+}({\mathpzc e}{\scriptstyle\wedge} e'_{2})({\mathpzc i}{\scriptstyle\wedge} i''_{2})
{+}({\mathpzc i}{\scriptstyle\wedge} i'_{3})({\mathpzc e}{\scriptstyle\wedge} e''_{3})
{+}({\mathpzc i}{\scriptstyle\wedge} i'_{4})({\mathpzc i}{\scriptstyle\wedge} i''_{4})\\
\hspace*{20mm}=\phi(e'_{1})\phi(e''_{1}){+}\phi(e'_{2})\phi(i''_{2})
{+}\phi(i'_{2})\phi(e''_{3}){+}\phi(i'_{4})\phi(i''_{4})\\
\hspace*{20mm}=(\phi{\otimes}\phi)({\vartriangle}_q^{\!\XX}(e))$.

So $(\widehat\psi\,\widehat\otimes\,{\vartriangle}_q)\phi=(\phi{\otimes}\phi){\vartriangle}_q^{\!\XX}$.

The support coalgebra of $(T_*^\XX\!,\widehat\psi)$ is itself and the support coalgebra of $(H_*^\XX\!,{\vartriangle}_q)$ is
$(H_*^\SS\!,{\vartriangle}_q)$. So by Theorem~6.11,
the support coalgebra of $(T_*^\XX\widehat\otimes H_*^\XX\!,\widehat\psi\widehat\otimes{\vartriangle}_q)$
is $(T_*^\SS\widehat\otimes H_*^\SS\!,\widehat\psi\widehat\otimes{\vartriangle}_q)$.
It is obvious that $(T_*^\SS\widehat\otimes H_*^\SS\!,\widehat\psi\widehat\otimes{\vartriangle}_q)=(T_*^\SS\widehat\otimes H_*^\SS\!,\psi_q\widehat\otimes{\vartriangle}_q)$. So we have the last commutative diagram of the theorem by Theorem~6.11.
\hfill$\Box$\vspace{3mm}

{\bf Definition~8.4} Let $(\underline{D_*},\underline{C_*})=\{(D_{k\,*},C_{k\,*})\}_{k=1}^m$
be a sequence of pairs such that each $(D_{k\,*},C_{k\,*})$ is a split coalgebra pair by Definition~8.1.
Write $\Lambda_k,\theta_k,q_k,\cdots$ for the $\Lambda,\theta,q,\cdots$ in Definition~8.2 for $(D_*,C_*)=(D_{k\,*},C_{k\,*})$.
Then all the groups and complexes in Definition~3.7 and Definition~3.11 are (co)algebras when all the $q_k$ are chosen.
Let $\underline{q}=q_1{\otimes}{\cdots}{\otimes}q_m$.

The {\it character chain coalgebra} of $(\underline{D_*},\underline{C_*})$ with respect to $\underline{q}$
is the chain $(\XX_m{\times}\underline{\Lambda})$-coalgebra ($\underline{\Lambda}=\Lambda_1{\times}{\cdots}{\times}\Lambda_m$)
$$(C_*^{\XX_m}(\underline{D_*},\underline{C_*}),{\vartriangle}_{\underline{q}}^{\!\XX_m})=
(\otimes_{k=1}^m\,C_*^\XX\!(D_{k\,*},C_{k\,*}),\otimes_{k=1}^m\,{\vartriangle}_{q_k}^{\!\XX}).$$

The {\it indexed homology coalgebra} and {\it indexed cohomology algebra}
of $(\underline{D_*},\underline{C_*})$ with respect to $\underline{q}$
are the  $(\XX_m{\times}\underline{\Lambda})$-(co)algebras
$$(H_*^{\XX_m}(\underline{D_*},\underline{C_*}),{\vartriangle}_{\underline{q}})
=(\otimes_{k=1}^m\,H_*^{\XX}(D_{k\,*},C_{k\,*}),\otimes_{k=1}^m\,{\vartriangle}_{q_k}),$$
$$(H^{\,*}_{\!\XX_m}(\underline{D_*},\underline{C_*}),{\triangledown}\!_{\underline{q}})
=(\otimes_{k=1}^m\,H^{\,*}_{\!\XX}(D_{k\,*},C_{k\,*}),\otimes_{k=1}^m\,{\triangledown}\!_{q_k}).$$

The {\it normal homology coalgebra} and {\it normal cohomology algebra}
of $(\underline{D_*},\underline{C_*})$ (irrelevant to $\underline{q}$) are the following $(\XX_m{\times}\underline{\Lambda})$-(co)algebras
$$(H_*^{\XX_m}(\underline{D_*},\underline{C_*}),\vartriangle_{(\underline{D},\underline{C})})
=(\otimes_{k=1}^m\,H_*^{\XX}(D_{k\,*},C_{k\,*}),\otimes_{k=1}^m\,\vartriangle_{(D_k,C_k)}),$$
$$(H^{\,*}_{\!\XX_m}(\underline{D_*},\underline{C_*}),{\triangledown}_{(\underline{D},\underline{C})})
=(\otimes_{k=1}^m\,H^{\,*}_{\!\XX}(D_{k\,*},C_{k\,*}),\otimes_{k=1}^m\,\triangledown_{(D_k,C_k)}).$$

For an index set $\DD\subset\XX_m$, the analogue coalgebra $H_*^\DD(\underline{X},\underline{A})$ of $(\underline{X},\underline{A})$ on $\DD$
is the restriction coalgebra of $H_*^{\XX_m}(\underline{X},\underline{A})$ on $\DD$ (the dual case is similar).

The {\it total coproduct} $\underline{\psi_q}$ with respect to $\underline{q}$ is $\psi_{q_1}{\otimes}{\cdots}{\otimes}\psi_{q_m}$.
The {\it total chain coalgebra} $(T_*^{\XX_m},{\underline{\psi_q}})$ with respect to $\underline{q}$ is
$(T_*^{\XX}\!{\otimes}{\cdots}{\otimes}T_*^{\XX}\!,\,\psi_{q_1}{\otimes}{\cdots}{\otimes}\psi_{q_m})$.

Let $K$ be a simplicial complex on $[m]$.

The polyhedral product chain complex ${\cal Z}(K;\underline{D_*},\underline{C_*})$
in Definition~3.1 is a chain subcoalgebra of $(\otimes_kD_{k\,*},\otimes_k\psi_{D_k})$.
The product of its cohomology algebra is denoted by
$\cup_{{\cal Z}(K;\underline{D},\underline{C})}$.
The polyhedral product character chain complex ${\cal Z}^{\XX_m}(K;\underline{D_*},\underline{C_*})$ in Definition~3.11
is a chain subcoalgebra of the character chain coalgebra
$(C_*^{\XX_m}(\underline{D_*},\underline{C_*}),{\vartriangle}_{\underline{q}}^{\XX_m})$.

The total chain complex of $K$ with respect to $\underline{\psi_q}$ by Definition~7.4 is denoted by
$(T_*^{\XX_m}(K),\psi_{(K;\underline{q})})$.
The total cohomology algebra of $K$ with respect to $\underline{\psi_q}$ by Definition~7.4 is denoted by
$(H^{\,*}_{\!\XX_m}(K),\cup_{(K;\underline{q})})$.

For an index set $\DD\subset\XX_m$, we also have the total objects of $K$ on $\DD$ with respect to ${\underline{\psi_q}}$.
When $\DD=\RR_m$, we have the {\it right total objects} of $K$.
Denote them as in the following table.
\begin{center}
{\rm \begin{tabular}{|c|c|c|}
\hline
{\rule[-3mm]{0mm}{8mm}}
total object\,\,&\,\,right total object\,\,&\,\, object on $\DD$\\
\hline
{\rule[-3mm]{0mm}{8mm}}
\,\,$(T_*^{\XX_m}(K),\psi_{(K;\underline{q})})$\,\,&\,\,$(T_*^{\RR_m}(K),\psi_{(K;\underline{q})})$\,\,&\,\,$(T_*^\DD(K),\psi_{(K;\underline{q})})$\,\,\\
\hline
{\rule[-3mm]{0mm}{8mm}}
\,\,$(T^{\,*}_{\!\XX_m}(K),\pi_{(K;\underline{q})})$\,\,&\,\,$(T^{\,*}_{\!\RR_m}(K),\pi_{(K;\underline{q})})$
\,\,&\,\,$(T^{\,*}_{\!\DD}(K),\pi_{(K;\underline{q})})$\,\,\\
\hline
{\rule[-3mm]{0mm}{8mm}}
\,\,$(H_*^{\XX_m}(K),\vartriangle_{(K;\underline{q})})$\,\,&\,\,$(H_*^{\RR_m}(K),\vartriangle_{(K;\underline{q})})$
\,\,&\,\,$(H_*^{\DD}(K),\vartriangle_{(K;\underline{q})})$\,\,\\
\hline
{\rule[-3mm]{0mm}{8mm}}
\,\,$(H^{\,*}_{\!\XX_m}(K),\triangledown_{(K;\underline{q})})$\,\,&\,\,$(H^{\,*}_{\!\RR_m}(K),\triangledown_{(K;\underline{q})})$
\,\,&\,\,$(H^{\,*}_{\!\DD}(K),\triangledown_{(K;\underline{q})})$\,\,\\
\hline
{\rule[-3mm]{0mm}{8mm}}
\,\,$(H^{\,*}_{\!\XX_m}(K),\cup_{(K;\underline{q})})$\,\,&\,\,$(H^{\,*}_{\!\RR_m}(K),\cup_{(K;\underline{q})})$\,\,&\,\,$(H^{\,*}_{\!\DD}(K),\cup_{(K;\underline{q})})$\,\,\\
\hline
\end{tabular}}
\end{center}
\vspace{5mm}

{\bf Theorem~8.5} {\it If ${\cal Z}_*(K;\underline{D_*},\underline{C_*})$ satisfies that each $(D_{k\,*},C_{k\,*})$ is a split coalgebra pair,
then the chain homomorphisms in Theorem~3.12 are indexed chain coalgebra homomorphisms as follows
(${\cal Z}^{\XX_m}={\cal Z}^{\XX_m}(K;\underline{D_*},\underline{C_*})$).
$$q_{(K;\underline{D},\underline{C})}\colon({\cal Z}(K;\underline{D_*},\underline{C_*}),\psi_{(K;\underline{D_*},\underline{C_*})})
\stackrel{\simeq}{\longrightarrow}
({\cal Z}^{\XX_m},{\vartriangle}_{\underline{q}}^{\!\XX_m})$$
is a chain $\underline{\Lambda}$-coalgebra weak homomorphism (index $\XX_m$ neglected) and
$$\phi_{(K;\underline{D},\underline{C})}\colon({\cal Z}^{\XX_m},{\vartriangle}_{\underline{q}}^{\!\XX_m})
\stackrel{\cong}{\longrightarrow}
(T_*^{\XX_m}(K)\,\widehat\otimes\,H_*^{\XX_m}(\underline{D_*},\underline{C_*}),\cup_{(K;\underline{q})}\widehat\otimes\cup_{\underline{q}}).$$
is a chain $(\XX_m{\times}\underline{\Lambda})$-coalgebra isomorphism.

The cohomology group isomorphisms in Theorem~3.12 are $\underline{\Lambda}$-algebra (all other index neglected) isomorphisms as follows.
We always have
$$(H^*({\cal Z}(K;\underline{D_*},\underline{C_*})),\cup_{{\cal Z}(K;\underline{D},\underline{C})})
\cong(H^{\,*}_{\!\XX_m}(K)\,\widehat\otimes\, H^{\,*}_{\!\XX_m}(\underline{D_*},\underline{C_*}),
\widehat\cup_K\,\widehat\otimes\,{\triangledown}\!_{\underline{q}}),
$$
where $(H_*^{\XX_m}(K),\widehat\cup_K)$ is the universal cohomology algebra in Definition~7.4.
Suppose $(T_*^{\XX_m},\underline{\psi})=(T_*^\XX{\otimes}{\cdots}{\otimes}T_*^\XX,\psi_1{\otimes}{\cdots}{\otimes}\psi_m)$
is a total chain coalgebra such that each $\psi_{q_k}\prec\psi_k$ (by Definition~6.10), then we have
$$\begin{array}{l}
\quad (H^{\,*}_{\!\XX_m}(K)\,\widehat\otimes\, H^{\,*}_{\!\XX_m}(\underline{D_*},\underline{C_*}),
\widehat\cup_K\,\widehat\otimes\,{\triangledown}\!_{\underline{q}})\vspace{2mm}\\
\cong (H^{\,*}_{\!\DD}(K)\,\widehat\otimes\, H^{\,*}_{\!\DD}(\underline{D_*},\underline{C_*}),
\cup_K\,\widehat\otimes\,{\triangledown}\!_{\underline{q}})\vspace{2mm}\\
\cong (H^{\,*}_{\!\underline{\SS}}(K)\,\widehat\otimes\, H^{\,*}_{\!\underline{\SS}}(\underline{D_*},\underline{C_*}),
\cup_{(K;\underline{q})}\,\widehat\otimes\,{\triangledown}\!_{\underline{q}}),
\end{array}$$
where $(H^{\,*}_{\!\DD}(K),\cup_K)$ is the total cohomology algebra of $K$ on $\DD$ with respect to $\underline{\psi}$
in Definition~7.4 and $(H^{\,*}_{\!\underline{\SS}}(K),\cup_{(K;\underline{q})})$ is as in Definition~8.4.
\vspace{2mm}

\it Proof}\, By Theorem~8.3,
the $q_\tau=p_1{\otimes}{\cdots}{\otimes}p_m$ in the proof of Theorem~3.12 is a chain $\underline{\Lambda}$-coalgebra weak homomorphism.
So $q_{(K;\underline{D_*},\underline{C_*})}=+_{\tau\in K}\,q_\tau$ is a chain $\underline{\Lambda}$-coalgebra weak homomorphism.
Similarly, $\phi_\tau=\lambda_1{\otimes}{\cdots}{\otimes}\lambda_m$ is a chain $(\XX_m{\times}\underline{\Lambda})$-coalgebra isomorphism.
So $\phi_{(K;\underline{D_*},\underline{C_*})}=+_{\tau\in K}\,\phi_\tau$ is a chain $(\XX_m{\times}\underline{\Lambda})$-coalgebra isomorphism.

The last equalities of the theorem holds by Theorem~6.11.
\hfill$\Box$\vspace{3mm}

{\bf Example~8.6} We check Theorem~8.5 for $K=\Delta\!^S$ with $S\subset[m]$.
By K\"{u}nneth Theorem,
$$(H_*({\cal Z}(K;\underline{D_*},\underline{C_*})),{\vartriangle}_{{\cal Z}(K;\underline{D_*},\underline{C_*})})
\cong(E_1{\otimes}{\cdots}{\otimes}E_m,{\vartriangle}_1{\otimes}{\cdots}{\otimes}{\vartriangle}_m)=(\w E,\w{\vartriangle}),$$
where $(E_k,{\vartriangle}_k)=(H_*(D_{k\,*}),{\vartriangle}_{D_k})$ if $k\in S$ and $(E_k,{\vartriangle}_k)=(H_*(C_{k\,*}),{\vartriangle}_{C_k})$ if $k\notin S$.
Denote by $H_*^{\XX_m}(\underline{D_*},\underline{C_*})=\oplus_{(\sigma,\omega)\in\XX_m}\,H_*^{\sigma,\omega}$.
Then $\w E$ is naturally the subgroup $\oplus_{\sigma\subset S,\,\omega{\cap}S=\emptyset}\,H_*^{\sigma,\omega}$
of $H_*^{\XX_m}(\underline{D_*},\underline{C_*})$.
From the computation of Example~3.13 we have a group isomorphism $\phi$ from $\w E$ to
$H_*^{\XX_m}(K)\widehat\otimes H_*^{\XX_m}(D_*,C_*)$ defined by
$\phi(x_{\sigma,\omega})=1_{\sigma,\omega}\widehat\otimes(x_{\sigma,\omega})$ for all $x_{\sigma,\omega}\in H_*^{\sigma,\omega}$,
where $1_{\sigma,\omega}$ is the generator of $H_0^{\sigma,\omega}(K)$ as in Example~7.8.
Now we check that $\phi$ is an algebra isomorphism.
By definition, for $x_{\sigma,\omega}=x_1{\otimes}{\cdots}{\otimes}x_m\in H_*^{\sigma,\omega}\subset\w E$,
$${\vartriangle}_{(\underline{D},\underline{C})}(x_1{\otimes}{\cdots}{\otimes}x_m)
=\w{\vartriangle}(x_1{\otimes}{\cdots}{\otimes}x_m)
+\Sigma(y'_1{\otimes}{\cdots}{\otimes}y'_m){\otimes}(y''_1{\otimes}{\cdots}{\otimes}y''_m)
$$
where for some $k\in S$, at least one of $y'_k$ and $y''_k$ is in ${\rm ker}\,\theta_k$.
This implies that at least one of $y_{\mu'\!,\nu'}=y'_1{\otimes}{\cdots}{\otimes}y'_m$
and $y_{\mu''\!,\nu''}=y''_1{\otimes}{\cdots}{\otimes}y''_m$ is not in $\w E$.
Suppose $\w{\vartriangle}(x_{\sigma,\omega})=\Sigma x_{\sigma'\!,\omega'}{\otimes}x_{\sigma''\!,\omega''}$
with $x_{-,-}\in H_*^{-,-}$.
Then $\sigma',\sigma''\subset S$, $\omega'{\cap}S=\omega''{\cap}S=\emptyset$.
So by the computation of Example~7.8,
$$\begin{array}{l}
\quad(\w{\vartriangle}_K\widehat\otimes{\vartriangle}_{(\underline{D},\underline{C})})(\phi(x_{\sigma,\omega}))\vspace{1mm}\\
=(\w{\vartriangle}_K(1_{\sigma,\omega}))\,\widehat\otimes\,({\vartriangle}_{(\underline{D},\underline{C})}(x_{\sigma,\omega}))\vspace{1mm}\\
=(\Sigma\,1_{\sigma'\!,\omega'}{\otimes}1_{\sigma''\!,\omega''})\,\widehat\otimes\,
(\Sigma\,x_{\sigma'\!,\omega'}{\otimes}x_{\sigma''\!,\omega''}{+}\Sigma\,y_{\mu'\!,\nu'}{\otimes}y_{\mu'',\nu''})\vspace{1mm}\\
=\Sigma\,(1_{\sigma'\!,\omega'}{\otimes}x_{\sigma'\!,\omega'})\,\widehat\otimes\,
(1_{\sigma''\!,\omega''}{\otimes}x_{\sigma''\!,\omega''})\vspace{1mm}\\
=(\phi{\otimes}\phi)(\w{\vartriangle}(x_{\sigma,\omega})).
\end{array}$$
So $(\w{\vartriangle}_K\widehat\otimes{\vartriangle}_{(\underline{D},\underline{C})})\phi
=(\phi{\otimes}\phi)(\w{\vartriangle})$. But $\w{\vartriangle}_K=\widehat{\vartriangle}_K$ when $K=\Delta\!^S$.
This implies $\w{\vartriangle}_K\widehat\otimes{\vartriangle}_{(\underline{D},\underline{C})}
=\widehat{\vartriangle}_K\widehat\otimes{\vartriangle}_q$. So
$(\widehat{\vartriangle}_K\widehat\otimes{\vartriangle}_{q})\phi
=(\phi{\otimes}\phi)\w{\vartriangle}$.
\vspace{3mm}

The local product of the product $\cup_{(K;\underline{q})}$ in Theorem~8.5 may be very complicated,
so we have to use $\cup_K$ with simpler local product instead of $\cup_{(K;\underline{q})}$ in actual computation.
The following theorem lists three cases that are often used.
\vspace{2mm}

{\bf Theorem~8.7} {\it Suppose ${\cal Z}_*(K;\underline{D_*},\underline{C_*})$ satisfies that each $(D_{k\,*},C_{k\,*})$ is a split coalgebra pair.

If each pair $(D_{k\,*},C_{k\,*})$ is normal, then we have
$$(H^*({\cal Z}(K;\underline{D_*},\underline{C_*})),\cup_{(K;\underline{D},\underline{C})})
\cong(H^{\,*}_{\!\XX_m}(K)\,\widehat\otimes\, H^{\,*}_{\!\XX_m}(\underline{D_*},\underline{C_*}),
\w\cup_K\,\widehat\otimes\,{\triangledown}\!_{(\underline{D},\underline{C})}),
$$
where $(H_*^{\XX_m}(K),\w\cup_K)$ is the normal cohomology algebra of $K$.

If each $\theta_k\colon H_*(C_{k\,*})\to H_*(D_{k\,*})$ is an epimorphism, then we have
$$(H^*({\cal Z}(K;\underline{D_*},\underline{C_*})),\cup_{(K;\underline{D},\underline{C})})
\cong(H^{\,*}_{\!\RR_m}(K)\,\widehat\otimes\, H^{\,*}_{\!\RR_m}(\underline{D_*},\underline{C_*}),
\widehat\cup_K\,\widehat\otimes\,{\triangledown}\!_{\underline{q}}),
$$
where $(H_*^{\RR_m}(K),\widehat\cup_K)$ is the right universal cohomology algebra of $K$.

If each pair $(D_{k\,*},C_{k\,*})$ is normal and each $\theta_k$ is an epimorphism, then
$(H^*(C_{k\,*}),{\triangledown}_{C_k})=(H^{\XX}(D_{k\,*},C_{k\,*}),{\triangledown}_{\!(D,C)})$ and we have
$$(H^*({\cal Z}(K;\underline{D_*},\underline{C_*})),\cup_{(K;\underline{D},\underline{C})})
\cong(H^{\,*}_{\!\RR_m}(K)\,\widehat\otimes\,\big(\otimes_k H^*(C_{k\,*})\big),
\w\cup_K\,\widehat\otimes\,(\otimes_k{\triangledown}_{\!C_k})),
$$
where the $\RR_m{\times}\underline{\Lambda}$ index of $\otimes_k H^*(C_{k\,*})=H^*_{\RR_m}(\underline{D_*},\underline{C_*})$ is defined as follows.
For $a_k\in H^*(C_{k\,*})$, $a_1{\otimes}{\cdots}{\otimes}a_m\in H^*_{\emptyset,\omega}(\underline{D_*},\underline{C_*})$
with $\omega=\{k\,|\,a_k\in{\rm ker}\,\theta_k\}$.
\vspace{2mm}

\it Proof}\, If each pair $(D_{k\,*},C_{k\,*})$ is normal,
then we may take all $\psi_k$ in Theorem~8.5 to be the normal coproduct $\w\psi$ in Definition~7.1 and take $\DD$ to be $\XX_m$.

If each $\theta_k\colon H_*(C_{k\,*})\to H_*(D_{k\,*})$ is an epimorphism,
then we may take all $\psi_k$ in Theorem~8.5 to be the right universal coproduct $\widehat\psi'$ in Definition~7.1 and take $\DD$ to be $\RR_m$.
\hfill$\Box$\vspace{3mm}

{\bf Theorem~8.8} {\it Let $(D_{i\,*},C_{i\,*})$ be split coalgebra pairs with indexed homology coalgebra
$(H_*^\XX(D_{i\,*},C_{i\,*}),{\vartriangle}_{q_i})$ for $i=1,2$.
Then the pair $(D_*,C_*)=(D_{1\,*}{\otimes}D_{2\,*},C_{1\,*}{\otimes}C_{2\,*})$ is also a split coalgebra pair
with indexed homology coalgebra $(H_*^\XX(D_*,C_*),{\vartriangle}_{q})$ defined as follows.
There is a group isomorphism (index $\XX$ neglected)
$H_*^\XX(D_*,C_*)\cong H_*^\XX(D_{1\,*},C_{1\,*}){\otimes}H_*^\XX(D_{2\,*},C_{2\,*})$ with

(1) ${\mathpzc i}_{\,*}={\mathpzc i}_{\,1*}{\scriptstyle\otimes}{\mathpzc i}_{\,2*}$.

(2) ${\mathpzc n}_{\,*}={\mathpzc n}_{\,1*}{\scriptstyle\otimes}{\mathpzc n}_{\,2*}\oplus{\mathpzc n}_{\,1*}{\scriptstyle\otimes}{\mathpzc i}_{\,2*}\oplus
{\mathpzc n}_{\,1*}{\scriptstyle\otimes}{\mathpzc e}_{\,2*}\oplus{\mathpzc i}_{\,1*}{\scriptstyle\otimes}{\mathpzc n}_{\,2*}\oplus{\mathpzc e}_{\,1*}{\scriptstyle\otimes}{\mathpzc n}_{\,2*}$.

(3) ${\mathpzc e}_{\,*}={\mathpzc e}_{\,1*}{\scriptstyle\otimes}{\mathpzc e}_{\,2*}\oplus
{\mathpzc i}_{\,1*}{\scriptstyle\otimes}{\mathpzc e}_{\,2*}\oplus{\mathpzc e}_{\,1*}{\scriptstyle\otimes}{\mathpzc i}_{\,2*}$.

For $x_1\in H_*^\XX(D_{1\,*},C_{1\,*})$ and $x_2\in H_*^\XX(D_{2\,*},C_{2\,*})$,
$${\vartriangle}_q(x_1{\otimes}x_2)=
\left\{\begin{array}{ll}
({\vartriangle}_{q_1}{\otimes}{\vartriangle}_{(D_2,C_2)})(x_1{\otimes}x_2)&{\rm if}\,\,x_1\in{\mathpzc n}_{\,1*}\vspace{1mm}\\
({\vartriangle}_{(D_1,C_1)}{\otimes}{\vartriangle}_{q_2})(x_1{\otimes}x_2)
&{\rm if}\,\,x_1\notin{\mathpzc n}_{\,1*}\,\,{\rm but}\,\,x_2\in{\mathpzc n}_{\,2*}\vspace{1mm}\\
({\vartriangle}_{(D_1,C_1)}{\otimes}{\vartriangle}_{(D_2,C_2)})(x_1{\otimes}x_2)&{\rm otherwise}
\end{array}\right.$$

So if both $(D_{1\,*},C_{1\,*})$ and $(D_{2\,*},C_{2\,*})$ are normal, then $(D_{1\,*}{\otimes}D_{2\,*},C_{1\,*}{\otimes}C_{2\,*})$ is normal
and we have (co)algebra isomorphism
$$\begin{array}{l}
\quad(H_*^\XX(D_{1\,*}{\otimes}D_{2\,*},C_{1\,*}{\otimes}C_{2\,*}),{\vartriangle}_{(D_1{\otimes}D_2,C_1{\otimes}C_2)})\vspace{2mm}\\
\cong (H_*^\XX(D_{1\,*},C_{1\,*}){\otimes}H_*^\XX(D_{2\,*},C_{2\,*}),{\vartriangle}_{(D_1,C_1)}{\otimes}{\vartriangle}_{(D_2,C_2)}),
\end{array}$$
$$\begin{array}{l}
\quad(H^{\,*}_{\!\XX}(D_{1\,*}{\otimes}D_{2\,*},C_{1\,*}{\otimes}C_{2\,*}),{\triangledown}_{(D_1{\otimes}D_2,C_1{\otimes}C_2)})\vspace{2mm}\\
\cong (H^{\,*}_{\!\XX}(D_{1\,*},C_{1\,*}){\otimes}H^{\,*}_{\!\XX}(D_{2\,*},C_{2\,*}),{\triangledown}_{(D_1,C_1)}{\otimes}{\triangledown}_{(D_2,C_2)}).
\end{array}\vspace{2mm}$$

Proof} By definition.
\hfill$\Box$\vspace{3mm}

{\bf Theorem~8.9} {\it Let $(D_{i\,*}^\TT,C_{i\,*}^\TT)$ be split coalgebra pairs with indexed homology coalgebra
$(H_*^\XX(D_{i\,*}^\TT,C_{i\,*}^\TT),{\vartriangle}_{q_i})$ for $i=1,2$.
Then the pair $(D_*,C_*)=(D_{1\,*}^\TT{\widehat\otimes}D_{2\,*}^\TT,C_{1\,*}^\TT{\widehat\otimes}C_{2\,*}^\TT)$
(the diagonal tensor product is with respect to $\TT$) is also a split coalgebra pair
with indexed homology coalgebra $(H_*^\XX(D_*,C_*),{\vartriangle}_{q})$ defined as follows.
We have group isomorphism ($\XX$ neglected)
$H_*^\XX(D_*,C_*)\cong H_*^\XX(D_{1\,*}^\TT,C_{1\,*}^\TT){\widehat\otimes}H_*^\XX(D_{2\,*}^\TT,C_{2\,*}^\TT)$ with

(1) ${\mathpzc i}_{\,*}={\mathpzc i}_{\,1*}{\scriptstyle\widehat\otimes}{\mathpzc i}_{\,2*}$.

(2) ${\mathpzc n}_{\,*}={\mathpzc n}_{\,1*}{\scriptstyle\widehat\otimes}{\mathpzc n}_{\,2*}\oplus{\mathpzc n}_{\,1*}{\scriptstyle\widehat\otimes}{\mathpzc i}_{\,2*}\oplus
{\mathpzc n}_{\,1*}{\scriptstyle\widehat\otimes}{\mathpzc e}_{\,2*}\oplus{\mathpzc i}_{\,1*}{\scriptstyle\widehat\otimes}{\mathpzc n}_{\,2*}\oplus{\mathpzc e}_{\,1*}{\scriptstyle\widehat\otimes}{\mathpzc n}_{\,2*}$.

(3) ${\mathpzc e}_{\,*}={\mathpzc e}_{\,1*}{\scriptstyle\widehat\otimes}{\mathpzc e}_{\,2*}\oplus
{\mathpzc i}_{\,1*}{\scriptstyle\widehat\otimes}{\mathpzc e}_{\,2*}\oplus
{\mathpzc e}_{\,1*}{\scriptstyle\widehat\otimes}{\mathpzc i}_{\,2*}$.

For $x_1\in H_*^\XX(D_{1\,*}^\TT,C_{1\,*}^\TT)$ and $x_2\in H_*^\XX(D_{2\,*}^\TT,C_{2\,*}^\TT)$,
$${\vartriangle}_q(x_1{\widehat\otimes}x_2)=
\left\{\begin{array}{ll}
({\vartriangle}_{q_1}{\widehat\otimes}{\vartriangle}_{(D_2,C_2)})(x_1{\widehat\otimes}x_2)&{\rm if}\,\,x_1\in{\mathpzc n}_{\,1*}\vspace{1mm}\\
({\vartriangle}_{(D_1,C_1)}{\widehat\otimes}{\vartriangle}_{q_2})(x_1{\widehat\otimes}x_2)
&{\rm if}\,\,x_1\notin{\mathpzc n}_{\,1*}\,\,{\rm but}\,\,x_2\in{\mathpzc n}_{\,2*}\vspace{1mm}\\
({\vartriangle}_{(D_1,C_1)}{\widehat\otimes}{\vartriangle}_{(D_2,C_2)})(x_1{\widehat\otimes}x_2)&{\rm otherwise}
\end{array}\right.$$

So if both $(D_{1\,*}^\TT,C_{1\,*}^\TT)$ and $(D_{2\,*}^\TT,C_{2\,*}^\TT)$ are normal, then $(D_{1\,*}^\TT{\widehat\otimes}D_{2\,*}^\TT,C_{1\,*}^\TT{\widehat\otimes}C_{2\,*}^\TT)$ is normal
and we have (co)algebra isomorphism
$$\begin{array}{l}
\quad(H_*^\XX(D_{1\,*}^\TT{\widehat\otimes}D_{2\,*}^\TT,C_{1\,*}^\TT{\widehat\otimes}C_{2\,*}^\TT),
{\vartriangle}_{(D_1{\widehat\otimes}D_2,C_1{\widehat\otimes}C_2)})\vspace{2mm}\\
\cong (H_*^\XX(D_{1\,*}^\TT,C_{1\,*}^\TT){\widehat\otimes}H_*^\XX(D_{2\,*}^\TT,C_{2\,*}^\TT),
{\vartriangle}_{(D_1,C_1)}{\widehat\otimes}{\vartriangle}_{(D_2,C_2)}),
\end{array}$$
$$\begin{array}{l}
\quad(H^{\,*}_{\!\XX}(D_{1\,*}^\TT{\widehat\otimes}D_{2\,*}^\TT,C_{1\,*}^\TT{\widehat\otimes}C_{2\,*}^\TT),
{\triangledown}_{(D_1{\widehat\otimes}D_2,C_1{\widehat\otimes}C_2)})\vspace{2mm}\\
\cong (H^{\,*}_{\!\XX}(D_{1\,*}^\TT,C_{1\,*}^\TT){\widehat\otimes}H^{\,*}_{\!\XX}(D_{2\,*}^\TT,C_{2\,*}^\TT),
{\triangledown}_{(D_1,C_1)}{\widehat\otimes}{\triangledown}_{(D_2,C_2)}),
\end{array}\vspace{2mm}$$
where the diagonal tensor product is with respect to $\TT$ with $\XX$ neglected.

Proof} By definition.
\hfill$\Box$\vspace{3mm}

\section{Cohomology Algebra of Polyhedral Product Objects}\vspace{3mm}

\hspace*{5.5mm}{\bf Definition 9.1} For a topological space $X$,
let $\psi_X\colon S_*(X)\to S_*(X){\otimes}S_*(X)$ be the coproduct (unique up to homotopy) of the singular chain complex induced by the diagonal map  of $X$.
Then $(S_*(X),\psi_X,d)$ is a chain $\Lambda$-coalgebra such that $\Lambda$ is trivial.

For a homology split topological pair $(X,A)$ by Definition~4.3, the pair $\big((S_*(X),\psi_X,d),(S_*(A),\psi_A,d)\big)$
is a split coalgebra pair by Definition~8.1.
So the (co)homology groups (not reduced or suspension reduced groups) in Definition~4.4 are (co)algebras as defined in Definition~8.2
when the homotopy equivalence $q$ is chosen.
Replace the $(D_*,C_*)$ in Definition~8.2 by $(X,A)$ and we have all definitions for $(X,A)$ as follows.

Let $\theta\colon(H_*(A),{\vartriangle}_{A})\to(H_*(X),{\vartriangle}_X)$ be the homology coalgebra homomorphism induced by inclusion.
Denote by
$${\mathpzc i}_{\,*}={\rm coim}\,\theta,\quad{\mathpzc n}_{\,*}={\rm ker}\,\theta,\quad
\overline{\mathpzc n}_{\,*}=\Sigma{\rm ker}\,\theta,\quad{\mathpzc e}_{\,*}={\rm coker}\,\theta.$$
Then we have group equalities $H_*(A)={\mathpzc i}_{\,*}\oplus{\mathpzc n}_{\,*}$,
$H_*^\XX(X,A)={\mathpzc i}_{\,*}\oplus{\mathpzc n}_{\,*}\oplus{\mathpzc e}_{\,*}$ and
$C_*^\XX(X,A)={\mathpzc i}_{\,*}\oplus{\mathpzc n}_{\,*}\oplus\overline{\mathpzc n}_{\,*}\oplus{\mathpzc e}_{\,*}$.

The {\it normal homology coalgebra} $(H_*^\XX(X,A),\vartriangle_{(X,A)})$ of $(X,A)$ (irrelevant to $q$) is given by
$${\vartriangle}_{(X,A)}(x)=\left\{\begin{array}{cl}
{\vartriangle}_A(x)&{\rm if}\,\,x\in {\mathpzc i}_{\,*},\vspace{1mm}\\
{\vartriangle}_A(x)&{\rm if}\,\,x\in {\mathpzc n}_{\,*},\vspace{1mm}\\
{\vartriangle}_X(x)&{\rm if}\,\,x\in {\mathpzc e}_{\,*}.
\end{array}\right.$$

The {\it normal cohomology algebra} $(H^{\,*}_{\!\XX}(X,A),{\triangledown}_{(X,A)})$ of $(X,A)$ is
the dual algebra of $(H_*^\XX(X,A),\vartriangle_{(X,A)})$.

The {\it character chain coalgebra} $(C_*^\XX(X,A),{\vartriangle}_q^{\!\XX},d)$ of $(X,A)$ with respect to $q$
is the chain $\XX$-coalgebra defined as follows.
The restriction of ${\vartriangle}_q^{\!\XX}$ on $H_*^\XX(X,A)$ is ${\vartriangle}_{(X,A)}$.
For $x\in{\mathpzc n}_{\,*}$ with ${\vartriangle}_A(x)=\Sigma x'_i{\otimes}x''_i+\Sigma y'_j{\otimes}y''_j$,
where each $x'_i\in{\mathpzc n}_{\,*}$ and each $y'_j\notin{\mathpzc n}_{\,*}$ but $y''_j\in{\mathpzc n}_{\,*}$, define
${\vartriangle}_q^{\!\XX}(\overline x)=\Sigma\overline x'_i{\otimes}x''_i+\Sigma(-1)^{|y'_j|}y'_j{\otimes}\overline y''_j
+\xi(q{\otimes}q)\psi_X(\overline x)$, where $\xi$ is the projection from $C_*^\XX(X,A)\otimes C_*^\XX(X,A)$ to its homology.

The {\it indexed homology coalgebra} $(H_*^\XX(X,A),{\vartriangle}_q)$ of $(X,A)$ with respect to $q$ is
the $\XX$-coalgebra given by
$${\vartriangle}_q(x)=\left\{\begin{array}{ll}
{\vartriangle}_A(x)&{\rm if}\,\,x\in{\mathpzc i}_{\,*},\vspace{1mm}\\
{\vartriangle}_A(x){+}\xi(q{\otimes}q)\psi_X(\overline x)&{\rm if}\,\,x\in {\mathpzc n}_{\,*},\vspace{1mm}\\
{\vartriangle}_X(x)&{\rm if}\,\,x\in {\mathpzc e}_{\,*}.
\end{array}\right.$$

The {\it indexed cohomology algebra} $(H^{\,*}_{\!\XX}(X,A),{\triangledown}_{\!q})$ of $(X,A)$ with respect to $q$ is
the dual algebra of $(H_*^\XX(X,A),{\vartriangle}_q)$.

$(X,A)$ is called {\it normal} if there is a homotopy equivalence $q$ such that ${\vartriangle}_q={\vartriangle}_{(X,A)}$
and ${\triangledown}_{\!q}={\triangledown}_{\!(X,A)}$.

The {\it atom coproduct} $\psi_q$ and the {\it atom chain coalgebra $(T_*^\XX\!,\psi_q)$} with respect to $q$ are defined as follows.
For ${\mathpzc s}\in T_*^\XX$, suppose ${\vartriangle}_q^{\!\XX}({\mathpzc s}_*)\subset\oplus_i\,({\mathpzc s}'_i)_*{\otimes}({\mathpzc s}''_i)_*$
but ${\vartriangle}_q^{\!\XX}({\mathpzc s}_*){\cap}(({\mathpzc s}'_i)_*{\otimes}({\mathpzc s}''_i)_*)\neq 0$ for each $i$,
then define $\psi_q({\mathpzc s})=\Sigma_i\,{\mathpzc s}'_i{\otimes}{\mathpzc s}''_i$.

The {\it index set} of $(X,A)$ is $\SS=\{{\mathpzc s}\in T_*^\XX\,|\,{\mathpzc s}_{\,*}\neq 0\}$.
By definition, $\psi_q({\mathpzc s})=0$ if ${\mathpzc s}\notin\SS$.

We have all analogue definitions for simplicial complex pair $(X,A)$, where for a simplicial complex $K$,
the coproduct $\psi_K\colon C_*(K)\to C_*(K){\otimes}C_*(K)$ of simplicial chain complex is defined as follows.
For an ordered simplex $\{i_1,{\cdots},i_n\}\in K\subset C_*(K)$,
$$\psi_K(\{i_1,{\cdots},i_n\})=\Sigma_{k=1}^{n}\{i_1,{\cdots},i_{k-1},i_k\}\otimes\{i_{k},i_{k+1}{\cdots},i_n\}.\vspace{3mm}$$

{\bf Definition 9.2} Let $(\underline{X},\underline{A})=\{(X_k,A_k)\}_{k=1}^m$
be such that each $(X_k,A_k)$ is homology split.
Then all the (co)homology groups in Definition~4.5 are (co)algebras as defined in Definition~8.4
when all the homotopy equivalences $q_k$ for $(X_k,A_k)$ are chosen.
Precisely, replace the $(D_{k\,*},C_{k\,*})$ in Definition~8.4 by $(X_k,A_k)$
and we have all definitions for $(\underline{X},\underline{A})$.
We only list the definitions needed in Theorem~9.3.

The {\it indexed homology coalgebra} and {\it indexed cohomology algebra}
of $(\underline{X},\underline{A})$ with respect to $\underline{q}=q_1{\otimes}{\cdots}{\otimes}q_m$ are the  $\XX_m$-(co)algebras
$$(H_*^{\XX_m}(\underline{X},\underline{A}),{\vartriangle}_{\underline{q}})
=(\otimes_{k=1}^m\,H_*^{\XX}(X_k,A_k),\otimes_{k=1}^m\,{\vartriangle}_{q_k}),$$
$$(H^{\,*}_{\!\XX_m}(\underline{X},\underline{A}),{\triangledown}\!_{\underline{q}})
=(\otimes_{k=1}^m\,H^{\,*}_{\!\XX}(X_k,A_k),\otimes_{k=1}^m\,{\triangledown}\!_{q_k}).$$

The {\it normal homology coalgebra} and {\it normal cohomology algebra}
of $(\underline{X},\underline{A})$ (irrelevant to $\underline{q}$) are the following $\XX_m$-(co)algebras
$$(H_*^{\XX_m}(\underline{X},\underline{A}),\vartriangle_{(\underline{X},\underline{A})})
=(\otimes_{k=1}^m\,H_*^{\XX}(X_k,A_k),\otimes_{k=1}^m\,\vartriangle_{(X_k,A_k)}),$$
$$(H^{\,*}_{\!\XX_m}(\underline{X},\underline{A}),{\triangledown}\!_{(\underline{X},\underline{A})})
=(\otimes_{k=1}^m\,H^{\,*}_{\!\XX}(X_k,A_k),\otimes_{k=1}^m\,\triangledown\!_{(X_k,A_k)}).$$

For an index set $\DD\subset\XX_m$, the analogue coalgebra $H_*^\DD(\underline{X},\underline{A})$ of $(\underline{X},\underline{A})$ on $\DD$
is the restriction coalgebra of $H_*^{\XX_m}(\underline{X},\underline{A})$ on $\DD$ (the dual case is similar).

The {\it total coproduct} $\underline{\psi_q}$ with respect to $\underline{q}$ is $\psi_{q_1}{\otimes}{\cdots}{\otimes}\psi_{q_m}$.
The {\it total chain coalgebra} $(T_*^{\XX_m},{\underline{\psi_q}})$ with respect to $\underline{q}$ is
$(T_*^{\XX}\!{\otimes}{\cdots}{\otimes}T_*^{\XX}\!,\,\psi_{q_1}{\otimes}{\cdots}{\otimes}\psi_{q_m})$.

For a simplicial complex $K$ on $[m]$, the total chain complex of $K$ with respect to $\underline{\psi_q}$ by Definition~7.4 is denoted by
$(T_*^{\XX_m}(K),\psi_{(K;\underline{q})})$.
The total cohomology algebra of $K$ with respect to $\underline{\psi_q}$ by Definition~7.4 is denoted by
$(H^{\,*}_{\!\XX_m}(K),\cup_{(K;\underline{q})})$.
For an index set $\DD\subset\XX_m$, the total cohomology algebra of $K$ on $\DD$ with respect to $\underline{\psi_q}$ is denoted by
$(H^{\,*}_{\!\DD}(K),\cup_{(K;\underline{q})})$.
\vspace{3mm}

{\bf Theorem 9.3} {\it The cohomology group isomorphisms in Theorem~4.6 for topological space and simplicial complex pairs
(not reduced or suspension reduced) are algebra isomorphisms by Theorem~8.4 and Theorem~8.7 as follows.
Suppose $(T_*^{\XX_m},\underline{\psi})=(T_*^\XX{\otimes}{\cdots}{\otimes}T_*^\XX,\psi_1{\otimes}{\cdots}{\otimes}\psi_m)$
is a total chain coalgebra such that each atom coproduct $\psi_{q_k}\prec\psi_k$.
Let $\underline{\SS}=\SS_1{\times}{\cdots}{\times}\SS_m$, where each $\SS_k$ is the support index set of $(X_k,A_k)$.
Then for any index set $\DD$ such that $\underline{\SS}\subset\DD\subset\XX_m$,
we have algebra isomorphisms (index $\XX_m$, $\DD$, $\underline{\SS}$ neglected),
$$\begin{array}{l}
\quad (H^{\,*}({\cal Z}(K;\underline{X},\underline{A})),\cup_{{\cal Z}(K;\underline{X},\underline{A})})\vspace{2mm}\\
\cong (H^{\,*}_{\!\XX_m}(K)\,\widehat\otimes\, H^{\,*}_{\!\XX_m}(\underline{X},\underline{A}),
\widehat\cup_K\,\widehat\otimes\,{\triangledown}\!_{\underline{q}})\vspace{2mm}\\
\cong (H^{\,*}_{\!\DD}(K)\,\widehat\otimes\, H^{\,*}_{\!\DD}(\underline{X},\underline{A}),
\cup_K\,\widehat\otimes\,{\triangledown}\!_{\underline{q}})\vspace{2mm}\\
\cong (H^{\,*}_{\!\underline{\SS}}(K)\,\widehat\otimes\, H^{\,*}_{\!\underline{\SS}}(\underline{X},\underline{A}),
\cup_{(K;\underline{q})}\,\widehat\otimes\,{\triangledown}\!_{\underline{q}}),
\end{array}$$
where $(H^{\,*}({\cal Z}(K;\underline{X},\underline{A})),\cup_{{\cal Z}(K;\underline{X},\underline{A})})$
is the singular (simplicial) cohomology algebra of ${\cal Z}(K;\underline{X},\underline{A})$,
$(H^{\,*}_{\!\XX_m}(K),\widehat\cup_K)$ is the universal cohomology algebra,
$(H^{\,*}_{\!\DD}(K),\cup_K)$ is the total cohomology algebra on $\DD$ with respect to $\underline{\psi}$
and $(H^{\,*}_{\!\SS}(K),\cup_{(K;\underline{q})})$ is as defined in Definition~9.2.

If each pair $(X_k,A_k)$ is normal, then all the ${\triangledown}\!_{\underline{q}}$ in the above equalities
can be replaced by ${\triangledown}\!_{(\underline{X},\underline{A})}$ and
$(H_*^{\XX_m}(K),\widehat\cup_K)$ can be replaced by the normal cohomology algebra $(H_*^{\XX_m}(K),\w\cup_K)$.

If each $\theta_k\colon H_*(A_k)\to H_*(X_k)$ is an epimorphism, then we have
$$(H^*({\cal Z}(K;\underline{X},\underline{A})),\cup_{{\cal Z}(K;\underline{X},\underline{A})})
\cong(H^{\,*}_{\!\RR_m}(K)\,\widehat\otimes\, H^{\,*}_{\!\RR_m}(\underline{X},\underline{A}),
\widehat\cup_K\,\widehat\otimes\,{\triangledown}\!_{\underline{q}}).
$$

If each pair $(X_k,A_k)$ is normal and each $\theta_k$ is an epimorphism, then
$(H^*(A_k),{\triangledown}_{A_k})=(H^{\XX}(X_k,A_k),{\triangledown}_{\!(X_k,A_k)})$ and we have
$$(H^*({\cal Z}(K;\underline{X},\underline{A})),\cup_{{\cal Z}(K;\underline{X},\underline{A})})
\cong(H^{\,*}_{\!\RR_m}(K)\,\widehat\otimes\,\big(\otimes_k H^*(A_k)\big),
\w\cup_K\,\widehat\otimes\,(\otimes_k{\triangledown}_{\!A_k})),
$$
where the $\RR_m$ index of $\otimes_k H^*(A_k)=H^*_{\RR_m}(\underline{X},\underline{A})$ is defined as follows.
For $a_k\in H^*(A_k)$, $a_1{\otimes}{\cdots}{\otimes}a_m\in H^*_{\emptyset,\omega}(\underline{X},\underline{A})$
with $\omega=\{k\,|\,a_k\in{\rm ker}\,\theta_k\}$.
\vspace{2mm}

Proof}\, Corollary of Theorem~8.5 and Theorem~8.7 by taking $(D_{k\,*},C_{k\,*})$ to be
the singular (simplicial) chain coalgebra pair of $(X_k,A_k)$.
\hfill$\Box$\vspace{3mm}

{\bf Example 9.4} Let $M={\cal Z}(K;\underline{X},\underline{A})$ be the polyhedral product space in Theorem~9.3.
We compute the cohomology algebra $(H^*(M),\cup_M)$ in the following cases.

(1) Suppose each $\theta_k\colon H_*(A_k)\to H_*(X_k)$ is a monomorphism.
Let $J_k={\rm coker}\,\theta_k$ and $I_k={\rm ker}\,\theta_k^\circ$.
This is the only case when we may neglect the diagonal tensor product structure.
By K\"{u}nneth Theorem, $H_*(M)$ is the subcoalgebra $\oplus_{\sigma\in K}\,(\otimes_{k\in\sigma}J_k){\otimes}(\otimes_{k\notin\sigma}H_*(A_k))$
of $\otimes_{k=1}^m H_*(X_k)$.
$H^*(M)$ is the quotient algebra of $\otimes_{k=1}^m H^*(X_k)$ over the ideal
$\oplus_{\tau\notin K}\,(\otimes_{k\in\tau}I_k){\otimes}(\otimes_{k\notin\tau}H^*(A_k))$.

Now we compute the cohomology algebra by Theorem~9.3.
The support index set of $(X_k,A_k)$ is $\{{\mathpzc i},{\mathpzc e}\}$.
So $H^{\,*}_{\sigma\!,\omega}(\underline{X},\underline{A})=0$ if $\omega\neq\emptyset$ or $\sigma\notin K$.
Since $H_*^{\sigma\!,\,\emptyset}(K)=\Bbb Z$ for $\sigma\in K$,
we may identify $H^*_{\sigma\!,\emptyset}(K){\otimes}H^*_{\sigma\!,\emptyset}(\underline{X},\underline{A})$
with $H^{\,*}_{\sigma\!,\emptyset}(\underline{X},\underline{A})$.
Then we have
$$H_*(M)=\oplus_{\sigma\in K}H_*^{\sigma\!,\emptyset}(\underline{X},\underline{A})
=\oplus_{\sigma\in K}\,(\otimes_{k\in\sigma}J_k){\otimes}(\otimes_{k\notin\sigma}H_*(A_k)),$$
$$H^*(M)=\oplus_{\sigma\in K}H^*_{\sigma\!,\emptyset}(\underline{X},\underline{A})
=\oplus_{\sigma\in K}\,(\otimes_{k\in\sigma}I_k){\otimes}(\otimes_{k\notin\sigma}H^*(A_k)).$$
Since $H_*(M)$ is a subcoalgebra of $\otimes_{k=1}^m H_*(X_k)$, $H^*(M)$ is the quotient algebra of
$\otimes_{k=1}^m H^*(X_k)$ over the ideal
$\oplus_{\tau\notin K}\,(\otimes_{k\in\tau}I_k){\otimes}(\otimes_{k\notin\tau}H^*(A_k))$.

Take each $(X_k,A_k)=(\Bbb CP^{\infty}\!,\{*\})$, then
$(H^*(M),\cup_M)$ is the Stanley-Reisner face ring $F(K)=\mak[m]/I_K$, where $I_K$ is as in Definition~5.8.
Dually, $H_*(M)$ is the dual coalgebra of $F(K)$.

(2) Suppose each $X_k$ is contractible.
Let $\{*,a_{k,1},{\cdots}\}$ be the base points of all path-connected components of $A_k$
and $*$ is the base point of $A_k$.
Apply Definition~9.1 for the pair $(X_k,A_k)$.
We have ${\mathpzc e}_{\,k*}=0$, ${\mathpzc i}_{\,k*}\cong H_0(X_k)=\Bbb Z(1)$ with $1=[*]$.
Then ${\mathpzc n}_{\,k0}=\w H_0(A_k)$ is freely generated by all $\varepsilon_{k,i}=[a_{k,i}{-}*]$
with $i>0$.
By definition,
${\vartriangle_{A_k}}(1)=1{\otimes}1$,
${\vartriangle_{A_k}}(\varepsilon_{k,i})=[a_{k,i}{\otimes}a_{k,i}-*{\otimes}*]
=1{\otimes}\varepsilon_{k,i}+\varepsilon_{k,i}{\otimes}1+\varepsilon_{k,i}{\otimes}\varepsilon_{k,i}$.
For $t>0$ and $a\in{\mathpzc n}_{\,kt}=\w H_t(A_k)$, suppose
${\vartriangle_{A_k}}(a)=1{\otimes}a+a{\otimes}1+\Sigma a'_j{\otimes}a''_j$ with $a'_j,a''_j\in \w H_*(A_k)$,
then ${\vartriangle}_{q_k}^{\!\XX}(\overline a)=1{\otimes}\overline a+\overline a{\otimes}1+\Sigma\overline a'_j{\otimes}a''_j$,
for $1{\otimes}1$ can not appear in ${\vartriangle}_{q_k}^{\!\XX}(\overline a)$ by degree comparing.
So we have

${\vartriangle}_{q_k}^{\!\XX}({\mathpzc i}_{\,k*})\subset{\mathpzc i}_{\,k*}{\otimes}{\mathpzc i}_{\,k*}$.

${\vartriangle}_{q_k}^{\!\XX}({\mathpzc n}_{\,k*})\subset{\mathpzc n}_{\,k*}{\otimes}{\mathpzc n}_{\,k*}
\oplus {\mathpzc i}_{\,k*}{\otimes}{\mathpzc n}_{\,k*}\oplus {\mathpzc n}_{\,k*}{\otimes}{\mathpzc i}_{\,k*}$.

${\vartriangle}_{q_k}^{\!\XX}({\bar{\mathpzc n}_{\,k*}})\subset{\bar{\mathpzc n}_{\,k*}}{\otimes}{\mathpzc n}_{\,k*}
\oplus {\mathpzc i}_{\,k*}{\otimes}{\bar{\mathpzc n}_{\,k*}}\oplus {\bar{\mathpzc n}_{\,k*}}{\otimes}{\mathpzc i}_{\,k*}$.

So $(X_k,A_k)$ is normal and the atom coproduct $\psi_{q_k}$ is the right coproduct $\tilde\psi'$ in Definition~7.1.
Take the $\DD$ and $\underline{\psi}$ in Theorem~9.3 respectively to be $\RR_m$ and $\tilde\psi'^{(m)}$, we have
$$(H^*(M),\cup_M)\cong\big(H^{\,*}_{\!\RR_m}(K)\widehat\otimes(\otimes_{k=1}^m H^{*}(A_k)),
\tilde\cup_K\widehat\otimes(\otimes_{k=1}^m{\cup}\!_{A_k})\big),$$
where $(H^*_{\!\RR_m}(K),\tilde\cup_K)$ is the right strictly normal product in Definition~7.4.

Take each $(X_k,A_k)=(D^1\!,S^0)$ and identify $a\widehat\otimes x$ with $a$ as in Example~4.8.
Then we have
$$(H^*({\cal Z}(K;D^1\!,S^0)),\cup_{{\cal Z}(K;D^1\!,S^0)})\cong\big(H^{\,*}_{\!\RR_m}(K),\tilde\cup_K\big).$$
This shows that the right strictly normal cohomology algebra $H^{\,*}_{\!\RR_m}(K)$ is an associative, commutative algebra with unit.

(3) Suppose for each $k$, there is a contractible space $Y_k$ such that $A_k\subset Y_k\subset X_k$.
Apply Definition~9.1 for the pair $(X_k,A_k)$.
We have ${\mathpzc i}_{\,k*}={\mathpzc i}_{\,k0}=\Bbb Z(1)$ with $1=[*]$,
${\mathpzc n}_{\,k*}=\w H_*(A_k)$ and ${\mathpzc e}_{\,k*}=\w H_*(X_k)$.

The restriction of ${\vartriangle}_{q_k}^{\!\XX}$ on ${\mathpzc i}_{\,k*}\oplus{\mathpzc n}_{\,k*}{\oplus}\bar{\mathpzc n}_{\,k*}$
is the same as (2).
Let $\{b_{k,1},b_{k,2},{\cdots}\}$ be the base points of all path-connected components of $X_k$ other than $Y_k$.
Then ${\mathpzc e}_{\,k0}$ is freely generated by all $\varepsilon_{k,i}=[b_{k,i}{-}*]$.
By definition,
${\vartriangle_{X_k}}(\varepsilon_{k,i})=[b_{k,i}{\otimes}b_{k,i}-*{\otimes}*]
=1{\otimes}\varepsilon_{k,i}+\varepsilon_{k,i}{\otimes}1+\varepsilon_{k,i}{\otimes}\varepsilon_{k,i}$.
For $t>0$, $i_{\,k*}{\otimes}i_{\,k*}=i_{\,k0}{\otimes}i_{\,k0}$ can not be a summand of ${\vartriangle_{X_k}}({\mathpzc e}_{\,kt})$ by degree comparing.
So we have

${\vartriangle}_{q_k}^{\!\XX}({\mathpzc i}_{\,k*})\subset{\mathpzc i}_{\,k*}{\otimes}{\mathpzc i}_{\,k*}$.

${\vartriangle}_{q_k}^{\!\XX}({\mathpzc n}_{\,k*})\subset{\mathpzc n}_{\,k*}{\otimes}{\mathpzc n}_{\,k*}
\oplus {\mathpzc n}_{\,k*}{\otimes}{\mathpzc i}_{\,k*}\oplus {\mathpzc i}_{\,k*}{\otimes}{\mathpzc n}_{\,k*}$.

${\vartriangle}_{q_k}^{\!\XX}(\bar{\mathpzc n}_{\,k*})\subset{\bar{\mathpzc n}}_{\,k*}{\otimes}{\mathpzc n}_{\,k*}
\oplus {\mathpzc i}_{\,k*}{\otimes}{\bar{\mathpzc n}}_{\,k*}\oplus {\bar{\mathpzc n}}_{\,k*}{\otimes}{\mathpzc i}_{\,k*}$.

${\vartriangle}_{q_k}^{\!\XX}({\mathpzc e}_{\,k*})\subset{\mathpzc e}_{\,k*}{\otimes}{\mathpzc e}_{\,k*}
\oplus {\mathpzc i}_{\,k*}{\otimes}{\mathpzc e}_{\,k*}\oplus {\mathpzc e}_{\,k*}{\otimes}{\mathpzc i}_{\,k*}$.

This implies the pair $(X_k,A_k)$ is normal and the atom coproduct $\psi_{q_k}$ is the coproduct $\tilde\psi$ in Definition~7.1.
Note that $(H_*^\XX\!(X_k,A_k),\vartriangle_{(X_k,A_k)})$ happens to be the homology coalgebra $(H_*(X_k{\vee}A_k),{\vartriangle}_{X_k\vee A_k})$.
We identify the two coalgebras and their dual algebras and have
$$(H^*(M),\cup_M)\cong\big(H^{\,*}_{\!\XX_m}(K)\widehat\otimes(\otimes_{k=1}^m H^{*}(X_k{\vee}A_k)),
\tilde\cup_K\widehat\otimes(\otimes_{k=1}^m{\cup}_{\!X_k{\vee}A_k})\big),$$
where $\tilde\cup_K$ is the strictly normal product in Definition~7.4.

Take each $(X_k,A_k)=(S^2,S^0)$ and identify $a\widehat\otimes x$ with $a$ as in Example~4.7.
Then we have
$$(H^*({\cal Z}(K;S^2\!,S^0)),\cup_{{\cal Z}(K;S^2\!,S^0)})\cong\big(H^{\,*}_{\!\XX_m}(K),\tilde\cup_K\big).$$
This shows that the strictly normal cohomology algebra $(H^{\,*}_{\!\XX_m}(K),\tilde\cup_K)$ is an associative, commutative algebra with unit.

(4) Suppose each $(X_k,A_k)=(\Sigma Y_k,\Sigma B_k)$, where $\Sigma$ means the suspension as defined before Definition~4.1
and the pair is regarded as unpointed.
Apply Definition~9.1 for the pair $(X_k,A_k)$.
We have ${\mathpzc i}_{\,k\,*}={\mathpzc i}_{\,k\,0}=\Bbb Z(1)$ with $1=[*]$ ( $*$ the base point).
${\vartriangle_{A_k}}(1)=1{\otimes}1$.
All non-unit generators $x$ of $H_*(X_k)$ and $H_*(A_k)$ are primitive, i.e., ${\vartriangle_{-}}(x)=1{\otimes}x+x{\otimes}1$.
So by degree comparing, we have

${\vartriangle}_{q_k}^{\!\XX}({\mathpzc i}_{\,k*})\subset{\mathpzc i}_{\,k*}{\otimes}{\mathpzc i}_{\,k*}$.

${\vartriangle}_{q_k}^{\!\XX}({\mathpzc n}_{\,k*})\subset {\mathpzc i}_{\,k*}{\otimes}{\mathpzc n}_{\,k*}\oplus {\mathpzc n}_{\,k*}{\otimes}{\mathpzc i}_{\,k*}$.

${\vartriangle}_{q_k}^{\!\XX}({\bar{\mathpzc n}}_{\,k*})\subset
{\mathpzc i}_{\,k*}{\otimes}{\bar{\mathpzc n}}_{\,k*}\oplus {\bar{\mathpzc n}}_{\,k*}{\otimes}{\mathpzc i}_{\,k*}$.

${\vartriangle}_{q_k}^{\!\XX}({\mathpzc e}_{\,k*})\subset{\mathpzc i}_{\,k*}{\otimes}{\mathpzc e}_{\,k*}\oplus {\mathpzc e}_{\,k*}{\otimes}{\mathpzc i}_{\,k*}$.

This implies the pair $(X_k,A_k)$ is normal and the atom coproduct $\psi_{q_k}$ is
the coproduct $\overline\psi$ in Definition~7.1.
$(H^{\,*}_{\!\XX}(X_k,A_k),{\triangledown}_{(X_k,A_k)})$
is a trivial algebra (a free algebra with unit satisfying $xy=0$ for all non-unit generators $x,y$).
Denote by $0$ the product of a trivial algebra. The tensor product of trivial algebras is also a trivial algebra.
So we have
$$(H^*(M),\cup_M)\cong\big(H^{\,*}_{\!\XX_m}(K)\widehat\otimes (\otimes_{k=1}^mH^{\,*}_{\!\XX}(X_k,A_k)),
\overline\cup_K\widehat\otimes 0\big),$$
where $(H^{\,*}_{\!\XX_m}(K),\overline\cup_K)$ is the special cohomology algebra in Definition~7.4.
If each $\theta_k$ is an epimorphism, then we have
$$(H^*(M),\cup_M)\cong\big(H^{\,*}_{\!\RR_m}(K)\widehat\otimes (\otimes_{k=1}^m H^*(A_k)),
\overline\cup_K\widehat\otimes 0\big).$$

Take each $(X_k,A_k)=(S^4,S^2)$ and identify $a\widehat\otimes x$ with $a$ as in Example~4.7.
Then we have
$$(H^*({\cal Z}(K;S^4\!,S^2)),\cup_{{\cal Z}(K;S^4\!,S^2)})\cong\big(H^{\,*}_{\!\XX_m}(K),\overline\cup_K\big).$$
This shows that the special cohomology algebra $(H^{\,*}_{\!\XX_m}(K),\overline\cup_K)$ is an associative, commutative algebra with unit.

Take each $(X_k,A_k)=(S^r,S^p)$ with $r>p>0$ and identify $a\widehat\otimes x$ with $a$ as in Example~4.7.
Then we have
$$(H^*({\cal Z}(K;S^r\!,S^p)),\cup_{{\cal Z}(K;S^r\!,S^p)})\cong\big(H^{\,*}_{\!\XX_m}(K),\overline\cup_K^{r\!,p}\big),$$
where the product $\overline\cup_K^{r\!,p}$ is defined as follows.
For $a\in H^s_{\sigma'\!,\,\omega'}(K)$ and $b\in H^t_{\sigma''\!,\,\omega''}(K)$,
$a\,\overline\cup_K^{r\!,p} b=(-1)^{t(r|\sigma'|+p|\omega'|)}a\,\overline\cup_K b$.

Take each $(X_k,A_k)=(D^{3}\!,S^2)$ and identify $a\widehat\otimes x$ with $a$ as in Example~4.8.
Then we have
$$(H^*({\cal Z}(K;D^{3}\!,S^2)),\cup_{{\cal Z}(K;D^{3}\!,S^2)})\cong\big(H^{\,*}_{\!\RR_m}(K),\overline\cup_K\big).$$
This shows that the right special cohomology algebra $(H^{\,*}_{\!\RR_m}(K),\overline\cup_K)$ is an associative, commutative algebra with unit.

Take each $(X_k,A_k)=(D^{n+1}\!,S^n)$ with $n>0$ and identify $a\widehat\otimes x$ with $a$ as in Example~4.8.
Then we have
$$(H^*({\cal Z}(K;D^{n+1}\!,S^n)),\cup_{{\cal Z}(K;D^{n+1}\!,S^n)})\cong\big(H^{\,*}_{\!\RR_m}(K),\overline\cup_K^{\,n}\big),$$
where the product $\overline\cup_K^{\,n}$ is defined as follows.
For $a\in H^s_{\emptyset,\omega'}(K)$ and $b\in H^t_{\emptyset,\omega''}(K)$,
$a\,\overline\cup_K^{\,n} b=(-1)^{tn|\omega'|}a\,\overline\cup_K b$.

(5) Consider the map
$$f\colon S^3\stackrel{\mu'}{-\!\!\!-\!\!\!\longrightarrow}S^3\vee S^3
\stackrel{g\,\vee\, 1}{-\!\!\!-\!\!\!\longrightarrow} S^2\vee S^3,\vspace{-2mm}$$
where $\mu'$ is the coproduct of the co-$H$-space $S^3$, $g\colon S^3\to S^2$ is the Hopf bundle and $1$ is the identity map of $S^3$.
Let $A_1= A_2= S^2\vee S^3$, $X_1= S^2\vee CS^3$, $X_2= C_f$, where $C$ means the cone of
a space and $C_f$ means the mapping cone of $f$.
The equality $H_*(X_i)\cong H_*(S^2)$ implies $X_1\simeq X_2\simeq S^2$.
Apply Definition~9.1 for the pair $(X,A)=(X_i,A_i)$ for $i=1,2$.
By definition, $(C_*^\XX\!(X_1{,}A_1),d)=(C_*^\XX\!(X_2{,}A_2),d)$ with
${\mathpzc i}_{\,*}=\Bbb Z(1,x)$, ${\mathpzc n}_{\,*}=\Bbb Z(y)$, ${\bar{\mathpzc n}}_{\,*}=\Bbb Z(\overline y)$, ${\bar{\mathpzc e}}_{\,*}=0$,
where $1$ is the generator of $H_0(A_i)$, $x$ is the generator of $H_2(A_i)$, $y$ is the generator of $H_3(A_i)$.
But $(C_*^\XX\!(X_1{,}A_1),{\vartriangle}_{q_1}^{\!\XX})\neq(C_*^\XX\!(X_2{,}A_2),{\vartriangle}_{q_2}^{\!\XX})$.

${\vartriangle}_{q_1}^{\!\XX}(1)=1{\otimes}1$,
${\vartriangle}_{q_2}^{\!\XX}(1)=1{\otimes}1$.

${\vartriangle}_{q_1}^{\!\XX}(x)=1{\otimes}x+x{\otimes}1$,
${\vartriangle}_{q_2}^{\!\XX}(x)=1{\otimes}x+x{\otimes}1$.

${\vartriangle}_{q_1}^{\!\XX}(y)=1{\otimes}y+y{\otimes}1$,
${\vartriangle}_{q_2}^{\!\XX}(y)=1{\otimes}y+y{\otimes}1$.

${\vartriangle}_{q_1}^{\!\XX}(\overline y)=1{\otimes}\overline y+\overline y{\otimes}1$,
${\vartriangle}_{q_2}^{\!\XX}(\overline y)=1{\otimes}\overline y+\overline y{\otimes}1+x{\otimes}x$.

So $(X_1,A_1)$ is normal but $(X_2,A_2)$ is not and
the atom coproducts $\psi_{q_1}$ and $\psi_{q_2}$ are respectively the right coproduct $\overline\psi'_r$ and $\bar\psi'_r$ in Definition~7.1.
Dually,
$$(H^{\,*}_\RR(X_1,A_1),{\triangledown}_{\!(X_1,A_1)})=\Bbb Z[x,y]/(x^2,y^2,xy),$$
$$(H^{\,*}_\RR(X_2,A_2),{\triangledown}_{\!q_2})=\Bbb Z[x,y]/(x^2{-}y,y^2,xy),$$
where $\Bbb Z[x,y]$ is the polynomial algebra.
So we have
$$(H^*({\cal Z}(K;X_1,A_1)),\cup_{{\cal Z}(K;X_1,A_1)})
\cong\Big(H^*_{\RR_m}(K)\,\widehat\otimes\,(H^*(A_1)^{\otimes m}),
\overline\cup_K\,\widehat\otimes\,\cup_1\Big),$$
$$(H^*({\cal Z}(K;X_2,A_2)),\cup_{{\cal Z}(K;X_2,A_2)})
\cong\Big(H^*_{\RR_m}(K)\widehat\otimes\,(H^*(A_2)^{\otimes m}),
\bar\cup_K\widehat\otimes\,\cup_2\Big),$$
where $(H^{\,*}_{\!\RR_m}(K),\overline\cup_K)$ and $(H^{\,*}_{\!\RR_m}(K),\bar\cup_K)$ are respectively
the right special and right weakly special cohomology algebra in Definition~7.4 and
$$(H^*(A_1)^{\otimes m},\cup_1)=\Bbb Z[x_1,{\cdots},x_m,y_1,{\cdots},y_m]/(x_i^2,y_i^2,x_iy_i),$$
$$(H^*(A_2)^{\otimes m},\cup_2)=\Bbb Z[x_1,{\cdots},x_m,y_1,{\cdots},y_m]/(x_i^2{-}y_i,y_i^2,x_iy_i),$$
where $x_{i_1}{\cdots}x_{i_s}y_{j_1}{\cdots}y_{j_t}\in
H^*_{\emptyset,\{j_1,{\cdots},j_t\}}(\underline{X_i},\underline{A_i})$.
Note that the product $\cup_1$
keeps degree but the product $\cup_2$
does not keep degree.

(6) Let everything be as in (5) and $Y_1=X_1{\vee}S^4$, $Y_2=X_2{\vee}S^4$.
Then $(C_*^\XX\!(Y_1{,}A_1),d)=(C_*^\XX\!(Y_2{,}A_2),d)$ with
${\mathpzc i}_{\,*}=\Bbb Z(1,x)$, ${\mathpzc n}_{\,*}=\Bbb Z(y)$, ${\bar{\mathpzc n}}_{\,*}=\Bbb Z(\overline y)$,
${\mathpzc e}_{\,*}=\Bbb Z(z)$,
where $1$ is the generator of $H_0(A_i)$,  $x$ is the generator of $H_2(A_i)$, $y$ is the generator of $H_3(A_i)$, $z$ is the generator of $H_4(Y_i)$.

${\vartriangle}_{q_1}^{\!\XX}(1)=1{\otimes}1$,
${\vartriangle}_{q_2}^{\!\XX}(1)=1{\otimes}1$.

${\vartriangle}_{q_1}^{\!\XX}(x)=1{\otimes}x+x{\otimes}1$,
${\vartriangle}_{q_2}^{\!\XX}(x)=1{\otimes}x+x{\otimes}1$.

${\vartriangle}_{q_1}^{\!\XX}(y)=1{\otimes}y+y{\otimes}1$,
${\vartriangle}_{q_2}^{\!\XX}(y)=1{\otimes}y+y{\otimes}1$.

${\vartriangle}_{q_1}^{\!\XX}(\overline y)=1{\otimes}\overline y+\overline y{\otimes}1$,
${\vartriangle}_{q_2}^{\!\XX}(\overline y)=1{\otimes}\overline y+\overline y{\otimes}1+x{\otimes}x$.

${\vartriangle}_{q_1}^{\!\XX}(z)=1{\otimes}z+z{\otimes}1$,
${\vartriangle}_{q_2}^{\!\XX}(z)=1{\otimes}z+z{\otimes}1$.

So $(Y_1,A_1)$ is normal but $(Y_2,A_2)$ is not and the atom coproducts $\psi_{q_1}$ and $\psi_{q_2}$
are respectively the coproduct $\overline\psi$ and $\bar\psi$ in Definition~7.1. Dually,
$$(H^{\,*}_\XX(Y_1,A_1),{\triangledown}\!_{(Y_1,A_1)})=\Bbb Z[x,y,z]/(x^2,y^2,z^2,xy,yz,xz),$$
$$(H^{\,*}_\XX(Y_2,A_2),{\triangledown}\!_{q_2})=\Bbb Z[x,y,z]/(x^2{-}y,y^2,z^2,xy,yz,xz).$$
So we have
$$(H^*({\cal Z}(K;Y_1,A_1)),\cup)\cong\Big(H^{\,*}_{\!\XX_m}(K)\,\widehat\otimes\,R_1,
\overline\cup_K\,\widehat\otimes\,{\cup}_1\Big),$$
$$(H^*({\cal Z}(K;Y_2,A_2)),\cup)\cong\Big(H^{\,*}_{\!\XX_m}(K)\widehat\otimes\,R_2,
\bar\cup_K\widehat\otimes\,{\cup}_2\Big),$$
where $(H^{\,*}_{\XX_m}(K),\overline\cup_K)$ and $(H^{\,*}_{\XX_m}(K),\bar\cup_K)$ are respectively
the special and weakly special cohomology algebra in Definition~7.4 and
$$(R_1,{\cup}_1)=\Bbb Z[x_1,{\cdots},x_m,y_1,{\cdots},y_m,z_1,{\cdots},z_m]/(x_i^2,y_i^2,z_i^2,x_iy_i,x_iz_i,y_iz_i),$$
$$(R_2,{\cup}_2)=\Bbb Z[x_1,{\cdots},x_m,y_1,{\cdots},y_m,z_1,{\cdots},z_m]/(x_i^2{-}y_i,y_i^2,z_i^2,x_iy_i,x_iz_i,y_iz_i),$$
where $x_{i_1}{\cdots}x_{i_s}y_{j_1}{\cdots}y_{j_t}z_{k_1}{\cdots}z_{k_u}\in
H^*_{\{k_1,\cdots,k_u\},\{j_1,{\cdots},j_t\}}(\underline{Y_i},\underline{A_i})=R_i$.
\vspace{3mm}

{\bf Definition 9.5} Let $(T_*^{\XX_m},\underline{\psi})$ be a total chain coalgebra (Definition~7.2)
and $\TT\subset\XX_m$ is an index set.
Then for a simplicial complex $K$ on $[m]$, the total chain coalgebra $(T_*^\TT(K),\psi_K,d)$
of $K$ on $\TT$ with respect to $\underline{\psi}$ in Definition~7.4 is a chain $\TT$-coalgebra.

For a simplicial pair $(X,A)$ on $[m]$ densely split on $\TT$ by Definition~4.11,
the chain complex pair $\big((T_*^\TT(X),\psi_X,d),(T_*^\TT(A),\psi_A,d)\big)$ is a split coalgebra pair by Definition~8.1 with $\Lambda=\TT$.
So the (co)homology groups in Definition~4.11 are (co)algebras as defined in Definition~8.2
when the homotopy equivalence $q$ is chosen.
Precisely, replace the $(D_*,C_*)$ in Definition~8.2 by $(X,A)$ and we have all definitions for $(X,A)$ as follows.

Let $\theta\colon(H_*^\TT(A),{\vartriangle}_{A})\to(H_*^\TT(X),{\vartriangle}_X)$ be the total homology coalgebra
(with respect to $\underline{\psi}$) homomorphism induced by inclusion.
Denote by
$${\mathpzc i}_{\,*}={\rm coim}\,\theta,\quad{\mathpzc n}_{\,*}={\rm ker}\,\theta,\quad
\overline{\mathpzc n}_{\,*}=\Sigma{\rm ker}\,\theta,\quad{\mathpzc e}_{\,*}={\rm coker}\,\theta.$$
Then we have group equalities $H_*^\TT(A)={\mathpzc i}_{\,*}\oplus{\mathpzc n}_{\,*}$,
$H_*^{\XX;\TT}(X,A)={\mathpzc i}_{\,*}\oplus{\mathpzc n}_{\,*}\oplus{\mathpzc e}_{\,*}$ and
$C_*^{\XX;\TT}(X,A)={\mathpzc i}_{\,*}\oplus{\mathpzc n}_{\,*}\oplus\overline{\mathpzc n}_{\,*}\oplus{\mathpzc e}_{\,*}$.

The index set $\TT$ is called the {\it reference index set} of $(X,A)$ as in Definition~4.11.
The coproduct $\underline{\psi}$ is called the {\it reference coproduct} of $(X,A)$.

The {\it densely normal homology coalgebra} $(H_*^{\XX;\TT}(X,A),\vartriangle_{(X,A;\underline{\psi})})$
of $(X,A)$ (irrelevant to $q$) on $\TT$ with respect to $\underline{\psi}$ is given by
$${\vartriangle}_{(X,A;\underline{\psi})}(x)=\left\{\begin{array}{cl}
{\vartriangle}_A(x)&{\rm if}\,\,x\in {\mathpzc i}_{\,*},\vspace{1mm}\\
{\vartriangle}_A(x)&{\rm if}\,\,x\in {\mathpzc n}_{\,*},\vspace{1mm}\\
{\vartriangle}_X(x)&{\rm if}\,\,x\in {\mathpzc e}_{\,*}.
\end{array}\right.$$

The {\it densely normal cohomology algebra} $(H^{\,*}_{\!\XX;\TT}(X,A),{\triangledown}_{(X,A;\underline{\psi})})$ of $(X,A)$
on $\TT$ with respect to $\underline{\psi}$ is the dual algebra of $(H_*^{\XX;\TT}(X,A),\vartriangle_{(X,A;\underline{\psi})})$.

The {\it dense character chain coalgebra} $(C_*^{\XX;\TT}(X,A),{\vartriangle}_{(q;\underline{\psi})}^{\!\XX})$
of $(X,A)$ on $\TT$ with respect to $q$ and $\underline{\psi}$ is the chain $(\XX\!{\times}\TT)$-coalgebra defined as follows.
The restriction of ${\vartriangle}_{(q;\underline{\psi})}^{\!\XX}$ on $H_*^{\XX;\TT}(X,A)$ is ${\vartriangle}_{(X,A;\underline{\psi})}$.
For $x\in{\mathpzc n}_{\,*}$ with ${\vartriangle}_A(x)=\Sigma x'_i{\otimes}x''_i+\Sigma y'_j{\otimes}y''_j$,
where each $x'_i\in{\mathpzc n}_{\,*}$ and each $y'_j\notin{\mathpzc n}_{\,*}$ but $y''_j\in{\mathpzc n}_{\,*}$, define
${\vartriangle}_{(q;\underline{\psi})}^{\!\XX}(\overline x)=\Sigma\overline x'_i{\otimes}x''_i+\Sigma(-1)^{|y'_j|}y'_j{\otimes}\overline y''_j
+\xi(q{\otimes}q)\psi_X(\overline x)$ with $\xi$ the projection from $C_*^{\XX;\TT}(X,A){\otimes}C_*^{\XX;\TT}(X,A)$ to its homology.

The {\it densely indexed homology coalgebra} $(H_*^{\XX;\TT}(X,A),{\vartriangle}_{(q;\underline{\psi})})$
of $(X,A)$ on $\TT$ with respect to $q$ and $\underline{\psi}$ is the $(\XX\!{\times}\TT)$-coalgebra given by
$${\vartriangle}_{(q;\underline{\psi})}(x)=\left\{\begin{array}{ll}
{\vartriangle}_A(x)&{\rm if}\,\,x\in{\mathpzc i}_{\,*},\vspace{1mm}\\
{\vartriangle}_A(x){+}\xi(q{\otimes}q)\psi_X(\overline x)&{\rm if}\,\,x\in {\mathpzc n}_{\,*},\vspace{1mm}\\
{\vartriangle}_X(x)&{\rm if}\,\,x\in {\mathpzc e}_{\,*}.
\end{array}\right.$$

The {\it densely indexed cohomology algebra} $(H^{\,*}_{\!\XX;\TT}(X,A),{\triangledown}_{\!(q;\underline{\psi})})$
of $(X,A)$ on $\TT$ with respect to $q$ and $\underline{\psi}$ is the dual algebra of $(H_*^{\XX;\TT}(X,A),{\vartriangle}_{(q;\underline{\psi})})$.

The densely split pair $(X,A)$ is called {\it densely normal} on $\TT$ with respect to $\underline{\psi}$
if there is a homotopy equivalence $q$ such that ${\vartriangle}_{(q;\underline{\psi})}={\vartriangle}_{(X,A;\underline{\psi})}$
and ${\triangledown}_{\!(q;\underline{\psi})}={\triangledown}_{\!(X,A;\underline{\psi})}$.

The {\it atom coproduct} $\psi_{(q;\underline{\psi})}$ and
the {\it atom chain coalgebra $(T_*^\XX\!,\psi_{(q;\underline{\psi})})$} with respect to $q$ and $\underline{\psi}$ are defined as follows.
For ${\mathpzc s}\in T_*^\XX$, suppose ${\vartriangle}_q^{\!\XX}({\mathpzc s}_*)\subset\oplus_i\,({\mathpzc s}'_i)_*{\otimes}({\mathpzc s}''_i)_*$
but ${\vartriangle}_q^{\!\XX}({\mathpzc s}_*){\cap}(({\mathpzc s}'_i)_*{\otimes}({\mathpzc s}''_i)_*)\neq 0$ for each $i$,
then define $\psi_{(q;\underline{\psi})}({\mathpzc s})=\Sigma_i\,{\mathpzc s}'_i{\otimes}{\mathpzc s}''_i$.

The {\it dense index set} of $(X,A)$ on $\TT$ with respect to $\underline{\psi}$ is $\SS=\{{\mathpzc s}\in T_*^\XX\,|\,{\mathpzc s}_{\,*}\neq 0\}$.
By definition, $\psi_{(q;\underline{\psi})}({\mathpzc s})=0$ if ${\mathpzc s}\notin\SS$.
\vspace{3mm}

{\bf Definition 9.6} Let $(\underline{X},\underline{A})=\{(X_k,A_k)\}_{k=1}^m$
be such that each simplicial complex pair $(X_k,A_k)$ on $[n_k]$ is densely split on $\TT_k$ with respect to $\underline{\psi_k}$.
Then the (co)homology groups in Definition~4.12 are (co)algebras as defined in Definition~8.4 when
all the homotopy equivalences $q_k$ for $(X_k,A_k)$ are chosen.
Precisely, replace the $(D_{k\,*},C_{k\,*})$ in Definition~8.4 by $(X_k,A_k)$ and we have all the definitions for $(\underline{X},\underline{A})$.
We only list the definitions needed in Theorem~9.7.

The {\it reference index set} of $(\underline{X},\underline{A})$ is $\underline{\TT}=\TT_1{\times}{\cdots}{\times}\TT_m$,
where each $\TT_k$ is the reference index set of $(X_k,A_k)$.
The {\it reference coproduct} of $(\underline{X},\underline{A})$ is
$\underline{\psi_n}=\underline{\psi_1}{\otimes}{\cdots}{\otimes}\underline{\psi_m}$,
where each $\underline{\psi_k}$ is the reference coproduct of $(X_k,A_k)$ and $[n]=[n_1]\sqcup{\cdots}\sqcup[n_m]$.

The {\it densely indexed homology coalgebra} and {\it densely indexed cohomology algebra}
of $(\underline{X},\underline{A})$ on $\underline{\TT}$ with respect to $\underline{q}$ and $\underline{\psi_n}$
are the  $(\XX_m{\times}\underline{\TT})$-(co)algebras
$$(H_*^{\XX_m;\underline{\TT}}(\underline{X},\underline{A}),{\vartriangle}_{(\underline{q};\underline{\psi_n})})
=(\otimes_{k=1}^m\,H_*^{\XX;\TT_k}(X_k,A_k),\otimes_{k=1}^m\,{\vartriangle}_{(q_k;\underline{\psi_k})}),$$
$$(H^{\,*}_{\!\XX_m;\underline{\TT}}(\underline{X},\underline{A}),{\triangledown}\!_{(\underline{q};\underline{\psi_n})})
=(\otimes_{k=1}^m\,H^{\,*}_{\!\XX;\TT_k}(X_k,A_k),\otimes_{k=1}^m\,{\triangledown}\!_{(q_k;\underline{\psi_k})}).$$

The {\it densely normal homology coalgebra} and {\it densely normal cohomology algebra}
of $(\underline{X},\underline{A})$ on $\underline{\TT}$ with respect to $\underline{\psi_n}$ (irrelevant to $\underline{q}$)
are the following $(\XX_m{\times}\underline{\TT})$-(co)algebras
$$(H_*^{\XX_m;\underline{\TT}}(\underline{X},\underline{A}),\vartriangle_{(\underline{X},\underline{A};\underline{\psi_n})})
=(\otimes_{k=1}^m\,H_*^{\XX;\TT_k}(X_k,A_k),\otimes_{k=1}^m\,\vartriangle_{(X_k,A_k;\underline{\psi_k})}),$$
$$(H^{\,*}_{\!\XX_m;\underline{\TT}}(\underline{X},\underline{A}),{\triangledown}_{(\underline{X},\underline{A};\underline{\psi_n})})
=(\otimes_{k=1}^m\,H^{\,*}_{\!\XX;\TT_k}(X_k,A_k),\otimes_{k=1}^m\,\triangledown_{(X_k,A_k;\underline{\psi_k})}).$$

For an index set $\DD\subset\XX_m$, the analogue coalgebra $H_*^{\DD;\underline{\TT}}(\underline{X},\underline{A})$
of $(\underline{X},\underline{A})$ on $\DD{\times}\underline{\TT}$
is the restriction coalgebra of $H_*^{\XX_m;\underline{\TT}}(\underline{X},\underline{A})$ on $\DD{\times}\underline{\TT}$ (the dual case is similar).

The {\it total coproduct} $\psi_{(\underline{q};\underline{\psi_n})}$ with respect to $\underline{q}$ and $\underline{\psi_n}$
is $\psi_{(q_1;\underline{\psi_1})}{\otimes}{\cdots}{\otimes}\psi_{(q_m;\underline{\psi_m})}$.
The {\it total chain coalgebra} with respect to $\underline{q}$ and $\underline{\psi_n}$ is
$(T_*^{\XX_m},\psi_{(\underline{q};\underline{\psi_n})})$.

For a simplicial complex $K$ on $[m]$,
the total chain complex of $K$ with respect to $\psi_{(\underline{q};\underline{\psi_n})}$ by Definition~7.4 is denoted by
$(T_*^{\XX_m}(K),\psi_{(K;\underline{q};\psi_n)})$.
The total cohomology algebra of $K$ with respect to $\psi_{(\underline{q};\underline{\psi_n})}$ by Definition~7.4 is denoted by
$(H^{\,*}_{\!\XX_m}(K),\cup_{(K;\underline{q};\psi_n)})$.
For an index set $\DD\subset\XX_m$,
the total cohomology algebra of $K$ on $\DD$ with respect to $\psi_{(\underline{q};\underline{\psi_n})}$ is denoted by
$(H^{\,*}_{\!\DD}(K),\cup_{(K;\underline{q};\underline{\psi_n})})$.
\vspace{3mm}

{\bf Theorem 9.7} {\it Suppose the polyhedral join simplicial complex ${\cal Z}^*(K;\underline{X},\underline{A})$ satisfies
that each $(X_k,A_k)$ is densely split on $\TT_k$ with respect to $\underline{\psi_k}$,
then the cohomology group isomorphisms in Theorem~4.13 are algebra isomorphisms by Theorem~8.4 and Theorem~8.7 as follows.
Suppose $(T_*^{\XX_m},\underline{\psi})=(T_*^\XX{\otimes}{\cdots}{\otimes}T_*^\XX,\psi_1{\otimes}{\cdots}{\otimes}\psi_m)$
is a total chain coalgebra  such that each $\psi_{q_k}\prec\psi_k$.
Let $\underline{\SS}=\SS_1{\times}{\cdots}{\times}\SS_m$,
where each $\SS_k$ is the dense support index set of $(X_k,A_k)$ on $\TT_k$ with respect to $\underline{\psi_k}$.
Then for any index set $\DD$ such that $\underline{\SS}\subset\DD\subset\XX_m$,
we have algebra isomorphisms (index $\XX_m$, $\DD$, $\underline{\SS}$ neglected),

$$\begin{array}{l}
\quad (H^{\,*}_{\!\underline{\TT}}({\cal Z}^*(K;\underline{X},\underline{A})),\cup_{{\cal Z}^*(K;\underline{X},\underline{A})})\vspace{2mm}\\
\cong (H^{\,*}_{\!\XX_m}(K)\,\widehat\otimes\, H^{\,*}_{\!\XX_m;\underline{\TT}}(\underline{X},\underline{A}),
\widehat\cup_K\,\widehat\otimes\,{\triangledown}\!_{(\underline{q};\underline{\psi_n})})\vspace{2mm}\\
\cong (H^{\,*}_{\!\DD}(K)\,\widehat\otimes\, H^{\,*}_{\!\DD;\underline{\TT}}(\underline{X},\underline{A}),
\cup_K\,\widehat\otimes\,{\triangledown}\!_{(\underline{q};\underline{\psi_n})})\vspace{2mm}\\
\cong (H^{\,*}_{\!\underline{\SS}}(K)\,\widehat\otimes\, H^{\,*}_{\!\underline{\SS};\underline{\TT}}(\underline{X},\underline{A}),
\cup_{(K;\underline{q};\underline{\psi_n})}\,\widehat\otimes\,{\triangledown}\!_{(\underline{q};\underline{\psi_n})}),
\end{array}$$
where $(H^{\,*}_{\!\underline{\TT}}(-),\cup_{-})$
is the total cohomology algebra of ${\cal Z}^*(K;\underline{X},\underline{A})$ on $\underline{\TT}$ with respect to $\underline{\psi_n}$,
$(H^{\,*}_{\!\XX_m}(K),\widehat\cup_K)$ is the universal cohomology algebra,
$(H^{\,*}_{\!\DD}(K),\cup_K)$ is the total cohomology algebra on $\DD$ with respect to $\underline{\psi}$
and $(H^{\,*}_{\!\SS}(K),\cup_{(K;\underline{q};\underline{\psi_n})})$ is as defined in Definition~9.6.

If each pair $(X_k,A_k)$ is densely normal on $\TT_k$ with respect to $\underline{\psi_k}$, then
all the ${\triangledown}\!_{(\underline{q};\underline{\psi_n})}$ in the above equalities
can be replaced by ${\triangledown}\!_{(\underline{X},\underline{A};\underline{\psi_n})}$
and $(H_*^{\XX_m}(K),\widehat\cup_K)$ can be replaced by the normal cohomology algebra $(H_*^{\XX_m}(K),\w\cup_K)$.

If each $\theta_k\colon H_*^{\TT_k}(A_k)\to H_*^{\TT_k}(X_k)$ is an epimorphism, then we have
$$\begin{array}{l}
\quad (H^{\,*}_{\!\underline{\TT}}({\cal Z}^*(K;\underline{X},\underline{A})),\cup_{{\cal Z}^*(K;\underline{X},\underline{A})})\vspace{2mm}\\
\cong (H^{\,*}_{\!\RR_m}(K)\,\widehat\otimes\, H^{\,*}_{\!\RR_m;\underline{\TT}}(\underline{X},\underline{A}),
\widehat\cup_K\,\widehat\otimes\,{\triangledown}\!_{(\underline{q};\underline{\psi_n})}),
\end{array}$$
where $(H_*^{\RR_m}(K),\widehat\cup_K)$ is the right universal cohomology algebra.

If each pair $(X_k,A_k)$ is densely normal on $\TT_k$ with respect to $\underline{\psi_k}$ and each $\theta_k$ is an epimorphism, then
$(H^{\,*}_{\!\TT_k}(A_k),{\triangledown}_{\!A_k})=(H^{\,*}_{\!\XX;\TT_k}(X_k,A_k),{\triangledown}_{\!(X_k,A_k;\underline{\psi_k})})$ and we have
$$\begin{array}{l}
\quad (H^{\,*}_{\!\underline{\TT}}({\cal Z}^*(K;\underline{X},\underline{A})),\cup_{{\cal Z}^*(K;\underline{X},\underline{A})})\vspace{2mm}\\
\cong (H^{\,*}_{\!\RR_m}(K)\,\widehat\otimes\,\big(\otimes_k H^{\,*}_{\!\TT_k}(A_k)\big),
\w\cup_K\,\widehat\otimes\,\big(\otimes_k {\triangledown}_{\!A_k}\big),
\end{array}$$
where the $\RR_m$ index of $\otimes_k H^{\,*}_{\!\TT_k}(A_k)=H^*_{\RR_m;\underline{\TT}}(\underline{X},\underline{A})$ is defined as follows.
For $a_k\in H^*_{\TT_k}(A_k)$, $a_1{\otimes}{\cdots}{\otimes}a_m\in H^*_{\emptyset,\omega;\underline{\TT}}(\underline{X},\underline{A})$
with $\omega=\{k\,|\,a_k\in{\rm ker}\,\theta_k\}$.
\vspace{2mm}

Proof}\, Take the $(D_{k\,*},C_{k\,*})$ in Theorem~8.4 to be $(T_*^{\TT_k}(X_k),T_*^{\TT_k}(A_k))$
with reference coproduct $\underline{\psi_k}$.
\hfill$\Box$\vspace{3mm}

{\bf Example 9.8} We compute the right universal cohomology algebra of the composition
complex ${\cal Z}^*(K;L_1,{\cdots},L_m)$ in Example~4.15.
Suppose each $H_*^{\RR_{n_k}}(L_k)$ is a free group or the indexed (co)homology group is taken over a field.

Apply Definition~9.5 to the pair $(\Delta\!^{[n_k]},L_k)$
with reference index set $\RR_{n_k}$ and reference coproduct $\underline{\psi_k}=\widehat\psi^{(n_k)}$,
where $\widehat\psi$ is the universal coproduct in Definition~7.1.
Then $\theta_k\colon H_*^{\RR_{n_k}}(L_k)\to H_*^{\RR_{n_k}}(\Delta\!^{[n_k]})$ is an epimorphism.
So ${\mathpzc i}_{\,k*}={\mathpzc i}_{\,k0}=\Bbb Z(1)$ with $1=[\emptyset]\in H_0^{\emptyset,\emptyset}(L_k)$.
${\mathpzc n}_{\,k*}=H_*^{\overline\RR_{n_k}}(L_k)$ with $\overline\RR_{n_k}=\{(\emptyset,\omega)\in\RR_{n_k}\,|\,\omega\neq\emptyset\}$ as in Example~4.15.
${\mathpzc e}_{\,k*}=0$. Since ${\mathpzc i}_{\,k\,t}=0$ for $t>0$,
the homomorphism
$$\xi(q{\otimes}q)\widehat\psi_{\Delta\!^{[n_k]}}
\colon \overline{\mathpzc n}_{\,k*}\to {\mathpzc i}_{\,k*}{\otimes}{\mathpzc i}_{\,k*}
=H_0^{\emptyset,\emptyset}(\Delta\!^{[n_k]}){\otimes}H_0^{\emptyset,\emptyset}(\Delta\!^{[n_k]})$$
is $0$ by degree compairing.
So the pair $(\Delta\!^{[n_k]},L_k)$ is densely normal on $\RR_{n_k}$ with respect to $\widehat\psi^{(n_k)}$.
By Theorem~9.7, we have
\begin{eqnarray*}&&(H^{\,*}_{\!\RR_n}({\cal Z}^*(K;L_1,{\cdots},L_m)),
\widehat\cup_{{\cal Z}^*(K;L_1,{\cdots},L_m)}))\\
&\cong&\big(H^{\,*}_{\!\RR_m}(K)\,\widehat\otimes\,(\otimes_{k=1}^m H^{\,*}_{\RR_{n_k}}(L_k)),
\w\cup_{K}\,\widehat\otimes\,(\otimes_{k=1}^m{\widehat\triangledown}_{\!L_k})\big).
\end{eqnarray*}

Specifically, take each $L_k=\partial\Delta\!^{[n_k]}$ with $n_k$ an odd number $>2$.
Denote by $g_k$ the generator of $H^{n_k{-}1}_{\emptyset,[n_k]}(L_k)$
and $H_{n_k{-}1}^{\emptyset,[n_k]}(L_k)$.
Then $H_*^{\RR_{n_k}}(L_k)$ is generated by $1$ and $g_k$ with coproduct given by
$\widehat\psi_{L_k}(1)=1{\otimes}1$ and $\widehat\psi_{L_k}(g_k)=1{\otimes}g_k{+}g_k{\otimes}1$.
This implies the atom coproduct $\psi_{(q_k;\underline{\psi_k})}$ is the right special coproduct $\overline\psi'$ in Definition~7.1.
Dually, $(H^{\,*}_{\!\RR_{n_k}}(L_k),\widehat\triangledown_{L_k})$ is the truncated polynomial $\Bbb Z[g_k]/(g_k^2)$.
So we have
\begin{eqnarray*}&&(H^{\,*}_{\!\RR_n}({\cal Z}^*(K;\partial\Delta\!^{[n_1]},{\cdots},\partial\Delta\!^{[n_m]})),
\widehat\cup_{{\cal Z}^*(K;L_1,{\cdots},L_m)}))\\
&\cong&\big(H^{\,*}_{\!\RR_m}(K)\,\widehat\otimes\,\Bbb Z(g_1,{\cdots},g_m)/(g_k^2),
\overline\cup_{K}\,\widehat\otimes\,\times\big),
\end{eqnarray*}
where $(H^{\,*}_{\!\RR_m}(K),\overline\cup_{K})$ is the right special cohomology algebra
and $\times$ is the usual product of polynomials.
For $a\in H^{\,*}_{\emptyset,\omega}(K)$ ($\omega=\{i_1,{\cdots},i_s\}$), identify $a\widehat\otimes(g_{i_1}{\cdots}g_{i_s})$ with $a$.
With this identification,
$$(H^{\,*}_{\!\RR_n}({\cal Z}^*(K;\partial\Delta\!^{[n_1]},{\cdots},\partial\Delta\!^{[n_m]})),
\widehat\cup_{{\cal Z}^*(K;\partial\Delta\!^{[n_1]},{\cdots},\partial\Delta\!^{[n_m]})}))
\cong (H^{\,*}_{\!\RR_m}(K),\overline\cup_{K}).
$$

Similarly, take all $L_k=\partial\Delta\!^{[2]}$. Denote by $g$ the generator of $H^{\,1}_{\emptyset,[2]}(\partial\Delta\!^{[2]})$
and $H_{1}^{\emptyset,[2]}(\partial\Delta\!^{[2]})$.
Then $H_*^{\RR_2}(\partial\Delta\!^{[2]})$ is generated by $1$ and $g$ with coproduct given by
$\widehat\psi_{\partial\Delta\!^{[2]}}(1)=1{\otimes}1$ and $\widehat\psi_{\partial\Delta\!^{[2]}}(g)=1{\otimes}g{+}g{\otimes}1{+}g{\otimes}g$.
So the atom coproduct $\psi_{(q_k;\underline{\psi_k})}$ is the right strictly normal coproduct $\tilde\psi'$ in Definition~7.1.
With the same identification as above,
$$(H^{\,*}_{\!\RR_n}({\cal Z}^*(K;\partial\Delta\!^{[2]},{\cdots},\partial\Delta\!^{[2]})),
\widehat\cup_{{\cal Z}^*(K;\partial\Delta\!^{[2]},{\cdots},\partial\Delta\!^{[2]})}))
\cong (H^{\,*}_{\!\RR_m}(K),\tilde\cup^1_{K}),
$$
where $(H^{\,*}_{\!\RR_m}(K),\tilde\cup^1_{K})$ is determined by the right strictly normal cohomology algebra $(H^{\,*}_{\!\RR_m}(K),\tilde\cup_{K})$ as follows. For $a\in H^{\,s}_{\emptyset,\omega'}(K)$ and $b\in H^{\,t}_{\emptyset,\omega''}(K)$,
$a\tilde\cup^1_K b=(-1)^{t|\omega'|}a\tilde\cup_K b$.
\vspace{3mm}

{\bf Example 9.9} We compute the right universal cohomology algebra of the polyhedral join
${\cal Z}^*(K;\underline{CA},\underline{A})$, where $(\underline{CA},\underline{A})=\{(CA_k,A_k)\}_{k=1}^m$
and $C$ means the cone complex.

Apply Definition~9.5 to the pair $(CA_k,A_k)$ with reference index set $\RR_{n_k}$ and reference coproduct $\underline{\psi_k}=\widehat\psi^{(n_k)}$,
where $\widehat\psi$ is the universal coproduct in Definition~7.1.
Then $\theta_k\colon H_*^{\RR_{n_k}}(A_k)\to H_*^{\RR_{n_k}}(CA_k)$ is an epimorphism.
Suppose the vertex set of $A_k$ is a subset of $[n_k{-}1]$ and the new vertex of $CA_k$ is $n_k$.
Denote by $\hat A_k$ the same $A_k$ regarded as a simplicial complex on $[n_k{-}1]$.
For subsets $\tau,\tau_1,\cdots$ of $[n_k{-}1]$, we denote by $\tau',\tau'_1,\cdots$
the subsets $\tau\cup\{n_k\},\tau_1\cup\{n_k\},\cdots$ of $[n_k]$.
Then we have the following three group isomorphisms
$$\quad H_*^{\RR_{n_k-1}}(\hat A_k)={\mathpzc i}_{\,k*}=\oplus_{\omega\subset[n_k-1]}H_*^{\emptyset,\omega}(A_k),$$
$$n\colon H_*^{\RR_{n_k-1}}(\hat A_k)\to {\mathpzc n}_{\,k*}=\oplus_{\omega\subset[n_k-1]}H_*^{\emptyset,\omega'}(A_k),$$
$$\overline n\colon H_*^{\RR_{n_k-1}}(\hat A_k)\to \overline{\mathpzc n}_{\,k*}\subset\oplus_{\omega\subset[n_k-1]}T_*^{\emptyset,\omega'}(CA_k)$$
defined as follows.
The isomorphism $n$ is induced by the chain isomorphism from $\oplus_{\omega\subset[n_k-1]}\, T_*^{\emptyset,\omega}(\hat A_k)$
to $\oplus_{\omega\subset[n_k-1]}\, T_*^{\emptyset,\omega'}(A_k)$.
For $x=\Sigma\,k_i\tau_i\in H_*^{\RR_{n_k-1}}(\hat A_k)$ with $k_i\in\Bbb Z$ and $\tau_i\in (A_k)_{\emptyset,\omega_i}$,
$\overline n(x)=\Sigma\,k_i\tau'_i$.

Denote by $(\widehat\psi_K)^{\omega}_{\omega_1,\omega_2}\colon
T_*^{\emptyset,\omega}(K)\to T_*^{\emptyset,\omega_1}(K)\otimes T_*^{\emptyset,\omega_2}(K)$ the local coproduct in Theorem~7.6.
Suppose for $\omega,\omega_1,\omega_2\subset[n_k{-}1]$ and $\tau\in T_*^{\emptyset,\omega}(\hat A_k)$,
we have $(\widehat\psi_{\hat A_k})^\omega_{\omega_1,\omega_2}(\tau)=\pm\tau_1{\otimes}\tau_2$.
Then
$$(\widehat\psi_{A_k})^{\omega}_{\omega_1,\omega_2}(\tau)=\pm\tau_1{\otimes}\tau_2,\,\,
(\widehat\psi_{A_k})^{\omega}_{\omega'_1,\omega_2}(\tau)=\pm\tau_1{\otimes}\tau_2,$$
$$(\widehat\psi_{A_k})^{\omega}_{\omega_1,\omega'_2}(\tau)=\pm\tau_1{\otimes}\tau_2,\,\,
(\widehat\psi_{A_k})^{\omega}_{\omega'_1,\omega'_2}(\tau)=\pm\tau_1{\otimes}\tau_2,$$
$$(\widehat\psi_{A_k})^{\omega'}_{\omega_1,\omega_2}(\tau)=0\,({\scriptstyle\omega'{\setminus}(\omega_1{\cup}\omega_2)\notin A_k}),\,\,
(\widehat\psi_{A_k})^{\omega'}_{\omega'_1,\omega_2}(\tau)=\pm\tau_1{\otimes}\tau_2,$$
$$(\widehat\psi_{A_k})^{\omega'}_{\omega_1,\omega'_2}(\tau)=\pm\tau_1{\otimes}\tau_2,\,\,
(\widehat\psi_{A_k})^{\omega'}_{\omega'_1,\omega'_2}(\tau)=\pm\tau_1{\otimes}\tau_2,$$
$$(\widehat\psi_{CA_k})^{\omega'}_{\omega_1,\omega_2}(\tau')=\pm\tau_1{\otimes}\tau_2,\,\,
(\widehat\psi_{CA_k})^{\omega'}_{\omega'_1,\omega_2}(\tau')=\pm\tau'_1{\otimes}\tau_2,$$
$$(\widehat\psi_{CA_k})^{\omega'}_{\omega_1,\omega'_2}(\tau')=\pm\tau_1{\otimes}\tau'_2,\,\,
(\widehat\psi_{CA_k})^{\omega'}_{\omega'_1,\omega'_2}(\tau')=\pm\tau'_1{\otimes}\tau_2.$$
This implies that if ${\vartriangle_{\hat A_k}}(x)=\Sigma_k\,x_{1,k}{\otimes}x_{2,k}$ for $x\in H_*^{\RR_{n_k-1}}(\hat A_k)$, then
$${\vartriangle_{A_k}}(x)=\Sigma_k\,\big(x_{1,k}{\otimes}x_{2,k}+n(x_{1,k}){\otimes}x_{2,k}+x_{1,k}{\otimes}n(x_{2,k})+n(x_{1,k}){\otimes}n(x_{2,k})\big),$$
$${\vartriangle_{A_k}}(n(x))=\Sigma_k\,\big(n(x_{1,k}){\otimes}x_{2,k}+x_{1,k}{\otimes}n(x_{2,k})+n(x_{1,k}){\otimes}n(x_{2,k})\big)$$
$$\xi(q{\otimes}q)\widehat\psi_{CA_k}(\overline n(x))=\Sigma_k\,x_{1,k}{\otimes}x_{2,k}$$
So the coproduct $\vartriangle_{q_k}^{\!\XX}$ of the dense character chain coalgebra $C_{*}^{\XX\!;\XX_{n_k}}(CA_k,A_k)$
satisfies

${\vartriangle}_{q_k}^{\!\XX}({\mathpzc i}_{\,k*})\subset{\mathpzc i}_{\,k*}{\otimes}{\mathpzc i}_{\,k*}\oplus
{\mathpzc n}_{\,k*}{\otimes}{\mathpzc i}_{\,k*}\oplus
{\mathpzc i}_{\,k*}{\otimes}{\mathpzc n}_{\,k*}\oplus{\mathpzc n}_{\,k*}{\otimes}{\mathpzc n}_{\,k*}$.

${\vartriangle}_{q_k}^{\!\XX}({\mathpzc n}_{\,k*})\subset {\mathpzc n}_{\,k*}{\otimes}{\mathpzc n}_{\,k*}\oplus
{\mathpzc i}_{\,k*}{\otimes}{\mathpzc n}_{\,k*}\oplus {\mathpzc n}_{\,k*}{\otimes}{\mathpzc i}_{\,k*}$.

${\vartriangle}_{q_k}^{\!\XX}({\bar{\mathpzc n}}_{\,k*})\subset \overline{\mathpzc n}_{\,k*}{\otimes}{\mathpzc n}_{\,k*}\oplus
{\mathpzc i}_{\,k*}{\otimes}{\bar{\mathpzc n}}_{\,k*}\oplus {\bar{\mathpzc n}}_{\,k*}{\otimes}{\mathpzc i}_{\,k*}
\oplus {{\mathpzc i}}_{\,k*}{\otimes}{\mathpzc i}_{\,k*}$.

So the atom coproduct $\psi_{(q_k;\underline{\psi_k})}$ is the right universal coproduct $\widehat\psi'$ in Definition~7.1.
The coproduct $\vartriangle_{(q_k;\underline{\psi_k})}$ of the densely indexed homology coalgebra $H_*^{\XX\!;\RR_{n_k}}(CA_k,A_k)$
($={\mathpzc i}_{\,k*}\oplus{\mathpzc n}_{\,k*}$)
satisfies
${\vartriangle_{(q_k;\underline{\psi_k})}}(x)={\vartriangle_{(q_k;\underline{\psi_k})}}(n(x))={\vartriangle_{A_k}}(x)$.
So we have chain coalgebra isomorphism
\begin{eqnarray*}&\phi\colon\big(T_*^{\RR_m}(K)\,\widehat\otimes\,(\otimes_{k=1}^m H_*^{\RR;\RR_{n_k}}(CA_k,A_k)),
\widehat\psi_{K}\,\widehat\otimes\,(\otimes_{k=1}^m{\vartriangle_{q_k}})\big)\\
&\quad\quad\longrightarrow\,\,\big(T_*^{\RR_m}(K)\otimes(\otimes_{k=1}^m H_*^{\RR_{n_k-1}}(\hat A_k)),
\widehat\psi_{K}\otimes(\otimes_{k=1}^m{\widehat\vartriangle_{\hat A_k}})\big)
\end{eqnarray*}
defined as follows. For $a\in T_*^{\emptyset,\omega}(K)$ and $b_i\in H_*^{\RR_{n_k-1}}(\hat A_k)$, define
$$\phi(a\widehat\otimes(\otimes_{k\in\omega}n(b_k))\otimes(\otimes_{j\notin\omega}b_j))=
a{\otimes}b_1\otimes\cdots\otimes b_m.$$

Dually, we have
\begin{eqnarray*}&&(H^{\,*}_{\RR_{n}}({\cal Z}^*(K;\underline{CA},\underline{A})),\widehat\cup_{{\cal Z}^*(K;\underline{CA},\underline{A})})\\
&\cong&\big(H^{\,*}_{\RR_m}(K)\,\widehat\otimes\,(\otimes_{k=1}^m H^{\,*}_{\RR;\RR_{n_k}}(CA_k,A_k)),
\widehat\cup_{K}\,\widehat\otimes\,(\otimes_{k=1}^m{\triangledown_{\!q_k}})\big)\\
&\cong&\big(H^{\,*}_{\RR_m}(K)\otimes(\otimes_{k=1}^m H^{\,*}_{\RR_{n_k-1}}(\hat A_k)),
\widehat\cup_{K}\otimes(\otimes_{k=1}^m{\widehat\triangledown_{\!\hat A_k}})\big)
\end{eqnarray*}

From another point of view, $(CA_k,A_k)=(\Delta\!^{[1]}*\hat A_k,\{\emptyset\}*\hat A_k)$ and so
${\cal Z}^*(K;\underline{CA},\underline{A})\cong {\cal Z}^*(K;\Delta\!^{[1]},\{\emptyset\})*{\hat A_1}*{\cdots}*{\hat A_m}=K*{\hat A_1}*{\cdots}*{\hat A_m}$.
So the cohomology algebra equality holds by K\"{u}nneth theorem.
\vspace{3mm}

{\bf Theorem~9.10} {\it Suppose $(X_k,A_k)$ is a simplicial pair on $[n_k]$
densely normal on $\TT_k$ with respect to $\underline{\psi_k}$
for $k=1,{\cdots},m$ and
$(Y,B)$ is a simplicial pair on $[m]$ densely normal on $\underline{\SS}=\SS_1{\times}{\cdots}{\times}\SS_m$ with respect to $\w\psi^{(m)}$,
where $\w\psi$ is the normal coproduct in Definition~7.1 and
$\SS_k$ is the dense support index set of $(X_k,A_k)$ on $\TT_k$ with respect to $\underline{\psi_k}$.
Then the polyhedral join pair
$$({\cal Z}^*(Y;\underline{X},\underline{A}),{\cal Z}^*(B;\underline{X},\underline{A})),\quad
(\underline{X},\underline{A})=\{(X_k,A_k)\}_{k=1}^m$$
is densely normal on $\underline{\TT}=\TT_1{\times}{\cdots}{\times}\TT_m$ with respect to
$\underline{\psi_n}=\underline{\psi_1}{\otimes}{\cdots}{\otimes}\underline{\psi_m}$.
\vspace{2mm}

Proof}\, By Theorem~9.7, we have the commutative diagram
$$\begin{array}{ccc}
(T_*^{\underline{\TT}}({\cal Z}^*(B;\underline{X},\underline{A}),\psi_n,d)&\stackrel{\phi_B}{\longrightarrow}&
(T_*^{\underline{\SS}}(B)\widehat\otimes H_*^{\underline{\SS};\underline{\TT}}(\underline{X},\underline{A}),
\w\psi_B\widehat\otimes{\triangledown}_{\!(\underline{X},\underline{A};\underline{\psi_k})})\,\vspace{2mm}\\
\cap&&\cap\\
(T_*^{\underline{\TT}}({\cal Z}^*(Y;\underline{X},\underline{A}),\psi_n,d)&\stackrel{\phi_Y}{\longrightarrow}&
(T_*^{\underline{\SS}}(Y)\widehat\otimes H_*^{\underline{\SS};\underline{\TT}}(\underline{X},\underline{A}),
\w\psi_Y\widehat\otimes{\triangledown}_{\!(\underline{X},\underline{A};\underline{\psi_k})}),
\end{array}$$
where $\phi_B,\phi_Y$ are chain homotopy equivalences. So if the pair
$$(T_*^{\underline{\SS}}(Y)\widehat\otimes H_*^{\underline{\SS};\underline{\TT}}(\underline{X},\underline{A}),
T_*^{\underline{\SS}}(B)\widehat\otimes H_*^{\underline{\SS};\underline{\TT}}(\underline{X},\underline{A}))$$
is densely normal, then the theorem holds.
This is by Theorem~8.9 when we take $(D_{1\,*}^\TT,C_{1\,*}^\TT)=(T_*^{\underline{\SS}}(Y),T_*^{\underline{\SS}}(B))$
and $D_{2\,*}^\TT=C_{2\,*}^\TT=H_*^{\underline{\SS};\underline{\TT}}(\underline{X},\underline{A})$.
\hfill$\Box$\vspace{3mm}

Suppose the pair $(X,A)$ is densely split on $\TT\subset\XX_m$, then is $(X,A)$
densely normal on $\TT$ with respect to the normal coproduct $\w\psi^{(m)}$?
The answer should be negative but a counterexample is not found yet.

\end{document}